%% file: Thesis2.tex
\documentclass{article}
\input{headerarticle.tex}
\usepackage{amssymb}
\usepackage{young}
\usepackage[vcentermath]{youngtab}
\usepackage{graphicx}
\begin{document}
\large

\input{TitlePage}

\newpage
\begin{center} \textbf{Abstract} \end{center}
The fusion procedure provides a way to construct new solutions to the Yang-Baxter equation. In the case of the symmetric group the fusion procedure has been used to construct diagonal matrix elements using a decomposition of the Young diagram into its rows or columns. We present a new construction which decomposes the diagram into hooks, the great advantage of this is that it minimises the number of auxiliary parameters needed in the procedure. We go on to use the hook fusion procedure to find diagonal matrix elements computationally and calculate supporting evidence to a previous conjecture.\\
We are motivated by the construction of certain elements that allow us to generate representations of the symmetric group and single out particular irreducible components. In this way we may construct higher representations of the symmetric group from elementary ones. We go some way to generalising the hook fusion procedure by considering other decompositions of Young diagrams, specifically into ribbons. Finally, we adapt our construction to the quantum deformation of the symmetric group algebra known as the Hecke algebra.

\newpage
\tableofcontents
\newpage
\input{Acknowledgments}
\newpage
\input{Introduction}

\newpage
\input{TheHookFusionProcedure}
\newpage
\input{FindingEigenvalues}
\newpage
\input{TheRibbonFusionProcedure}
\newpage
\input{TheHookFusionProcedureForHeckeAlgebras}
\newpage
\begin{appendix}
\input{Appendix}

\end{appendix}
\newpage
\input{References}

\end{document}

%% file: headerarticle.tex
\newenvironment{proof}{\begin{trivlist} \item[]
{\bf Proof.}}{\nolinebreak
\hfill \rule{2mm}{2mm} \end{trivlist}}
\newcounter{ctr}

\newcounter{ctr1}
\newenvironment{listr}{\begin{list}{(\roman{ctr1})}%
{\usecounter{ctr1}
\setlength{\itemsep}{0mm} \setlength{\topsep}{-2mm}
\setlength{\leftmargin}{3em}}}{\end{list}}

\newcounter{ctr2}

\newenvironment{example}{\begin{trivlist} \item[]
{\addtocounter{definition}{1} \bf Example \thedefinition .\phantom{x}}}{\nolinebreak
\hfill \rule{2mm}{2mm} \end{trivlist}}

\newcounter{ctr3}

\newtheorem{definition}{Definition}[section]    %for book style
\newtheorem{theorem}[definition]{Theorem}
\newtheorem{lemma}[definition]{Lemma}
\newtheorem{corollary}[definition]{Corollary}
\newtheorem{proposition}[definition]{Proposition}

\newtheorem{conjecture}[definition]{Conjecture}
\newtheorem{result}[definition]{Result}
\newcommand{\NN}{\hbox{{\sf I}\kern-.15em\hbox{\sf N}}}
\newcommand{\ZZ}{\hbox{{\sf Z}\kern-.77em\hbox{\sf Z}\kern.15em}}
\newcommand{\QQ}{\hbox{{\sf I}\kern-.4em\hbox{\sf Q}}}
\newcommand{\RR}{\hbox{{\sf I}\kern-.15em\hbox{\sf R}}}
\newcommand{\CC}{\hbox{{\sf I}\kern-.4em\hbox{\sf C}}}

\newcommand{\ud}{\, {\rm d} \kern-.015em }

\newcommand{\mod}[1]{\left| \kern.05em #1 \kern.05em \right|}
\newcommand{\norm}[1]{\left\| \kern.05em #1 \kern.05em \right\|}
\newcommand{\inner}[1]{\left\langle \kern.05em #1 \kern.05em \right\rangle }

\newcommand{\pick}[2]{\renewcommand{\arraystretch}{0.6}
\left( \kern-.4em \begin{array}{c} #1 \\ #2 \end{array} \kern-.4em \right) }

%% file: TitlePage.tex
%\begin{titlepage}
\thispagestyle{empty}

\vspace{3cm}

\begin{center}
 \Large\textbf{The hook fusion procedure and its generalisations.}
  \end{center}
  \vspace{2cm}
  \begin{center}
   \textbf{James Grime}
    \end{center}
    \vspace{1cm}
   \begin{center}
   \textbf{A Thesis Submitted for the Degree of PhD.}
   \end{center}
   \vspace{1cm}
   \begin{center}
   \textbf{Department of Mathematics\\ University of York\\ Heslington\\ 
    York~Y010~5DD}
   \end{center}
   \vspace{1cm}
   \begin{center}
   \textbf{July 2007}
\end{center}

%\end{titlepage}

%% file: Acknowledgments.tex
\begin{center} \textbf{Acknowledgments} \end{center}

There are many people I wish to thank. First of all I would like to thank everyone in the department of mathematics at the University of York, especially Stephen Donkin, Niall MacKay and my supervisor Maxim Nazarov for all their help and guidance. Also I would like to give special thanks to the PhD students of the department for all their friendship and support. I thank EPSRC for funding my research, and all those of the GAP support mailing list, especially Alexander Konovalov and Alexander Hulpke, for their help and their amazing patience. \\And finally, but most of all, I'd like to thank my parents, Rob and Liz Grime, without whose advice, love and support this would not have been possible.

\quad

\quad

\quad

\quad

\quad

\quad

\quad

\quad

\quad

\quad

\quad

\quad

\begin{flushright}
\emph{Moriarty: How are you at Mathematics? \\
Seagoon: I speak it like a native.}
\end{flushright}

\newpage
\textbf{Declaration}

I, James Grime, confirm that the work presented in this thesis is my own. Where information has been derived from other sources, I confirm that this has been indicated in the thesis.

%% file: Introduction.tex
\section{Introduction}

There has been much interest in the Yang-Baxter equation since it first appeared in the works of McGuire in 1964 and Yang in 1967, \cite{MG}, \cite{Ya}. They considered a quantum mechanical many body problem and found the scattering matrix factorised to that of the two body problem. Here the Yang-Baxter equation arises as the consistency condition for the factorisation. It has since continued to appear in many disparate branches of mathematics including solvable lattice models, soliton theory and quantum integrable models, and quantum groups. Therefore there has been much study into the Yang-Baxter equation and the construction of solutions.\\ In this thesis we present the construction of certain solutions in the symmetric group algebra and deformations of this algebra. We use a method of construction known as `the fusion procedure' which requires introducing a number of spectral parameters before taking a certain product of operators. The result is an operator which can be used to generate a representation of the symmetric group algebra. This representation decomposes into irreducible representations, and for certain values of the spectral parameters, the fusion procedure gives us a way to single out one of these components.\\ We present a new fusion procedure which minimises the number of auxiliary parameters necessary in the procedure. We will call this new fusion procedure the `hook fusion procedure' and go on to use it to construct elements of the symmetric group computationally, and thus provide supporting evidence to a previous conjecture. We will also generalise our results in two ways. We will construct a hook fusion procedure for Hecke algebras, which is a deformation of the symmetric group algebra, and we will go some way to generalising the fusion procedure beyond our hook version, which we will call the `ribbon fusion procedure'.

In chapter 2 we give some of the background to the problem, this includes a description and formulation of the Yang-Baxter equation itself, and its variant known as the classical Yang-Baxter equation. We state several well known definitions relating to the representation theory of the symmetric group, including the definition of a Young diagram and its Young symmetrizer. The Young symmetrizers are primitive idempotents so may be used to form simple left ideals. We then go on to describe diagonal matrix elements, related to the Young symmetrizers by a certain multiplication of invertible elements, and so may use these elements to generate irreducible representations. We also give a description of the fusion procedure for constructing solutions to the Yang-Baxter equation, first in a general setting and then specifically for the symmetric group algebra due to Cherednik \cite{C1} and Nazarov \cite{N3}. We show how the fusion procedure may be used to construct the diagonal matrix elements in the symmetric group.

In Chapter 3 we present one of our main results, which appears in the author's first paper \cite{Gr1}. Up to now the fusion procedure for the symmetric group required a Young diagram to be decomposed into either its rows or its columns. We consider a different decomposition of the Young diagram, namely we decompose the diagram into its `principal hooks'. Not only does this reduce the number of auxiliary parameters needed in the fusion procedure, but it actually minimises the number. This will prove to be useful in later computational calculations. We introduce an element $F_\Lambda(z_1, \dots, z_n)$ whose coefficients are rational polynomials in the parameters $z_1, \dots ,z_n$, which may contain a singularity on the line $z_1= \cdots = z_n$. After restriction to a certain subspace we will prove the function $F_\Lambda(z_1, \dots, z_n)$ is regular on the line, with any singularities being removable. By showing certain divisibilities we will prove that the value of this function on the line coincides with the diagonal matrix element. We call Theorem \ref{jimtheorem1} of Chapter 3 the hook fusion procedure. We are motivated by the construction of certain induced representations and their decomposition into irreducible representations.

We use our main result of Chapter 3 in Chapter 4 to calculate diagonal matrix elements using the hook fusion procedure. Since the hook fusion procedure minimises the number of necessary spectral parameters we may use it to calculate the diagonal matrix elements computationally. We are motivated by the irreducibility criterion of representations of degenerate affine Hecke algebras, which we define. It has been conjecture that this criterion can be defined solely in terms of certain eigenvalues. We provide supporting evidence that this is the case. We also go on to calculate more eigenvalues of a type not previously calculated. It is interesting to note that these eigenvalues all have integral roots.

In Chapter 5 we propose a possible ribbon fusion procedure which may be considered a generalisation of the hook fusion procedure and those procedures that relied on the decomposition of a Young diagram into its rows or columns. In this chapter we present many new ideas, including a new combinatorial construction of ribbons from $2 \times 2$ squares, and a new approach to proving such a ribbon fusion procedure based on chains of diagrams that resolve singularities of $F_\Lambda(z_1, \dots, z_n)$ using a common ordering. We also provide evidence supporting the opposite implication, that if the decomposition of the Young diagram is not as prescribed then the singularities of $F_\Lambda(z_1, \dots, z_n)$ are not removable and are in fact poles. We also mention Schur functions, $s_\lambda(\textbf{x})$, and how the decomposition of Young diagrams into row, columns, hooks and ribbons can give determinantal expressions for $s_\lambda(\textbf{x})$.

Finally, in Chapter 6 we adapt the hook fusion procedure of Chapter 3 to find representations of Hecke algebras, denoted $H_n$. The Hecke algebra may be viewed as a deformation of the symmetric group algebra, and like the symmetric group algebra, other fusion procedures for Hecke algebras have already been considered. However, the hook fusion procedure again minimises the now multiplicative spectral parameters. As an adaptation the proof both built upon the ideas developed in Chapter 3 while in detail contains several necessary changes. We finish by describing how the element defined in the limit of this procedure generates irreducible representations of $H_n$. We also touch on how this relates to the affine Hecke algebra and is the subject of the author's second paper \cite{Gr2}.

\newpage
\section{Background}
\subsection{The Yang-Baxter equation}

%The Yang-Baxter equation first manifested itself in the work of McGuire in 1964 and Yang in 1967. They considered a quantum mechanical many body problem and found the scattering matrix factorised to that of the two body problem. Here the Yang-Baxter equation arises as the consistency condition for the factorisation. \\
%The hunt for solvable lattice models culminated in Baxter's solution of the eight vertex model in 1972. This contributed to Zamolodchikov's work of 1979 on factorised $S$-matrices in two dimensional quantum theory, where it was pointed out that the mechanism at work was the same as that in Baxter's and other's works.\\
%In 1978-79 Faddeev, Sklyanin and Takhtajan proposed the quantum inverse method as a unification of the classical integrable models, known as soliton theory, and the quantum integrable models. In their theory the basic commutation relation of operators is described by a solution of the Yang-Baxter equation, and the terminology itself is due to them.\\
%These works led to the idea  of introducing certain deformations of groups or Lie algebras, called quantum groups by Drinfeld. At about the same time there appeared the discovery of new invariants of links, and subsequently the aspect of the Yang-Baxter equation as the braid-type relation has been brought to attention. Closely related structures have also been revealed in conformal field theory.

We may view the Yang-Baxter equation as a relation involving operators of a certain vector space. This property turns out to be an important condition arising in several areas of mathematics, including solvable lattice models, quantum theory and conformal field theory. 

The Yang-Baxter equation first appeared in 1964 in the study of the quantum mechanical many body problem, \cite{MG}, \cite{Ya}. McGuire and Yang found the scattering matrix factorised to that of the two body problem, see \cite{C4}. The consistency condition of that factorisation manifests itself as the Yang-Baxter equation.\\
The term `Yang-Baxter equation' was given by Faddeev, Skylanin and Takhtajan in their papers of 1978-79 which describe the quantum inverse method. This work was a unification of the classical integrable models, known as soliton theory, and the quantum integrable models. Here the commutation relation between operators is described by a solution of the Yang-Baxter equation, \cite{FST1}, \cite{FST2}.
These works led to the idea of quantum groups as certain deformations of groups, or Lie algebras, and consequently the study of Yang-Baxter equation as a type of braid relation, \cite{FRT}, \cite{Sk}, \cite{D1}, \cite{D3}, \cite{Ji1}.

%\subsubsection{The Yang-Baxter equation} 

Let $V$ be a complex vector space, and $R(u)$ a function of $u \in \mathbb{C}$ taking values in $\textrm{End}_{\mathbb{C}}(V \otimes V)$. Furthermore, let $R_{ij}$ denote the operator on $V^{\otimes 3}$, acting as $R(u)$ on the $i$-th and $j$-th components and as the identity on the other component, for example $R_{23}(u) = I \otimes R(u)$. Then we call the following equation for $R(u)$ the \emph{Yang-Baxter equation}: \begin{equation}\label{YBE} R_{12}(u)R_{13}(u+v)R_{23}(v) = R_{23}(v)R_{13}(u+v)R_{12}(u). \end{equation} The variable $u$ is called the spectral parameter. In most cases we assume that $N = \textrm{dim} V < \infty$, hence a solution of (\ref{YBE}) is often referred to as an $R$ matrix. For example, a typical solution of (\ref{YBE}) in the case $V = \mathbb{C}^2$ might be 

\[ R(u) = 
\left[ \begin{array}{cccc}
	1+u & 0 & 0 & 0\\
	0&u&1&0\\
	0&1&u&0\\
	0&0&0&1+u
\end{array} \right] = P + uI, \] 

where $P \in \textrm{End}_\mathbb{C} (V \otimes V)$ denotes the transposition $P(x\otimes y) = y \otimes x$.

We may extend the idea to matrices acting on $V^{\otimes n}$, $n \geqslant 2$, by defining $R_{ij}(u)$ as acting as $R(u)$ on the $i$-th and $j$-th components, such that \[ R_{ij}(u)R_{ik}(u+v)R_{jk}(v) = R_{jk}(v)R_{ik}(u+v)R_{ij}(u), \] then for pairwise distinct indices \[R_{ij}(u)R_{kl}(v) = R_{kl}(v)R_{ij}(u). \] 

%We may generalise the Yang-Baxter equation by considering instead a family of vector spaces $\mathcal{F} = \{V\}$ and operators $\{ R_{V_i V_j(u)} \in \textrm{End}_\mathbb{C}(V_i \otimes V_j) : V_i , V_j  \in \mathcal{F} \}$. The Yang-Baxter equation (\ref{YBE}) is then an equation in $\textrm{End}_\mathbb{C}(V_1 \otimes V_2 \otimes V_3)$, where $R_{ij}(u) = R_{V_i V_j}$. \\

Further, one may consider the Yang-Baxter equation for a function of two variables $R(u,u')$: \begin{equation}\label{YBE3}R_{12}(u_1,u_2)R_{13}(u_1,u_3)R_{23}(u_2,u_3) = R_{23}(u_2,u_3)R_{13}(u_1,u_3)R_{12}(u_1, u_2). \end{equation} Equation (\ref{YBE}) is a special case of (\ref{YBE3}) where the $(u,u')$-dependence enters only through the difference \[ R(u,u') \equiv R(u-u'). \]

We could rewrite the Yang-Baxter equation in terms of the matrix $\check{R}(u) = PR(u)$. Let $n \geqslant 2$, and define matrices on $V^{\otimes n}$ by $\check{R}_i(u) = I \otimes \cdots \otimes \check{R}(u) \otimes \cdots \otimes I$, where $\check{R}(u)$ is in the $(i, i+1)$-th slot, $i = 1, \dots, n-1$. One then has \[ \check{R}_i(u) \check{R}_j(v) = \check{R}_j(v)\check{R}_i(u) \qquad \textrm{ if } |i-j| > 1, \]\[\check{R}_{i+1}(u)\check{R}_i(u+v)\check{R}_{i+1}(v) = \check{R}_i(v)\check{R}_{i+1}(u+v)\check{R}_i(u). \] Notice that, without the spectral parameters $u, v$, these relations are simply Artin's braid relations. For an introduction to the Yang-Baxter equation see \cite{Ji}.

%Suppose $V_1=V_2=V$. Regarding $\textrm{End}_\mathbb{C}(V \otimes V_3) = \textrm{End}_\mathbb{C}(V) \otimes \mathcal{A}$, $\mathcal{A} = \textrm{End}_\mathbb{C}(V_3)$, let us write $R_{V V_3}(u)$ as $T(u) = \sum t_{ij}(u)E_{ij}$ with $t_{ij}(u) \in \mathcal{A}$ and $E_{ij} = (\delta_{ia} \delta_{jb})_{a,b= 1, \dots, N}$. In this notation (\ref{YBE}) becomes \begin{equation}\label{YBE2} \check{R}(u-v)T(u)\otimes T(v) = T(v) \otimes T(u)\check{R}(u-v). \end{equation} Here $T(u) \otimes T(v) = \sum t_{ij}t_{kl} E_{ij} \otimes E_{kl}$. Equation (\ref{YBE2}) can be viewed as giving commutation relations among the generators $t_{ij}(u)$ of an abstract algebra $\mathcal{A}$. The Yang-Baxter equation for $\check{R}(u)$ guarantees the associativity of $\mathcal{A}$ thus defined. An important feature of (\ref{YBE2}) is that the comultiplication \[ \Delta : \mathcal{A} \to \mathcal{A} \otimes \mathcal{A}, \qquad \Delta(t_{ij}(u)) = \sum_k t_{ik}(u) \otimes t_{kj}(u) \] preserves the relations (\ref{YBE2}). The formulas (\ref{YBE}) and (\ref{YBE2}) are the basic constituents in the quantum inverse method.

\subsubsection{The classical Yang-Baxter equation} If solution of the Yang-Baxter equation contains an extra parameter $\hbar$ such that as $\hbar \rightarrow 0$ it has the expansion \[ R(u, \hbar) = (\textrm{scalar}) \times (I + \hbar r(u) + O(\hbar^2)), \] then we call $R(u, \hbar)$ quasi-classical, and $r(u) \in \textrm{End}_\mathbb{C}(V \otimes V)$ the classical limit of  $R(u, \hbar)$. For quasi-classical $R$-matrices, the Yang-Baxter equation (\ref{YBE}) implies the following \emph{classical Yang-Baxter equation} for $r(u)$ (\cite{BD}, \cite{KS}):
\begin{equation}\label{CYBE} [r_{12}(u), r_{13}(u+v)] + [r_{12}(u), r_{23}(v)] + [r_{13}(u+v), r_{23}(v)] = 0. \end{equation} 

Not all solutions of the Yang-Baxter equation are quasi-classical. However, quasi-classical solutions are a significant subset. In particular, one can see that the classical Yang-Baxter equation is formulated using only the Lie algebra structure of $\textrm{End}(V)$. We continue by giving some important results for the classical case.

Let $\mathfrak{g}$ be a Lie algebra, and let $r(u)$ be a $\mathfrak{g} \otimes \mathfrak{g}$-valued function. In terms of a basis $\{X_\mu\}$ of $\mathfrak{g}$, write  \[ r(u) = \sum_{\mu, \nu} r^{\mu \nu}(u) X_\mu \otimes X_\nu \] with $\mathbb{C}$-valued functions $r^{\mu \nu}(u)$. Further, let $r_{12}(u) = \sum r^{\mu \nu}(u) X_\mu \otimes X_\nu \otimes I \in (U \mathfrak{g})^{\otimes 3}$ and so on, where $U \mathfrak{g}$ denotes the universal enveloping algebra. One then has \[ [r_{12}(u), r_{23}(v)] = \sum r^{\mu \nu}(u)r^{\sigma \rho}(u) X_\mu \otimes [X_\nu, X_\rho] \otimes X_\sigma, etc \] so that each term in (\ref{CYBE}) actually lies inside $\mathfrak{g}^{\otimes 3}$. For each triplet of representations $(\pi_i, V_i), i = 1, 2, 3$ of $\mathfrak{g}$, $(\pi_i \otimes \pi_j)(r_{ij}(u))$ gives a matrix solution of the classical Yang-Baxter equation in $V_1 \otimes V_2 \otimes V_3$.

We say $r(u)$ is \emph{regular} on a subset of the complex plane if it is a holomorphic function everywhere except on a set of isolated points, which are poles for the function, and call $r(u)$ \emph{non-degenerate} if $\det (r^{\mu \nu}(u)) \neq 0$. Then we have the following result.

\begin{theorem}[Belavin-Drinfeld] If $\mathfrak{g}$ is a finite dimensional complex simple Lie algebra, $\{X_\mu\}$ an orthonormal basis of $\mathfrak{g}$, and $r(u)$ is a regular, non-degenerate $\mathfrak{g} \otimes \mathfrak{g}$-valued solution of (\ref{CYBE}) defined in a neighbourhood of $0 \in \mathbb{C}$, then 

\begin{listr}
\item $r(u)$ extends regularly to the whole complex plane $\mathbb{C}$, with all its poles being simple.
\item $\Gamma = \textrm{ the set of poles of } r(u)$ is a discrete subgroup relative to the addition of $\mathbb{C}$.
\item As a function of $u$ there are the following possibilities for the $r^{\mu \nu}(u)$:
\begin{trivlist}
	\item rank $\Gamma = 2$ : elliptic function;
	\item rank $\Gamma = 1$ : trigonometric function (i.e. a rational function in the variable $e^{\textrm{\scriptsize const. } u}$);
	\item rank $\Gamma = 0$ : rational function.
\end{trivlist}
\end{listr}
\end{theorem}

Belavin-Drinfeld further show that elliptic solutions exist only for $\mathfrak{g} = \mathfrak{sl}(n)$, in which case it is unique up to a certain equivalence of solutions. They also show that trigonometric solutions exist for each type, and can be classified using the Dynkin diagram for affine Lie algebras. Although we will not use the Yang-Baxter equation in this form ourselves, these remain an important class of solutions.

\subsubsection{The fusion procedure} There are a few approaches to constructing $R$ matrices. For instance, for finite dimensional matrix solutions one may solve a system of linear equations that uniquely determine $R(u)$, up to a scalar factor. Here we describe another method known as the \emph{fusion procedure} developed by Kulish, Reshetikhin and Skylanin, \cite{KRS}. This method relates to that of a standard technique to get irreducible representations of Lie algebras, in which we form a tensor product of fundamental representations and decompose it.

We will describe the procedure in a more general setting than the ones we have so far considered. We begin by considering instead a family of vector spaces $\mathcal{F} = \{V_i\}$ and operators $\{ R_{V_i V_j}(u) \in \textrm{End}_\mathbb{C}(V_i \otimes V_j) : V_i , V_j  \in \mathcal{F} \}$. The Yang-Baxter equation (\ref{YBE}) is then an equation in $\textrm{End}_\mathbb{C}(V_1 \otimes V_2 \otimes V_3)$, where $R_{ij}(u) = R_{V_i V_j}(u)$. 

Fix $u_1, u_2$ and set \[ W = R_{V_1 V_2}(u_2-u_1)(V_1 \otimes V_2) \subseteq V_1 \otimes V_2, \] and define \[ R_{V_2 \otimes V_1, V_3}(u) = R_{V_2 V_3}(u + u_1) R_{V_1 V_3}(u+u_2).\] Then using the Yang-Baxter equation (\ref{YBE}) one finds 

\begin{normalsize}
\begin{eqnarray*} 
R_{V_2 \otimes V_1, V_3}(u) (W \otimes V_3) 
&=& R_{V_2 V_3}(u+u_1) R_{V_1 V_3}(u+u_2) R_{V_1 V_2}(u_2 - u_1) (V_1 \otimes V_2 \otimes V_3) \\
&=& R_{V_1 V_2}(u_2 - u_1)R_{V_1 V_3}(u+u_2) R_{V_2 V_3}(u+u_1)  (V_1 \otimes V_2 \otimes V_3) \\
&\subseteq & W \otimes V_3. 
\end{eqnarray*}
\end{normalsize}

And so, with appropriate choices of $u_1$ and $u_2$, $W$ becomes a proper, non-trivial, subspace giving new $R$ matrices \[ R_{W V_3}(u) = \left. R_{V_2 \otimes V_1, V_3}(u) \right|_{W \otimes V_3}, \] and similarly \[ R_{V_3 \otimes W}(u) = \left. R_{V_3, V_2 \otimes V_1}(u) \right|_{V_3 \otimes W}. \]

\subsection{Young symmetrizers}

Soon we will apply the fusion procedure to representations of the symmetric group, but first we give some well known classical results.

Let $K$ be a field, and $e_1, e_2, \dots, e_k$ the set of orthogonal idempotents in a $K$-algebra $\mathcal{A}$, in other words we have $e_i^2 = e_i$ and $e_ie_j = 0 = e_je_i$ for all $1 \leqslant i \neq j \leqslant k$. Any idempotent $e$ in $\mathcal{A}$ is \emph{primitive} (or minimal) if $e$ cannot be written as a sum of two non-zero orthogonal idempotents. The importance of primitive idempotents derives from the following proposition.

\begin{proposition} Let $e \in \mathcal{A}$ be an idempotent element, then $e$ is primitive if and only if the left ideal $\mathcal{A} \cdot e$ is indecomposable. \end{proposition}

From now on we will work exclusively in the field $K = \mathbb{C}$. No direct construction of the element $e$ is known for a general $\mathcal{A}$, but if $\mathcal{A}$ is the symmetric group algebra $\mathbb{C}S_n$ there exists the following construction due to Alfred Young, where the left ideals are not only indecomposable but also irreducible, \cite{Y1}.

A partition $\lambda = (\lambda_1, \lambda_2, \dots, \lambda_k)$ of a positive integer $n$ is a sequence of weakly decreasing integers whose sum is equal to $n$. This may be represented by its \emph{Young diagram} which consists of $k$ rows with $\lambda_i$ boxes in the $i^{\mbox{\scriptsize th}}$ row, and with each row being left justified. The coordinates of the boxes in such a diagram are labeled using matrix notation, that is to say we let the first coordinate $i$ increase as one goes
downwards, and the second coordinate $j$ increase from left to
right. For example the partition $\lambda = (3,3,2)$ is a partition of $8$ and gives the
diagram 

\[ \yng(3,3,2)\]

Let us identify a partition $\lambda$ with its corresponding diagram. If $\lambda_1, \lambda_2, \dots , \lambda_k$ are the row lengths of $\lambda$ then we denote its columns by $\lambda_1', \lambda_2', \dots, \lambda_l'$. The conjugate of a diagram $\lambda$ is then a reflection in the main diagonal and is denoted $\lambda' = (\lambda_1', \lambda_2', \dots, \lambda_l')$. Furthermore, if two diagrams correspond to partitions $\lambda = (\lambda_1, \dots, \lambda_k)$ and $\mu = (\mu_1, \dots, \mu_l)$, we write $\mu \subset \lambda$ if  $\mu_i \leqslant \lambda_i$ for all $i$, and say `$\lambda$ contains $\mu$'. A skew diagram is obtained by removing a Young diagram from one that contains it. The resulting skew shape is denoted $\lambda/\mu$. For example, if $\lambda = (4,3,3,1)$ and $\mu = (2,1)$, we get the following skew diagram:

\begin{center}
\begin{picture}(40,40)
\put(0,12){\framebox(12,12)[r]{}}
\put(0,0){\framebox(12,12)[r]{}}
\put(12,24){\framebox(12,12)[r]{}}
\put(12,12){\framebox(12,12)[r]{}}
\put(24,36){\framebox(12,12)[r]{}}
\put(24,24){\framebox(12,12)[r]{}}
\put(24,12){\framebox(12,12)[r]{}}
\put(36,36){\framebox(12,12)[r]{}}
\end{picture}
\end{center}

Any way of putting a positive integer in each box of the Young
diagram is called a $filling$. A \emph{semistandard tableau}, is a
filling of a Young diagram $\lambda$ that is both weakly
increasing across each row and strictly increasing down each
column. A \emph{standard tableau} $\Lambda$ of shape $\lambda$ is a semistandard tableau of shape $\lambda$ in which the
entries are the numbers 1 to $n$, each occurring once. For example 

\[ \young(147,258,36)\] 

is a standard tableau of shape $(3,3,2)$. In this thesis we are only interested in standard tableaux and so may refer to one simply as a tableau. We will always denote tableau of a diagram using the corresponding capital Greek letter. 

Since a permutation $\sigma \in S_n$ act on a tableau $\Lambda$ by permutation of its entries, the row group of $\Lambda$, $R(\Lambda)$, is the set of all permutations in the symmetric group $S_n$ that preserve the rows of $\Lambda$.  For example, if $\Lambda$ is as above then (1 4 7)(2 8) $\in R(\Lambda)$. The group $R(\Lambda)$ is a subgroup of $S_n$ isomorphic to $S_{\lambda_1} \times S_{\lambda_2} \times \cdots \times S_{\lambda_k}$.\\
Similarly we have the column group of $\Lambda$, $C(\Lambda)$, as the set of all permutations that preserve the columns of $\Lambda$.

We denote sign homomorphism by sgn and define it to be the homomorphism from the symmetric group $S_n$ to the multiplicative group $\{1, -1\}$ such that sgn$t$ = -1, for any transposition $t \in S_n$. We now define two elements of $\mathbb{C}S_n$ \[ P_\Lambda = \sum_{p \in R(\Lambda)}p  \qquad \textrm{and} \qquad Q_\Lambda = \sum_{q \in C(\Lambda)} (\textrm{sgn} q) \cdot q . \] Then $Y_\Lambda = P_\Lambda Q_\Lambda$ is the \emph{Young symmetrizer} of $\Lambda$. The element $Y_\Lambda$ is a scalar multiple of a primitive idempotent, in particular \begin{equation}\label{youngidempotent}Y_\Lambda^2 = \alpha_\Lambda \cdot Y_\Lambda,\end{equation} for some rational number $\alpha_\Lambda \in \mathbb{Q}$. Consequently we have the following theorem, \cite{Y1}.

\newpage
\begin{theorem}[Young] For any standard tableau $\Lambda$ of shape $\lambda$, the ideal $\mathbb{C}S_n \cdot Y_\Lambda$ of $\mathbb{C}S_n$ is an irreducible left regular $\mathbb{C}S_n$-module, denoted $V_\Lambda$. Two
modules, $V_\Lambda$ and $V_{\Lambda'}$ obtained from tableaux
$\Lambda$ and $\Lambda'$ in this way are isomorphic if and only if
$\Lambda$ and $\Lambda'$ are of the same shape. Furthermore, by choosing one standard tableau for each partition we obtain a full set of irreducible left $\mathbb{C}S_n$-modules. \end{theorem}

Therefore irreducible representations of the group algebra $\mathbb{C}S_n$ are parameterized by partitions of $n$. We will use $V_\lambda$ to denote $V_\Lambda$ if $\Lambda$ is an unspecified standard tableau of shape $\lambda$. The coefficient $\alpha_\Lambda$ in (\ref{youngidempotent}) is equal to $\frac{n!}{\textrm{\small dim}V_\lambda}$.

In the following example, for any two distinct numbers $1 \leqslant p < q \leqslant n$, we denote the transposition of $p$ and $q$ by $(pq) \in S_n$.

\begin{example}\label{youngexample} The partitions of $n=3$ are
\[\lambda = \yng(3) \phantom{x}, \qquad\qquad\mu=\yng(1,1,1) \qquad \textrm{and} \qquad \nu=\yng(2,1) \phantom{x}.\qquad\] 

If $N = \young(13,2)$ \phantom{x} then \begin{eqnarray*} 
Y_N &=& \left( 1 + (1\phantom{x} 3)\right)\left(1-(1\phantom{x}2)\right) \\
		&=& 1 + (1\phantom{x}3) - (1\phantom{x}2) - (1\phantom{x}2\phantom{x}3). \end{eqnarray*}
		
And $\mathbb{C}S_3 \cdot Y_N$ has basis $Y_N$, (1 2)$Y_N$, and is a 2-dimensional irreducible $\mathbb{C}S_3$-module.
\end{example}

\subsection{The fusion procedure for the symmetric group}

The symmetric group $S_n$ acts on $V^{\otimes n}$ by permuting its componants.

Let \[ \{1\} = S_0 \subset S_1 \subset S_2 \subset \dots \] be a
chain of symmetric groups with the standard embedding. Let $V_\mu$ be an irreducible $S_{n-1}$-module, $V_\lambda$ an irreducible $S_n$-module and $k$ the multiplicity of $V_\mu$ in the restriction
of $V_\lambda$ to the group $S_{n-1}$. For symmetric groups we always have $k \in \{0,1\}$. By induction
we obtain a decomposition of a module in $V_\lambda$ into
irreducible $S_0$-modules. These are simply 1-dimensional subspaces indexed by all
standard tableaux $\Lambda$ of shape $\lambda$, thus we have the decomposition
\[ V_\lambda \cong \bigoplus_\Lambda \mathbb{C} v_\Lambda . \]   Fix an $S_n$-invariant
inner product $(\phantom{x} ,\phantom{x} )$ such that $(v_\Lambda,
v_\Lambda) = 1$ for every $\Lambda$. The vectors $v_\Lambda$ are pairwise orthogonal and are known as the \emph{Young basis} of $V_\lambda$, \cite{Y2}. 

In the algebra of $n \times n$ matrices $\textrm{Mat}_n(\mathbb{C})$, let $E_{ij}$ denote the elementary matrix whose $(i, j)^{\mbox{\scriptsize th}}$ entry is 1 with zeros elsewhere. The elementary matrices have the following properties \[ E_{ij}E_{kl} = \delta_{jk} E_{il}, \quad \textrm{ and } \]\[ E_{11} + E_{22} + \cdots + E_{nn} = I_n, \] where $\delta_{jk}$ is the Kronecker delta product and $I_n$ is the identity matrix. Note that the $E_{ii}$ form a complete set of pairwise orthogonal idempotents. Any set of matrices satisfying the two properties above is called a set of \emph{matrix units}.

The group algebra $\mathbb{C}S_n$ is isomorphic to the direct sum of matrix algebras \[ \mathbb{C}S_n \cong \bigoplus_{\lambda} \textrm{Mat}_{f_\lambda}(\mathbb{C}), \] where the sum is over all partitions of $n$ and $f_\lambda$ is equal to the number of standard tableaux of shape $\lambda$. The matrix units $E_{\Lambda \Lambda'} \in \textrm{Mat}_{f_\lambda}(\mathbb{C})$ are parameterised by pairs of standard tableaux, however when considering diagonal matrix units we shall write $E_\Lambda = E_{\Lambda \Lambda}$. 

Let $F_{\Lambda \Lambda'}$ be the element \[ F_{\Lambda \Lambda'} = \sum_{\sigma \in S_n} (v_\Lambda, \sigma v_{\Lambda'}) \sigma \in \mathbb{C}S_n.\] And denote the element $\frac{f_\lambda}{n!} F_{\Lambda \Lambda'}$ by $\overline{E}_{\Lambda \Lambda'}$. Then there exists an isomorphism from the matrix algebras to the symmetric group algebras sending $E_{\Lambda \Lambda'}$ to $\overline{E}_{\Lambda \Lambda'}$ such that the elements $\overline{E}_\Lambda$ act diagonally on the Young basis, that is to say $\overline{E}_\Lambda v_{\Lambda'} = \delta_{\Lambda \Lambda'}v_\Lambda$.

Murphy's construction of the diagonal matrix units, \cite{M1}, used certain elements $X_1, \dots, X_n \in \mathbb{C}S_n$ known as Jucys-Murphy elements where \[ X_1 = 0 \phantom{xxxxxxxxxxxxxxxxxxxxxxxxxxxxxxxxxxxx} \] \begin{equation}\label{JMelements} X_k = (1 \;\; k) + (2 \;\; k) + \cdots + (k-1 \;\; k), \quad k = 2, \dots, n. \end{equation} If $M$ is a standard tableaux with $n-1$ boxes then \[ \overline{E}_M = \sum_{M \to \Lambda} \overline{E}_\Lambda, \] where $M \to \Lambda$ means that the standard tableaux $\Lambda$ is obtained from $M$ by adding one box containing the entry $n$. We then have the action \begin{equation}\label{content} X_k \overline{E}_\Lambda = \overline{E}_\Lambda X_k = c_k(\Lambda) \overline{E}_\Lambda,\end{equation} where $c_k(\Lambda) = j-i$ if the box $(i,j) \in \lambda$ is filled with the number $k$ in
$\Lambda$. The difference $j-i$ is the \emph{content} of box
$(i,j)$.

Let $z$ be a complex variable, and $M$ the standard tableau obtained from $\Lambda$ after removing the box containing $n$. Murphy's recurrence relation for diagonal matrix units may then be written in terms of Jucys-Murphy elements as the following limit; \[ \overline{E}_\Lambda = \left. \overline{E}_M \frac{z}{z+c_n(\Lambda) - X_n} \right|_{z=0} . \]

We shall now move from diagonal matrix units and instead proceed using the \emph{diagonal
matrix element}, $F_\Lambda$, defined by,
\begin{equation}\label{dme} F_\Lambda = \sum_{\sigma \in S_n} (v_\Lambda, \sigma v_\Lambda)
\sigma \in \mathbb{C}S_n. \end{equation}

We naturally inherit the following equalities \[F_\Lambda^2 = \frac{n!}{f_\lambda} \cdot F_\Lambda, \qquad F_\Lambda v_{\Lambda'} = \frac{n!}{f_\lambda} \cdot \delta_{\Lambda \Lambda'} v_\Lambda \quad \textrm{and} \quad X_k F_\Lambda = F_\Lambda X_k = \frac{f_\lambda}{n!} \cdot c_k(\Lambda) F_\Lambda.\] 

We define the \emph{row tableau} of shape $\lambda$, $\Lambda^r$, to be the standard tableau obtain by filling the first row of $\lambda$ consecutively with the numbers 1 to $\lambda_1$, then filling the second row with $\lambda_1 +1$ to $\lambda_2$, and so on. We may similarly define a \emph{column tableau} of shape $\lambda$, denoted $\Lambda^c$, by consecutively filling the columns of $\lambda$. There is an explicit formula for the element $F_\Lambda$ in the group ring $\mathbb{C}S_n$ which is particularly simple when $\Lambda$ is the row, or column, tableau.

\begin{proposition}\label{dmeandyoung} (i) Diagonal matrix elements may be expressed as follows;
\begin{equation}\label{rowbigfexpression} F_{\Lambda^r} = P_{\Lambda^r}Q_{\Lambda^r}P_{\Lambda^r} / \lambda_1! \lambda_2! \dots,\end{equation}
\begin{equation}\label{columnbigfexpression}F_{\Lambda^c} = Q_{\Lambda^c}P_{\Lambda^c}Q_{\Lambda^c} / \lambda'_1! \lambda'_2! \dots, \end{equation}
where $\lambda'_i$ is the length of the $i^\textrm{\tiny th}$ column of $\lambda$. \\
(ii) There exist invertible elements $P, Q \in \mathbb{C}S_n$ such that \[P_{\Lambda^r}Q_{\Lambda^r}P_{\Lambda^r} = P_{\Lambda^r}Q_{\Lambda^r}P \qquad \textrm{and} \qquad Q_{\Lambda^c}P_{\Lambda^c}Q_{\Lambda^c} = Q_{\Lambda^c}P_{\Lambda^c}Q. \] \end{proposition}

For any $k = 1, \dots, n-1$ let $\sigma_k \in S_n$ be the transposition of $k$ and $k+1$. The following expression for the matrix element $F_\Lambda$ for an arbitrary $\Lambda$ can be obtained from either of (\ref{rowbigfexpression}) or (\ref{columnbigfexpression}) by using the formulas for the action of the generators $\sigma_1, \dots, \sigma_{n-1}$ of the group $S_n$ on the vectors of the Young basis.

Fix any standard tableau $\Lambda$. Consider the tableau $\sigma_k\Lambda$ obtained from $\Lambda$ by exchanging the numbers $k$ and $k+1$. Since the numbers $k$ and $k+1$ cannot lie on the same diagonal in $\Lambda$ we have $c_k(\Lambda) \neq c_{k+1}(\Lambda)$ always, and so we may set $d_k(\Lambda)=(c_{k+1}(\Lambda) - c_k(\Lambda))^{-1}$. If $\sigma_k\Lambda$ is non-standard then we must have $k$ and $k+1$ standing next to each other in the same row or column of $\Lambda$, and hence $d_k(\Lambda)=1$ or $d_k(\Lambda)=-1$.\\
Due to \cite{Y2} all the vectors of the Young basis may be normalised so that for any standard tableau $\Lambda$ and $k= 1, \dots, n-1$ \begin{equation}\label{sigmaonbasis} \begin{array}{lllll} \sigma_k \cdot v_\Lambda &=& d_k(\Lambda) v_\Lambda + \sqrt{1-d_k(\Lambda)^2}v_{\sigma_k\Lambda} &\quad& \textrm{for $\sigma_k \Lambda$ standard} \\
 &=& d_k(\Lambda) v_\Lambda  &\quad& \textrm{for $\sigma_k \Lambda$ non-standard}. \end{array} \end{equation}
This normalisation determines all the vectors of the Young basis up to a common multiplier $\omega \in \mathbb{C}$ with $|\omega| =1$. If the tableau $\sigma_k\Lambda$ is standard, then by (\ref{sigmaonbasis}) \begin{equation} \sqrt{1-d_k(\Lambda)^2}v_{\sigma_k\Lambda} = (\sigma_k - d_k(\Lambda))v_\Lambda,\end{equation}
and by definition (\ref{dme}) we have the relation \begin{equation}\label{sigmaondme} (1-d_k(\Lambda)^2)F_{\sigma_k\Lambda} = (\sigma_k - d_k(\Lambda))F_\Lambda(\sigma_k - d_k(\Lambda)). \end{equation}
For any standard tableau $\Lambda$ there is a sequence of transpositions $\sigma_{k_1}, \dots, \sigma_{k_b}$ such that $\Lambda = \sigma_{k_b} \cdots \sigma_{k_1}\Lambda^r$. Set $\Lambda^a = \sigma_{k_a} \cdots \sigma_{k_1}\Lambda^r$. Then $\Lambda^{\circ} = \Lambda^r$, $\Lambda^b = \Lambda$ and the tableaux $\Lambda^a$ are standard for all $a=0,\dots,b$. Using (\ref{sigmaondme}) repeatedly, one can derive from (\ref{rowbigfexpression}) the explicit formula for the element $F_\Lambda \in \mathbb{C}S_n$;
\begin{eqnarray*} \prod_{a=1}^b (1-d_{k_a}(\Lambda^{a-1})^2) F_\Lambda &=& (\sigma_{k_b} - d_{k_b}(\Lambda^{b-1})) \cdots (\sigma_{k_1} - d_{k_1}(\Lambda^r)) \\ &\phantom{=}& \times F_{\Lambda^r} \times (\sigma_{k_1} - d_{k_1}(\Lambda^r)) \cdots (\sigma_{k_b} - d_{k_b}(\Lambda^{b-1})).\end{eqnarray*}

So multiplying $F_\Lambda$ on the left and right by certain
invertible elements returns the Young symmetrizer. So one can see that the diagonal matrix element
$F_\Lambda$ may also be used to generate the irreducible module
$V_\Lambda$, \cite[section 2]{N5}.

\begin{example}
Let $N$ be the column tableau of shape $\nu = (2,1)$ as in Example 2.4, and $N' = \sigma_2 N$ the row tableau.

\[ N = \young(13,2) \qquad \qquad N' = \young(12,3) \]

Then $v_N, v_{N'}$ is the complete orthonormal Young basis of $V_\nu$. Using (\ref{sigmaonbasis}) we calculate the action of $S_3$ on $v_N$ to be;
\begin{eqnarray*}
\sigma_1 v_N &=& \phantom{\;}-v_N; \\
\sigma_2 v_N &=& \phantom{+}\frac{1}{2}v_N + \omega_{N'}; \\
\sigma_2\sigma_1 v_N &=& -\frac{1}{2}v_N - \omega v_{N'};\\
\sigma_1\sigma_2 v_N &=& -\frac{1}{2}v_N + \omega v_{N'};\\
\sigma_1\sigma_2\sigma_1 v_N &=& \phantom{+}\frac{1}{2}v_N - \omega v_{N'};
\end{eqnarray*}
where $\omega = \sqrt{3/4}$.

Hence the diagonal matrix element may be calculated as

\begin{eqnarray*}
F_N &=& \sum_{\sigma \in S_3} (v_N, \sigma v_N) \sigma \\ && \\
&=& (v_N,v_N)1 + (v_N, \sigma_1 v_N)(1\phantom{x}2) + (v_N, \sigma_1\sigma_2\sigma_1 v_N)(1 \phantom{x}3) + (v_N, \sigma_2 v_N)(2\phantom{x}3) \\ && \phantom{x} + (v_N, \sigma_1\sigma_2 v_N)(1\phantom{x}2\phantom{x}3) + (v_N, \sigma_2\sigma_1 v_N)(1\phantom{x}3\phantom{x}2) \\ &&\\
&=& 1 - (1\phantom{x}2) + \frac{(1\phantom{x}3)}{2} + \frac{(2\phantom{x}3)}{2} - \frac{(1\phantom{x}2\phantom{x}3)}{2} - \frac{(1\phantom{x}3\phantom{x}2)}{2}
\end{eqnarray*}

Proposition \ref{dmeandyoung} states the diagonal matrix element of the row tableau may be obtained from the Young symmetrizer after multiplication on the right by some invertible element, indeed in this case we have $F_{N'} = Y_{N'} (1 - \frac{(2\phantom{x}3)}{2})$. Furthermore, by (\ref{sigmaondme}), $F_N = \frac{4}{3} ((2\phantom{x}3) + \frac{1}{2}) F_{N'} ((2\phantom{x}3) + \frac{1}{2}) = (1 + \frac{(2\phantom{x}3)}{2}) Y_N$ and both elements generate modules isomorphic to $V_\nu$.
\end{example}

In 1986 Ivan Cherednik proposed another description of the diagonal matrix element, \cite{C1}. For any two distinct numbers $1 \leqslant p < q \leqslant n$, let
$(pq)$ be the transposition in the symmetric group $S_n$. Consider
the rational function of two complex variables $u$, $v$, $u \neq v$, with
values in the group ring $\mathbb{C}S_n$:
\begin{equation}\label{smallf} f_{pq}(u,v) = 1 -
\frac{(pq)}{u-v}.\end{equation}

Now introduce $n$ complex variables $z_1, \dots , z_n$. Consider
the set of pairs $(p,q)$ with $1 \leqslant p < q \leqslant n$.
Ordering the pairs lexicographically we form the product
\begin{equation}\label{bigf} F_\Lambda(z_1, \dots, z_n) = \prod_{(p,q)}^\rightarrow
f_{pq}(z_p + c_p(\Lambda), z_q + c_q(\Lambda)).\end{equation} Notice, if $p$ and $q$ sit on
the same diagonal in the tableau $\Lambda$, then $f_{pq}(z_p +
c_p(\Lambda), z_q + c_q(\Lambda))$ has a pole at $z_p = z_q$. However, using appropriate limiting procedure, such singularities prove to be removable.

\newpage
\begin{theorem}[Cherednik] The consecutive evaluations \[ \left. \left. \left. F_\Lambda(z_1, \dots, z_n) \right|_{z_1=0} \right|_{z_2=0} \cdots \right|_{z_n=0} \] of the rational function $F_\Lambda(z_1, \dots, z_n)$ are well defined. The corresponding value coincides with the matrix element $F_\Lambda$. \end{theorem}

However, Cherednik's paper does not contain complete proofs. A complete proof was given by Nazarov, \cite{N3}, while a relatively simple proof was recently found by Molev \cite{Mo}. Cheredik's limiting procedure, however, requires $n$ auxiliary parameters. Nazarov describes a different limiting procedure that reduces the number of such parameters. The following version of Cherednik's
description can be found in \cite{N4}. Let $\mathcal{R}_\Lambda$ be the vector subspace in $\mathbb{C}^n$
consisting of all tuples $(z_1, \dots , z_n)$ such that $z_p =
z_q$ whenever the numbers $p$ and $q$ appear in the same row of
the tableau $\Lambda$.

By direct calculation we can show the following identities are
true; \begin{equation}\label{triple}
f_{pq}(u,v)f_{pr}(u,w)f_{qr}(v,w) =
f_{qr}(v,w)f_{pr}(u,w)f_{pq}(u,v) \end{equation} for all pairwise
distinct indices $p$, $q$, $r$, and
\begin{equation}\label{commute} f_{pq}(u,v)f_{st}(z,w) =
f_{st}(z,w)f_{pq}(u,v) \end{equation} for all pairwise distinct
$p$, $q$, $s$, $t$.

Let's refer to (\ref{triple}) and (\ref{commute}) as the Yang-Baxter
relations. Using these relations we obtain reduced expressions for
$F_\Lambda(z_1, \dots, z_n)$ different from the right hand side of
(\ref{bigf}). This situation is analogous to the case of the
symmetric group where two reduced decompositions of the same
permutation are always connected by a sequence of Coxeter moves.
For more details on different expressions for $F_\Lambda(z_1,
\dots, z_n)$ see \cite{GP}.

Using (\ref{triple}) and (\ref{commute}) we may reorder the
product $F_\Lambda(z_1, \dots , z_n)$ such that each singularity
is contained in an expression known to be regular at $z_1 = z_2 =
\dots = z_n$, \cite{N4}. It is by this method that it was shown
that the restriction of the rational function $F_\Lambda(z_1,
\dots ,z_n)$ to the subspace $\mathcal{R}_\Lambda$ is regular at
$z_1 = z_2 = \cdots = z_n$. Furthermore, by showing divisibility on
the left and on the right by certain elements of $\mathbb{C}S_n$
it was shown that the value of $F_\Lambda(z_1, \dots ,z_n)$ at
$z_1 = z_1 = \cdots = z_n$ is the diagonal matrix element
$F_\Lambda$. Thus we have the following theorem,

\begin{theorem}[Nazarov] \label{MNtheorem} Restriction to $\mathcal{R}_\Lambda$ of the
rational function $F_\Lambda(z_1, \dots , z_n)$ is regular at $z_1
= z_2 = \cdots = z_n$. The value of this restriction at $z_1 = z_2
= \cdots = z_n$ coincides with the element $F_\Lambda \in
\mathbb{C}S_n$.
\end{theorem}

We now give a typical example for a diagonal matrix element obtained by such a limiting procedure, and obtain the corresponding irreducible representation by using it in the of the Young symmetrizer.

\begin{example}\label{hookexample} Again let $N=\young(13,2)$ \phantom{x} as in Example 2.6. Then on the subspace $\mathcal{R}_N$ we have $z_1 = z_3$. By setting $z = z_1 - z_2$ and applying Theorem \ref{MNtheorem} we have,

%\newpage
\begin{eqnarray*} F_N =  \left. F_N(z)\right|_{z = 0} &=&  \left.  \left(1 - \frac{(1 \phantom{x}2)}{z+1}\right)\left(1 + (1\phantom{x}3)\right)\left(1 + \frac{(2\phantom{x}3)}{z+2}\right)\right|_{z=0} \\ &\phantom{=}&\\ &=& \left( 1 - (1\phantom{x}2) \right)\left(1 + (1\phantom{x}3)\right)\left(1 + \frac{(2\phantom{x}3)}{2}\right) \\ &\phantom{=}& \\
&=& 1 - (1 \phantom{x}2) + \frac{(1\phantom{x}3)}{2} + \frac{(2\phantom{x}3)}{2} - \frac{(1\phantom{x}2\phantom{x}3)}{2} - \frac{(1\phantom{x}3\phantom{x}2)}{2} \end{eqnarray*}

Notice we may reverse the second line of the above product using (\ref{triple}). In which case we have an invertible element of $\mathbb{C}S_3$ followed by the Young symmetrizer $Y_N$, and $\mathbb{C}S_3 \cdot F_N$ is the 2-dimensional irreducible $\mathbb{C}S_3$-module.

The smallest examples containing a singularity are for standard tableau of shape $\gamma = (2,2)$. Let $\Gamma$ be the column tableau of shape $\gamma$, 
\[\Gamma = \young(13,24)\] 
Then on the subspace $\mathcal{R}_\Gamma$ we have $z_1=z_3$ and $z_2=z_4$, and by setting $z=z_1-z_2$ we have,

\begin{small}
\begin{eqnarray*} F_\Gamma &=&  \left.  \left(1 - \frac{(1 \phantom{x}2)}{z+1}\right)\left(1 + (1\phantom{x}3)\right)\left(1 - \frac{(1 \phantom{x}4)}{z}\right)\left(1 + \frac{(2\phantom{x}3)}{z+2}\right)\left(1 + (2 \phantom{x}4)\right) \left(1 - \frac{(3 \phantom{x}4)}{z+1}\right)\right|_{z=0} \end{eqnarray*}
\end{small}

Using (\ref{triple}) we may reorder the product as follows. This ordering quickly shows the singularity in the limit is indeed a removable singularity.

\begin{small}
\begin{eqnarray*} \phantom{F_\Gamma} &=& \left.  \left(1 - \frac{(1 \phantom{x}2)}{z+1}\right)\left[\left(1 + (1\phantom{x}3)\right)\left(1 - \frac{(1 \phantom{x}4)}{z}\right)\left(1 - \frac{(3 \phantom{x}4)}{z+1}\right)\right]\left(1 + (2 \phantom{x}4)\right)\left(1 + \frac{(2\phantom{x}3)}{z+2}\right) \right|_{z=0} \\ &\phantom{=}& \\
&=& \left.  \left(1 - \frac{(1 \phantom{x}2)}{z+1}\right)\left[ 1 + (1\phantom{x}3) - \frac{(1\phantom{x}4)}{z+1} - \frac{(3\phantom{x}4)}{z+1} - \frac{(1\phantom{x}3\phantom{x}4)}{z+1} - \frac{(1\phantom{x}4\phantom{x}3)}{z+1}\right]\left(1 + (2 \phantom{x}4)\right)\left(1 + \frac{(2\phantom{x}3)}{z+2}\right) \right|_{z=0} \\ &\phantom{=}& \\
&=& \left(1 - (1 \phantom{x}2)\right)\left[ 1 + (1\phantom{x}3) - (1\phantom{x}4) - (3\phantom{x}4) - (1\phantom{x}3\phantom{x}4) - (1\phantom{x}4\phantom{x}3)\right]\left(1 + (2 \phantom{x}4)\right)\left(1 + \frac{(2\phantom{x}3)}{2}\right) \\ &\phantom{=}& \\
&=& 1 - (1 \phantom{x}2) - (3 \phantom{x}4) + (1 \phantom{x}2)(3 \phantom{x}4) + (1 \phantom{x}3)(2 \phantom{x}4) + (1 \phantom{x}4)(2 \phantom{x}3) - (1\phantom{x}3\phantom{x}2 \phantom{x}4) - (1\phantom{x}4\phantom{x}2 \phantom{x}3) \phantom{x} + \frac{(1\phantom{x}3)}{2}\\ &\phantom{=}& \phantom{x} + \frac{(1\phantom{x}4)}{2} + \frac{(2\phantom{x}3)}{2} + \frac{(2\phantom{x}4)}{2} - \frac{(1\phantom{x}2\phantom{x}3)}{2} - \frac{(1\phantom{x}2\phantom{x}4)}{2} - \frac{(1\phantom{x}3\phantom{x}2)}{2} - \frac{(1\phantom{x}3 \phantom{x}4)}{2}  - \frac{(1\phantom{x}4\phantom{x}2)}{2} - \frac{(1\phantom{x}4\phantom{x}3)}{2} \\ &\phantom{=}& \phantom{x} - \frac{(2\phantom{x}3 \phantom{x}4)}{2} - \frac{(2\phantom{x}4\phantom{x}3)}{2} + \frac{(1\phantom{x}2\phantom{x}3 \phantom{x}4)}{2} + \frac{(1\phantom{x}2\phantom{x}4 \phantom{x}3)}{2} + \frac{(1\phantom{x}3\phantom{x}4 \phantom{x}2)}{2} + \frac{(1\phantom{x}4\phantom{x}3 \phantom{x}2)}{2}
\end{eqnarray*}
\end{small}
\end{example}

In particular, if $\lambda$ is a partition with only one part, that is to say $\lambda$ has a Young diagram consisting of a single row, then the fusion procedure uses only one parameter, $z$. However, this means all parameters in $F_\Lambda(z)$ cancel and the diagonal matrix element may be written without parameters as \[F_\Lambda = \prod_{(p,q)}^\rightarrow f_{pq}(c_p(\Lambda), c_q(\Lambda)). \] 

Similarly, we may form another expression for $F_\Lambda$ by
considering the subspace in $\mathbb{C}^n$ consisting of all
tuples $(z_1, \dots , z_n)$ such that $z_p = z_q$ whenever the
numbers $p$ and $q$ appear in the same column of the tableau
$\Lambda$ \cite{NT}.

\subsection{Constructing representations}

The fusion procedure described above required the decomposition of a Young diagram $\lambda$ into its rows, with a similar limiting procedure corresponding to a decomposition of $\lambda$ into its columns. We may use the same decompositions of $\lambda$ to construct certain symmetric functions.

Let $\Lambda(i,j)$ be the entry in box $(i,j)$ of a
Young tableau $\Lambda$ of shape $\lambda$. Then we can define a
product $x^\Lambda$ as follows,
\[ x^\Lambda = \prod_{(i,j) \in \lambda} x_{\Lambda(i,j)} .\]

Given a partition $\lambda$, the associated \emph{Schur function}
is \[ s_\lambda (x_1, \dots , x_m) = \sum_{\Lambda} x^\Lambda , \]
where the sum is over all semistandard tableaux $\Lambda$ of shape
$\lambda$ with entries $1, \dots , m$. The Schur functions are
symmetric functions, \cite{MD}.

The \emph{$n^{\mbox{th}}$ complete homogeneous symmetric
function}, $h_n(x_1, \dots , x_m)$, is the corresponding Schur
function of a tableau consisting of a single row of $n$ boxes.
That is to say,
\[ h_n(x_1, \dots , x_m) = s_{(n)}(x_1, \dots , x_m).\]
Similarly, the \emph{$n^{\mbox{th}}$ elementary symmetric
function}, $e_n(x_1, \dots , x_m)$, is the corresponding Schur
function of a tableau consisting of a single column of $n$ boxes,
\[ e_n(x_1, \dots , x_m) = s_{(1^n)}(x_1, \dots , x_m).\]

Define $\textbf{x} = (x_1, \dots , x_m)$. We can obtain the Schur
function of any shape from Young diagrams consisting of single
rows, or alternatively single columns, by means of the
\emph{Jacobi-Trudi identities} \cite[Chapter I3]{MD}. Let $\lambda
= (\lambda_1, \lambda_2, \dots , \lambda_k)$, then we have the
determinantal expressions \begin{equation}\label{jt1}
s_\lambda(\textbf{x}) = \textrm{det} [ h_{\lambda_i - i +
j}(\textbf{x}) ]_{i, j =1}^k
\end{equation} and
\begin{equation}\label{jt2} s_{\lambda}(\textbf{x}) = \textrm{det} [ e_{\lambda'_i - i +
j}(\textbf{x}) ]_{i,j=1}^l, \end{equation} where $\lambda'$ is the
conjugate of $\lambda$, $l$ is the number of columns of $\lambda$,
and any symmetric function with negative subscript is defined to
be zero.

We can think of the following identities as related to the
Jacobi-Trudi identities. Given $\lambda = (\lambda_1, \dots ,
\lambda_k)$ such that $n = \lambda_1 + \dots + \lambda_k$ then the
following identity holds,
\begin{equation}\label{dual1}  h_{\lambda_1}(\textbf{x})
h_{\lambda_2}(\textbf{x})\dots h_{\lambda_k}(\textbf{x}) =
 \sum_{\mu} K_{\mu \lambda}
s_\mu(\textbf{x}),\end{equation} where the sum is over all
partitions of $n$. The coefficients $K_{\mu \lambda}$ are
non-negative integers known as \emph{Kostka numbers}, \cite{MD}.
Importantly, we have $K_{\lambda \lambda} = 1.$

There exists an isomorphism from the graded ring of symmetric functions to the ring of equivalence classes of representations of all $\mathbb{C}S_n$, which sends the schur function of shape $\lambda$ to the equivalence class of shape $\lambda$. This translates the identity (\ref{dual1}) into the following identity of $S_n$-modules;
\begin{equation}\label{repdual1} \textrm{Ind}_{S_{\lambda_1} \times S_{\lambda_2} \times \cdots
\times S_{\lambda_k}}^{S_n } V_{(\lambda_1)} \otimes
V_{(\lambda_2)} \otimes \cdots \otimes V_{(\lambda_k)} \cong
\bigoplus_{\mu} (V_\mu)^{\oplus K_{\mu \lambda}},\end{equation} where the sum
is over all partitions of $n$. Note that $V_{(\lambda_i)}$ is the
trivial representation of $S_{\lambda_i}$.

Let $\mu$ be a
partition of $m$ and $M$ is a standard tableaux of shape $\mu$.
Denote by $\epsilon_x$ the embedding $\mathbb{C}S_{n} \to
\mathbb{C}S_{n+x}$ defined by
\begin{equation}\label{epsilon} \epsilon_x : (pq) \mapsto (p+x, q+x)
\end{equation} for all distinct $p, q = 1, \dots, n$. Then the
$S_{n+m}$-module induced from the $(S_n \times S_m)$-module
$V_\lambda \otimes V_\mu$ can be realised as the left ideal in
$\mathbb{C}S_{n+m}$ generated by the product $F_\Lambda \cdot
\epsilon_n(F_M)$.

Therefore on the subspace $\mathcal{R}_\Lambda$ we may write (\ref{bigf}) as
the product of diagonal matrix elements corresponding to the rows
of $\lambda$ and an element depending on parameters $z_1, z_2,
\dots, z_n$. Explicitly we have,
\[ F_\Lambda(z_1, \dots , z_n) = F_{(\lambda_1)} \cdot \epsilon_{\lambda_1}(F_{(\lambda_2)})
\cdot \epsilon_{\lambda_1 + \lambda_2}(F_{(\lambda_3)}) \cdots
\epsilon_{\lambda_1 + \dots + \lambda_{k-1}}(F_{(\lambda_k)}) \cdot
A(z_1, \dots , z_n).\]

\newpage
If $z_i - z_j \notin \mathbb{Z}$ when $i$ and $j$ are in different
rows of $\Lambda$ then $A(z_1, \dots , z_n)$ is invertible in
$\mathbb{C}S_n$. Hence, under this condition, (\ref{repdual1}) may be realised as the left ideal in
$\mathbb{C}S_n$ generated by $F_\Lambda(z_1, \dots, z_n)$.

The irreducible representation $V_\lambda$ appears in the
decomposition of this induced module with coefficient 1, and is
the ideal of $\mathbb{C}S_n$ generated by $F_\Lambda(z_1, \dots ,
z_n)$ when $z_1 = z_2 = \cdots = z_n$. The fusion procedure of
Theorem \ref{MNtheorem} provides a way of singling out this
irreducible component.

\begin{example}
Let $\Lambda$ be a standard tableau of shape $\lambda = (3,3,2)$. The first Jacobi-Trudi identity, (\ref{jt1}), gives us a way to write the Schur
function $s_\lambda(\textbf{x})$ as a determinant of complete homogeneous symmetric functions.
\[
s_{(3,3,2)}(\textbf{x}) = \left|
\begin{array}{ccc}
  h_3(\textbf{x}) & h_4(\textbf{x}) &h_5(\textbf{x}) \\
  h_2(\textbf{x}) & h_3(\textbf{x}) &h_4(\textbf{x}) \\
  h_0(\textbf{x}) & h_1(\textbf{x}) &h_2(\textbf{x}) \\
\end{array}
\right| =
 \left|
\begin{array}{ccc}
  s_{\textrm{\tiny \yng(3)}} & s_{\textrm{\tiny \yng(4)}} & s_{\textrm{\tiny \yng(5)}}\\
  s_{\textrm{\tiny \yng(2)}} & s_{\textrm{\tiny \yng(3)}} & s_{\textrm{\tiny \yng(4)}}\\
  1 & s_{\textrm{\tiny \yng(1)}} & s_{\textrm{\tiny \yng(2)}}\\
\end{array}
\right|,
\]
or dually,
\[ \begin{array}{rl} h_3(\textbf{x})h_3(\textbf{x})h_2(\textbf{x}) = & s_{(3,3,2)}(\textbf{x}) + s_{(4,4)}(\textbf{x}) + 2 s_{(4,3,1)}(\textbf{x})
+ s_{(4,2,2)}(\textbf{x}) + 3 s_{(5,3)}(\textbf{x})\\& +
2 s_{(5,2,1)}(\textbf{x}) + 3 s_{(6,2)}(\textbf{x}) +
2 s_{(7,1)}(\textbf{x}) + s_{(8)}(\textbf{x}).\\
\end{array}\]

In terms of modules we therefore obtain the following identity;
\[ \begin{array}{rl} \textrm{Ind}_{S_3 \times S_3 \times S_2}^{S_{8}} V_{(3)} \otimes
V_{(3)} \otimes V_{(2)} =&  V_{(3,3,2)} \oplus V_{(4,4)} \oplus 2 V_{(4,3,1)} \oplus  V_{(4,2,2)} \oplus 3 V_{(5,3)} \\& \oplus 2 V_{(5,2,1)} \oplus 3 V_{(6,2)} \oplus 2 V_{(7,1)} \oplus
V_{(8)}.
\end{array} \]

Here, the induced module may be realised as
the left ideal of $\mathbb{C}S_8$ generated by $F_\Lambda(z_1,
\dots, z_{8})$ which has three parameters on the subspace $\mathcal{R}_\Lambda$, where $z_i - z_j
\notin \mathbb{Z}$ when $i$ and $j$ are in different rows of $\Lambda$. The irreducible
representation corresponding to $\lambda$, $V_{(3,3,2)}$, appears in the decomposition above with coefficient 1, and is generated by $F_\Lambda(z_1,
\dots, z_{8})$ in the limit where the three parameters are equal.
\end{example}

\newpage
Similarly we have the equivalent identity for columns, \[
\textrm{Ind}_{S_{\lambda'_1} \times S_{\lambda'_2} \times \cdots
\times S_{\lambda'_l}}^{S_n } V_{(1^{\lambda'_1})} \otimes
V_{(1^{\lambda'_2})} \otimes \cdots \otimes V_{(1^{\lambda'_l})}
\cong \bigoplus_{\mu} (V_{\mu})^{\oplus K_{\mu' \lambda'}}.\] In
this case $V_{(1^{\lambda'_i})}$ is the alternating representation
of $S_{\lambda'_i}$. This induced module is isomorphic to the left
ideal of $\mathbb{C}S_n$ generated by $F_{\Lambda}(z_1, \dots,
z_n)$ considered on the subspace $\mathcal{R}_{\Lambda'}$, with
$z_i - z_j \notin \mathbb{Z}$ when $i, j$ are in different columns
of $\Lambda$. Here $\Lambda'$ denotes the transpose of the
standard tableau $\Lambda$ of shape $\lambda$. Again the
irreducible representation $V_{\lambda}$ appears in the
decomposition of this induced module with coefficient 1, and is
the ideal of $\mathbb{C}S_n$ generated by $F_{\Lambda}(z_1, \dots,
z_n)$ when $z_1 = z_2 = \cdots = z_n$.

%% file: TheHookFusionProcedure.tex
\section{The hook fusion procedure}\label{hookfusionprocedure}

In this chapter we present a new expression for the diagonal
matrix elements which minimises the number of auxiliary parameters
needed in the fusion procedure. We do this by considering diagonal hooks of
standard tableaux rather than their rows or columns. This minimises the number of auxiliary parameters since using fewer parameters would require at least one pair $(p,q)$ where $z_p =z_q$ and $p$ and $q$ lie on the same diagonal, leading to a pole in the function (\ref{bigf}) in the limit $z_1= \cdots = z_n$.  This is an expanded version of the author's first paper, \cite{Gr1}.

If $(i,j)$ is a box in the diagram of $\lambda$, then its $(i,j)$-\emph{hook} is the set of boxes in $\lambda$
\[  \{ (i, j') : j' \geqslant j \} \cup \{ (i', j)
: i' \geqslant i \}. \]
For example, if $\lambda = (3, 3, 2)$ then the (2,2)-hook is
given by the dotted cells in the following diagram;
\[ \begin{Young}
& &  \cr & $\bullet$ &  $\bullet$ \cr
 & $\bullet$ \cr
\end{Young} \]
We call the $(i,i)$-hook in $\lambda$ the \emph{$i^{\mbox{\tiny th}}$
principal hook}.

Let $\mathcal{H}_\Lambda$ be the vector subspace in $\mathbb{C}^n$
consisting of all tuples $(z_1, \dots , z_n)$ such that $z_p =
z_q$ whenever the numbers $p$ and $q$ appear in the same principal
hook of the tableau $\Lambda$. We will prove the following
theorem.

\begin{theorem}\label{fulltheorem} Restriction to $\mathcal{H}_\Lambda$ of the
rational function $F_\Lambda(z_1, \dots , z_n)$ is regular at $z_1
= z_2 = \cdots = z_n$. The value of this restriction at $z_1 = z_2
=  \dots = z_n$ coincides with the element $F_\Lambda \in
\mathbb{C}S_n$.
\end{theorem}

In particular, this hook fusion procedure can be used to form
irreducible representations of $\mathbb{C}S_n$ corresponding to Young
diagrams of hook shape using only one auxiliary parameter, $z$. If $\nu$ is a
partition of hook shape, this means we may calculate the diagonal matrix element without parameters as
\begin{equation}\label{bigfhook} F_N = \prod_{(p,q)}^\rightarrow f_{pq} (c_p(N), c_q(N)), \end{equation} with
the pairs $(p, q)$ in the product ordered lexicographically. For example, if $\nu = (2,1)$ as in Example 2.4, then the hook fusion procedure directly gives us the diagonal matrix element $F_N$.

%%%%%%%%%%%%%%%%%%%%%%%%%%55

There is another expression for Schur functions known as the
\emph{Giambelli identity} \cite{Gi}. This identity expresses any
Schur function $s_\lambda(\textbf{x})$ as a determinant of Schur
functions of hook shape. The following description along with a
combinatorial proof of the Giambelli identity can be found in
\cite{ER}.

Divide a Young diagram $\lambda$ into boxes with positive and
non-positive content. We may illustrate this on the Young diagram
by drawing `steps' above the main diagonal. Denote the boxes above
the steps by $\alpha(\lambda)$ and the rest by $\beta(\lambda)$.
For example, the following figure illustrates $\lambda$,
$\alpha(\lambda)$ and $\beta(\lambda)$ for $\lambda = (3,3,2)$.

\begin{figure}[h]
\label{steps}
\begin{center}
\begin{picture}(150, 60)
\put(0,30){\framebox(15,15)[r]{  }}
\put(15,30){\framebox(15,15)[r]{ }}
\put(30,30){\framebox(15,15)[r]{  }}

\put(0,15){\framebox(15,15)[r]{  }}
\put(15,15){\framebox(15,15)[r]{  }}
\put(30,15){\framebox(15,15)[r]{  }}

\put(0,0){\framebox(15,15)[r]{  }}
\put(15,0){\framebox(15,15)[r]{}}

\put(72.5,30){\framebox(15,15)[r]{  }}
\put(87.5,30){\framebox(15,15)[r]{ }}
\put(125,30){\framebox(15,15)[r]{ }}

\put(87.5,15){\framebox(15,15)[r]{  }}
\put(125,15){\framebox(15,15)[r]{  }}
\put(140,15){\framebox(15,15)[r]{  }}

\put(125,0){\framebox(15,15)[r]{ }}
\put(140,0){\framebox(15,15)[r]{}}

\put(20,50){$\lambda$} \put(77.5,50){$\alpha(\lambda)$}
\put(130,50){$\beta(\lambda)$}

\linethickness{2pt} \put(0,45){\line(1,0){15}}
\put(15,30){\line(1,0){15}} \put(30,15){\line(1,0){15}}

\put(15,30){\line(0,1){15}} \put(30,15){\line(0,1){15}}

\end{picture}
\end{center}
 \caption{The Young diagram $(3,3,2)$ divided into boxes with positive content and non-positive content}
 \end{figure}
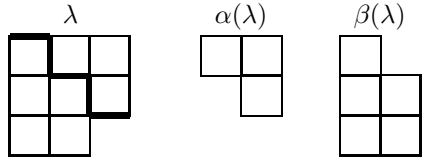

If we denote the rows of $\alpha(\lambda)$ by $\alpha_1 > \alpha_2
> \dots > \alpha_d > 0$ and the columns of $\beta(\lambda)$
by $\beta_1 > \beta_2 > \dots > \beta_d > 0$, then we have the
following alternative notation for $\lambda$;
\[ \lambda = ( \alpha | \beta ), \]
where $\alpha = (\alpha_1, \dots , \alpha_d)$ and $\beta =
(\beta_1, \dots , \beta_d)$. \footnote{A related notation is known as the \emph{Frobenius notation}, where $\beta' = (\beta_1 - 1 , \dots , \beta_d - 1)$ and $\lambda$ is written as $(\alpha | \beta')$.}

Here $d$ denotes the length of the side of the \emph{Durfee
square} of the shape $\lambda$, which is the set of boxes
corresponding to the largest square that fits inside $\lambda$,
and is equal to the number of principal hooks in $\lambda$. In our
example $d=2$ and $\lambda = (2, 1 | 3, 2)$. 

Using this notation the Giambelli identity may be written as
follows;
\begin{equation}\label{giambelli}  s_{(\alpha|\beta)}(\textbf{x}) = \det [ s_{(\alpha_i | \beta_j)}
(\textbf{x})]_{i,j=1}^d, \end{equation} where the determinant is $d \times
d$.

Dually, for a partition $\lambda = (\alpha_1, \alpha_2, \dots,
\alpha_d | \beta_1, \beta_2, \dots ,\beta_d)$ of $n$, we have the
identity \begin{equation}\label{dualgiambelli} s_{(\alpha_1 |
\beta_1)}(\textbf{x})s_{(\alpha_2 | \beta_2)}(\textbf{x}) \dots
s_{(\alpha_d | \beta_d)}(\textbf{x}) = \sum_{\mu}
D_{\mu \lambda} s_\mu(\textbf{x}),
\end{equation} where the sum is over all partitions of $n$.
The coefficients, $D_{\mu \lambda}$, are non-negative integers,
and in particular $D_{\lambda \lambda} =1$, see Proposition \ref{Dlambdalambda}. This translates to a
decomposition of the induced representation of the tensor product
of modules of hook shape. Further these hooks are the principal
hooks of $\lambda$.
\[ \textrm{Ind}_{S_{h_1} \times S_{h_2} \times \cdots
\times S_{h_d}}^{S_n} V_{(\alpha_1 | \beta_1)} \otimes
V_{(\alpha_2 | \beta_2)} \otimes \cdots \otimes V_{(\alpha_d |
\beta_d)} \cong \bigoplus_{\mu} (V_{\mu})^{\oplus D_{\mu
\lambda}},\] where $h_i$ is the length of the $i^{\mbox{\scriptsize th}}$
principal hook.

Now, on the subspace $\mathcal{H}_\Lambda$, we may write
(\ref{bigf}) as the product of diagonal matrix elements
corresponding to the principal hooks of $\lambda$ and an element
depending on parameters $z_1, z_2, \dots, z_n$. Explicitly we
have,
\[  F_\Lambda(z_1, \dots ,
z_n) =  F_{(\alpha_1 | \beta_1)} \cdot \epsilon_{h_{1}}(F_{(\alpha_2
| \beta_2)}) \cdot \epsilon_{h_{1} + h_{2}}(F_{(\alpha_3 |
\beta_3)}) \cdot \cdots \phantom{xxxxxxxxxxxxxx} \]\[ \phantom{xxxxxxxxxxxxxxxxxx}
\cdot \cdots \cdot \epsilon_{h_{1} + \dots +
h_{d-1}}(F_{(\alpha_d | \beta_d)}) \cdot B(z_1, \dots, z_n), \] where $\epsilon$ is the embedding described in (\ref{epsilon}).

If $z_i - z_j \notin \mathbb{Z}$ when $i$ and $j$ are in different
principal hooks of $\Lambda$ then $B(z_1, \dots , z_n)$ is
invertible in $\mathbb{C}S_n$. Therefore, under this condition,
the above induced module may be realised as the left ideal in
$\mathbb{C}S_n$ generated by $F_\Lambda(z_1, \dots, z_n)$.\\
The irreducible representation $V_\lambda$ appears in the
decomposition of this induced module with coefficient 1, and is
the ideal of $\mathbb{C}S_n$ generated by $F_\Lambda(z_1, \dots ,
z_n)$ when $z_1 = z_2 = \cdots = z_n$.

Hence, in this way,  our hook fusion procedure relates to the
Giambelli identity in the same way that Cherednik's original
fusion procedure relates to the Jacobi-Trudi identity. Namely, it
provides a way of singling out the irreducible component
$V_\lambda$ from the above induced module.

%%%%%%%%%%%%%%%%%%%%%%%%%555

\subsection{Proof of the hook fusion procedure}

%The diagonal matrix element $F_\Lambda$ determines the irreducible
%module $V_\lambda$ of $S_n$, up to isomorphism, for any tableau
%$\Lambda$ of shape $\lambda$. Therefore in the sequel we will only
%use one particular tableau, the \emph{hook tableau}. In which case
%we may denote the diagonal matrix element $F_\Lambda$ by
%$F_\lambda$, and the space $\mathcal{H}_\Lambda$ by
%$\mathcal{H}_\lambda$.

We now begin our proof of Theorem \ref{fulltheorem}, the hook fusion procedure. 
We fill a diagram $\lambda$ by hooks to form a tableau $\Lambda^{\circ}$
in the following way: For the first principal hook we fill the
column with entries $1$, $2$, \dots , $r$ and then fill the row
with entries $r+1$, $r+2$, \dots , $s$. We then fill the column of
the second principal hook with $s+1$, $s+2$, \dots , $t$ and fill
the row with $t+1$, $t+2$, \dots , $x$. Continuing in this way we
form the \emph{hook tableau}.

\begin{example} On the
left is the hook tableau, $\Lambda^\circ$, of the diagram $\lambda = (3,3,2)$, and
on the right the same diagram with the content of each box.

\begin{normalsize}

\begin{center}
\begin{picture}(50,50)
\put(0,30){\framebox(15,15)[r]{ 1 }}
\put(15,30){\framebox(15,15)[r]{ 4 }}
\put(30,30){\framebox(15,15)[r]{ 5 }}
\put(0,15){\framebox(15,15)[r]{ 2 }}
\put(15,15){\framebox(15,15)[r]{ 6 }}
\put(30,15){\framebox(15,15)[r]{ 8 }}
\put(0,0){\framebox(15,15)[r]{ 3 }}
\put(15,0){\framebox(15,15)[r]{ 7 }}
\end{picture}
\qquad \qquad \qquad
\begin{picture}(50,50)
\put(0,30){\framebox(15,15)[r]{ 0 }}
\put(15,30){\framebox(15,15)[r]{ 1 }}
\put(30,30){\framebox(15,15)[r]{ 2 }}
\put(0,15){\framebox(15,15)[r]{ -1 }}
\put(15,15){\framebox(15,15)[r]{ 0 }}
\put(30,15){\framebox(15,15)[r]{ 1 }}
\put(0,0){\framebox(15,15)[r]{ -2 }}
\put(15,0){\framebox(15,15)[r]{ -1 }}
\end{picture}
\end{center}
\end{normalsize}
Therefore the sequence $(c_1(\Lambda^{\circ}), c_2(\Lambda^{\circ}), \dots , c_8(\Lambda^{\circ}))$ is given by \\ $(0,
-1, -2, 1, 2, 0 , -1, 1)$. \end{example}

Consider (\ref{bigf}) as a rational function of the variables
$z_1, \dots , z_n$ with values in $\mathbb{C}S_n$. The factor
$f_{pq}(z_p + c_p(\Lambda), z_q + c_q(\Lambda))$ has a pole at $z_p = z_q$ if and
only if the numbers $p$ and $q$ stand on the same diagonal of the
tableau $\Lambda$. We then call the pair $(p, q)$ a
\emph{singularity}. And we call the corresponding term $f_{pq}(z_p
+ c_p(\Lambda), z_q + c_q(\Lambda))$ a \emph{singularity term}, or
simply a singularity.

Let $p$ and $q$ be in the same hook of $\Lambda$. If $p$ and $q$
are next to one another in the column of the hook then, on
$\mathcal{H}_\Lambda$, $f_{pq}(z_p + c_p(\Lambda), z_q + c_q(\Lambda)) = 1 - (pq)$.
Similarly, if $p$ and $q$ are next to one another in
the same row of the hook then $f_{pq}(z_p + c_p(\Lambda), z_q + c_q(\Lambda)) = 1 +
(pq)$. In both cases
$\frac{1}{2} f_{pq}(z_p + c_p(\Lambda), z_q + c_q(\Lambda))$ is an idempotent.

Also, for distinct $p$, $q$, we have the identity
\begin{equation}\label{inverses} f_{pq}(u,v)f_{qp}(v,u) = 1 - \frac{1}{(u-v)^2}.
\end{equation}
Therefore, if the contents $c_p$ and $c_q$ differ by a number
greater than one, then the factor $f_{pq}(c_p, c_q)$
is invertible in $\mathbb{C}S_n$ with inverse
$\frac{(c_p - c_q)^2}{(c_p - c_q)^2 - 1} f_{qp}(c_q, c_p)$.

The presence of singularity terms in the product $F_\Lambda(z_1,
\dots , z_n)$ means this product may or may not be regular on the
vector subspace of $\mathcal{H}_\Lambda$ in the limit $z_1 = z_2 = \cdots = z_n$. Using
the following lemma, we will be able to show that $F_\Lambda(z_1,
\dots , z_n)$ is in fact regular on this subspace, with all singularities being removable.

\begin{lemma}\label{regular} The restriction of $f_{pq}(u, v)f_{pr}(u, w)f_{qr}(v,
w)$ to the set $(u,v,w)$ such that $u = v \pm 1$ is regular at $u
= w$. \end{lemma}
\begin{proof} Under the condition $u = v \pm 1$, the product $f_{pq}(u, v)f_{pr}(u, w)f_{qr}(v,
w)$ can be written as \[ (1 \mp (pq)) \cdot \left(1 - \frac{(pr) +
(qr)}{v-w}\right). \] And this rational function of $v$, $w$ is
clearly regular at $w = v \pm 1$. \end{proof}

Notice that the three term product, or \emph{triple}, in the
statement of the lemma can be written in reverse order due to
(\ref{triple}). In particular, if the middle term is a singularity
and the other two terms are an appropriate idempotent and
\emph{triple term}, then the triple is regular at $z_1 = z_2 =
\dots = z_n$. We may now prove the first statement of Theorem
\ref{fulltheorem}.

\begin{proposition}\label{jimtheorem1} The restriction of the rational function
$F_\Lambda (z_1, \dots , z_n)$ to the subspace
$\mathcal{H}_\Lambda$ is regular at $z_1 = z_2 = \cdots = z_n$.
\end{proposition}

\newpage 
\begin{proof}
Consider any standard tableau $\Lambda'$ obtained from the tableau $\Lambda$
by an adjacent transposition of its entries, say by $\sigma_k\in S_n$.
Using the relations (\ref{triple}), (\ref{commute}) and (\ref{inverses}), we derive
the equality of rational functions in the variables $z_1, \ldots , z_n$
\[
F_\Lambda(z_1, \ldots , z_n) \cdot f_{k+1, k}\left(z_{k+1} + c_{k+1}(\Lambda), z_k + c_k(\Lambda)
\right) \sigma_k
=
\]
\begin{equation}\label{jimtheorem1equation}
\sigma_k f_{k, k+1}\left( z_k + c_k(\Lambda), z_{k+1} + c_{k+1}(\Lambda)
\right) \cdot F_{\Lambda'}(z'_1, \ldots ,z'_n),
\end{equation}
where the sequence of variables $(z'_1, \ldots ,z'_n)$\ is obtained from
the sequence $(z_1, \ldots , z_n)$ by exchanging the terms $z_k$ and 
$z_{k+1}$. Observe that
\[
(z'_1, \ldots ,z'_n)\in\mathcal{H}_{\Lambda'}
\quad\Leftrightarrow\quad
(z_1, \ldots , z_n)\in\mathcal{H}_\Lambda.
\]
Also observe that here $| c_k(\Lambda)-c_{k+1}(\Lambda)|\geqslant2$
because the tableaux $\Lambda$ and $\Lambda'$ are standard.
Therefore the functions 
\[
f_{k,k+n}\left( z_k + c_k(\Lambda), z_{k+1} + c_{k+1}(\Lambda)
\right)
\ \quad\textrm{and}\ \quad
f_{k+1,k}\left( z_{k+1} + c_{k+1}(\Lambda), z_k + c_k(\Lambda)
\right)
\]
appearing in the equality (\ref{jimtheorem1equation}),
are regular at $z_k=z_{k+1}$.
Moreover, their values at $z_k=z_{k+1}$ are invertible
in $\mathbb{C}S_n$. 
Due to these two observations, the equality (\ref{jimtheorem1equation})
shows that Proposition \ref{jimtheorem1} is equivalent to its counterpart for
the tableau $\Lambda'$ instead of $\Lambda$.

Let us take the hook tableau $\Lambda^\circ$ of shape $\lambda$.  
There is a chain $\Lambda,\Lambda', \ldots ,\Lambda^\circ$ of standard tableaux
of the same shape $\lambda$, such that each subsequent tableau in the 
chain is
obtained from the previous one by an adjacent transposition of the 
entries.
Due to the above argument, it now suffices to prove Proposition \ref{jimtheorem1} 
only in the case $\Lambda=\Lambda^\circ$.

We will prove the statement by reordering the factors of the
product \\ $F_{\Lambda^{\circ}} (z_1, \dots , z_n)$, using relations
(\ref{triple}) and (\ref{commute}), in such a way that each
singularity is part of a triple which is regular at $z_1 = z_2 =
\dots = z_n$, and hence the whole of $F_{\Lambda^{\circ}} (z_1, \dots ,
z_n)$ will be manifestly regular.

Define $g_{pq}$ to be the following; \[ g_{pq} = \left\{
\begin{array}{ccc}
  f_{pq} (z_p + c_p(\Lambda^{\circ}), z_q + c_q(\Lambda^{\circ})) & \textrm{if} & p<q \\
  1 & \textrm{if} & p \geqslant q
\end{array} \right. \]

Now, let us divide the diagram $\lambda$ into two parts,
consisting of those boxes with positive contents and those with
non-positive contents as in Figure \ref{steps}. Consider the
entries of the $i^{\mbox{\tiny th}}$ column of the hook tableau
$\Lambda$ of shape $\lambda$ that lie below the steps. If $u_1,
u_2, \dots , u_k$ are the entries of the $i^{\mbox{\tiny th}}$ column
below the steps, we define
\begin{equation}\label{cproduct} C_i = \prod_{q=1}^n g_{u_1 , q}
g_{u_2 , q} \dots g_{u_k , q}.
\end{equation} Now consider the entries of the $i^{\mbox{\tiny th}}$ row of $\Lambda$ that lie
above the steps. If $v_1, v_2, \dots , v_l$ are the entries of the
$i^{\mbox{\tiny th}}$ row above the steps, we define
\begin{equation}\label{rproduct} R_i = \prod_{q=1}^n g_{v_1 , q}
g_{v_2 , q} \dots g_{v_l , q}.
\end{equation}

Our choice of the hook tableau was such that the following is
true; if $d$ is the number of principal hooks of $\lambda$ then by
relations (\ref{triple}) and (\ref{commute}) we may reorder the
factors of $F_{\Lambda^{\circ}} (z_1, \dots , z_n)$ such that \[ F_{\Lambda^{\circ}}
(z_1, \dots , z_n) = \prod_{i=1}^d C_i R_i . \]

Now each singularity $(p,q)$ has its corresponding term $f_{pq}
(z_p + c_p(\Lambda^{\circ}), z_q + c_q(\Lambda^{\circ}))$ contained in some product $C_i$ or $R_i$.
This singularity term will be on the immediate left of the
 term $f_{p+1,q} (z_{p+1} +
c_{p+1}(\Lambda^{\circ}), z_q + c_q(\Lambda^{\circ}))$. Also, this ordering has been chosen such
that the product of factors to the left of any such singularity in
$C_i$ or $R_i$ is divisible on the right by $f_{p, p+1} (z_p +
c_p(\Lambda^{\circ}), z_{p+1} + c_{p+1}(\Lambda^{\circ}))$.

Therefore we can replace the pair \begin{small} \[ f_{pq}(z_p + c_p(\Lambda^{\circ}), z_q + c_q(\Lambda^{\circ}))
f_{p+1,q}(z_p + c_p(\Lambda^{\circ}), z_q + c_q(\Lambda^{\circ})) \] \end{small} in the product by the triple
\begin{small} \[ \frac{1}{2} f_{p, p+1} (z_p + c_p(\Lambda^{\circ}), z_{p+1} + c_{p+1}(\Lambda^{\circ})) f_{pq}(z_p
+c_p(\Lambda^{\circ}), z_q + c_q(\Lambda^{\circ})) f_{p+1, q} (z_{p+1} +c_{p+1}(\Lambda^{\circ}), z_q + c_q(\Lambda^{\circ})),\]\end{small}where
$\frac{1}{2} f_{p, p+1} (z_p + c_p(\Lambda^{\circ}), z_{p+1} + c_{p+1}(\Lambda^{\circ}))$ is an
idempotent. Divisibility on the right by $f_{p, p+1}(z_p + c_p(\Lambda^{\circ}),
z_{p+1} + c_{p+1}(\Lambda^{\circ}))$ means the addition of the idempotent has no
effect on the value of the product $C_i$ or
$R_i$.\\
By Lemma \ref{regular}, the above triples are regular at $z_1 =
z_2 = \cdots = z_n$, and therefore, so too are the products $C_i$
and $R_i$, for all $1 \leqslant i \leqslant d$. Moreover, this
means $F_{\Lambda^{\circ}} (z_1, \dots , z_n)$ is regular at $z_1 = z_2 =
\dots = z_n$.
\end{proof}

{\addtocounter{definition}{1} \bf Example \thedefinition .} As an
example consider the hook tableau of the Young diagram $\lambda =
(3,3,2)$.

In the original lexicographic ordering the product $F_{\Lambda^{\circ}}(z_1,
\dots , z_n)$ is written as follows, for simplicity we will write
$f_{pq}$ in place of the term $f_{pq}(z_p + c_p(\Lambda^{\circ}), z_q + c_q(\Lambda^{\circ}))$.
\[
\begin{array}{rl}
   F_{\Lambda^{\circ}} (z_1, \dots , z_n) = & f_{12}f_{13}f_{14}f_{15}\textbf{$f_{16}$}f_{17}f_{18}f_{23}f_{24}f_{25}f_{26}\textbf{$f_{27}$}f_{28}f_{34}f_{35}f_{36} f_{37}f_{38}\\&
   f_{45}f_{46}f_{47}\textbf{$f_{48}$}f_{56}f_{57}f_{58}f_{67}f_{68}f_{78}\\
\end{array}
\]
We may now reorder this product into the form below using
relations (\ref{triple}) and (\ref{commute}) as described in the
above proposition. The terms bracketed are the singularity terms
with their appropriate triple terms.
\[
\begin{array}{rl}
   F_{\Lambda^{\circ}} (z_1, \dots , z_n) = &
   f_{12}f_{13}f_{23}f_{14}f_{24}f_{34}f_{15}f_{25}f_{35}(f_{16}f_{26})f_{36}f_{17}(f_{27}f_{37})f_{18}f_{28}f_{38}\\&
   \cdot f_{45}f_{46}f_{56}f_{47}f_{57}(f_{48}f_{58}) \cdot f_{67}f_{68}f_{78}\\
\end{array}
\]
And so, for each singularity $f_{pq}$, we can replace
$f_{pq}f_{p+1,q}$ in the product by the triple $\frac{1}{2} f_{p,
p+1}f_{pq}f_{p+1, q}$, where $\frac{1}{2} f_{p, p+1}$ is an
idempotent, without changing the value of $F_{\Lambda^{\circ}} (z_1, \dots ,
z_n)$.
\[ \begin{array}{rl}
   F_{\Lambda^{\circ}} (z_1, \dots , z_n) = & f_{12}f_{13}f_{23}f_{14}f_{24}f_{34}f_{15}f_{25}f_{35}(\frac{1}{2}f_{12}f_{16}f_{26})f_{36}f_{17}(\frac{1}{2}f_{23}f_{27}f_{37})f_{18}\\&f_{28}f_{38}
   \cdot f_{45}f_{46}f_{56}f_{47}f_{57}(\frac{1}{2}f_{45}f_{48}f_{58}) \cdot f_{67}f_{68}f_{78}\\
\end{array}\]
And since each of these triples are regular at $z_1 = z_2 = \dots
= z_n$ then so too is the whole of $F_{\Lambda^{\circ}} (z_1, \dots , z_n)$.
{\nolinebreak \hfill \rule{2mm}{2mm}

%\quad

Therefore, due to the above proposition an element $F_\Lambda \in
\mathbb{C}S_n$ can now be defined as the value of $F_\Lambda (z_1,
\dots , z_n)$ at $z_1 = z_2 = \cdots = z_n$. We proceed by showing
this $F_\Lambda$ is indeed the diagonal matrix element. To this
end we will need the following propositions.

Note that for $n=1$, we have $F_\Lambda =1$. For any $n \geqslant
1$, we have the following fact.

\begin{proposition}\label{coeff1} The coefficient of $F_\Lambda \in
\mathbb{C}S_n$ at the unit element of $S_n$ is 1.
\end{proposition}
\begin{proof} For each $r = 1, \dots , n-1$ let $\sigma_r \in S_n$ be
the adjacent transposition $(r, r+1)$. Let $\sigma_0 \in S_n$ be the
element of maximal length. Multiply the ordered product
(\ref{bigf}) by the element $\sigma_0$ on the right. Using the reduced
decomposition
\begin{equation}\label{longestelement} \sigma_0 = \prod_{(p,q)}^\rightarrow \sigma_{q-p}
\end{equation} we get the product \[ \prod_{(p,q)}^\rightarrow
\left(\sigma_{q-p} - \frac{1}{z_p - z_q + c_p(\Lambda) - c_q(\Lambda)}\right) . \] For
the derivation of this formula see \cite[(2.4)]{GP}. This formula
expands as a sum of the elements $\sigma \in S_n$ with coefficients
from the field of rational functions of $z_1, \dots , z_n$ valued
in $\mathbb{C}$. Since the decomposition (\ref{longestelement}) is
reduced, the coefficient at $\sigma_0$ is 1. By the definition of
$F_\Lambda$, this implies that the coefficient of $F_\Lambda \sigma_0
\in \mathbb{C}S_n$ at $\sigma_0$ is also 1.
\end{proof}

In particular this shows that $F_\Lambda \neq 0$ for any nonempty
diagram $\lambda$. Let us now denote by $\varphi$ the involutive
antiautomorphism of the group ring $\mathbb{C}S_n$ defined by
$\varphi (\sigma) = \sigma^{-1}$ for every $\sigma \in S_n$.

\begin{proposition}\label{varphi} The element $F_\Lambda \in \mathbb{C}S_n$ is
$\varphi$-invariant. \end{proposition}
\begin{proof} Any element of the group ring $\mathbb{C}S_n$ of
the form $f_{pq}(u,v)$ is $\varphi$-invariant. Applying the
antiautomorphism $\varphi$ to an element of the form (\ref{bigf})
just reverses the ordering of the factors corresponding to the
pairs (p,q). However, the initial ordering can then be restored
using relations (\ref{triple}) and (\ref{commute}), for more
detail see \cite{GP}. Therefore, any value of the function
$F_\Lambda (z_1, \dots , z_n)$ is $\varphi$-invariant. Therefore,
so too is the element $F_\Lambda$.
\end{proof}

\begin{proposition}\label{stripping} Let $x$ be last entry in the $k^{\mbox{\scriptsize th}}$ row
of the hook tableau of shape $\lambda$. If $\lambda= (\alpha_1,
\alpha_2, \dots, \alpha_d | \beta_1, \beta_2, \dots , \beta_d)$
\\and $\mu = (\alpha_{k+1}, \alpha_{k+2}, \dots , \alpha_d |
\beta_{k+1}, \beta_{k+2}, \dots , \beta_d)$, then $F_{\Lambda^{\circ}} = P
\cdot \epsilon_x(F_{M^{\circ}})$, for some element $P \in \mathbb{C}S_n$,
where $\epsilon_x$ is the embedding in (\ref{epsilon}).
\end{proposition}
\begin{proof}
Here the shape $\mu$ is obtained by removing the first $k$
principal hooks of $\lambda$. By the ordering described in
Proposition \ref{jimtheorem1},
\[ F_{\Lambda^{\circ}}(z_1, \dots ,
z_n) = \prod_{i=1}^k C_iR_i \cdot \epsilon_x(F_{M^{\circ}}(z_{x+1}, \dots,
z_{n})),\]
where $C_i$ and $R_i$ are defined by (\ref{cproduct}) and (\ref{rproduct}).\\
Since all products $C_i$ and $R_i$ are regular at $z_1 = z_2 =
\dots = z_n$, Proposition \ref{jimtheorem1} then gives us the
required statement. \end{proof}

In any given ordering of $F_\Lambda(z_1, \dots , z_n)$, we want a
singularity term to be placed next to an appropriate triple term
such that we may then form a regular triple. In that case we will
say these two terms are `tied'.\\
If $u$,
$v$ stand next to each other in the same row, or same column, of
$\Lambda^{\circ}$ the following results show that $F_{\Lambda^{\circ}}$ is divisible
on the left (and on the right) by $1 - (uv)$, or $1+(uv)$
respectively, or divisible by these terms preceded (followed) by
some invertible elements of $\mathbb{C}S_n$. \\%Hence $F_\lambda$ is
%the diagonal matrix element. \\
However, proving the divisibilities described requires some pairs
to be `untied', in which case we must form a new ordering. This is
the content of the following proofs. Some explicit examples will
then given in Example 3.11 below.

\begin{proposition}\label{jimtheorem2} Suppose the numbers $u < v$ stand next to each
other in the same column of the hook tableau $\Lambda^{\circ}$ of shape
$\lambda$. First, let $s$ be the last entry in the row containing
$u$. If $c_v(\Lambda^{\circ}) < 0$ then the element $F_{\Lambda^{\circ}} \in \mathbb{C}S_n$
is divisible on the left and on the right by $f_{u,v}(c_u(\Lambda^{\circ}), c_v(\Lambda^{\circ})) =
1 - (uv)$. If $c_v(\Lambda^{\circ}) \geqslant 0$ then the element $F_{\Lambda^{\circ}} \in
\mathbb{C}S_n$ is divisible on the left by the product
\[ \prod_{i = u, \dots, s}^\leftarrow \left( \prod_{j= s+1, \dots,
v}^\rightarrow f_{ij}(c_i(\Lambda^{\circ}), c_j(\Lambda^{\circ})) \right) \]
\end{proposition}

\begin{proof}
By Proposition \ref{varphi}, the divisibility of $F_{\Lambda^{\circ}}$ by
the element $1- (uv)$ on the left is equivalent to the
divisibility by the same element on the right. Let us prove
divisibility on the left.

By Proposition \ref{stripping}, if $\epsilon_x(F_{M^{\circ}})$ is divisible
on the right by $f_{uv}(c_u(\Lambda^{\circ}), c_v(\Lambda^{\circ}))$, or $f_{uv}(c_u(\Lambda^{\circ}), c_v(\Lambda^{\circ}))$ followed
by some invertible terms, then so too will $F_{\Lambda^{\circ}}$. If
$\epsilon_x(F_{M^{\circ}})$ is divisible on the left by $f_{uv}(c_u(\Lambda^{\circ}), c_v(\Lambda^{\circ}))$,
or $f_{uv}(c_u(\Lambda^{\circ}), c_v(\Lambda^{\circ}))$ preceded by some invertible terms, then by
using Proposition \ref{varphi} twice, first on $\epsilon_x(F_{M^{\circ}})$
and then on the product $F_{\Lambda^{\circ}} = P \cdot \epsilon_x(F_{M^{\circ}})$, so too
will $F_{\Lambda^{\circ}}$. Hence we only need to prove the statement for
$(u,v)$ such that $u$ is in the first row or first column of
$\Lambda^{\circ}$.

Let $r$ be the last entry in the first column of $\Lambda^{\circ}$, $s$
the last entry in the first row of $\Lambda^{\circ}$, and $t$ the last
entry in the second column of $\Lambda^{\circ}$, as shown in Figure
\ref{jimtheorem2fig}.

\begin{figure}[h] \label{jimtheorem2fig}

\begin{normalsize}
\begin{center}
\begin{picture}(275,200)
\begin{small}
\put(0,175){\framebox(25,25)[c]{ 1 }}
\put(25,175){\framebox(25,25)[c]{ $r + 1$ }}
\put(50,175){\framebox(25,25)[c]{ $r+2$ }}
\put(0,150){\framebox(25,25)[c]{ 2 }}
\put(25,150){\framebox(25,25)[c]{ $s+1$ }}
\put(50,150){\framebox(25,25)[c]{ $t+1$ }}

\put(0,25){\line(0,1){125}} \put(25,25){\line(0,1){125}}
\put(0,0){\framebox(25,25)[c]{ $r$ }}

\put(75,200){\line(1,0){50}} \put(75,175){\line(1,0){50}}
\put(125,175){\framebox(25,25)[c]{ $u$ }}
\put(150,200){\line(1,0){100}} \put(150,175){\line(1,0){100}}
\put(250,175){\framebox(25,25)[c]{ $s$ }}

\put(50,100){\line(0,1){50}} \put(25,75){\framebox(25,25)[c]{ $t$
}}

\put(75,150){\line(1,0){50}} \put(125,150){\framebox(25,25)[c]{
$v$ }} \put(150,150){\line(1,0){50}} \put(200,150){\line(0,1){25}}

\linethickness{1.5pt} \put(0,200){\line(1,0){25}}
\put(25,175){\line(1,0){25}} \put(50,150){\line(1,0){25}}
\put(25,175){\line(0,1){25}} \put(50,150){\line(0,1){25}}

\put(125,150){\line(0,1){50}}\put(150,150){\line(0,1){50}}
\put(125,200){\line(1,0){25}}\put(125,150){\line(1,0){25}}

\put(95,185){$\dots$} \put(170,185){$\dots$}\put(220,185){$\dots$}
\put(95,160){$\dots$} \put(170,160){$\dots$}

\put(10,120){$\vdots$} \put(35,
120){$\vdots$}\put(10,80){$\vdots$} \put(10,40){$\vdots$}

\end{small}
\end{picture}
\end{center}
\end{normalsize}
 \caption{The first two principal hooks of the hook tableau $\Lambda^{\circ}$}
 \end{figure}

We now continue this proof by considering three cases and showing
the appropriate divisibility in each.

\emph{(i)} \quad  If $c_v(\Lambda^{\circ}) < 0$ (i.e. $u$ and $v$ are in the first
column of $\Lambda^{\circ}$) then $v = u+1$ and $F_{\Lambda^{\circ}}(z_1, \dots ,
z_n)$ can be written as $F_{\Lambda^{\circ}}(z_1, \dots, z_n) = f_{uv}(z_u
+c_u(\Lambda^{\circ}), z_v
+c_v(\Lambda^{\circ})) \cdot F$. \\
Starting with $F_{\Lambda^{\circ}}(z_1, \dots , z_n)$ written in the
ordering described in Proposition \ref{jimtheorem1} and simply
moving the term $f_{u,v}(z_u + c_u(\Lambda^{\circ}), z_v + c_v(\Lambda^{\circ}))$ to the left
results in all the singularity terms in the product $F$ remaining
tied to the same triple terms as originally described. Therefore we may still form regular triples for
each singularity in $F$, and hence $F$ is regular at $z_1 = z_2 = \cdots = z_n$.\\
So by considering this expression for $F_{\Lambda^{\circ}}(z_1, \dots ,
z_n)$ at $z_1 = z_2 = \cdots = z_n$ we see that $F_{\Lambda^{\circ}}$ will be
divisible on the left or right by $f_{uv}(c_u(\Lambda^{\circ}), c_v(\Lambda^{\circ})) = (1 - (uv))$.

\emph{(ii)} \quad  If $c_v(\Lambda^{\circ}) = 0$ then $v=s+1$, and $F_{\Lambda^{\circ}}(z_1,
\dots , z_n)$ can be written as \[ F_{\Lambda^{\circ}}(z_1, \dots , z_n) =
\prod_{i = u, \dots , s}^\leftarrow f_{i, s+1}(z_i + c_i(\Lambda^{\circ}), z_{s+1}
+ c_{s+1}(\Lambda^{\circ})) \cdot F' . \] Again, starting with the ordering
described in Proposition \ref{jimtheorem1}, this results in  all
the singularity terms in the product $F'$ remaining tied to the
same triple terms as originally described in that ordering. Hence
$F'$ is regular at $z_1 = z_2 = \cdots = z_n$. And so $F_{\Lambda^{\circ}}$
is divisible on the left by
\[ \prod_{i = u, \dots , s}^\leftarrow f_{i, s+1}(c_i(\Lambda^{\circ}), c_{s+1}(\Lambda^{\circ})).
\]

\emph{(iii)} \quad  If $c_v(\Lambda^{\circ}) > 0$ (i.e. $v$ is above the steps)
then $f_{uv}(z_u + c_u(\Lambda^{\circ}), z_v + c_v(\Lambda^{\circ}))$ is tied to the singularity
$f_{u-1, v}(z_{u-1} + c_{u-1}, z_v + c_v(\Lambda^{\circ}))$ as a triple term. To
show divisibility by $f_{uv}(z_u + c_u(\Lambda^{\circ}), z_v + c_v(\Lambda^{\circ}))$ in this case
we need an alternative expression for $F_{\Lambda^{\circ}}(z_1, \dots ,
z_n)$ that is regular when $z_1 = z_2 = \cdots = z_n$. Define a
permutation $\tau$ as follows,

\[ \tau = \prod_{i = u, \dots, s}^\rightarrow \left( \prod_{j= s+1, \dots,
v}^\leftarrow (i j) \right)
\phantom{XXXXXXXXXXXXXXXXXXXXXXXXXXXXXXX}
\]
\[ = \left(
\textrm{ \scriptsize $\begin{array}{ccccccccccccccccc}
                   1 & 2 & \dots & u-1 & u & u+1 & \dots &  & \dots & \dots  &  & \dots & v-1 & v & v+1 & \dots & n \\
                   1 & 2 & \dots & u-1 & s+1 & s+2 & \dots & v-1 & v & u & u+1 & \dots & s-1 & s & v+1 & \dots & n \\
                 \end{array}$ } \right) \]

From the definition of $C_1$ in (\ref{cproduct}) we now define
$C'_1 = \tau C_1$, where $\tau$ acts on the indices of the product
$C_1$ such that
\[ \tau C_1 = \prod_{j=1}^n \left( \prod_{i=1}^r g_{i, \tau
\cdot j} \right). \]

%\newpage
For the rest of this proof we will simply write $f_{ij}$
instead of $f_{ij}(z_i + c_i(\Lambda^{\circ}), z_j + c_j(\Lambda^{\circ}))$. Define $R'_1$ as,
\newpage
\begin{eqnarray*}
R'_1 & = & \prod_{i= r+2, \dots, u-1}^\leftarrow \left(
\prod_{j=s+1, \dots ,v}^\rightarrow f_{ij} \right) \cdot \prod_{i=
s+1, \dots, t-1}^\rightarrow \left( \prod_{j=i+1, \dots
,t}^\rightarrow f_{ij} \right) \cdot \left( \prod_{j=s+1, \dots
,t}^\leftarrow f_{r+1, j} \right)  \\
   && \times \left( \prod_{j=t+1, \dots
,v}^\rightarrow f_{r+1, j} \right) \cdot \prod_{i= r+1,
\dots, s-1}^\rightarrow \left( \prod_{j=i+1, \dots ,s}^\rightarrow
f_{ij} \right) \cdot \prod_{j= v+1, \dots, n}^\rightarrow
\left( \prod_{i=r+1, \dots ,s}^\rightarrow f_{ij} \right). 
\end{eqnarray*}

Finally, define $C'_2$ as, \[ C'_2 = \prod_{j= t+1, \dots,
n}^\rightarrow \left( \prod_{i=s+1, \dots ,t}^\rightarrow f_{ij}
\right). \]

Then, \[ F_{\Lambda^{\circ}}(z_1, \dots , z_n) = \prod_{i = u, \dots ,
s}^\leftarrow \left( \prod_{j=s+1, \dots, v}^\rightarrow f_{ij}
\right) \cdot C'_1 R'_1 C'_2 R_2 \cdot \prod_{i=3}^d C_iR_i, \]
where $d$ is the number of principal hooks of $\lambda$.

The product $C'_1 R'_1 C'_2 R_2 \cdot \prod C_iR_i$ is regular at
$z_1 = z_2 = \cdots = z_n$ since, as before, for any singularity
$(p,q)$ the terms $f_{pq}f_{p+1 ,q}$ can be replaced by the triple
$\frac{1}{2}f_{p, p+1}f_{pq}f_{p+1, q}$ -- except in the
expression $R'_1$ where the terms $f_{pl}f_{pq}$ are replaced by
$\frac{1}{2} f_{pl}f_{pq}f_{lq}$, where $l$ is the entry to the
immediate left of $q$. Note that $l = q-1$ when $c_q(\Lambda^{\circ}) > 1$ and $l =
s+1$ when $c_q(\Lambda^{\circ}) = 1$. \\
And so by letting $z_1 = z_2 = \cdots =z_n$ we see that $F_{\Lambda^{\circ}}$
is divisible on the left by \[ \prod_{i = u, \dots , s}^\leftarrow
\left( \prod_{j=s+1, \dots, v}^\rightarrow f_{ij}(c_i(\Lambda^{\circ}), c_j(\Lambda^{\circ}))
\right).
\]
\end{proof}

\begin{proposition}\label{jimtheorem3} Suppose the numbers $u < v$ stand next to each
other in the same row of the hook tableau $\Lambda^{\circ}$ of shape
$\lambda$. Let $r$ be the last entry in the column containing $u$.
If $c_u(\Lambda^{\circ}) > 0$ then the element $F_{\Lambda^{\circ}} \in \mathbb{C}S_n$ is
divisible on the left and on the right by $f_{u,v}(c_u(\Lambda^{\circ}), c_v(\Lambda^{\circ})) = 1 +
(uv)$. If $c_u(\Lambda^{\circ}) \leqslant 0$ then the element $F_{\Lambda^{\circ}} \in
\mathbb{C}S_n$ is divisible on the left by the product
\[ \prod_{i = u, \dots, r}^\leftarrow \left( \prod_{j= r+1, \dots,
v}^\rightarrow f_{ij}(c_i(\Lambda^{\circ}), c_j(\Lambda^{\circ})) \right) \]
\end{proposition}
\begin{proof}
Suppose $u$ is to the immediate left of $v$ in some row of
$\Lambda^{\circ}$. As in the proof of Proposition \ref{jimtheorem2}, we
need only consider $(u,v)$ such that $u$ is in the first row or
first column of $\Lambda^{\circ}$. Let $r$ be the last entry in the first
column of $\Lambda^{\circ}$ and $s$ the last entry in the first row of
$\Lambda^{\circ}$.

As in the proof of Proposition \ref{jimtheorem2}, we consider
three cases.

\emph{(i)} \quad  If $c_u(\Lambda^{\circ}) > 0$ (i.e. $u$ and $v$ are in the first
row of $\Lambda^{\circ}$) then $v = u+1$ and $F_{\Lambda^{\circ}}(z_1, \dots , z_n)$
can be written as $F_{\Lambda^{\circ}}(z_1, \dots, z_n) = f_{uv}(z_u +c_u(\Lambda^{\circ}),
z_v +c_v(\Lambda^{\circ})) \cdot F$. Singularity terms in the product $F$ remaining
tied to the same triple terms as originally described in
Proposition \ref{jimtheorem1}, hence $F$ is regular at $z_1 = z_2
= \cdots = z_n$. By considering this expression for $F_{\Lambda^{\circ}}(z_1,
\dots , z_n)$ at $z_1 = z_2 = \cdots = z_n$ we see that $F_{\Lambda^{\circ}}$
will be divisible on the left/right by $f_{uv}(c_u(\Lambda^{\circ}), c_v(\Lambda^{\circ})) = (1 +
(uv))$.

\emph{(ii)} \quad  If $c_u(\Lambda^{\circ}) = 0$ then $u=1$, $v=r+1$, and
$F_{\Lambda^{\circ}}(z_1, \dots , z_n)$ can be written as \[ F_{\Lambda^{\circ}}(z_1,
\dots , z_n) = \prod_{i = 1, \dots , r}^\leftarrow f_{i, r+1}(z_i
+ c_i(\Lambda^{\circ}), z_{r+1} + c_{r+1}(\Lambda^{\circ})) \cdot F' .
\] Again, singularities in $F'$ remain tied to the same triple terms as
originally described in Proposition \ref{jimtheorem1}. And so
$F_{\Lambda^{\circ}}$ is divisible on the left by
\[ \prod_{i = 1, \dots , r}^\leftarrow f_{i, r+1}(c_i(\Lambda^{\circ}), c_{r+1}(\Lambda^{\circ})).
\]

\emph{(iii)} \quad  If $c_u(\Lambda^{\circ}) < 0$ (i.e. $v$ is below the steps)
then $f_{uv}(z_u + c_u(\Lambda^{\circ}), z_v + c_v(\Lambda^{\circ}))$ is tied to the singularity
$f_{u-1, v}(z_{u-1} + c_{u-1}, z_v + c_v(\Lambda^{\circ}))$ as a triple term.

For the rest of this proof we will simply write $f_{ij}$ instead
of $f_{ij}(z_i + c_i(\Lambda^{\circ}), z_j + c_j(\Lambda^{\circ}))$. Define $C''_1$ as,
\begin{eqnarray*}
C''_1 & = & \prod_{i= 2, \dots, u-1}^\leftarrow \left(
\prod_{j=r+1, \dots ,v}^\rightarrow f_{ij} \right) \cdot \prod_{i=
r+1, \dots, s-1}^\rightarrow \left( \prod_{j=i+1, \dots
,s}^\rightarrow f_{ij} \right) \cdot \left( \prod_{j=r+1, \dots
,s}^\leftarrow f_{1, j} \right)  \\
   && \times f_{1, s+1} \cdot \left(\prod_{i=r+1, \dots ,s}^\rightarrow f_{i, s+1} \right) \cdot
\prod_{j= s+2, \dots, v}^\rightarrow f_{1, j} \cdot \prod_{i= 1,
\dots, r-1}^\rightarrow \left( \prod_{j=i+1, \dots ,r}^\rightarrow
f_{ij} \right)  \\
   && \times \prod_{j= v+1, \dots, n}^\rightarrow
\left( \prod_{i=1, \dots ,r}^\rightarrow f_{ij} \right),
\end{eqnarray*}

and define $R''_1$ as, \[ R''_1 = \prod_{j= s+2, \dots,
n}^\rightarrow \left( \prod_{i=r+1, \dots ,s}^\rightarrow f_{ij}
\right). \]

Then, \[ F_{\Lambda^{\circ}}(z_1, \dots , z_n) = \prod_{i = u, \dots ,
r}^\leftarrow \left( \prod_{j=r+1, \dots, v}^\rightarrow f_{ij}
\right) \cdot C''_1 R''_1 \cdot \prod_{i=2}^d C_iR_i, \] where $d$
is the number of principal hooks of $\lambda$.

The product $C''_1 R''_1 \cdot \prod C_iR_i$ is regular at $z_1 =
z_2 = \cdots = z_n$ since, as before, for any singularity $(p,q)$ the
terms $f_{pq}f_{p+1 ,q}$ can be replaced by the triple
$\frac{1}{2}f_{p, p+1}f_{pq}f_{p+1, q}$ -- except in the
expression $C''_1$ where $f_{1, s+1}f_{r+1, s+1}$ is replaced by
$\frac{1}{2} f_{1, r+1}f_{1, s+1}f_{r+1, s+1}$, and $f_{p,
q-1}f_{pq}$ is replaced by
$\frac{1}{2} f_{p, q-1}f_{pq}f_{q-1, q}$ for all other singularities $(p, q)$ in $C''_1$.\\
And so by letting $z_1 = z_2 = \dots =z_n$ we see that $F_{\Lambda^{\circ}}$
is divisible on the left by \[ \prod_{i = u, \dots , r}^\leftarrow
\left( \prod_{j=r+1, \dots, v}^\rightarrow f_{ij}(c_i(\Lambda^{\circ}), c_j(\Lambda^{\circ}))
\right).
\]
\end{proof}

Let us now
consider an example that allows us to see how the product
\\ $F_{\Lambda^{\circ}}(z_1, \dots, z_n)$ is broken down in the proof of
Proposition \ref{jimtheorem2}.

\begin{example} We
again consider the hook tableau of the Young diagram $\lambda =
(3,3,2)$, as shown in Example 3.2.

We begin with the product $F_{\Lambda^{\circ}}(z_1, \dots, z_n)$ in the
ordering described in Proposition \ref{jimtheorem1}. For
simplicity we again write $f_{pq}$ in place of the term
$f_{pq}(z_p + c_p(\Lambda^{\circ}), z_q + c_q(\Lambda^{\circ}))$. We have also marked out the
singularities in this expansion along with their triple terms, but
no idempotents have yet been added which would form regular
triples.

\begin{equation}\label{example1}
\begin{array}{rl}
   F_{\Lambda^{\circ}} (z_1, \dots , z_n) = &
   f_{12}f_{13}f_{23}f_{14}f_{24}f_{34}f_{15}f_{25}f_{35}(f_{16}f_{26})f_{36}f_{17}(f_{27}f_{37})f_{18}f_{28}f_{38}\\&
   \cdot f_{45}f_{46}f_{56}f_{47}f_{57}(f_{48}f_{58}) \cdot f_{67}f_{68}f_{78}\\
\end{array}
\end{equation}

Let $u=4$ and $v=6$ in Proposition \ref{jimtheorem2}. Then by that
proposition we may arrange the above product as follows. Notice
since $c_v(\Lambda^{\circ}) = 0$ all singularity-triple term pairs remain the same.

\[
\begin{array}{rl}
   F_{\Lambda^{\circ}} (z_1, \dots , z_n) = &
   f_{56}f_{46} \cdot f_{12}f_{13}f_{23}(f_{16}f_{26})f_{36}f_{14}f_{24}f_{34}f_{15}f_{25}f_{35}f_{17}(f_{27}f_{37})\\&
   f_{18}f_{28}f_{38} \cdot f_{45}f_{47}f_{57}(f_{48}f_{58}) \cdot f_{67}f_{68}f_{78}\\
\end{array}
\]
We may now add the appropriate idempotents so that all
singularities remain in regular triples. And so by considering the
product at $z_1 = z_2 = \cdots = z_n$ we have that $F_{\Lambda^{\circ}}$ is
divisible on the left by $(1 - (46))$, preceded only by invertible
terms, as desired.

Now let $u=3$ and $v=7$ in Proposition \ref{jimtheorem3}. Then by
that proposition we may arrange (\ref{example1}) as follows.
Singularities in $C_1$ have been marked out with their alternative
triple terms, while all other singularity-triple term pairs remain
the same.
\[
\begin{array}{rl}
   F_{\Lambda^{\circ}} (z_1, \dots , z_n) = &
   f_{34}f_{35}f_{36}f_{37} \cdot f_{24}f_{25}(f_{26}f_{27}) \cdot f_{45} \cdot f_{13}(f_{14} \cdot f_{16}) \cdot f_{46}f_{56} \cdot f_{17} \\&
   \cdot f_{12}f_{13}f_{23} \cdot f_{18}f_{28}f_{38} \cdot f_{47}f_{57}(f_{48}f_{58}) \cdot f_{67}f_{68}f_{78}\\
\end{array}
\]
In moving $f_{37}$ to the left it is untied from the singularity
$f_{27}$. So we must form new triples which are regular at $z_1 =
z_2 = \cdots = z_n$ by the method described in Proposition
\ref{jimtheorem3}. Therefore, by considering the product at $z_1 =
z_2 = \cdots = z_n$, we have that $F_{\Lambda^{\circ}}$ is divisible on the
left by $(1 + (37))$, again preceded only by invertible terms, as
desired. \end{example}

\begin{lemma}\label{divisibilitybyadjacenttransposition} Let $\sigma_k$ be the adjacent transposition of $k$ and $k+1$. Let $\Lambda$ and $\tilde{\Lambda}$ be tableaux of the same shape such that $k= \Lambda(a,b) = \Lambda(a+1,b)-1$ and $\tilde{k}=\tilde{\Lambda}(a,b) = \tilde{\Lambda}(a+1,b)-1$. Then $F_\Lambda \in \mathbb{C}S_n$ is divisible on the left by $1 - \sigma_k$ if and only if $F_{\tilde{\Lambda}} \in \mathbb{C}S_n$ is divisible on the left by $1- \sigma_{\tilde{k}}$. \end{lemma}

\begin{proof}
Let $\sigma$ be the permutation such that $\tilde{\Lambda}=\sigma\cdot\Lambda$.
% Note that $\sigma$ is not necessarily an adjacent transposition. 
There is a decomposition
$\sigma=\sigma_{i_N}\ldots\sigma_{i_1}$ such that for each $M=1, \ldots , N-1$
the tableau $\Lambda_{ M}=\sigma_{i_M}\ldots\sigma_{i_1}\cdot\Lambda$ is
standard. Note that this decomposition is not necessarily reduced.

Denote $F_\Lambda$ by $f_{k,k+1}(c_k(\Lambda),c_{k+1}(\Lambda)) \cdot F$. Then, by using the relations
(\ref{triple}) and (\ref{commute}) and  Proposition \ref{jimtheorem1}, we have the following chain of equalities:

\[
f_{\tilde{k}, \tilde{k}+1}
\bigl( c_{\tilde{k}}(\tilde{\Lambda}), c_{\tilde{k}+1}(\tilde{\Lambda})
\bigr) \sigma_{\tilde{k}}
\ \cdot\ 
\prod_{M=1, \ldots , N}^{\longleftarrow}
f_{i_M+1, i_M}
\bigl( c_{i_M+1}(\Lambda_M)
, c_{i_M}(\Lambda_M)
\bigr) \sigma_{i_M} \cdot F   \]
\begin{eqnarray*} &=&\prod_{M=1, \ldots , N}^{\longleftarrow}
\sigma_{i_M} f_{i_M, i_M+1}
\bigl( c_{i_M}(\Lambda_M)
, c_{i_M+1}(\Lambda_M)
\bigr)
\ \cdot\ 
\sigma_k f_{k, k+1}
\bigl( c_k(\Lambda), c_{k+1}(\Lambda)\bigr) \cdot F \\
&=&
\prod_{M=1, \ldots , N}^{\longleftarrow}
\sigma_{i_M} f_{i_M, i_M+1}
\bigl( c_{i_M}(\Lambda_M)
, c_{i_M+1}(\Lambda_M)
\bigr)
\ \cdot\ F_\Lambda\
\\ 
&=&
F_{\tilde{\Lambda}}\ \cdot
\prod_{M=1, \ldots , N}^{\longleftarrow}
f_{i_M+1, i_M}
\bigl( c_{i_M+1}(\Lambda_M)
, c_{i_M}(\Lambda_M)
\bigr) \sigma_{i_M}
\end{eqnarray*}

Hence divisibility by $1-\sigma_k$ for $F_\Lambda$ implies its counterpart for the tableau
$\tilde{\Lambda}$ and the index $\tilde{k}$, and vice versa.
Here we also use the equalities
\[
f_{k,k+1}
\bigl( c_k(\Lambda), c_{k+1}(\Lambda)\bigr)
=1-\sigma_k,
\]\[
f_{\tilde{k}, \tilde{k}+1}
\bigl( c_{\tilde{k}}(\tilde{\Lambda}), c_{\tilde{k}+1}(\tilde{\Lambda})
\bigr)
=1-\sigma_{\tilde{k}}.
\]
\end{proof}

\begin{corollary}\label{divisibilitycorollary}
If $k=\Lambda(a, b)$ and $k+1=\Lambda(a+1, b)$
then the element $F_\Lambda\in \mathbb{C}S_n$ is divisible on the left and right by $1-\sigma_k$. If $k=\Lambda(a, b)$ and $k+1=\Lambda(a, b+1)$
then the element $F_\Lambda\in \mathbb{C}S_n$ is divisible on the left and right by 
$1+\sigma_k$.
\end{corollary}

\begin{proof}
By Proposition \ref{varphi} divisibility on the left is equivalent to divisibility on the right. Let us show divisibility on the left.

Due to Lemma \ref{divisibilitybyadjacenttransposition}
it suffices to prove the first part of Corollary \ref{divisibilitycorollary} for only one 
standard tableau $\Lambda$ of shape $\lambda$. Therefore, using Proposition \ref{jimtheorem2} and taking $\tilde{\Lambda}$ to be the hook tableau $\Lambda^\circ$ of shape $\lambda$ we have shown the first part of Corollary \ref{divisibilitycorollary} in the case $c_v(\Lambda) < 0$.

Next let $\Lambda^{\circ}(a,b)= u$, $\Lambda^{\circ}(a+1,b) = v$ and $s$ be the last entry in the row containing $u$. Then for $c_v({\Lambda^\circ}) \geqslant 0$ Proposition \ref{jimtheorem2} showed that $F_{\Lambda^\circ} \in \mathbb{C}S_n$ is divisible on the left by 
\begin{equation}\label{divisibilitycorollaryequation2} \prod_{i = u, \dots, s}^\leftarrow \left( \prod_{j= s+1, \dots,
v}^\rightarrow f_{i,j}(c_i(\Lambda^\circ), c_j(\Lambda^\circ)) \right) \end{equation}

Put $k=u+v-s-1$ and let $\Lambda$ be the tableau
such that $\Lambda^\circ$ is obtained from 
the tableau $\sigma_k\cdot\Lambda$ by the permutation
\[ \sigma = 
\prod_{i =  u, \ldots, s}^{\longleftarrow}\,
\biggl(\ 
\prod_{j = s+1, \ldots,  v}^{\longrightarrow}
\sigma_{ i+j-s-1} \biggr).
\]
The tableau $\Lambda$ is standard. Moreover, then
$\Lambda(a, b)=k$ and $\Lambda(a+1, b)=k+1$.
Note that the rightmost factor in the product (\ref{divisibilitycorollaryequation2}),
corresponding to $i=u$ and $j=v$, is
\[
f_{u,v}
\bigl( c_u(\Lambda^\circ), c_{v}(\Lambda^\circ)\bigr)
=
1 - (uv).
\]

Denote by $F$ the product of all factors in (\ref{divisibilitycorollaryequation2})
but the rightmost one. Further, denote by $F'$ the product obtained by replacing each factor in $F$ 
\[
f_{i,j}
\bigl(c_i(\Lambda^\circ), c_j(\Lambda^\circ)\bigr)
\]
respectively by
\[
\sigma_{i+j-s-1}f_{i+j-s-1,i+j-s}
\bigl(c_i(\Lambda^\circ), c_j(\Lambda^\circ)\bigr),
\]
and denote by $G'$ the product obtained
by replacing each factor in $F$ by \[
f_{i+j-s,i+j-s-1}
\bigl(c_j(\Lambda^\circ), c_i(\Lambda^\circ)\bigr) \sigma_{i+j-s-1}.
\]

The element $F \in \mathbb{C}S_n$ is invertible, and we have
\[ F \sigma \sigma_k \cdot F_\Lambda = F' \cdot F_\Lambda = F_{\Lambda^\circ} \cdot G' = F \cdot (1- (uv)) \cdot (C'_1R'_1C'_2R_2 \prod C_iR_i) \cdot G',\] 
where the final equality is as described in Proposition \ref{jimtheorem2}. Finally we use the identity $(uv)\sigma= \sigma\sigma_k$. Therefore the divisibility of 
the element $F_{\Lambda^\circ}$ on the left by the product (\ref{divisibilitycorollaryequation2})
will imply the divisibility of the element $F_\Lambda$ on the left by 
$1-\sigma_k$.

This shows the required divisibility for the tableau $\Lambda = \sigma \sigma_k \cdot \Lambda^{\circ}$. Using Lemma \ref{divisibilitybyadjacenttransposition} again concludes the proof of the first part of Corollary \ref{divisibilitycorollary}. 

The second part of Corollary \ref{divisibilitycorollary} may be shown similarly.
\end{proof}

We may now complete the proof of Theorem \ref{fulltheorem}.

\begin{proposition}The element $F_\Lambda$ is the diagonal matrix element for $\Lambda$.\end{proposition}

\begin{proof} Let $F_\Lambda = \sum_{\sigma \in S_n} \alpha_\sigma(\Lambda) \cdot \sigma$, where the coefficients $\alpha_\sigma(\Lambda)$ are in $\mathbb{C}$.

Let $\sigma_k$ be the adjacent transposition of $k$ and $k+1$. For $\sigma_k \Lambda$ nonstandard we have, by Proposition, \ref{divisibilitycorollary} \[\sigma_kF_\Lambda = F_\Lambda = F_\Lambda\sigma_k, \textrm{ for all } \sigma_k \in R(\Lambda), \textrm{ and, }\]\[ \sigma_kF_\Lambda= -F_\Lambda = F_\Lambda\sigma_k, \textrm{ for all } \sigma_k \in C(\Lambda).\]

Comparing coefficients gives us;
\begin{eqnarray*}\alpha_{\sigma_k\sigma}(\Lambda) &=& \alpha_\sigma(\Lambda), \textrm{for all } \sigma_k \in R(\Lambda), \textrm{and} \\
\alpha_{\sigma_k\sigma}(\Lambda) &=& - \alpha_\sigma(\Lambda), \textrm{for all } \sigma_k \in C(\Lambda). \end{eqnarray*}

For $\Lambda' = \sigma_k\Lambda$ standard we have, by Proposition \ref{jimtheorem1}, the identity \begin{equation}\label{FlambdaFlambdadash} \left(1- d_k(\Lambda')^2\right) F_\Lambda = \left( \sigma_k - d_k(\Lambda') \right) F_{\Lambda'} \left( \sigma_k - d_k(\Lambda') \right), \end{equation} where $d_k(\Lambda) = (c_{k+1}(\Lambda) - c_k(\Lambda))^{-1}$.

Comparing coefficients here gives us;
\begin{eqnarray*} \left(1-d_k(\Lambda')^2\right) \alpha_{\sigma_k\sigma}(\Lambda) \phantom{xxxxxxxxxxxxxxxxxxxxxxxxxxxxxxxxxxxxxx}\\ = \alpha_{\sigma\sigma_k}(\Lambda') - d_k(\Lambda')\alpha_{\sigma_k\sigma\sigma_k}(\Lambda') - d_k(\Lambda')\alpha_\sigma(\Lambda') + d_k(\Lambda')^2 \alpha_{\sigma_k\sigma}(\Lambda').\end{eqnarray*}

By Proposition \ref{coeff1} $\alpha_1(\Lambda) = 1$, therefore these identities completely determine the coefficients of $F_\Lambda$. In particular we see $\alpha_{\sigma_k} = d_k(\Lambda)$ for all adjacent transpositions $\sigma_k$.

On the other hand, consider the inner product of (\ref{dme}). By definition $(v_\Lambda, v_\Lambda) = 1$ and $(v_\Lambda, v_{\Lambda'}) = 0$ for $\Lambda \neq \Lambda'$.

\newpage 
Then by (\ref{sigmaonbasis}) we have for $\sigma_k\Lambda$ nonstandard; \[ (v_\Lambda, \sigma_k v_\Lambda) = 1 \textrm{ for all } \sigma_k \in R(\Lambda), \textrm{ and, } \]\[ (v_\Lambda, \sigma_k v_\Lambda) = -1 \textrm{ for all } \sigma_k \in C(\Lambda).\]

And if $\sigma_k\Lambda$ standard; \begin{eqnarray*}(v_\Lambda, \sigma_k v_\Lambda) &=& (v_\Lambda, d_k(\Lambda) v_\Lambda + \sqrt{1-d_k(\Lambda)^2} v_{\sigma_k \Lambda}) \\
&=& d_k(\Lambda). \end{eqnarray*}

We also have the equivalent identity for (\ref{FlambdaFlambdadash}) for diagonal matrix elements, (see (\ref{sigmaondme})). Hence the coefficient of $\sigma$ in the expression for the diagonal matrix element is the same as the coefficient $\alpha_\sigma(\Lambda)$ above for all $\sigma \in S_n$.
\end{proof}

Finally, let us consider the Giambelli identity,
(\ref{giambelli}), and its dual, (\ref{dualgiambelli}), for a
diagram $\lambda$. We obtain a decomposition of the induced
representation of the tensor product of representations of hook
shapes, where these hook shapes are the principal hooks of
$\lambda$.

\begin{proposition}\label{Dlambdalambda}The coefficient $D_{\lambda \lambda}$ in (\ref{dualgiambelli}) is equal to 1.\end{proposition}
\begin{proof} Consider the Littlewood-Richardson coefficient $c_{\mu \nu}^\lambda$. This is the coefficient of $s_\lambda(\textbf{x})$ in the decomposition of $s_\mu(\textbf{x})s_\nu(\textbf{x})$ and is equal to the number of ways to fill the skew tableau $\lambda/\mu$ with content $\nu$ such that the filling is a reverse lattice permutation.

Let $\lambda=(\alpha_1, \dots, \alpha_d|\beta_1, \dots, \beta_d)$ and let $\mu$ be a partition of $\sum_{i=1}^k h_i$, where $h_i$ is the length of the $i^\textrm{\tiny th}$ principal hook. If $\mu$ is not contained in \\ $(\alpha_1, \dots, \alpha_{k+1}|\beta_1, \dots, \beta_{k+1})$ then $c_{\mu (\alpha_{k+1}|\beta_{k+1})}^{(\alpha_1, \dots, \alpha_{k+1}|\beta_1, \dots, \beta_{k+1})} =0$.\\
If $\mu$ is contained in $(\alpha_1, \dots, \alpha_{k+1})$ then either $\mu = (\alpha_1, \dots, \alpha_k|\beta_1, \dots, \beta_k)$ and \\ $c_{\mu(\alpha_{k+1}|\beta_{k+1})}^{(\alpha_1, \dots, \alpha_{k+1}|\beta_1, \dots, \beta_{k+1})} = 1$, or $s_\mu(\textbf{x})$ has coefficient zero in the decomposition of \\ $s_{(\alpha_1|\beta_1)}(\textbf{x})s_{(\alpha_2|\beta_2)}(\textbf{x})\cdots s_{(\alpha_k|\beta_k)}(\textbf{x})$.\\
The same argument holds for the base case so, by induction, $s_\lambda$(\textbf{x}) appears with coefficient 1 in the decomposition of $s_{(\alpha_1|\beta_1)}(\textbf{x})s_{(\alpha_2|\beta_2)}(\textbf{x})\cdots s_{(\alpha_d|\beta_d)}(\textbf{x})$.
\end{proof}

\begin{example} Once
more let $\Lambda$ be a standard tableau of shape $\lambda = (3,3,2)$, or in the alternative notation used
in the Giambelli identity, $\lambda = (2,1|3,2)$. As described
before, the Giambelli identity gives us a way to write the Schur
function $s_\lambda(\textbf{x})$ as a determinant of Schur
functions of hook shape,
\[
s_{(2,1|3,2)}(\textbf{x}) = \left|
\begin{array}{cc}
  s_{(2|3)}(\textbf{x}) & s_{(2|2)}(\textbf{x}) \\
  s_{(1|4)}(\textbf{x}) & s_{(1|2)}(\textbf{x}) \\
\end{array}
\right| =
 \left|
\begin{array}{cc}
  s_{\textrm{\tiny \yng(3,1,1)}} & s_{\textrm{\tiny \yng(3,1)}} \\
  s_{\textrm{\tiny \yng(2,1,1)}} & s_{\textrm{\tiny \yng(2,1)}} \\
\end{array}
\right|,
\]
or dually,
\[ \begin{array}{rl} s_{(2|3)}(\textbf{x})s_{(1|2)}(\textbf{x}) = & s_{(2,1|3,2)}(\textbf{x}) + s_{(2,1|4,1)}(\textbf{x}) + s_{(2 |5,1)}(\textbf{x})
+ s_{(2 |4,2)}(\textbf{x}) + s_{(3,1|3,1)}(\textbf{x})\\& +
s_{(3|3,2)}(\textbf{x}) + 2 s_{(3|4,1)}(\textbf{x}) +
s_{(3|5)}(\textbf{x}) + s_{(4|3,1)}(\textbf{x}) +
s_{(4|4)}(\textbf{x}).
\end{array}\]

And therefore, in terms of modules, we obtain the following identity;
\[ \begin{array}{rl} \textrm{Ind}_{S_5 \times S_3}^{S_{8}} V_{(2|3)} \otimes
V_{(1|2)} =& V_{(2,1|3,2)} \oplus V_{(2,1|4,1)} \oplus V_{(2
|5,1)} \oplus V_{(2 |4,2)} \oplus V_{(3,1|3,1)} \\& \oplus
V_{(3|3,2)} \oplus 2 V_{(3|4,1)} \oplus V_{(3|5)} \oplus
V_{(4|3,1)} \oplus V_{(4|4)}.
\end{array} \]

Here, the induced module on the left hand side may be realised as
the left ideal of $\mathbb{C}S_n$ generated by $F_\Lambda(z_1,
\dots, z_{8})$ on $\mathcal{H}_\Lambda$ such that $z_i - z_j
\notin \mathbb{Z}$ when $i$ and $j$ are in different principal
hooks of $\lambda$. Since $\lambda$ only has two principal hooks
we only need two auxiliary parameters in the hook fusion
procedure. Letting these parameters be equal gives us the diagonal
matrix element $F_\Lambda$ which generates the irreducible
representation corresponding to $\lambda$, $V_{(2,1|3,2)}$, which
appears in the decomposition above with coefficient 1.
\end{example}

%% file: FindingEigenvalues.tex
\section{Representations of degenerate affine Hecke algebras}

We may now use our hook fusion procedure of Chapter 3 to calculate diagonal matrix elements computationally. This is only possible thanks to the minimisation of the parameters, making such calculations easier. We then use these elements in the construction of the eigenvalues of a certain operator which have been conjectured to provide the irreducibility criterion of certain representations. We provide supporting evidence for this conjecture, which includes calculating eigenvalues from a more general construction which have never been calculated before.

\subsection{The mixed hook length formula}\label{mhlf}

In this chapter we work with the degenerate affine Hecke algebra $H_n'$ of type $A$, for a description of degenerate affine Hecke algebras of other types see \cite{C3}, and was introduced by V. Drinfeld in [D]. As before, let $\sigma_k$ be the transposition of $k$ and $k+1$. Then the complex associative algebra $H_n'$ is generated by the symmetric group algebra $\mathbb{C}S_n$ and the pairwise commuting elements $y_1, \dots, y_n$ with the cross relations
\begin{eqnarray*} \sigma_i y_j &=& y_j \sigma_i, \quad j\neq i, i+1, \\
 \sigma_iy_i &=& y_{i+1}\sigma_i -1. \end{eqnarray*}
 
Let $V_\lambda$ be the irreducible $\mathbb{C}S_n$-module corresponding to the partition $\lambda$. We may also regard $V_\lambda$ as a $H_n'$-module by sending the element $y_k$ to the Jucys-Murphy element $X_k = (1 \;\; k) + (2 \;\; k) + \cdots + (k-1 \;\; k) \in \mathbb{C}S_n$ for each $k=1, \dots , n-1$.  For any number $z \in \mathbb{C}$ there is also an automorphism of $H_n'$ which acts as the identity on the subalgebra $\mathbb{C}S_n$, but such that $y_k \mapsto y_k + z$ for each $k=1, \dots, n$. We will denote by $V_\lambda(z)$ the $H_n'$-module obtained by pulling $V_\lambda$ back through this automorphism. The module $V_\lambda(z)$ is irreducible by construction and can be realised as the left ideal in $\mathbb{C}S_n$ generated by the diagonal matrix element $F_\Lambda$, for some standard tableau $\Lambda$ of shape $\lambda$. The subalgebra $\mathbb{C}S_n \subset H_n'$ acts here via left multiplication, and the action of the generators $y_1, \dots, y_n$ is $y_k \cdot F_\Lambda = (c_k(\Lambda) + z)F_\Lambda$ by (\ref{content}).

Consider the algebra $H_n' \otimes H_m'$. It is isomorphic to the subalgebra in $H_{n+m}'$ generated the transpositions $(p q)$, where $1 \leqslant p < q \leqslant n$ or $n+1 \leqslant p < q \leqslant n+m$, along with the elements $y_1, \dots, y_{n+m}$. For any partition $\mu$ of $m$ and any number $w \in \mathbb{C}$ take the corresponding $H_m'$-module $V_\mu(w)$. Now consider the $H_{n+m}'$-module $W$ induced from the $H_n' \otimes H_m'$-module $V_\lambda(z) \otimes V_\mu(w)$. Let $M$ be a standard tableaux of shape $\mu$ and set $\overline{F}_M = \epsilon_n(F_M)$, where $\epsilon$ is as defined in (\ref{epsilon}). Then the module $W$ may be realised as the left ideal in $\mathbb{C}S_{n+m}$ generated by the product $F_\Lambda \overline{F}_M$. The action of the generators $y_1, \dots y_{n+m}$ in this ideal can be determined using 
\[ y_k \cdot F_\Lambda \overline{F}_M = (c_k(\Lambda) + z)F_\Lambda \overline{F}_M \quad \textrm{ for each } k=1, \dots, n;\]\[ y_{n+k} \cdot F_\Lambda \overline{F}_M = (c_k(M) + w)F_\Lambda \overline{F}_M \quad \textrm{ for each } k=1, \dots, m.\] 
If $z-w \notin \mathbb{Z}$, then the module $W$ is irreducible, see \cite{C1}.

Now introduce the ordered products in the symmetric group algebra $\mathbb{C}S_{n+m}$
\[ R_{\Lambda M}(z,w) = \prod_{i=1, \dots, n}^\rightarrow \left( \prod_{j=1, \dots, m}^\leftarrow f_{i, n+j}(z+ c_i(\Lambda) , w + c_j(M)) \right), \]
\[ R_{\Lambda M}'(z,w) = \prod_{i=1, \dots, n}^\leftarrow \left( \prod_{j=1, \dots, m}^\rightarrow f_{i, n+j}(z+ c_i(\Lambda) , w + c_j(M)) \right). \] We keep to the assumption $z-w \notin \mathbb{Z}$. Applying Theorem \ref{fulltheorem} and using relations (\ref{triple}) and (\ref{commute}) repeatedly, we get \begin{equation}\label{preserves} F_\Lambda \overline{F}_M R_{\Lambda M}(z,w) = R_{\Lambda M}'(z,w) F_\Lambda \overline{F}_M.\end{equation} Hence the right multiplication in $\mathbb{C}S_{n+m}$ by $R_{\Lambda M}(z,w)$ preserves the left ideal $W$. 

Consider the operator of the right multiplication in $\mathbb{C}S_{n+m}$ by $R_{\Lambda M}(z,w)$. This operator preserves the subspace $W$ due to (\ref{preserves}). The restriction of this operator to $W$ will be denoted $J$. Furthermore, $J$ commutes with the action of $H_{n+m}'$, in other words $J$ is a $H_{n+m}'$-module homomorphism, see \cite{N2}. Now regard $W$ as a $\mathbb{C}S_{n+m}$-module only. Then the subalgebra $\mathbb{C}S_{n+m} \subset H_{n+m}'$ acts on $W \subset \mathbb{C}S_{n+m}$ via left multiplication and, under this action, $W$ splits into irreducible components according to the Littlewood-Richardson rule \cite{MD}. Let $\nu$ be any partition of $n+m$ such that the irreducible $\mathbb{C}S_{n+m}$-module $V_\nu$ appears in $W$ with multiplicity one. The homomorphism $J: W \to W$ preserves the subspace $V_\nu \subset W$ and therefore, by Schur's Lemma, acts there as multiplication by a certain number from $\mathbb{C}$. Denote this number $r_{\lambda \mu}^\nu(z,w)$, it depends on $z$ and $w$ as a rational function of $z-w$, but does not depend on the choice of tableaux $\Lambda$ and $M$. We will now compute the eigenvalues $r_{\lambda \mu}^\nu(z,w)$ of $J$ for certain $\nu$.

First let us observe the following general property of the eigenvalues $r_{\lambda \mu}^\nu(z,w)$. Consider the $H_{n+m}'$-module $W'$ induced from the $H_m' \otimes H_n'$-module $V_\mu(w) \otimes V_\lambda(z)$. Similar to the definition of $J$, one can define a $H_{n+m}'$-module homomorphism $J' : W' \to W'$ as the restriction to $W'$ of the operator of right multiplication in $\mathbb{C}S_{n+m}$ by $R_{M \Lambda}(w,z)$. There is a unique irreducible $\mathbb{C}S_{n+m}$-submodule  $V_\nu' \subset W'$ equivalent to $V_\nu \subset W$. By considering the corresponding eigenvalue $r_{\mu \lambda}^\nu(w,z)$ of the operator $J'$ we have the following identity derived in \cite{N2}.
\begin{proposition}\label{mixedhookproposition} \[ r_{\lambda \mu}^\nu(z,w) r_{\mu \lambda}^\nu(w,z) = \prod_{i=1}^{\lambda_1'} \prod_{k=1}^{\mu_1'} \frac{(z-w+\lambda_i - i +k)(z-w - \mu_k - i +k)}{(z-w+\lambda_i - \mu_k - i +k)(z-w-i+k)} . \] \end{proposition}
Furthermore, by the definition of $J$, we have the second identity $r_{\lambda \mu}^\nu(z,w) = r_{\lambda' \mu'}^{\nu'} (w,z)$, where $\lambda'$, $\mu'$ and $\nu'$ denote the conjugates of $\lambda$, $\mu$ and $\nu$.
We now continue by constructing certain partitions $\nu$ such that their corresponding irreducible $\mathbb{C}S_{n+m}$-modules $V_\nu$ appear in $W$ with multiplicity one.

Choose any sequence $a_1, \dots, a_{\lambda_1'} \in \mathbb{N}$ of pairwise distinct natural numbers. Note that the sequence does not need to be increasing. Consider the partition $\mu$ as an infinite sequence with finitely many non-zero terms. Define an infinite sequence $\gamma = (\gamma_1, \gamma_2, \dots )$ by 
\[ \gamma_{a_i} = \mu_{a_i} +\lambda_i, \quad i = 1, \dots , \lambda_1' ; \]
\[ \gamma_a = \mu_a, \quad a \neq a_1, \dots, a_{\lambda_1'} .\] If $\gamma_1 \geqslant \gamma_2 \geqslant \cdots$ then $\gamma$ is a partition of $n+m$. One may think of this partition as being constructed from the Young diagram of $\mu$ after adding the rows of $\lambda$ to the rows of $\mu$ in some order.\\
Similarly, if $b_1, \dots, b_{\lambda_1} \in \mathbb{N}$ is a sequence of pairwise distinct natural numbers, we may construct an infinite sequence $\delta' = (\delta_1', \delta_2', \dots)$ by 
\[ \delta_{b_j}' = \mu_{b_j}' + \lambda_j' , \quad j = 1, \dots, \lambda_1 ; \]
\[ \delta_b' = \mu_b', \quad b \neq b_1, \dots , b_{\lambda_1} .\] If $\delta_1' \geqslant \delta_2' \geqslant \cdots$ then $\delta'$ is a partition of $n+m$, as is its conjugate $\delta$. One may think of the partition $\delta$ as being constructed from the Young diagram of $\mu$ after adding the columns of $\lambda$ to the columns of $\mu$ in some order.

For such partitions, the irreducible $\mathbb{C}S_{n+m}$-modules $V_\gamma$ and $V_\delta$ appear in $W$ with multiplicity one. And the corresponding eigenvalues $r_{\lambda \mu}^\gamma(z,w)$ and $r_{\lambda \mu}^\delta(z,w)$ of the operator $J : W \to W$ have the following combinatorial description due to \cite{N2}.

\begin{theorem}[Nazarov] \label{mixedhooktheorem}
\[ r_{\lambda \mu}^\gamma(z,w) = \prod_{(i,j)} \frac{z-w-\lambda_j' - \mu_{a_i} + a_i + j -1}{z-w-i+j} ; \]
\[ r_{\lambda \mu}^\delta(z,w) = \prod_{(i,j)} \frac{z-w + \lambda_i + \mu_{b_j}' - i - b_j +1}{z-w-i+j} \] where the products are taken over all boxes $(i,j)$ of the Young diagram $\lambda$. \end{theorem}

So far we have assumed $z-w \notin \mathbb{Z}$, in which case the $H_{n+m}'$-module $W$ is irreducible. Now let us consider $z-w \in \mathbb{Z}$. If we multiply the terms of $R_{\Lambda M}(z,w)$ by their common denominator we will get a non-zero operator such that, for certain values of $z-w$, acts on the subspace $V_\nu$ with an eigenvalue of zero. Therefore the kernel of this operator is a proper non-trivial $H_{n+m}'$-submodule of $W$, and hence $W$ is reducible.

For example, there are two distinguished irreducible components of the $\mathbb{C}S_{n+m}$-module $W$ which always have multiplicity one. They correspond to partitions \[ \lambda + \mu = (\lambda_1+\mu_1, \lambda_2 +\mu_2, \dots ) \quad \textrm{ and } \quad (\lambda' + \mu')' = (\lambda_1' + \mu_1', \lambda_2' + \mu_2', \dots )'.\]
Denote by $h_{\lambda \mu}(z,w)$ the ratio of the corresponding eigenvalues of $J$, \\$r_{\lambda \mu}^{\lambda+\mu}(z,w) / r_{\lambda \mu}^{(\lambda'+\mu')'}(z,w)$. This ratio does not depend on the normalisation of the operator $J$. We then have the following corollary to Theorem \ref{mixedhooktheorem}.

\begin{corollary}\label{mixedhookcorollary} \[ h_{\lambda \mu}(z,w) = \prod_{(i,j)} \frac{z-w - \lambda_j' - \mu_i + i + j -1}{z-w+\lambda_i + \mu_j - i - j+1} \] where the product is taken over all boxes $(i,j)$ belonging to the intersection of the diagrams $\lambda$ and $\mu$. \end{corollary}

If $\lambda =\mu$ the numbers $\lambda_i + \lambda_j' - i - j +1$ are called the \emph{hook lengths} of the diagram $\lambda$, \cite{MD}. If $\lambda \neq \mu$ the numbers \[ \lambda_i + \mu_j' - i - j + 1 \quad \textrm{ and } \quad \lambda_j' + \mu_i - i - j +1 \] in the above fraction may be called the \emph{mixed hook lengths} of the first and second kind respectively. Both these numbers are positive for any box $(i,j)$ in the intersection of $\lambda$ and $\mu$. 

Then, by Corollary \ref{mixedhookcorollary}, $W$ is reducible if $z-w$ is either a mixed hook length of the first or second kind relative to $\lambda$ and $\mu$. When $\lambda = \mu$, the corollary implies that the module $W$ is reducible if $|z-w|$ is a hook length of $\lambda$. Furthermore, the module $W$ is irreducible for all remaining values $|z-w|$, \cite{LNT}.

The irreducibility criterion of the module $W$ for arbitrary $\lambda$ and $\mu$ has also been given in \cite{LNT}. This work shows that the module $W$ is reducible if and only if the difference $z-w$ belongs to a certain finite subset $\mathcal{S}_{\lambda \mu} \subset \mathbb{Z}$ determined in \cite{LZ}. %This subset satisfies the property $\mathcal{S}_{\lambda \mu} = - \mathcal{S}_{\mu \lambda}$. 
Denote by $\mathcal{D}_{\lambda \mu}$ the %union of the 
set of all zeros and poles of the rational functions $r_{\lambda \mu}^{\lambda+\mu}(z,w) / r_{\lambda \mu}^\nu(z,w)$ in $z-w$, where $\nu$ ranges over all partitions $\gamma$ and $\delta$ described above. Then $\mathcal{D}_{\lambda \mu} \subseteq \mathcal{S}_{\lambda \mu}$. %, and $-\mathcal{D}_{\mu \lambda} \subset \mathcal{S}_{\lambda \mu}$ also. 

At the end of \cite{N2} it was conjectured that the set $\mathcal{S}_{\lambda \mu}$ is equal to the set of all zeros and poles of the rational functions $r_{\lambda \mu}^{\lambda+\mu}(z,w) / r_{\lambda \mu}^\nu(z,w)$ in $z-w$, where $\nu$ ranges over all partitions of $n+m$ such that the $\mathbb{C}S_{n+m}$-module $V_\nu$ appears in $W$ with multiplicity one. The next section provides supporting evidence for this conjecture using a combination of the results in this section and computational methods.

\subsection{Non mixed hook calculations}

Given a partition $\lambda$ of $n$ and a partition $\mu$ of $m$, we call a partition $\nu$ of $n+m$ of \emph{mixed hook type} if it can be obtained either \begin{listr}
\item by adding rows of $\lambda$ to $\mu$;
\item by adding columns of $\lambda$ to $\mu$;
\item by adding rows of $\mu$ to $\lambda$;
\item by adding columns of $\mu$ to $\lambda$.
\end{listr}

We may calculated the eigenvalue $r_{\lambda \mu}^\nu(z,w)$ for partitions $\nu$ of mixed hook type (i) and (ii) by using Theorem \ref{mixedhooktheorem}, and those of type (iii) and (iv) by using a combination of Theorem \ref{mixedhooktheorem} and Proposition \ref{mixedhookproposition}.

Let $u = z-w$. The smallest examples of a partition $\nu$ that is \emph{not} of mixed hook type, while having Littlewood-Richardson coefficient $c_{\lambda \mu}^\nu =1$, occurs when $\lambda$ and $\mu$ are certain partitions of 4. There are only six such partitions, for which we have calculated the corresponding eigenvalues $r_{\lambda \mu}^\nu(u)$.

\begin{result}\label{nonmixedhookresults}

\[
\begin{tabular}{lllllllll}
1. &\qquad \qquad&  $\lambda =$ \small $\yng(3,1)$ &\qquad&  $\mu =$ \small $\yng(3,1)$ &\qquad&  $\nu =$ \small $\yng(5,3)$ &\qquad& $r_{\lambda \mu}^\nu(u) = \frac{(u-4)(u-1)}{(u)(u+1)}$ \\
&&&&&&&& \\
2. &\qquad \qquad&  $\lambda =$ \small $\yng(3,1)$ &\qquad&  $\mu =$ \small $\yng(2,2)$ &\qquad&  $\nu =$ \small $\yng(4,2,2)$ &\qquad& $r_{\lambda \mu}^\nu(u) = \frac{(u-3)(u+4)}{(u-1)(u+2)}$ \\
&&&&&&&& \\
3. &\qquad \qquad&  $\lambda =$ \small $\yng(2,2)$ &\qquad&  $\mu =$ \small $\yng(2,2)$ &\qquad&  $\nu =$ \small $\yng(4,3,1)$ &\qquad& $r_{\lambda \mu}^\nu(u) = \frac{(u-3)(u-2)(u+2)}{(u)(u)(u+1)}$ \\
&&&&&&&& \\
4. &\qquad \qquad&  $\lambda =$ \small $\yng(2,2)$ &\qquad&  $\mu =$ \small $\yng(2,2)$ &\qquad&  $\nu =$ \small $\yng(3,2,2,1)$ &\qquad& $r_{\lambda \mu}^\nu(u) = \frac{(u-2)(u+2)(u+3)}{(u-1)(u)(u)}$ \\
&&&&&&&& \\
5. &\qquad \qquad&  $\lambda =$ \small $\yng(2,2)$ &\qquad&  $\mu =$ \small $\yng(2,1,1)$ &\qquad&  $\nu =$ \small $\yng(3,3,1,1)$ &\qquad& $r_{\lambda \mu}^\nu(u) = \frac{(u-2)(u+3)}{u(u+1)}$ \\
&&&&&&&& \\
6. &\qquad \qquad&  $\lambda =$ \small $\yng(2,1,1)$ &\qquad&  $\mu =$ \small $\yng(2,1,1)$ &\qquad&  $\nu =$ \small $\yng(2,2,2,1,1)$ &\qquad& $r_{\lambda \mu}^\nu(u) = \frac{(u+1)(u+4)}{(u-1)(u)}$ \\
\end{tabular}
\]
\end{result}

%\newpage

Before we describe how we computed these eigenvalues, let us examine how these calculations support the conjecture at the end of section \ref{mhlf}. 

\begin{example}
Let $\lambda = (3,1)$ and $\mu =(2,2)$, then the following is a complete list of partitions $\nu$, such that the Littlewood-Richardson coefficient, $c_{\lambda \mu}^\nu$,  equals 1. We then give the corresponding eigenvalues $r_{\lambda \mu}^\nu(z,w)$ and $r_{\mu \lambda}^\nu(w,z)$.
\[\]

\newpage 

$\lambda =$ \small $\yng(3,1)$ \qquad  $\mu =$ \small $\yng(2,2)$

\[\]

\begin{tabular}{lllllll}
1. &\qquad&   $\nu =$ \small $\yng(5,3)$ &\qquad& \large $r_{\lambda \mu}^\nu(z,w) = \frac{(z-w-3)(z-w-2)}{(z-w+1)(z-w+2)}$ &\qquad& \large $r_{\mu \lambda}^\nu(w,z) = \frac{(z-w+3)(z-w+4)}{(z-w-1)(z-w)}$ \\
&&&&&& \\
2. &\qquad&   $\nu =$ \small $\yng(5,2,1)$ &\qquad& \large $r_{\lambda \mu}^\nu(z,w) = \frac{(z-w-3)}{(z-w+2)}$ &\qquad& \large $r_{\mu \lambda}^\nu(w,z) = \frac{(z-w-2)(z-w+3)(z-w+4)}{(z-w-1)(z-w)(z-w+1)}$ \\
&&&&&& \\
3. &\qquad&   $\nu =$ \small $\yng(4,3,1)$ &\qquad& \large $r_{\lambda \mu}^\nu(z,w) = \frac{(z-w-3)(z-w+4)}{(z-w+1)(z-w+2)}$ &\qquad& \large $r_{\mu \lambda}^\nu(w,z) = \frac{(z-w-2)(z-w+3)}{(z-w-1)(z-w)}$ \\
&&&&&& \\
$\dagger$4. &\qquad&   $\nu =$ \small $\yng(4,2,2)$ &\qquad& \large $r_{\lambda \mu}^\nu(z,w) = \frac{(z-w-3)(z-w+4)}{(z-w-1)(z-w+2)}$ &\qquad& \large $r_{\mu \lambda}^\nu(w,z) = \frac{(z-w-2)(z-w+3)}{(z-w)(z-w+1)}$ \\
&&&&&& \\
5. &\qquad&   $\nu =$ \small $\yng(4,2,1,1)$ &\qquad& \large $r_{\lambda \mu}^\nu(z,w) = \frac{(z-w+4)}{(z-w+2)}$ &\qquad& \large $r_{\mu \lambda}^\nu(w,z) = \frac{(z-w-3)(z-w-2)(z-w+3)}{(z-w-1)(z-w)(z-w+1)}$ \\
&&&&&& \\
6. &\qquad&   $\nu =$ \small $\yng(3,3,2)$ &\qquad& \large $r_{\lambda \mu}^\nu(z,w) = \frac{(z-w-3)(z-w+3)(z-w+4)}{(z-w-1)(z-w+1)(z-w+2)}$ &\qquad& \large $r_{\mu \lambda}^\nu(w,z) = \frac{(z-w-2)}{(z-w)}$ \\
&&&&&& \\
7. &\qquad&   $\nu =$ \small $\yng(3,2,2,1)$ &\qquad& \large $r_{\lambda \mu}^\nu(z,w) = \frac{(z-w+3)(z-w+4)}{(z-w-1)(z-w+2)}$ &\qquad& \large $r_{\mu \lambda}^\nu(w,z) = \frac{(z-w-3)(z-w-2)}{(z-w)(z-w+1)}$ \\
&&&&&& \\
\end{tabular}

\large 
The first result is the partition $\nu = \lambda +\mu$. The fourth result is not of mixed hook type and its eigenvalues were calaculated using GAP. All other eigenvalues were calculated using Theorem \ref{mixedhooktheorem} and Proposition \ref{mixedhookproposition}. Let $\mathcal{Z}_{\lambda \mu}$ denote the set of zeros and poles of $r_{\lambda \mu}^{\lambda+\mu}(z,w) / r_{\lambda \mu}^\nu(z,w)$, then for the $\nu$ listed above we have $\mathcal{Z}_{\lambda \mu} = \{\pm1, +2, \pm3, -4\}$.

Using \cite{LZ} one can demonstrate that if $\lambda_1', \mu_1' \leqslant 3$ then $\mathcal{S}_{\lambda \mu} = \mathcal{D}_{\lambda \mu} \cup (-\mathcal{D}_{\mu \lambda})$. Then we can see, using the above calculations, that $\mathcal{S}_{\lambda \mu}$ does indeed equal $\mathcal{Z}_{\lambda \mu}$ as expected. 
\end{example}

Similar calculations, where $\lambda$ and $\mu$ are partitions of 4, also support the conjecture. Although these are the smallest examples of partitions $\nu$ that are not of mixed hook type they are also the largest examples calculated so far. Results \ref{nonmixedhookresults} were calculated using GAP, see Appendix \ref{appendix a}. The results were calculated by constructing a test vector in $V_\nu$. This test vector is constructed by multiplying the primitive central idempotent, $Z_\nu$, corresponding to the representation $\nu$, on the right by the product $F_\Lambda \overline{F}_M$. The diagonal matrix elements $F_\Lambda$ and $F_M$ were computed using the hook fusion procedure of Section \ref{hookfusionprocedure}, which is a more efficient way to calculate the diagonal matrix elements since it minimises the number of auxiliary parameters needed. Finally, the eigenvalue $r_{\lambda \mu}^\nu(u)$ may be determined by multiplying the test vector on the right by the operator $R_{\Lambda M}(u)$.

We illustrate the method with a simple example below. This example is much smaller than the ones above. Here $\nu = \lambda +\mu$ and is the trivial representation of $\mathbb{C}S_4$.

\begin{example}
Let $\lambda = \yng(2)$ and $\mu = \yng(2)$. There is only one standard tableau, $\Lambda$, of shape $\lambda$. Similarly for $\mu$.

Take $\nu = \lambda + \mu = \yng(4)$, which has Littlewood-Richardson coefficient $c_{\lambda \mu}^\nu = 1$. Then the corresponding diagonal matrix elements are simply; \[ F_\Lambda = 1 + (1 \phantom{x} 2) \qquad \textrm{ and } \qquad \overline{F}_M = 1 + (3 \phantom{x} 4). \]

The primitive central idempotent corresponding to $\nu$ is given by the formula \[ Z_\nu = \frac{\dim V_\nu}{(n+m)!} \sum_{\sigma \in S_{n+m}} \chi^\nu(\sigma) \cdot \sigma = \frac{1}{4!} \sum_{\sigma \in S_4} \sigma , \] where $\chi^\nu$ is the character of the representation $\nu$. The test vector is then \[ v_\nu = Z_\nu F_\Lambda \overline{F}_M = 4 Z_\nu \in V_\nu.\]

The operator $J$ is then
\newpage
\begin{eqnarray*} 
R_{\Lambda M}(u) &=& f_{14} (z + c_1(\Lambda), w + c_2(M))f_{13}(z + c_1(\Lambda), w + c_1(M))\\ &\phantom{=}& \qquad \qquad \times f_{24}(z + c_2(\Lambda), w + c_2(M))f_{23}(z + c_2(\Lambda), w + c_1(M)) \\
&\phantom{=}& \\
&=& \left( 1 - \frac{(14)}{u-1} \right) \left( 1 - \frac{(13)}{u} \right) \left( 1 - \frac{(24)}{u} \right) \left( 1 - \frac{(23)}{u+1} \right) \\
&\phantom{=}&\\
&=& 1 - \frac{(14)}{u-1} - \frac{(13)}{u} + \frac{(134)}{u(u-1)} - \frac{(24)}{u} + \frac{(142)}{u(u-1)} + \frac{(13)(24)}{(u+1)(u-1)}\\ &\phantom{=}& \qquad  - \frac{2 \cdot (1342)}{u(u+1)(u-1)} - \frac{(23)}{u+1} + \frac{(14)(23)}{(u+1)(u-1)} + \frac{(132)}{u(u+1)} \\ &\phantom{=}& \qquad   - \frac{(1324)}{u(u+1)(u-1)} + \frac{(234)}{u(u+1)} - \frac{(1423)}{u(u+1)(u-1)} 
\end{eqnarray*}

And, \[ v_\nu \cdot R_{\Lambda M}(u) = \frac{u(u-1)(u-1)(u-2)}{u^2(u+1)(u-1)} \cdot v_\nu = r_{\lambda \mu}^\nu(u) v_\nu .\]
\end{example}

In particular, the numerators of the eigenvalues of \ref{nonmixedhookresults} have integral roots as a polynomial of $u$. This suggests there may be a combinatorial description for all such eigenvalues.

%% file: TheRibbonFusionProcedure.tex
\section{The ribbon fusion procedure}\label{ribbonfp}

We now consider a further generalisation of the fusion procedure. So far the fusion procedure has required us to decompose Young diagrams either into their rows, columns or hooks before applying a certain limiting process to obtain the diagonal matrix element. Therefore it is natural to consider how many other ways we may decompose the Young diagram in this procedure. So let us now consider ribbons of a  Young diagram, of which rows, columns and hooks are special cases. It is worth noting ribbons too have a Jacobi-Trudi or Giambelli style determinant to calculate Schur functions. We go on to introduce new notation which we may use to build certain types of ribbons, with the aim of making a ribbon fusion procedure.

\subsection{Ribbon Schur functions}

We have already seen some determinantal expressions for the Schur function $s_\lambda(\textbf{x})$ obtained either by decomposing the Young diagram $\lambda$ into rows or columns (Jacobi-Trudi identities (\ref{jt1}), (\ref{jt2})) or hooks (Giambelli (\ref{giambelli})). Since these decompositions also gave us a limiting procedure to find the diagonal matrix element $F_\Lambda$, let us consider some of the other ways to decompose a diagram $\lambda$. First we define a skew Schur function in the variables $\textbf{x} = (x_1, \dots, x_m)$ to be \[ s_{\lambda / \mu} (\textbf{x}) = \sum_{T} \prod_{(i,j) \in \lambda / \mu} x_{T(i,j)},\] where the sum is over all semistandard tableaux of shape $\lambda / \mu$ with entries $1, \dots, m$. We now state an important definition:

\begin{definition} A ribbon is a skew diagram with an edgewise connected set of boxes that contain no $2 \times 2$ block of boxes. \end{definition}

\newpage
For example, the skew diagram below is a ribbon.

\begin{center}
\begin{picture}(80,80)
%\put(0,60){\framebox(20,20)[r]{}}
%\put(0,40){\framebox(20,20)[r]{}}
\put(0,20){\framebox(20,20)[r]{}}
\put(0,0){\framebox(20,20)[r]{}}
%\put(20,60){\framebox(20,20)[r]{}}
%\put(20,40){\framebox(20,20)[r]{}}
\put(20,20){\framebox(20,20)[r]{}}
%\put(20,0){\framebox(20,20)[r]{}}
\put(40,60){\framebox(20,20)[r]{}}
\put(40,40){\framebox(20,20)[r]{}}
\put(40,20){\framebox(20,20)[r]{}}
%\put(40,0){\framebox(20,20)[r]{}}
\put(60,60){\framebox(20,20)[r]{}}
%\put(60,40){\framebox(20,20)[r]{}}
%\put(60,20){\framebox(20,20)[r]{}}
%\put(60,0){\framebox(20,20)[r]{}}
\end{picture}
\end{center}

These objects are also known as `strips' or `skew hooks' amongst other names. Note that rows, columns and hooks are all special cases of ribbon which are of non-skew shape. The term \emph{rim ribbon} denotes the maximal outer ribbon of a diagram. Below we give the successive rim ribbons of the partition  $(4, 3, 3, 1)$.

\begin{center}
\begin{picture}(80,80)
\put(0,60){\framebox(20,20)[r]{$\cdot \;\;\,$}}
\put(0,40){\framebox(20,20)[r]{}}
\put(0,20){\framebox(20,20)[r]{}}
\put(0,0){\framebox(20,20)[r]{}}
\put(20,60){\framebox(20,20)[r]{}}
\put(20,40){\framebox(20,20)[r]{}}
\put(20,20){\framebox(20,20)[r]{}}
%\put(20,0){\framebox(20,20)[r]{$8\;\;\,$}}
\put(40,60){\framebox(20,20)[r]{}}
\put(40,40){\framebox(20,20)[r]{}}
\put(40,20){\framebox(20,20)[r]{}}
%\put(40,0){\framebox(20,20)[r]{$12\;\;$}}
\put(60,60){\framebox(20,20)[r]{}}
%\put(60,40){\framebox(20,20)[r]{$14\;\;$}}
%\put(60,20){\framebox(20,20)[r]{$15\;\;$}}
%\put(60,0){\framebox(20,20)[r]{$16\;\;$}}

%vertical lines up
%\put(10,50){\line(0,1){20}}%2
%\put(10,30){\line(0,1){20}}%3
\put(10,10){\line(0,1){20}}%4

\put(30,50){\line(0,1){20}}%6
%\put(30,30){\line(0,1){20}}%7
%\put(30,10){\line(0,1){20}}%8

\put(50,50){\line(0,1){20}}%10
\put(50,30){\line(0,1){20}}%11
%\put(50,10){\line(0,1){20}}%12

%\put(70,50){\line(0,1){20}}%14
%\put(70,30){\line(0,1){20}}%15
%\put(70,10){\line(0,1){20}}%16

%Horizontal lines to the right
%\put(10,70){\line(1,0){20}}%1
\put(10,50){\line(1,0){20}}%2
\put(10,30){\line(1,0){20}}%3
%\put(10,10){\line(1,0){20}}%4

%\put(30,70){\line(1,0){20}}%5
%\put(30,50){\line(1,0){20}}%6
\put(30,30){\line(1,0){20}}%7
%\put(30,10){\line(1,0){20}}%8

\put(50,70){\line(1,0){20}}%9
%\put(50,50){\line(1,0){20}}%10
%\put(50,30){\line(1,0){20}}%11
%\put(50,10){\line(1,0){20}}%12
\end{picture}
\end{center}

In 1988, Lascoux and Pragacz obtained a new determinantal expression for the Schur function $s_\lambda(\textbf{x})$ obtained from a decomposition of $\lambda$ into rim ribbons, \cite{LP}. 

\begin{proposition}[Lascoux-Pragacz] Let $\lambda$ be a partition with $d$ boxes on the main diagonal, and let $(\theta_1, \theta_2, \dots, \theta_d)$ be its decomposition into rim ribbons, then \[ s_\lambda(\textbf{x}) = \det \left[ s_{\theta_i^- \& \theta_j^+ (\textbf{x})} \right]_{i,j=1}^d ,  \] where $\theta_i^- \& \theta_j^+$ denotes the skew diagram obtained by replacing the boxes below the diagonal of $\theta_j$ by that of $\theta_i$. \end{proposition}

Then, in 1995, Hamel and Goulden considered a much larger class of decompositions which gave determinantal expressions of skew Schur functions, and which also contain the Jacobi-Trudi, Giambelli and Lascoux-Pragacz identities as special cases, \cite{HG}. Suppose that $\theta_1, \dots , \theta_m$ are ribbons in  the skew diagram $\lambda / \mu$ such that each ribbon has a stating box on the left or bottom perimeter of the diagram and an ending box on the right or top perimeter of the diagram, then we say this set is a \emph{planar outside decomposition} of $\lambda / \mu$ if the disjoint union of these ribbons is the skew diagram $\lambda / \mu$. Below is the planar decomposition of $(4,3,3,1)$ into four ribbons.

\begin{center}
\begin{picture}(80,80)
\put(0,60){\framebox(20,20)[r]{}}
\put(0,40){\framebox(20,20)[r]{}}
\put(0,20){\framebox(20,20)[r]{}}
\put(0,0){\framebox(20,20)[r]{$\cdot \;\;\,$}}
\put(20,60){\framebox(20,20)[r]{}}
\put(20,40){\framebox(20,20)[r]{}}
\put(20,20){\framebox(20,20)[r]{}}
%\put(20,0){\framebox(20,20)[r]{$8\;\;\,$}}
\put(40,60){\framebox(20,20)[r]{}}
\put(40,40){\framebox(20,20)[r]{}}
\put(40,20){\framebox(20,20)[r]{}}
%\put(40,0){\framebox(20,20)[r]{$12\;\;$}}
\put(60,60){\framebox(20,20)[r]{}}
%\put(60,40){\framebox(20,20)[r]{$14\;\;$}}
%\put(60,20){\framebox(20,20)[r]{$15\;\;$}}
%\put(60,0){\framebox(20,20)[r]{$16\;\;$}}

%vertical lines up
%\put(10,50){\line(0,1){20}}%2
%\put(10,30){\line(0,1){20}}%3
%\put(10,10){\line(0,1){20}}%4

%\put(30,50){\line(0,1){20}}%6
%\put(30,30){\line(0,1){20}}%7
%\put(30,10){\line(0,1){20}}%8

\put(50,50){\line(0,1){20}}%10
%\put(50,30){\line(0,1){20}}%11
%\put(50,10){\line(0,1){20}}%12

%\put(70,50){\line(0,1){20}}%14
%\put(70,30){\line(0,1){20}}%15
%\put(70,10){\line(0,1){20}}%16

%Horizontal lines to the right
\put(10,70){\line(1,0){20}}%1
\put(10,50){\line(1,0){20}}%2
\put(10,30){\line(1,0){20}}%3
%\put(10,10){\line(1,0){20}}%4

%\put(30,70){\line(1,0){20}}%5
\put(30,50){\line(1,0){20}}%6
\put(30,30){\line(1,0){20}}%7
%\put(30,10){\line(1,0){20}}%8

\put(50,70){\line(1,0){20}}%9
%\put(50,50){\line(1,0){20}}%10
%\put(50,30){\line(1,0){20}}%11
%\put(50,10){\line(1,0){20}}%12
\end{picture}
\end{center}

Then we have the following remarkable identity;

\begin{theorem}[Hamel-Goulden] Let $\lambda / \mu$ be a skew shape partition. Then, for any planar outside decomposition $(\theta_1, \theta_2, \dots , \theta_m)$ of $\lambda / \mu$ we have \[ s_{\lambda / \mu} (\textbf{x}) = \left| s_{\theta_i \# \theta_j}(\textbf{x}) \right|_{m \times m} \] where $\theta_i \# \theta_j$ is the non-commutative superimposing of $\theta_i$ and $\theta_j$ described in \cite{HG}. \end{theorem}

However, in the same paper, Hamel and Goulden show that this determinantal expression does not necessarily hold for other decompositions that are not planar outside decompositions, in particular they give the following decomposition of $(4,4,4)$ as an example.

\begin{center}
\begin{picture}(80,60)
\put(0,40){\framebox(20,20)[r]{}}
\put(0,20){\framebox(20,20)[r]{}}
\put(0,0){\framebox(20,20)[r]{}}
\put(20,40){\framebox(20,20)[r]{}}
\put(20,20){\framebox(20,20)[r]{}}
\put(20,0){\framebox(20,20)[r]{}}
\put(40,40){\framebox(20,20)[r]{}}
\put(40,20){\framebox(20,20)[r]{}}
\put(40,0){\framebox(20,20)[r]{}}
\put(60,40){\framebox(20,20)[r]{}}
\put(60,20){\framebox(20,20)[r]{}}
\put(60,0){\framebox(20,20)[r]{}}

%vertical lines up
%\put(10,50){\line(0,1){20}}%2
\put(10,30){\line(0,1){20}}%3
\put(10,10){\line(0,1){20}}%4

%\put(30,50){\line(0,1){20}}%6
%\put(30,30){\line(0,1){20}}%7
\put(30,10){\line(0,1){20}}%8

%\put(50,50){\line(0,1){20}}%10
%\put(50,30){\line(0,1){20}}%11
%\put(50,10){\line(0,1){20}}%12

%\put(70,50){\line(0,1){20}}%14
\put(70,30){\line(0,1){20}}%15
\put(70,10){\line(0,1){20}}%16

%Horizontal lines to the right
%\put(10,70){\line(1,0){20}}%1
\put(10,50){\line(1,0){20}}%2
%\put(10,30){\line(1,0){20}}%3
%\put(10,10){\line(1,0){20}}%4

%\put(30,70){\line(1,0){20}}%5
\put(30,50){\line(1,0){20}}%6
\put(30,30){\line(1,0){20}}%7
%\put(30,10){\line(1,0){20}}%8

%\put(50,70){\line(1,0){20}}%9
%\put(50,50){\line(1,0){20}}%10
%\put(50,30){\line(1,0){20}}%11
\put(50,10){\line(1,0){20}}%12
\end{picture}
\end{center}

However, it turns out then that such a decomposition could be used in a ribbon fusion procedure, while other examples of planar outside decompositions fail. So instead we must construct our own class of decompositions which will form a basis for our ribbon fusion procedure.

\subsection{Towards a ribbon fusion procedure}

\newcommand{\type}[1]{\begin{picture}(10,10)\put(5,3.5){\circle{10}}\put(2.5,0.5){\textbf{\small #1}}\end{picture}\,}
\newcommand{\typeij}[1]{\begin{picture}(10,10)\put(5,3.5){\circle{10}}\put(3.5,0.5){\textbf{\small #1}}\end{picture}\,}
\newcommand{\bigtype}[1]{\begin{picture}(13,13)\put(6.5,4){\circle{13}}\put(1,1){\textbf{\small #1}}\end{picture}}
\newcommand{\bigtypecentre}[1]{\begin{picture}(13,13)\put(6.5,4){\circle{13}}\put(4,1){\textbf{\small #1}}\end{picture}}
\newcommand{\smalltype}[1]{\begin{picture}(7,7)\put(3.5,3.5){\circle{7}}\put(1,1.5){\textbf{\tiny #1}}\end{picture}}
\newcommand{\smalltypecentre}[1]{\begin{picture}(7,7)\put(3.5,3.5){\circle{7}}\put(2,1.5){\textbf{\tiny #1}}\end{picture}}

Consider a Young diagram of shape (2,2). That is a $2 \times 2$ square of four boxes. Each box has an associated parameter $z_i$. If $z_i = z_j$ connect the boxes containing $i$ and $j$ with a line. Below is a complete list of $2 \times 2$ squares connected using vertical and horizontal lines, we label them as types $\bigtypecentre{0}$ to $\bigtype{15}$.

\[
\begin{picture}(40,40)
\put(-20,25){\bigtypecentre{0}}
\put(0,0){\framebox(20,20)[r]{  }}
\put(0,20){\framebox(20,20)[r]{  }}
\put(20,0){\framebox(20,20)[r]{  }}
\put(20,20){\framebox(20,20)[r]{  }}
\put(10,10){\line(0,1){20}}
\put(30,10){\line(0,1){20}}
\end{picture}
\qquad \qquad
\begin{picture}(40,40)
\put(-20,25){\bigtypecentre{1}}
\put(0,0){\framebox(20,20)[r]{  }}
\put(0,20){\framebox(20,20)[r]{  }}
\put(20,0){\framebox(20,20)[r]{  }}
\put(20,20){\framebox(20,20)[r]{  }}
\put(10,10){\line(0,1){20}}
\put(29,9){.}
\put(29,29){.}
\end{picture}
\qquad \qquad
\begin{picture}(40,40)
\put(-20,25){\bigtypecentre{2}}
\put(0,0){\framebox(20,20)[r]{  }}
\put(0,20){\framebox(20,20)[r]{  }}
\put(20,0){\framebox(20,20)[r]{  }}
\put(20,20){\framebox(20,20)[r]{  }}
\put(10,10){\line(0,1){20}}
\put(10,30){\line(1,0){20}}
\put(29,9){.}
\end{picture}
\qquad \qquad
\begin{picture}(40,40)
\put(-20,25){\bigtypecentre{3}}
\put(0,0){\framebox(20,20)[r]{  }}
\put(0,20){\framebox(20,20)[r]{  }}
\put(20,0){\framebox(20,20)[r]{  }}
\put(20,20){\framebox(20,20)[r]{  }}
\put(10,30){\line(1,0){20}}
\put(9,9){.}
\put(29,9){.}
\end{picture}
\]
\[
\begin{picture}(40,40)
\put(-20,25){\bigtypecentre{4}}
\put(0,0){\framebox(20,20)[r]{  }}
\put(0,20){\framebox(20,20)[r]{  }}
\put(20,0){\framebox(20,20)[r]{  }}
\put(20,20){\framebox(20,20)[r]{  }}
\put(10,10){\line(1,0){20}}
\put(10,30){\line(1,0){20}}
\end{picture}
\qquad \qquad
\begin{picture}(40,40)
\put(-20,25){\bigtypecentre{5}}
\put(0,0){\framebox(20,20)[r]{  }}
\put(0,20){\framebox(20,20)[r]{  }}
\put(20,0){\framebox(20,20)[r]{  }}
\put(20,20){\framebox(20,20)[r]{  }}
\put(10,10){\line(1,0){20}}
\put(9,29){.}
\put(29,29){.}
\end{picture}
\qquad \qquad
\begin{picture}(40,40)
\put(-20,25){\bigtypecentre{6}}
\put(0,0){\framebox(20,20)[r]{  }}
\put(0,20){\framebox(20,20)[r]{  }}
\put(20,0){\framebox(20,20)[r]{  }}
\put(20,20){\framebox(20,20)[r]{  }}
\put(10,10){\line(1,0){20}}
\put(30,10){\line(0,1){20}}
\put(9,29){.}
\end{picture}
\qquad \qquad
\begin{picture}(40,40)
\put(-20,25){\bigtypecentre{7}}
\put(0,0){\framebox(20,20)[r]{  }}
\put(0,20){\framebox(20,20)[r]{  }}
\put(20,0){\framebox(20,20)[r]{  }}
\put(20,20){\framebox(20,20)[r]{  }}
\put(30,10){\line(0,1){20}}
\put(9,9){.}
\put(9,29){.}
\end{picture}
\]
\[
\begin{picture}(40,40)
\put(-20,25){\bigtypecentre{8}}
\put(0,0){\framebox(20,20)[r]{  }}
\put(0,20){\framebox(20,20)[r]{  }}
\put(20,0){\framebox(20,20)[r]{  }}
\put(20,20){\framebox(20,20)[r]{  }}
\put(10,30){\line(1,0){20}}
\put(30,10){\line(0,1){20}}
\put(9,9){.}
\end{picture}
\qquad \qquad
\begin{picture}(40,40)
\put(-20,25){\bigtypecentre{9}}
\put(0,0){\framebox(20,20)[r]{  }}
\put(0,20){\framebox(20,20)[r]{  }}
\put(20,0){\framebox(20,20)[r]{  }}
\put(20,20){\framebox(20,20)[r]{  }}
\put(10,10){\line(0,1){20}}
\put(10,10){\line(1,0){20}}
\put(29,29){.}
\end{picture}
\qquad \qquad
\begin{picture}(40,40)
\put(-20,25){\bigtype{10}}
\put(0,0){\framebox(20,20)[r]{  }}
\put(0,20){\framebox(20,20)[r]{  }}
\put(20,0){\framebox(20,20)[r]{  }}
\put(20,20){\framebox(20,20)[r]{  }}
\put(10,10){\line(0,1){20}}
\put(10,10){\line(1,0){20}}
\put(10,30){\line(1,0){20}}
\end{picture}
\qquad \qquad
\begin{picture}(40,40)
\put(-20,25){\bigtype{11}}
\put(0,0){\framebox(20,20)[r]{  }}
\put(0,20){\framebox(20,20)[r]{  }}
\put(20,0){\framebox(20,20)[r]{  }}
\put(20,20){\framebox(20,20)[r]{  }}
\put(10,10){\line(0,1){20}}
\put(30,10){\line(0,1){20}}
\put(10,30){\line(1,0){20}}
\end{picture}
\]
\[
\begin{picture}(40,40)
\put(-20,25){\bigtype{12}}
\put(0,0){\framebox(20,20)[r]{  }}
\put(0,20){\framebox(20,20)[r]{  }}
\put(20,0){\framebox(20,20)[r]{  }}
\put(20,20){\framebox(20,20)[r]{  }}
\put(30,10){\line(0,1){20}}
\put(10,10){\line(1,0){20}}
\put(10,30){\line(1,0){20}}
\end{picture}
\qquad \qquad
\begin{picture}(40,40)
\put(-20,25){\bigtype{13}}
\put(0,0){\framebox(20,20)[r]{  }}
\put(0,20){\framebox(20,20)[r]{  }}
\put(20,0){\framebox(20,20)[r]{  }}
\put(20,20){\framebox(20,20)[r]{  }}
\put(10,10){\line(0,1){20}}
\put(30,10){\line(0,1){20}}
\put(10,10){\line(1,0){20}}
\end{picture}
\qquad \qquad
\begin{picture}(40,40)
\put(-20,25){\bigtype{14}}
\put(0,0){\framebox(20,20)[r]{  }}
\put(0,20){\framebox(20,20)[r]{  }}
\put(20,0){\framebox(20,20)[r]{  }}
\put(20,20){\framebox(20,20)[r]{  }}
\put(10,10){\line(0,1){20}}
\put(30,10){\line(0,1){20}}
\put(10,10){\line(1,0){20}}
\put(10,30){\line(1,0){20}}
\end{picture}
\qquad \qquad
\begin{picture}(40,40)
\put(-20,25){\bigtype{15}}
\put(0,0){\framebox(20,20)[r]{  }}
\put(0,20){\framebox(20,20)[r]{  }}
\put(20,0){\framebox(20,20)[r]{  }}
\put(20,20){\framebox(20,20)[r]{  }}
\put(9,9){.}
\put(9,29){.}
\put(29,9){.}
\put(29,29){.}
\end{picture}
\]

If $p$ and $q$ sit on the same diagonal of a tableau $\Lambda$, we say the singularity $(p,q)$ is of degree $r$ if $q$ is $r$ steps away from $p$ on the diagonal. For example the $2\times2$ squares above contain a singularity of degree 1 only.

We say a singularity $(p,q)$ is of type $\typeij{i}$, $0 \leqslant i \leqslant 15$, if we can form a $2\times2$ square of type $\typeij{i}$ as follows: Take the boxes containing $p$ and $q$ and identify (or superimpose) the box below $p$ with the box to the left of $q$, and identify the box to the right of $p$ with the box above $q$ to form a $2\times2$ square, see the below figure.

\begin{figure}[h]
\label{identify}
\begin{center}
\begin{picture}(75, 75)

\put(0,0){\line(1,0){75}} \put(0,75){\line(1,0){75}}
\put(0,60){\line(1,0){30}} \put(45,15){\line(1,0){30}} 
\put(0,45){\line(1,0){15}} \put(60,30){\line(1,0){15}}

\put(0,0){\line(0,1){75}} \put(75,0){\line(0,1){75}}
\put(15,45){\line(0,1){30}} \put(60,0){\line(0,1){30}}
\put(30,60){\line(0,1){15}} \put(45,0){\line(0,1){15}}

\put(6,66){$p$} \put(66,6){$q$}
\put(21,65){$\dagger$} \put(66,20){$\dagger$}
\put(5,49){$\ast$} \put(51,5){$\ast$}
\put(23,40){$\ddots$} \put(42,25){$\ddots$}

\end{picture}
\end{center}
 \caption{We identify the boxes labeled $\dagger$ and the boxes labeled $\ast$ to form a $2 \times 2$ square.}
 \end{figure}
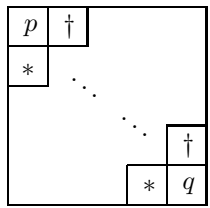

For example, the $3\times3$ square
\[\begin{picture}(60,60)
\put(0,0){\framebox(20,20)[r]{  }}
\put(0,20){\framebox(20,20)[r]{  }}
\put(0,40){\framebox(20,20)[r]{  }}
\put(20,0){\framebox(20,20)[r]{  }}
\put(20,20){\framebox(20,20)[r]{  }}
\put(20,40){\framebox(20,20)[r]{  }}
\put(40,0){\framebox(20,20)[r]{  }}
\put(40,20){\framebox(20,20)[r]{  }}
\put(40,40){\framebox(20,20)[r]{  }}
\put(10,10){\line(1,0){20}}
\put(30,10){\line(0,1){20}}
\put(50,10){\line(0,1){20}}
\put(10,30){\line(0,1){20}}
\put(30,30){\line(0,1){20}}
\put(49,49){.}
\end{picture}
\]
contains a singularity of degree 2 of type $\type{0}$.

If the singularity of degree $r$ in a $(r+1) \times (r+1)$ square is of type $\typeij{i}$, we denote the square by $\typeij{i}^r$. We omit $r$ when $r=1$.

We say $i$ union $j$, and write 
\[\begin{picture}(40,16)
\put(8,8){\circle{16}}%cicle 3
\put(32,8){\circle{16}}%circle 4
\put(16,8){\line(1,0){8}}%3 to 4
\put(16,14){$r$}%circle 3
\put(40,14){$r$}%circle 4
\put(7,5){i}
\put(30,5){j}
\end{picture}\]
to mean the second column of $\typeij{i}^r$ superimposed onto the first column of $\typeij{j}^r$, and write
\[\begin{picture}(16,40)
\put(8,32){\circle{16}}%circle 1
\put(8,8){\circle{16}}%cicle 3
\put(8,16){\line(0,1){8}}%1 to 3
\put(16,38){$r$}%circle 1
\put(16,14){$r$}%circle 3
\put(7,29){i}
\put(6,5){j}
\end{picture}\] to mean the second row of $\typeij{i}^r$ superimposed onto the first row of $\typeij{j}^r$.

So the $3\times3$ square above may be denoted $\type{0}^2$ for singularities of degree 2, and
\[ \begin{picture}(40,40)
\put(8,32){\circle{16}}%circle 1
\put(32,32){\circle{16}}%circle 2
\put(8,8){\circle{16}}%circle 3
\put(32,8){\circle{16}}%circle 4
\put(16,32){\line(1,0){8}}%1 to 2
\put(16,8){\line(1,0){8}}%3 to 4
\put(8,16){\line(0,1){8}}%1 to 3
\put(32,16){\line(0,1){8}}%2 to 4
\put(6,29){0}
\put(30,29){1}
\put(6,5){6}
\put(30,5){0}
\end{picture} \]
for singularities of degree 1.

\begin{definition}\label{validdiagrams}\emph{\textbf{: Valid Diagrams.}}\end{definition}
A diagram is a tableau $\Lambda$ in which we connect the boxes as described above. If $i$ is connected to $j$ in the diagram $\Lambda$ then we let $z_i =z_j$, and set $z_i-z_j \notin \mathbb{Z}$ otherwise.

It is our intention to re-order the product (\ref{bigf}) to form regular triples at $z_1 = \cdots = z_n$. Since they do not allow the formation of regular triples we exclude $2 \times 2$ squares of type $\bigtypecentre{8}$ to $\bigtype{15}$ in our diagrams. For the same reason we also exclude the unions

\[\begin{picture}(40,40)
\put(8,32){\circle{16}}%circle 1
\put(32,32){\circle{16}}%circle 2
\put(16,32){\line(1,0){8}}%1 to 2
\put(16,38){$r$}%circle 1
\put(40,38){$r$}%circle 2
\put(6,29){3}
\put(30,29){5}
\end{picture}
\qquad \qquad
\begin{picture}(16,40)
\put(8,32){\circle{16}}%circle 1
\put(8,8){\circle{16}}%cicle 3
\put(8,16){\line(0,1){8}}%1 to 3
\put(16,38){$r$}%circle 1
\put(16,14){$r$}%circle 3
\put(6,29){1}
\put(6,5){7}
\end{picture}
\qquad \qquad
\begin{picture}(40,40)
\put(8,32){\circle{16}}%circle 1
\put(32,32){\circle{16}}%circle 2
\put(8,8){\circle{16}}%cicle 3
\put(16,32){\line(1,0){8}}%1 to 2
\put(8,16){\line(0,1){8}}%1 to 3
\put(16,38){$r$}%circle 1
\put(40,38){$r$}%circle 2
\put(16,14){$r$}%circle 3
\put(6,29){2}
\put(30,29){5}
\put(6,5){7}
\end{picture}
\qquad \qquad
\begin{picture}(40,40)
\put(32,32){\circle{16}}%circle 2
\put(8,8){\circle{16}}%cicle 3
\put(32,8){\circle{16}}%circle 4
\put(16,8){\line(1,0){8}}%3 to 4
\put(32,16){\line(0,1){8}}%2 to 4
\put(40,38){$r$}%circle 2
\put(16,14){$r$}%circle 3
\put(40,14){$r$}%circle 4
\put(30,29){1}
\put(6,5){3}
\put(30,5){6}
\end{picture},\]

%\newpage
including these unions with any number of
\[\begin{picture}(16,32)
\put(8,16){\circle{16}}%circle 1
\put(8,0){\line(0,1){8}}%1 to 3
\put(8,24){\line(0,1){8}}
\put(16,22){$r$}%circle 1
\put(6,13){0}
\end{picture}
\qquad \qquad 
\begin{picture}(32,40)
\put(16,16){\circle{16}}%circle 1
\put(24,16){\line(1,0){8}}%1 to 2
\put(0,16){\line(1,0){8}}
\put(24,22){$r$}%circle 1
\put(14,13){4}
\end{picture}
\qquad \qquad
\begin{picture}(40,40)
\put(8,16){\circle{16}}%circle 1
\put(16,16){\line(1,0){8}}%1 to 2
\put(8,0){\line(0,1){8}}%1 to 3
\put(16,22){$r$}%circle 1
\put(6,13){2}
\end{picture}
\qquad \qquad
\begin{picture}(40,40)
\put(32,8){\circle{16}}%circle 4
\put(16,8){\line(1,0){8}}%3 to 4
\put(32,16){\line(0,1){8}}%2 to 4
\put(40,14){$r$}%circle 4
\put(30,5){6}
\end{picture}
\]
between the two ends of the chain, for all $r \geqslant 1$, where such unions are possible.

The direction of the connecting lines is important. For example,
\[\begin{picture}(64,40)
\put(56,32){\circle{16}}
\put(32,32){\circle{16}}%circle 2
\put(8,8){\circle{16}}%circle 3
\put(32,8){\circle{16}}%circle 4
\put(16,8){\line(1,0){8}}%3 to 4
\put(32,16){\line(0,1){8}}%2 to 4
\put(40,32){\line(1,0){8}}
\put(54,29){5}
\put(30,29){2}
\put(6,5){3}
\put(30,5){6}
\end{picture}\] is an excluded union, whereas
\[\begin{picture}(88,16)
\put(56,8){\circle{16}}
\put(32,8){\circle{16}}%circle 2
\put(8,8){\circle{16}}%circle 3
\put(80,8){\circle{16}}%circle 4
\put(16,8){\line(1,0){8}}%3 to 4
\put(64,8){\line(1,0){8}}%2 to 4
\put(40,8){\line(1,0){8}}
\put(54,5){2}
\put(30,5){6}
\put(6,5){3}
\put(77,5){5}
\end{picture}\]is not excluded.

If we list the possible unions to the right and below a singularity of type $\typeij{i}$, then the unions of degree 1 are;

\quad

\begin{picture}(44,58)
\put(46,46){$\left\{ \begin{array}{c} \type{0}, \type{1}\phantom{, \type{x}} \\ \type{2}\phantom{, \type{x}} \end{array} \right.$}
\put(0,0){$\overbrace{\begin{array}{c}\type{5}, \type{6}, \type{7} \\ \type{0}, \type{1} \end{array}}^{\phantom{x}}$}
\put(36,50){\line(1,0){12}}
\put(28,30){\line(0,1){12}}
\put(28,50){\circle{16}}
\put(26,47){0}
\end{picture}
\qquad \qquad \qquad \qquad
\begin{picture}(44,58)
\put(46,46){$\left\{ \begin{array}{c} \type{3}, \type{4}, \type{5} \\ \type{6}, \type{7} \end{array} \right.$}
\put(0,0){$\overbrace{\begin{array}{c}\type{5}, \type{6}\phantom{, \type{x}} \\ \type{0}, \type{1} \end{array}}^{\phantom{x}}$}
\put(36,50){\line(1,0){12}}
\put(28,30){\line(0,1){12}}
\put(28,50){\circle{16}}
\put(26,47){1}
\end{picture}
\qquad \qquad \qquad \qquad
\begin{picture}(44,58)
\put(46,46){$\left\{ \begin{array}{c} \type{3}, \type{4}, \type{5} \\ \type{6}, \type{7} \end{array} \right.$}
\put(0,0){$\overbrace{\begin{array}{c}\type{5}, \type{6}, \type{7} \\ \type{0}, \type{1} \end{array}}^{\phantom{x}}$}
\put(36,50){\line(1,0){12}}
\put(28,30){\line(0,1){12}}
\put(28,50){\circle{16}}
\put(26,47){2}
\end{picture}

\quad

\begin{picture}(44,58)
\put(46,46){$\left\{ \begin{array}{c} \type{3}, \type{4}\phantom{, \type{x}} \\ \type{6}, \type{7} \end{array} \right.$}
\put(0,0){$\overbrace{\begin{array}{c}\type{5}, \type{6}, \type{7} \\ \type{0}, \type{1} \end{array}}^{\phantom{x}}$}
\put(36,50){\line(1,0){12}}
\put(28,30){\line(0,1){12}}
\put(28,50){\circle{16}}
\put(26,47){3}
\end{picture}
\qquad \qquad \qquad \qquad
\begin{picture}(44,58)
\put(46,46){$\left\{ \begin{array}{c} \type{3}, \type{4}, \type{5} \\ \type{6}, \type{7} \end{array} \right.$}
\put(0,0){$\overbrace{\begin{array}{c}\type{2}, \type{3}, \type{4} \\ \phantom{\type{x}, \type{x}} \end{array}}^{\phantom{x}}$}
\put(36,50){\line(1,0){12}}
\put(28,30){\line(0,1){12}}
\put(28,50){\circle{16}}
\put(26,47){4}
\end{picture}
\qquad \qquad \qquad \qquad
\begin{picture}(44,58)
\put(46,46){$\left\{ \begin{array}{c} \type{3}, \type{4}, \type{5} \\ \type{6}, \type{7} \end{array} \right.$}
\put(0,0){$\overbrace{\begin{array}{c}\type{2}, \type{3}, \type{4} \\ \phantom{\type{x}, \type{x}} \end{array}}^{\phantom{x}}$}
\put(36,50){\line(1,0){12}}
\put(28,30){\line(0,1){12}}
\put(28,50){\circle{16}}
\put(26,47){5}
\end{picture}

\quad

\begin{picture}(44,58)
\put(46,46){$\left\{ \begin{array}{c} \type{0}, \type{1}\phantom{, \type{x}} \\ \type{2}\phantom{, \type{x}} \end{array} \right.$}
\put(0,0){$\overbrace{\begin{array}{c}\type{2}, \type{3}, \type{4} \\ \phantom{\type{x}, \type{x}} \end{array}}^{\phantom{x}}$}
\put(36,50){\line(1,0){12}}
\put(28,30){\line(0,1){12}}
\put(28,50){\circle{16}}
\put(26,47){6}
\end{picture}
\qquad \qquad \qquad \qquad
\begin{picture}(44,58)
\put(46,46){$\left\{ \begin{array}{c} \type{0}, \type{1}\phantom{, \type{x}} \\ \type{2}\phantom{, \type{x}} \end{array} \right.$}
\put(0,0){$\overbrace{\begin{array}{c}\type{5}, \type{6}, \type{7} \\ \type{0}, \type{1} \end{array}}^{\phantom{x}}$}
\put(36,50){\line(1,0){12}}
\put(28,30){\line(0,1){12}}
\put(28,50){\circle{16}}
\put(26,47){7}
\end{picture}

\quad

Similarly, a $2\times2$ square of type $\type{1}$ may not be diagonally below a $2\times2$ square of type $\type{7}$, and a $2\times2$ square of type $\type{3}$ may not be diagonally below a $2\times2$ square of type $\type{5}$. The following shows which types may be immediately diagonally below a $2\times2$ square of type $\typeij{i}$.

\quad

\begin{picture}(24,58)
\put(23,28){$\{ \type{5}, \type{6}, \type{7}, \type{0}, \type{1} \}$}
\put(15,46){\line(1,-1){10}}
\put(8,50){\circle{16}}
\put(6,47){0}
\end{picture}
\qquad \qquad \qquad \qquad \qquad
\begin{picture}(24,58)
\put(23,28){$\{ \type{0}, \type{1}, \type{2}, \type{3}, \type{4}, \type{5}, \type{6}, \type{7} \}$}
\put(15,46){\line(1,-1){10}}
\put(8,50){\circle{16}}
\put(6,47){1}
\end{picture}

\begin{picture}(24,58)
\put(23,28){$\{ \type{0}, \type{1}, \type{2}, \type{3}, \type{4}, \type{5}, \type{6}, \type{7} \}$}
\put(15,46){\line(1,-1){10}}
\put(8,50){\circle{16}}
\put(6,47){2}
\end{picture}
\qquad \qquad \qquad \qquad \qquad \qquad \qquad
\begin{picture}(24,58)
\put(23,28){$\{ \type{0}, \type{1}, \type{2}, \type{3}, \type{4}, \type{5}, \type{6}, \type{7} \}$}
\put(15,46){\line(1,-1){10}}
\put(8,50){\circle{16}}
\put(6,47){3}
\end{picture}

\begin{picture}(24,58)
\put(23,28){$\{ \type{3}, \type{4}, \type{5}, \type{6}, \type{7} \}$}
\put(15,46){\line(1,-1){10}}
\put(8,50){\circle{16}}
\put(6,47){4}
\end{picture}
\qquad \qquad \qquad \qquad \qquad
\begin{picture}(24,58)
\put(23,28){$\{ \type{4}, \type{5}, \type{6}, \type{7} \}$}
\put(15,46){\line(1,-1){10}}
\put(8,50){\circle{16}}
\put(6,47){5}
\end{picture}
\qquad \qquad \qquad \qquad \qquad
\begin{picture}(24,58)
\put(23,28){$\{ \type{5}, \type{6}, \type{7} \}$}
\put(15,46){\line(1,-1){10}}
\put(8,50){\circle{16}}
\put(6,47){6}
\end{picture}

\begin{picture}(24,58)
\put(23,28){$\{ \type{5}, \type{6}, \type{7}, \type{0} \}$}
\put(15,46){\line(1,-1){10}}
\put(8,50){\circle{16}}
\put(6,47){7}
\end{picture}

By considering diagonal chains we can see that, under the above conditions, we can never have a singularity of type $\bigtype{15}$, of any degree.

\newpage
Note, for a tableau $\Lambda$ of shape $\lambda = (\lambda_1, \dots, \lambda_k)$, we may connect the boxes $\Lambda(1, b)$ with $b > \lambda_2$ in the first row any way we like without consequence. Similarly we may connect the boxes $\Lambda(a,1)$ with $a> \lambda'_2$ in the first column in any way we wish.

We call a diagram \emph{invalid} if;
\begin{listr}
	\item It contains any $2\times2$ square of type $\bigtypecentre{8}$ to $\bigtype{15}$.
%	\item It breaks the rules of union, for example \type{4}--\type{1}.
	\item It contains excluded unions such as \bigtypecentre{3}--\bigtypecentre{5}.
	\item It contains a singularity of type $\bigtype{15}$, of any degree.
\end{listr}
We call all other diagrams \emph{valid}.

\quad

We may now state our conjecture for a possible ribbon fusion procedure. Let $\mathcal{P}_\Lambda$ be the vector subspace of $\mathbb{C}^n$ consisting of all tuples $(z_1, \dots, z_n)$ such that $z_i = z_j$ if $i$ is connected to $j$ in the diagram $\Lambda$, and $z_i-z_j \notin \mathbb{Z}$ otherwise.

\begin{conjecture}\label{ribbonconjecture}The value of $F_\Lambda(z_1, \dots, z_n)$ at $z_1 = \cdots = z_n$, after restriction to $\mathcal{P}_\Lambda$, is the diagonal matrix element $F_\Lambda \in \mathbb{C}S_n$ if and only if $\Lambda$ is a valid diagram.\end{conjecture}

Let us look at how valid diagrams can imply the element $F_\Lambda(z_1, \dots, z_n)$ is the diagonal matrix element. The valid diagrams are designed to be those diagrams  such that regular triples are possible at $z_1 = \cdots = z_n$. Consider the column tableau of the partition $(2,2)$ and then the following two diagrams;

\[ \begin{picture}(40,40)
\put(-25,25){$A =$}
\put(0,0){\framebox(20,20)[r]{$2\;\;$}}
\put(0,20){\framebox(20,20)[r]{$1\;\;$}}
\put(20,0){\framebox(20,20)[r]{$4\;\;$}}
\put(20,20){\framebox(20,20)[r]{$3\;\;$}}
\put(10,10){\line(0,1){20}}
\put(30,10){\line(0,1){20}}
\end{picture}
\qquad \qquad \qquad
\begin{picture}(40,40)
\put(-25,25){$B =$}
\put(0,0){\framebox(20,20)[r]{$2\;\;$}}
\put(0,20){\framebox(20,20)[r]{$1\;\;$}}
\put(20,0){\framebox(20,20)[r]{$4\;\;$}}
\put(20,20){\framebox(20,20)[r]{$3\;\;$}}
%\put(10,10){\line(0,1){20}}
\put(30,10){\line(0,1){20}}
\end{picture}\]

The ordering $f_{12}(f_{13}f_{14})f_{23}f_{24}f_{34}$, with $\frac{1}{2} f_{34}$ as an idempotent, shows regularity at $z_1 = \cdots z_n$ for both $F_A(z_1, \dots, z_n)$ and $F_B(z_1,\dots, z_n)$. In full we have;

\newpage
\begin{eqnarray*}F_A(z_1, \dots, z_n)= &\left(1 - (1\phantom{x}2)\right)\left[\left(1 - \frac{(1\phantom{x}3)}{z_1 -z_3-1}\right)\left(1 - \frac{(1\phantom{x}4)}{z_1 -z_3}\right)\frac{1}{2}\left(1 - (3\phantom{x}4)\right)\right]& \\ &\times \left(1 - \frac{(2\phantom{x}3)}{z_1 -z_3-2}\right)\left(1 - \frac{(2\phantom{x}4)}{z_1 -z_3-1}\right)\left(1 - (3\phantom{x}4)\right)&\end{eqnarray*} 
and
\begin{eqnarray*}F_B(z_1, \dots, z_n)= &\left(1 - \frac{(1\phantom{x}2)}{z_1-z_2+1}\right)\left[\left(1 - \frac{(1\phantom{x}3)}{z_1 -z_3-1}\right)\left(1 - \frac{(1\phantom{x}4)}{z_1 -z_3}\right)\frac{1}{2}\left(1 - (3\phantom{x}4)\right)\right]& \\ &\times \left(1 - \frac{(2\phantom{x}3)}{z_2 -z_3-2}\right)\left(1 - \frac{(2\phantom{x}4)}{z_2 -z_3-1}\right)\left(1 - (3\phantom{x}4)\right)&\end{eqnarray*}

Since we use the same tableau in both diagrams, they both have the same singularity. In both cases that singularity is resolved in exactly the same way, with a shared idempotent $\frac{1}{2}f_{34}$ and a shared triple term $f_{13}$. So the value of the triple is the same at $z_1 = \cdots = z_n$ in both cases. All other terms in the two products also have the same value at $z_1 = \cdots z_n$. So by comparing the two products, first comparing triples and then term by term, we see that $F_A(z_1, \dots, z_n)$ has the same value as $F_B(z_1, \dots, z_n)$ at $z_1 =\cdots =z_n$. We know, by the column fusion procedure, that $F_A(z_1, \dots, z_n)$ is the diagonal matrix element at $z_1= \cdots = z_n$, and hence so is $F_B(z_1, \dots z_n)$.

We cannot compare the diagrams 
\[ \begin{picture}(40,40)
\put(-25,25){$A =$}
\put(0,0){\framebox(20,20)[r]{$2\;\;$}}
\put(0,20){\framebox(20,20)[r]{$1\;\;$}}
\put(20,0){\framebox(20,20)[r]{$4\;\;$}}
\put(20,20){\framebox(20,20)[r]{$3\;\;$}}
\put(10,10){\line(0,1){20}}
\put(30,10){\line(0,1){20}}
\end{picture}
\qquad \quad \textrm{and} \quad \qquad
\begin{picture}(40,40)
\put(-25,25){$D =$}
\put(0,0){\framebox(20,20)[r]{$2\;\;$}}
\put(0,20){\framebox(20,20)[r]{$1\;\;$}}
\put(20,0){\framebox(20,20)[r]{$4\;\;$}}
\put(20,20){\framebox(20,20)[r]{$3\;\;$}}
\put(10,10){\line(1,0){20}}
\put(10,30){\line(1,0){20}}
\end{picture}\]
directly in this way. Instead we need an intermediate step such as
\[\begin{picture}(40,40)
\put(-25,25){$C =$}
\put(0,0){\framebox(20,20)[r]{$2\;\;$}}
\put(0,20){\framebox(20,20)[r]{$1\;\;$}}
\put(20,0){\framebox(20,20)[r]{$4\;\;$}}
\put(20,20){\framebox(20,20)[r]{$3\;\;$}}
\put(10,10){\line(0,1){20}}
\put(10,30){\line(1,0){20}}
\end{picture}\]
We may then form a chain of diagrams comparing $A$ and $C$ using $\frac{1}{2}f_{12}$ as the shared idempotent and $f_{24}$ as the shared triple term. Then rearranging the product of $C$ using the Yang-Baxter relations (\ref{triple}) and (\ref{commute}), which reorder the terms without changing the value of the product, such that we may then compare $C$ with $D$ using $\frac{1}{2}f_{13}$ as the shared idempotent and $f_{34}$ as the shared triple term.

This example illustrates our method of proof. We break the proof into two statements. The first states that the corresponding product $F_\Lambda(z_1, \dots z_n)$ of every valid diagram $\Lambda$ has at least one ordering such that each singularity may be resolved by regular triples. While the second statement shows there exists a chain of diagrams, ending with the diagram used in the column fusion procedure, such that we may compare each diagram in the chain to the one before in the way described above.

Given this plan of attack, we give the first statement formally as the following conjecture.

\begin{conjecture}\label{jimtheorem4}For every valid diagram $\Lambda$ there exists at least one ordering of $F_\Lambda(z_1, .., z_n)$ such that, after restriction to $\mathcal{P}_\Lambda$, we may resolve all singularities by regular triples at $z_1 = \cdots = z_n$. \end{conjecture}

We assert this to be a reasonable conjecture, which then has the immediate corollary;

\begin{corollary}For every valid diagram $\Lambda$, $F_\Lambda(z_1, .., z_n)$ is regular at $z_1 = \cdots = z_n$, after restriction to $\mathcal{P}_\Lambda$.\end{corollary}

Given Conjecture \ref{jimtheorem4} we now describe a chain of diagrams which allow us to compare any valid diagram to the diagram used in the column fusion procedure which is made entirely of $2\times2$ squares of type $\type{0}$.

\begin{proposition}\label{jimtheorem5}Let $\tilde{\Lambda}$ be the diagram made entirely of $2\times2$ squares of type \emph{$\type{0}$}. For any valid diagram $\Lambda$ there exist a chain of valid diagrams, $\Lambda = \Lambda^\circ, \Lambda^1, \Lambda^2, \dots, \Lambda^N = \tilde{\Lambda}$ such that, any two adjacent diagrams, $F_{\Lambda^k}(z_1, \dots, z_n)$ and $F_{\Lambda^{k+1}}(z_1, \dots, z_n)$ share an ordering which resolves all singularities in both products by regular triples at $z_1 = \cdots = z_n$.\end{proposition}

\begin{proof}
Let $u$ stand next to $v$ in the same row or column of $\Lambda$. We call a line connecting two such adjacent boxes an idempotent, since $\frac{1}{2}f_{uv}(c_u(\Lambda), c_v(\Lambda))$ is an idempotent in $\mathbb{C}S_n$. One diagram may be turned into another by addition or subtraction of an idempotent. Then, given Conjecture \ref{jimtheorem4}, the diagram with the fewest idempotents has an associated product with at least one ordering which resolves all singularities by regular triples at $z_1 = \cdots =z_n$. Any such ordering will also resolve the singularities of the product associated to the diagram with the additional idempotent. Note, this would not be the case if the additional idempotent was in the place of a triple term, but this can not happen since such an idempotent would connect the boxes of a singularity making an invalid diagram.

Let $j = i \pm 1 \textrm{ mod }8$. Then we may turn a $2\times2$ square of type $\typeij{i}$, $0\leqslant i \leqslant 7$, into a square of type $\typeij{j}$ by adding or removing an idempotent. The following figure shows which types may be turned into which other types\footnote{This figure is suggestive of modular arithmetic, mod 8. Indeed $\smalltype{-i}$ is the square $\smalltypecentre{i}$ rotated $180^\circ$. Also, if $j-i = 4\textrm{ mod }8$ then $\smalltypecentre{j}$ is opposite $\smalltypecentre{i}$ on the circle and is a reflection of $\smalltypecentre{i}$ in the line $y=x$.}.

\begin{figure}[h]
	\centering
		\includegraphics[width=5cm]{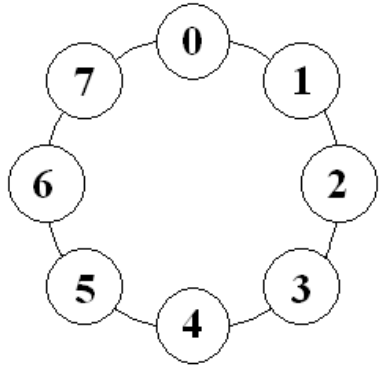}
	\caption{Changing types by adding and removing idempotents.}
	\label{fig:circle}
\end{figure}

We now provide an algorithm to turn any valid diagram into the diagram made entirely of $2\times2$ squares of type $\type{0}$. We do this by using Figure \ref{fig:circle} and moving away from $2\times2$ squares of type $\type{4}$ and towards $2\times2$ squares of type $\type{0}$. This will give us a valid diagram except in the following cases. In these cases we must perform the prescribed actions first. If this causes a loop of logic we perform all actions simultaneously. 

\newpage
\quad \textbf{Turning $\type{1}$ into $\type{0}$ is fine except when we have the following unions:}

\begin{list}{}{} 

\item \begin{picture}(80,16)
\put(8,8){\circle{16}}%circle 1
\put(32,8){\circle{16}}%circle 2
\put(16,8){\line(1,0){8}}%1 to 2
\put(6,5){1}
\put(30,5){3}
\put(74, 3){In which case first change $\type{3}  \mapsto \type{2}$}
\end{picture}

\item \begin{picture}(80,20)
\put(8,8){\circle{16}}%circle 1
\put(32,8){\circle{16}}%circle 2
\put(16,8){\line(1,0){8}}%1 to 2
\put(6,5){1}
\put(30,5){4}
\put(74, 3){First change $\type{4}  \mapsto \type{3}$ or $\type{4} \mapsto \type{5}$}
\end{picture}

\item \begin{picture}(80,20)
\put(8,8){\circle{16}}%circle 1
\put(32,8){\circle{16}}%circle 2
\put(16,8){\line(1,0){8}}%1 to 2
\put(6,5){1}
\put(30,5){5}
\put(74, 3){First change $\type{5}  \mapsto \type{6}$}
\end{picture}

\item \begin{picture}(80,20)
\put(8,8){\circle{16}}%circle 1
\put(32,8){\circle{16}}%circle 2
\put(16,8){\line(1,0){8}}%1 to 2
\put(6,5){1}
\put(30,5){6}
\put(74, 3){First change $\type{6}  \mapsto \type{7}$}
\end{picture}

\item \begin{picture}(80,20)
\put(8,8){\circle{16}}%circle 1
\put(32,8){\circle{16}}%circle 2
\put(16,8){\line(1,0){8}}%1 to 2
\put(6,5){1}
\put(30,5){7}
\put(74, 3){First change $\type{7}  \mapsto \type{0}$}
\end{picture}
\end{list}

%\newpage
\quad \textbf{Turning $\type{2}^r$ into $\type{1}^r$ is fine except when we have the following unions:}

\begin{list}{}{} 

\item \begin{picture}(80,40)
\put(32,8){\circle{16}}%circle 1
\put(32,32){\circle{16}}%circle 2
\put(32,16){\line(0,1){8}}%1 to 2
\put(30,5){2}
\put(30,30){4}
\put(74, 18){In which case first change $\type{4}  \mapsto \type{3}$ or $\type{4} \mapsto \type{5}$}
\end{picture}

\item \begin{picture}(80,40)
\put(32,8){\circle{16}}%circle 1
\put(32,32){\circle{16}}%circle 2
\put(32,16){\line(0,1){8}}%1 to 2
\put(30,5){2}
\put(30,30){5}
\put(74, 18){First change $\type{5}  \mapsto \type{6}$}
\end{picture}

\item \begin{picture}(80,40)
\put(32,8){\circle{16}}%circle 1
\put(32,32){\circle{16}}%circle 2
\put(32,16){\line(0,1){8}}%1 to 2
\put(30,5){2}
\put(30,30){6}
\put(74, 18){First change $\type{6}  \mapsto \type{7}$}
\end{picture}

\item \begin{picture}(80,40)
\put(32,8){\circle{16}}%circle 1
\put(32,32){\circle{16}}%circle 2
\put(32,16){\line(0,1){8}}%1 to 2
\put(30,5){7}
\put(30,30){2}
\put(74, 18){First change $\type{7}^r  \mapsto \type{0}^r$, for $r \geqslant 1$}
\put(40,38){$r$}%circle 1
\put(40,14){$r$}%circle 3
\put(-10,14){*}
\end{picture}

\item \begin{picture}(80,40)
\put(32,8){\circle{16}}%circle 1
\put(32,32){\circle{16}}%circle 2
\put(8,8){\circle{16}}
\put(32,16){\line(0,1){8}}%1 to 2
\put(16,8){\line(1,0){8}}
\put(6,5){3}
\put(30,5){6}
\put(30,30){2}
\put(40,38){$r$}%circle 1
\put(40,14){$r$}%circle 3
\put(16,14){$r$}%circle 3
\put(74, 18){First change $\type{5}^r  \mapsto \type{6}^r$, for $r \geqslant 1$}
\put(-10,14){*}
\end{picture}
\end{list}

\quad \textbf{Turning $\type{3}$ into $\type{2}$ is fine except when we have the following unions:}

\begin{list}{}{} 

\item \begin{picture}(80,16)
\put(8,8){\circle{16}}%circle 1
\put(32,8){\circle{16}}%circle 2
\put(16,8){\line(1,0){8}}%1 to 2
\put(6,5){1}
\put(30,5){3}
\put(74, 3){In which case first change $\type{1}  \mapsto \type{0}$}
\end{picture}

\item \begin{picture}(80,20)
\put(8,8){\circle{16}}%circle 1
\put(32,8){\circle{16}}%circle 2
\put(16,8){\line(1,0){8}}%1 to 2
\put(6,5){2}
\put(30,5){3}
\put(74, 3){First change $\type{2}  \mapsto \type{1}$}
\end{picture}

\item \begin{picture}(80,20)
\put(8,8){\circle{16}}%circle 1
\put(32,8){\circle{16}}%circle 2
\put(16,8){\line(1,0){8}}%1 to 2
\put(6,5){3}
\put(30,5){3}
\put(16,10){'}
\put(74, 3){First change $\type{3}'  \mapsto \type{2}$}
\end{picture}

\item \begin{picture}(80,20)
\put(8,8){\circle{16}}%circle 1
\put(32,8){\circle{16}}%circle 2
\put(16,8){\line(1,0){8}}%1 to 2
\put(6,5){4}
\put(30,5){3}
\put(74, 3){First change $\type{4}  \mapsto \type{3}$ or $\type{4} \mapsto \type{5}$}
\end{picture}

\item \begin{picture}(80,20)
\put(8,8){\circle{16}}%circle 1
\put(32,8){\circle{16}}%circle 2
\put(16,8){\line(1,0){8}}%1 to 2
\put(6,5){5}
\put(30,5){3}
\put(74, 3){First change $\type{5}  \mapsto \type{6}$}
\end{picture}
\end{list}

\quad \textbf{Turning $\type{4}^r$ into $\type{3}^r$ is fine except when we have the following unions:}

\begin{list}{}{}

\item \begin{picture}(80,40)
\put(32,8){\circle{16}}%circle 1
\put(32,32){\circle{16}}%circle 2
\put(32,16){\line(0,1){8}}%1 to 2
\put(30,5){2}
\put(30,30){4}
\put(74, 18){First change $\type{2}  \mapsto \type{1}$}
\end{picture}

\item \begin{picture}(80,40)
\put(32,8){\circle{16}}%circle 1
\put(32,32){\circle{16}}%circle 2
\put(32,16){\line(0,1){8}}%1 to 2
\put(30,5){3}
\put(30,30){4}
\put(74, 18){First change $\type{3}  \mapsto \type{2}$}
\end{picture}

\item \begin{picture}(80,40)
\put(32,8){\circle{16}}%circle 1
\put(32,32){\circle{16}}%circle 2
\put(32,16){\line(0,1){8}}%1 to 2
\put(30,5){4}
\put(30,30){4}
\put(40,10){'}
\put(74, 18){In which case first change $\type{4}'  \mapsto \type{3}$ or $\type{4}' \mapsto \type{5}$}
\end{picture}

\item \begin{picture}(80,40)
\put(32,8){\circle{16}}%circle 1
\put(8,32){\circle{16}}%circle 2
\put(16,16){\tiny $\ddots$}
\put(30,5){4}
\put(6,30){5}
\put(74, 18){i.e. If $\type{4}$ is diagonally below \type{5} then first change $\type{5}  \mapsto \type{6}$}
\end{picture}

\item \begin{picture}(80,20)
\put(8,8){\circle{16}}%circle 1
\put(32,8){\circle{16}}%circle 2
\put(16,8){\line(1,0){8}}%1 to 2
\put(6,5){4}
\put(30,5){5}
\put(74, 3){First change $\type{5}^r  \mapsto \type{6}^r$, for all $r \geqslant 1$}
\put(16,14){$r$}%circle 1
\put(40,14){$r$}%circle 3
\put(-10,3){*}
\end{picture}

\item \begin{picture}(80,40)
\put(32,8){\circle{16}}%circle 1
\put(32,32){\circle{16}}%circle 2
\put(8,8){\circle{16}}
\put(32,16){\line(0,1){8}}%1 to 2
\put(16,8){\line(1,0){8}}
\put(6,5){4}
\put(30,5){6}
\put(30,30){1}
\put(40,38){$r$}%circle 1
\put(40,14){$r$}%circle 3
\put(16,14){$r$}%circle 3
\put(74, 18){First change $\type{1}^r  \mapsto \type{0}^r$, for all $r \geqslant 1$}
\put(-10,14){*}
\end{picture}
\end{list}

\quad \textbf{Turning $\type{4}^r$ into $\type{5}^r$ is fine except when we have the following unions:}

\begin{list}{}{} 

\item \begin{picture}(80,40)
\put(32,8){\circle{16}}%circle 1
\put(32,32){\circle{16}}%circle 2
\put(32,16){\line(0,1){8}}%1 to 2
\put(30,5){4}
\put(30,30){4}
\put(40,34){'}
\put(74, 18){In which case first change $\type{4}'  \mapsto \type{3}$ or $\type{4}' \mapsto \type{5}$}
\end{picture}

\item \begin{picture}(80,40)
\put(32,8){\circle{16}}%circle 1
\put(32,32){\circle{16}}%circle 2
\put(32,16){\line(0,1){8}}%1 to 2
\put(30,5){4}
\put(30,30){5}
\put(74, 18){First change $\type{5}  \mapsto \type{6}$}
\end{picture}

\item \begin{picture}(80,40)
\put(32,8){\circle{16}}%circle 1
\put(32,32){\circle{16}}%circle 2
\put(32,16){\line(0,1){8}}%1 to 2
\put(30,5){4}
\put(30,30){6}
\put(74, 18){First change $\type{6}  \mapsto \type{7}$}
\end{picture}

\item \begin{picture}(80,40)
\put(32,8){\circle{16}}%circle 1
\put(8,32){\circle{16}}%circle 2
\put(16,16){\tiny $\ddots$}
\put(30,5){3}
\put(6,30){4}
\put(74, 18){i.e. If $\type{3}$ is diagonally below \type{4} then first change $\type{3}  \mapsto \type{2}$}
\end{picture}

\item \begin{picture}(80,20)
\put(8,8){\circle{16}}%circle 1
\put(32,8){\circle{16}}%circle 2
\put(16,8){\line(1,0){8}}%1 to 2
\put(6,5){3}
\put(30,5){4}
\put(74, 3){First change $\type{3}^r  \mapsto \type{2}^r$, for all $r \geqslant 1$}
\put(16,14){$r$}%circle 1
\put(40,14){$r$}%circle 3
\put(-10,3){*}
\end{picture}

\item \begin{picture}(80,40)
\put(8,32){\circle{16}}%circle 1
\put(32,32){\circle{16}}%circle 2
\put(8,8){\circle{16}}
\put(8,16){\line(0,1){8}}
\put(16,32){\line(1,0){8}}
\put(6,5){7}
\put(6,30){2}
\put(30,30){4}
\put(40,38){$r$}%circle 1
\put(16,38){$r$}%circle 3
\put(16,14){$r$}%circle 3
\put(74, 18){First change $\type{7}^r  \mapsto \type{0}^r$, for all $r \geqslant 1$}
\put(-10,14){*}
\end{picture}
\end{list}

\quad \textbf{Turning $\type{5}$ into $\type{6}$ is fine except when we have the following unions:}

\begin{list}{}{} 

\item \begin{picture}(80,16)
\put(8,8){\circle{16}}%circle 1
\put(32,8){\circle{16}}%circle 2
\put(16,8){\line(1,0){8}}%1 to 2
\put(6,5){5}
\put(30,5){3}
\put(74, 3){In which case first change $\type{3}  \mapsto \type{2}$}
\end{picture}

\item \begin{picture}(80,20)
\put(8,8){\circle{16}}%circle 1
\put(32,8){\circle{16}}%circle 2
\put(16,8){\line(1,0){8}}%1 to 2
\put(6,5){5}
\put(30,5){4}
\put(74, 3){First change $\type{4}  \mapsto \type{3}$ or $\type{4} \mapsto \type{5}$}
\end{picture}

\item \begin{picture}(80,20)
\put(8,8){\circle{16}}%circle 1
\put(32,8){\circle{16}}%circle 2
\put(16,8){\line(1,0){8}}%1 to 2
\put(6,5){5}
\put(30,5){5}
\put(40,10){'}
\put(74, 3){First change $\type{5}'  \mapsto \type{6}$}
\end{picture}

\item \begin{picture}(80,20)
\put(8,8){\circle{16}}%circle 1
\put(32,8){\circle{16}}%circle 2
\put(16,8){\line(1,0){8}}%1 to 2
\put(6,5){5}
\put(30,5){6}
\put(74, 3){First change $\type{6}  \mapsto \type{7}$}
\end{picture}

\item \begin{picture}(80,20)
\put(8,8){\circle{16}}%circle 1
\put(32,8){\circle{16}}%circle 2
\put(16,8){\line(1,0){8}}%1 to 2
\put(6,5){5}
\put(30,5){7}
\put(74, 3){First change $\type{7}  \mapsto \type{0}$}
\end{picture}
\end{list}

\quad \textbf{Turning $\type{6}^r$ into $\type{7}^r$ is fine except when we have the following unions:}

\begin{list}{}{} 

\item \begin{picture}(80,40)
\put(32,8){\circle{16}}%circle 1
\put(32,32){\circle{16}}%circle 2
\put(32,16){\line(0,1){8}}%1 to 2
\put(30,5){2}
\put(30,30){6}
\put(74, 18){First change $\type{2}  \mapsto \type{1}$}
\end{picture}

\item \begin{picture}(80,40)
\put(32,8){\circle{16}}%circle 1
\put(32,32){\circle{16}}%circle 2
\put(32,16){\line(0,1){8}}%1 to 2
\put(30,5){3}
\put(30,30){6}
\put(74, 18){First change $\type{3}  \mapsto \type{2}$}
\end{picture}

\item \begin{picture}(80,40)
\put(32,8){\circle{16}}%circle 1
\put(32,32){\circle{16}}%circle 2
\put(32,16){\line(0,1){8}}%1 to 2
\put(30,5){4}
\put(30,30){6}
\put(74, 18){In which case first change $\type{4}  \mapsto \type{3}$ or $\type{4} \mapsto \type{5}$}
\end{picture}

\item \begin{picture}(80,20)
\put(8,8){\circle{16}}%circle 1
\put(32,8){\circle{16}}%circle 2
\put(16,8){\line(1,0){8}}%1 to 2
\put(6,5){6}
\put(30,5){1}
\put(74, 3){First change $\type{1}^r  \mapsto \type{0}^r$, for all $r \geqslant 1$}
\put(16,14){$r$}%circle 1
\put(40,14){$r$}%circle 3
\put(-10,3){*}
\end{picture}

\item \begin{picture}(80,40)
\put(8,32){\circle{16}}%circle 1
\put(32,32){\circle{16}}%circle 2
\put(8,8){\circle{16}}
\put(8,16){\line(0,1){8}}
\put(16,32){\line(1,0){8}}
\put(6,5){6}
\put(6,30){2}
\put(30,30){5}
\put(40,38){$r$}%circle 1
\put(16,38){$r$}%circle 3
\put(16,14){$r$}%circle 3
\put(74, 18){First change $\type{8}^r  \mapsto \type{6}^r$, for all $r \geqslant 1$}
\put(-10,14){*}
\end{picture}
\end{list}

\quad \textbf{Turning $\type{7}$ into $\type{0}$ is fine except when we have the following unions:}

\begin{list}{}{} 

\item \begin{picture}(80,16)
\put(8,8){\circle{16}}%circle 1
\put(32,8){\circle{16}}%circle 2
\put(16,8){\line(1,0){8}}%1 to 2
\put(6,5){1}
\put(30,5){7}
\put(74, 3){In which case first change $\type{1}  \mapsto \type{0}$}
\end{picture}

\item \begin{picture}(80,20)
\put(8,8){\circle{16}}%circle 1
\put(32,8){\circle{16}}%circle 2
\put(16,8){\line(1,0){8}}%1 to 2
\put(6,5){2}
\put(30,5){7}
\put(74, 3){First change $\type{2}  \mapsto \type{1}$}
\end{picture}

\item \begin{picture}(80,20)
\put(8,8){\circle{16}}%circle 1
\put(32,8){\circle{16}}%circle 2
\put(16,8){\line(1,0){8}}%1 to 2
\put(6,5){3}
\put(30,5){7}
\put(74, 3){First change $\type{3}  \mapsto \type{2}$}
\end{picture}

\item \begin{picture}(80,20)
\put(8,8){\circle{16}}%circle 1
\put(32,8){\circle{16}}%circle 2
\put(16,8){\line(1,0){8}}%1 to 2
\put(6,5){4}
\put(30,5){7}
\put(74, 3){First change $\type{4}  \mapsto \type{3}$ or $\type{4} \mapsto \type{5}$}
\end{picture}

\item \begin{picture}(80,20)
\put(8,8){\circle{16}}%circle 1
\put(32,8){\circle{16}}%circle 2
\put(16,8){\line(1,0){8}}%1 to 2
\put(6,5){5}
\put(30,5){7}
\put(74, 3){First change $\type{5}  \mapsto \type{6}$}
\end{picture}
\end{list}

* with any number of
\[\begin{picture}(16,32)
\put(8,16){\circle{16}}%circle 1
\put(8,0){\line(0,1){8}}%1 to 3
\put(8,24){\line(0,1){8}}
\put(16,22){$r$}%circle 1
\put(6,13){0}
\end{picture}
\qquad \qquad 
\begin{picture}(32,40)
\put(16,16){\circle{16}}%circle 1
\put(24,16){\line(1,0){8}}%1 to 2
\put(0,16){\line(1,0){8}}
\put(24,22){$r$}%circle 1
\put(14,13){4}
\end{picture}
\qquad \qquad
\begin{picture}(40,40)
\put(8,16){\circle{16}}%circle 1
\put(16,16){\line(1,0){8}}%1 to 2
\put(8,0){\line(0,1){8}}%1 to 3
\put(16,22){$r$}%circle 1
\put(6,13){2}
\end{picture}
\qquad \qquad
\begin{picture}(40,40)
\put(32,8){\circle{16}}%circle 4
\put(16,8){\line(1,0){8}}%3 to 4
\put(32,16){\line(0,1){8}}%2 to 4
\put(40,14){$r$}%circle 4
\put(30,5){6}
\end{picture}
\]
between the two ends of the chain, for all $r \geqslant 1$.

These actions prevent us from making an invalid diagram from a diagram that was previously valid. We often have a choice of changing a square of type $\type{4}$ to either one of type $\type{3}$ or type $\type{5}$. For types of odd number it is enough to consider unions of degree 1, for types of even number we must also consider singularities of higher degree.
\end{proof}

Thanks to the above proposition, given any valid starting point, we can now form a chain of diagrams in which we know, for any two adjacent diagrams, $F_{\Lambda^k}(z_1, \dots, z_n)$ and $F_{\Lambda^{k+1}}(z_1, \dots, z_n)$ have the same value in $\mathbb{C}S_n$ at $z_1 = \cdots = z_n$. This chain ends with the diagram used in the column fusion procedure, for which the associated product is known to be the diagonal matrix element at $z_1 = \cdots = z_n$. Hence;

\begin{corollary}For any valid diagram $\Lambda$, the value of $F_\Lambda(z_1, \dots, z_n)$ at $z_1 = \cdots = z_n$, after restriction to $\mathcal{P}_\Lambda$, is the diagonal matrix element $F_\Lambda \in \mathbb{C}S_n$.\end{corollary}

Let us give an example of this algorithm in practice and how it implies the product associated to a particular diagram defines the diagonal matrix element.

\begin{example}
Let $\lambda = (4,4,4,4)$. Note $\lambda$ could be any partition, not necessarily square. Fix a standard tableau of shape $\lambda$, say the column tableau $\Lambda^c$. Consider the following valid diagram;

\[ \begin{picture}(80,80)
\put(0,60){\framebox(20,20)[r]{$1\;\;$}}
\put(0,40){\framebox(20,20)[r]{$2\;\;$}}
\put(0,20){\framebox(20,20)[r]{$3\;\;$}}
\put(0,0){\framebox(20,20)[r]{$4\;\;$}}
\put(20,60){\framebox(20,20)[r]{$5\;\;$}}
\put(20,40){\framebox(20,20)[r]{$6\;\;$}}
\put(20,20){\framebox(20,20)[r]{$7\;\;$}}
\put(20,0){\framebox(20,20)[r]{$8\;\;$}}
\put(40,60){\framebox(20,20)[r]{$9\;\;$}}
\put(40,40){\framebox(20,20)[r]{$10\;$}}
\put(40,20){\framebox(20,20)[r]{$11\;$}}
\put(40,0){\framebox(20,20)[r]{$12\;$}}
\put(60,60){\framebox(20,20)[r]{$13\;$}}
\put(60,40){\framebox(20,20)[r]{$14\;$}}
\put(60,20){\framebox(20,20)[r]{$15\;$}}
\put(60,0){\framebox(20,20)[r]{$16\;$}}

%vertical lines up
%\put(10,50){\line(0,1){20}}%2
%\put(10,30){\line(0,1){20}}%3
\put(10,10){\line(0,1){20}}%4
\put(30,50){\line(0,1){20}}%6
\put(30,30){\line(0,1){20}}%7
%\put(30,10){\line(0,1){20}}%8
%
%\put(50,50){\line(0,1){20}}%10
\put(50,30){\line(0,1){20}}%11
\put(50,10){\line(0,1){20}}%12
%
%\put(70,50){\line(0,1){20}}%14
%\put(70,30){\line(0,1){20}}%15
\put(70,10){\line(0,1){20}}%16

%Horizontal lines to the right
%\put(10,70){\line(1,0){20}}%1
%\put(10,50){\line(1,0){20}}%2
\put(10,30){\line(1,0){20}}%3
%\put(10,10){\line(1,0){20}}%4
%
%\put(30,70){\line(1,0){20}}%5
%\put(30,50){\line(1,0){20}}%6
%\put(30,30){\line(1,0){20}}%7
\put(30,10){\line(1,0){20}}%8
%
%\put(50,70){\line(1,0){20}}%9
\put(50,50){\line(1,0){20}}%10
%\put(50,30){\line(1,0){20}}%11
%\put(50,10){\line(1,0){20}}%12
\end{picture}\]

Starting with this diagram we then use the algorithm of Proposition \ref{jimtheorem5} to form the following chain of diagrams. We also write the diagrams in terms of their unions of both degree 1 and degree 2.

\[
\begin{picture}(80,80)
%\put(-40,36){$\mapsto$}
\put(0,60){\framebox(20,20)[r]{$1\;\;$}}
\put(0,40){\framebox(20,20)[r]{$2\;\;$}}
\put(0,20){\framebox(20,20)[r]{$3\;\;$}}
\put(0,0){\framebox(20,20)[r]{$4\;\;$}}
\put(20,60){\framebox(20,20)[r]{$5\;\;$}}
\put(20,40){\framebox(20,20)[r]{$6\;\;$}}
\put(20,20){\framebox(20,20)[r]{$7\;\;$}}
\put(20,0){\framebox(20,20)[r]{$8\;\;$}}
\put(40,60){\framebox(20,20)[r]{$9\;\;$}}
\put(40,40){\framebox(20,20)[r]{$10\;$}}
\put(40,20){\framebox(20,20)[r]{$11\;$}}
\put(40,0){\framebox(20,20)[r]{$12\;$}}
\put(60,60){\framebox(20,20)[r]{$13\;$}}
\put(60,40){\framebox(20,20)[r]{$14\;$}}
\put(60,20){\framebox(20,20)[r]{$15\;$}}
\put(60,0){\framebox(20,20)[r]{$16\;$}}

%vertical lines up
%\put(10,50){\line(0,1){20}}%2
%\put(10,30){\line(0,1){20}}%3
\put(10,10){\line(0,1){20}}%4
\put(30,50){\line(0,1){20}}%6
\put(30,30){\line(0,1){20}}%7
%\put(30,10){\line(0,1){20}}%8
%
%\put(50,50){\line(0,1){20}}%10
\put(50,30){\line(0,1){20}}%11
\put(50,10){\line(0,1){20}}%12
%
%\put(70,50){\line(0,1){20}}%14
%\put(70,30){\line(0,1){20}}%15
\put(70,10){\line(0,1){20}}%16

%Horizontal lines to the right
%\put(10,70){\line(1,0){20}}%1
%\put(10,50){\line(1,0){20}}%2
\put(10,30){\line(1,0){20}}%3
%\put(10,10){\line(1,0){20}}%4
%
%\put(30,70){\line(1,0){20}}%5
%\put(30,50){\line(1,0){20}}%6
%\put(30,30){\line(1,0){20}}%7
\put(30,10){\line(1,0){20}}%8
%
%\put(50,70){\line(1,0){20}}%9
\put(50,50){\line(1,0){20}}%10
%\put(50,30){\line(1,0){20}}%11
%\put(50,10){\line(1,0){20}}%12
\end{picture}
\qquad \qquad
\begin{picture}(64,80)
\put(8,64){\circle{16}}%circle 1
\put(32,64){\circle{16}}%circle 2
\put(56,64){\circle{16}}%circle 3
\put(8,40){\circle{16}}%circle 4
\put(32,40){\circle{16}}%circle 5
\put(56,40){\circle{16}}%circle 6
\put(8,16){\circle{16}}%circle 7
\put(32,16){\circle{16}}%circle 8
\put(56,16){\circle{16}}%circle 9

\put(16,64){\line(1,0){8}}%1 to 2
\put(40,64){\line(1,0){8}}%2 to 3
\put(16,40){\line(1,0){8}}%4 to 5
\put(40,40){\line(1,0){8}}%5 to 6
\put(16,16){\line(1,0){8}}%6 to 7
\put(40,16){\line(1,0){8}}%7 to 8

\put(8,48){\line(0,1){8}}%4 to 1
\put(32,48){\line(0,1){8}}%5 to 2
\put(56,48){\line(0,1){8}}%6 to 3
\put(8,24){\line(0,1){8}}%7 to 4
\put(32,24){\line(0,1){8}}%8 to 5
\put(56,24){\line(0,1){8}}%9 to 6

\put(6,61){7}%1
\put(30,61){1}%2
\put(54,61){5}%3
\put(6,37){6}%4
\put(30,37){0}%5
\put(54,37){2}%6
\put(6,13){2}%7
\put(30,13){6}%8
\put(54,13){0}%9
\end{picture}
\qquad \qquad
\begin{picture}(40,80)
\put(8,52){\circle{16}}%circle 1
\put(32,52){\circle{16}}%circle 2
\put(8,28){\circle{16}}%circle 3
\put(32,28){\circle{16}}%circle 4
\put(16,52){\line(1,0){8}}%1 to 2
\put(16,28){\line(1,0){8}}%3 to 4
\put(8,36){\line(0,1){8}}%1 to 3
\put(32,36){\line(0,1){8}}%2 to 4
%\put(16,38){$r$}%circle 1
\put(40,58){$2$}%circle 2
%\put(16,14){$r$}%circle 3
%\put(40,14){$r$}%circle 4
\put(6,49){7}%1
\put(30,49){1}%2
\put(6,25){6}%3
\put(30,25){0}%4
\end{picture}\]

\[
\begin{picture}(80,80)
\put(-40,36){$\mapsto$}
\put(0,60){\framebox(20,20)[r]{$1\;\;$}}
\put(0,40){\framebox(20,20)[r]{$2\;\;$}}
\put(0,20){\framebox(20,20)[r]{$3\;\;$}}
\put(0,0){\framebox(20,20)[r]{$4\;\;$}}
\put(20,60){\framebox(20,20)[r]{$5\;\;$}}
\put(20,40){\framebox(20,20)[r]{$6\;\;$}}
\put(20,20){\framebox(20,20)[r]{$7\;\;$}}
\put(20,0){\framebox(20,20)[r]{$8\;\;$}}
\put(40,60){\framebox(20,20)[r]{$9\;\;$}}
\put(40,40){\framebox(20,20)[r]{$10\;$}}
\put(40,20){\framebox(20,20)[r]{$11\;$}}
\put(40,0){\framebox(20,20)[r]{$12\;$}}
\put(60,60){\framebox(20,20)[r]{$13\;$}}
\put(60,40){\framebox(20,20)[r]{$14\;$}}
\put(60,20){\framebox(20,20)[r]{$15\;$}}
\put(60,0){\framebox(20,20)[r]{$16\;$}}

%vertical lines up
\put(10,50){\line(0,1){20}}%2
%\put(10,30){\line(0,1){20}}%3
\put(10,10){\line(0,1){20}}%4
\put(30,50){\line(0,1){20}}%6
\put(30,30){\line(0,1){20}}%7
%\put(30,10){\line(0,1){20}}%8
%
%\put(50,50){\line(0,1){20}}%10
\put(50,30){\line(0,1){20}}%11
\put(50,10){\line(0,1){20}}%12
%
%\put(70,50){\line(0,1){20}}%14
%\put(70,30){\line(0,1){20}}%15
\put(70,10){\line(0,1){20}}%16

%Horizontal lines to the right
%\put(10,70){\line(1,0){20}}%1
%\put(10,50){\line(1,0){20}}%2
\put(10,30){\line(1,0){20}}%3
%\put(10,10){\line(1,0){20}}%4
%
%\put(30,70){\line(1,0){20}}%5
%\put(30,50){\line(1,0){20}}%6
%\put(30,30){\line(1,0){20}}%7
\put(30,10){\line(1,0){20}}%8
%
%\put(50,70){\line(1,0){20}}%9
\put(50,50){\line(1,0){20}}%10
%\put(50,30){\line(1,0){20}}%11
%\put(50,10){\line(1,0){20}}%12
\end{picture}
\qquad \qquad
\begin{picture}(64,80)
\put(8,64){\circle{16}}%circle 1
\put(32,64){\circle{16}}%circle 2
\put(56,64){\circle{16}}%circle 3
\put(8,40){\circle{16}}%circle 4
\put(32,40){\circle{16}}%circle 5
\put(56,40){\circle{16}}%circle 6
\put(8,16){\circle{16}}%circle 7
\put(32,16){\circle{16}}%circle 8
\put(56,16){\circle{16}}%circle 9

\put(16,64){\line(1,0){8}}%1 to 2
\put(40,64){\line(1,0){8}}%2 to 3
\put(16,40){\line(1,0){8}}%4 to 5
\put(40,40){\line(1,0){8}}%5 to 6
\put(16,16){\line(1,0){8}}%6 to 7
\put(40,16){\line(1,0){8}}%7 to 8

\put(8,48){\line(0,1){8}}%4 to 1
\put(32,48){\line(0,1){8}}%5 to 2
\put(56,48){\line(0,1){8}}%6 to 3
\put(8,24){\line(0,1){8}}%7 to 4
\put(32,24){\line(0,1){8}}%8 to 5
\put(56,24){\line(0,1){8}}%9 to 6

\put(6,61){0}%1
\put(30,61){1}%2
\put(54,61){5}%3
\put(6,37){6}%4
\put(30,37){0}%5
\put(54,37){2}%6
\put(6,13){2}%7
\put(30,13){6}%8
\put(54,13){0}%9
\end{picture}
\qquad \qquad
\begin{picture}(40,80)
\put(8,52){\circle{16}}%circle 1
\put(32,52){\circle{16}}%circle 2
\put(8,28){\circle{16}}%circle 3
\put(32,28){\circle{16}}%circle 4
\put(16,52){\line(1,0){8}}%1 to 2
\put(16,28){\line(1,0){8}}%3 to 4
\put(8,36){\line(0,1){8}}%1 to 3
\put(32,36){\line(0,1){8}}%2 to 4
%\put(16,38){$r$}%circle 1
\put(40,58){$2$}%circle 2
%\put(16,14){$r$}%circle 3
%\put(40,14){$r$}%circle 4
\put(6,49){0}%1
\put(30,49){1}%2
\put(6,25){6}%3
\put(30,25){0}%4
\end{picture}\]

\[
\begin{picture}(80,80)
\put(-40,36){$\mapsto$}
\put(0,60){\framebox(20,20)[r]{$1\;\;$}}
\put(0,40){\framebox(20,20)[r]{$2\;\;$}}
\put(0,20){\framebox(20,20)[r]{$3\;\;$}}
\put(0,0){\framebox(20,20)[r]{$4\;\;$}}
\put(20,60){\framebox(20,20)[r]{$5\;\;$}}
\put(20,40){\framebox(20,20)[r]{$6\;\;$}}
\put(20,20){\framebox(20,20)[r]{$7\;\;$}}
\put(20,0){\framebox(20,20)[r]{$8\;\;$}}
\put(40,60){\framebox(20,20)[r]{$9\;\;$}}
\put(40,40){\framebox(20,20)[r]{$10\;$}}
\put(40,20){\framebox(20,20)[r]{$11\;$}}
\put(40,0){\framebox(20,20)[r]{$12\;$}}
\put(60,60){\framebox(20,20)[r]{$13\;$}}
\put(60,40){\framebox(20,20)[r]{$14\;$}}
\put(60,20){\framebox(20,20)[r]{$15\;$}}
\put(60,0){\framebox(20,20)[r]{$16\;$}}

%vertical lines up
\put(10,50){\line(0,1){20}}%2
%\put(10,30){\line(0,1){20}}%3
\put(10,10){\line(0,1){20}}%4
\put(30,50){\line(0,1){20}}%6
\put(30,30){\line(0,1){20}}%7
%\put(30,10){\line(0,1){20}}%8
%
%\put(50,50){\line(0,1){20}}%10
\put(50,30){\line(0,1){20}}%11
\put(50,10){\line(0,1){20}}%12
\put(70,50){\line(0,1){20}}%14
%\put(70,30){\line(0,1){20}}%15
\put(70,10){\line(0,1){20}}%16

%Horizontal lines to the right
%\put(10,70){\line(1,0){20}}%1
%\put(10,50){\line(1,0){20}}%2
\put(10,30){\line(1,0){20}}%3
%\put(10,10){\line(1,0){20}}%4
%
%\put(30,70){\line(1,0){20}}%5
%\put(30,50){\line(1,0){20}}%6
%\put(30,30){\line(1,0){20}}%7
\put(30,10){\line(1,0){20}}%8
%
%\put(50,70){\line(1,0){20}}%9
\put(50,50){\line(1,0){20}}%10
%\put(50,30){\line(1,0){20}}%11
%\put(50,10){\line(1,0){20}}%12
\end{picture}
\qquad \qquad
\begin{picture}(64,80)
\put(8,64){\circle{16}}%circle 1
\put(32,64){\circle{16}}%circle 2
\put(56,64){\circle{16}}%circle 3
\put(8,40){\circle{16}}%circle 4
\put(32,40){\circle{16}}%circle 5
\put(56,40){\circle{16}}%circle 6
\put(8,16){\circle{16}}%circle 7
\put(32,16){\circle{16}}%circle 8
\put(56,16){\circle{16}}%circle 9

\put(16,64){\line(1,0){8}}%1 to 2
\put(40,64){\line(1,0){8}}%2 to 3
\put(16,40){\line(1,0){8}}%4 to 5
\put(40,40){\line(1,0){8}}%5 to 6
\put(16,16){\line(1,0){8}}%6 to 7
\put(40,16){\line(1,0){8}}%7 to 8

\put(8,48){\line(0,1){8}}%4 to 1
\put(32,48){\line(0,1){8}}%5 to 2
\put(56,48){\line(0,1){8}}%6 to 3
\put(8,24){\line(0,1){8}}%7 to 4
\put(32,24){\line(0,1){8}}%8 to 5
\put(56,24){\line(0,1){8}}%9 to 6

\put(6,61){0}%1
\put(30,61){1}%2
\put(54,61){6}%3
\put(6,37){6}%4
\put(30,37){0}%5
\put(54,37){2}%6
\put(6,13){2}%7
\put(30,13){6}%8
\put(54,13){0}%9
\end{picture}
\qquad \qquad
\begin{picture}(40,80)
\put(8,52){\circle{16}}%circle 1
\put(32,52){\circle{16}}%circle 2
\put(8,28){\circle{16}}%circle 3
\put(32,28){\circle{16}}%circle 4
\put(16,52){\line(1,0){8}}%1 to 2
\put(16,28){\line(1,0){8}}%3 to 4
\put(8,36){\line(0,1){8}}%1 to 3
\put(32,36){\line(0,1){8}}%2 to 4
%\put(16,38){$r$}%circle 1
\put(40,58){$2$}%circle 2
%\put(16,14){$r$}%circle 3
%\put(40,14){$r$}%circle 4
\put(6,49){0}%1
\put(30,49){1}%2
\put(6,25){6}%3
\put(30,25){0}%4
\end{picture}\]

There is a product of the form (\ref{bigf}) associated to each diagram above. The following ordering is shared by all three diagrams, and resolves all the singularities of $\Lambda^c$ in the same way. For simplicity we write $f_{p,q}$ instead of $f_{p,q}(z_p + c_p(\Lambda^c), z_q + c_q(\Lambda^c))$. We indicate the singularities along with their triple terms, idempotents may then be added to form regular triples at $z_1 = \cdots = z_n$.

\small
\[\begin{array}{l} f_{1,2}f_{1,3}f_{2,3}f_{1,4}f_{2,4}f_{3,4}\left(f_{1,5}f_{1,6}\right)f_{1,7}f_{2,5}\left(f_{2,6}f_{2,7}\right)f_{3,5}f_{3,6}f_{3,7}f_{4,5}f_{4,6}f_{4,7}f_{5,6}f_{5,7}f_{6,7}f_{8,9}f_{1,9}f_{2,9}f_{3,9}\\
f_{4,9}f_{5,9}f_{6,9}f_{7,9}f_{8,10}f_{8,11}\left(f_{1,10}f_{1,11}\right)f_{2,10}f_{2,11}f_{3,10}f_{3,11}f_{4,10}f_{4,11}\left(f_{5,10}f_{6,10}\right)f_{7,10}f_{5,11}\left(f_{6,11}f_{7,11}\right)\\
f_{10,11}f_{9,11}f_{1,8}f_{1,12}\left(f_{2,8}f_{2,12}\right)\left(f_{3,8}f_{4,8}\right)f_{3,12}f_{4,12}f_{5,8}f_{5,12}f_{6,8}f_{6,12}\left(f_{7,8}f_{7,12}\right)f_{8,12}f_{9,12}f_{10,12}f_{11,12}\\
f_{1,13}f_{2,13}f_{3,13}f_{4,13}f_{5,13}f_{6,13}f_{7,13}f_{8,13}f_{10,13}f_{9,13}f_{11,13}f_{12,13}f_{1,14}f_{2,14}f_{3,14}f_{4,14}f_{5,14}f_{6,14}f_{7,14}f_{8,14}\\
\left(f_{9,10}f_{9,14}\right)f_{10,14}f_{11,14}f_{12,14}f_{13,14}\left(f_{1,15}f_{1,16}\right)f_{2,15}f_{2,16}f_{3,15}f_{3,16}f_{4,15}f_{4,16}\left(f_{5,15}f_{6,15}\right)f_{7,15}f_{5,16}\\
\left(f_{6,16}f_{7,16}\right)f_{8,15}f_{8,16}f_{9,15}f_{9,16}\left(f_{10,15}f_{11,15}\right)f_{12,15}f_{10,16}\left(f_{11,16}f_{12,16}\right)f_{13,15}f_{13,16}f_{14,15}f_{14,16}f_{15,16}
\end{array}\]
\large

And hence, by first comparing triples and then all other terms, these three products have the same value at $z_1 = \cdots = z_n$. We now continue the chain of diagrams, starting with our last diagram.

\[
\begin{picture}(80,80)
%\put(-40,36){$\mapsto$}
\put(0,60){\framebox(20,20)[r]{$1\;\;$}}
\put(0,40){\framebox(20,20)[r]{$2\;\;$}}
\put(0,20){\framebox(20,20)[r]{$3\;\;$}}
\put(0,0){\framebox(20,20)[r]{$4\;\;$}}
\put(20,60){\framebox(20,20)[r]{$5\;\;$}}
\put(20,40){\framebox(20,20)[r]{$6\;\;$}}
\put(20,20){\framebox(20,20)[r]{$7\;\;$}}
\put(20,0){\framebox(20,20)[r]{$8\;\;$}}
\put(40,60){\framebox(20,20)[r]{$9\;\;$}}
\put(40,40){\framebox(20,20)[r]{$10\;$}}
\put(40,20){\framebox(20,20)[r]{$11\;$}}
\put(40,0){\framebox(20,20)[r]{$12\;$}}
\put(60,60){\framebox(20,20)[r]{$13\;$}}
\put(60,40){\framebox(20,20)[r]{$14\;$}}
\put(60,20){\framebox(20,20)[r]{$15\;$}}
\put(60,0){\framebox(20,20)[r]{$16\;$}}

%vertical lines up
\put(10,50){\line(0,1){20}}%2
%\put(10,30){\line(0,1){20}}%3
\put(10,10){\line(0,1){20}}%4
\put(30,50){\line(0,1){20}}%6
\put(30,30){\line(0,1){20}}%7
%\put(30,10){\line(0,1){20}}%8
%
%\put(50,50){\line(0,1){20}}%10
\put(50,30){\line(0,1){20}}%11
\put(50,10){\line(0,1){20}}%12
\put(70,50){\line(0,1){20}}%14
%\put(70,30){\line(0,1){20}}%15
\put(70,10){\line(0,1){20}}%16

%Horizontal lines to the right
%\put(10,70){\line(1,0){20}}%1
%\put(10,50){\line(1,0){20}}%2
\put(10,30){\line(1,0){20}}%3
%\put(10,10){\line(1,0){20}}%4
%
%\put(30,70){\line(1,0){20}}%5
%\put(30,50){\line(1,0){20}}%6
%\put(30,30){\line(1,0){20}}%7
\put(30,10){\line(1,0){20}}%8
%
%\put(50,70){\line(1,0){20}}%9
\put(50,50){\line(1,0){20}}%10
%\put(50,30){\line(1,0){20}}%11
%\put(50,10){\line(1,0){20}}%12
\end{picture}
\qquad \qquad
\begin{picture}(64,80)
\put(8,64){\circle{16}}%circle 1
\put(32,64){\circle{16}}%circle 2
\put(56,64){\circle{16}}%circle 3
\put(8,40){\circle{16}}%circle 4
\put(32,40){\circle{16}}%circle 5
\put(56,40){\circle{16}}%circle 6
\put(8,16){\circle{16}}%circle 7
\put(32,16){\circle{16}}%circle 8
\put(56,16){\circle{16}}%circle 9

\put(16,64){\line(1,0){8}}%1 to 2
\put(40,64){\line(1,0){8}}%2 to 3
\put(16,40){\line(1,0){8}}%4 to 5
\put(40,40){\line(1,0){8}}%5 to 6
\put(16,16){\line(1,0){8}}%6 to 7
\put(40,16){\line(1,0){8}}%7 to 8

\put(8,48){\line(0,1){8}}%4 to 1
\put(32,48){\line(0,1){8}}%5 to 2
\put(56,48){\line(0,1){8}}%6 to 3
\put(8,24){\line(0,1){8}}%7 to 4
\put(32,24){\line(0,1){8}}%8 to 5
\put(56,24){\line(0,1){8}}%9 to 6

\put(6,61){0}%1
\put(30,61){1}%2
\put(54,61){6}%3
\put(6,37){6}%4
\put(30,37){0}%5
\put(54,37){2}%6
\put(6,13){2}%7
\put(30,13){6}%8
\put(54,13){0}%9
\end{picture}
\qquad \qquad
\begin{picture}(40,80)
\put(8,52){\circle{16}}%circle 1
\put(32,52){\circle{16}}%circle 2
\put(8,28){\circle{16}}%circle 3
\put(32,28){\circle{16}}%circle 4
\put(16,52){\line(1,0){8}}%1 to 2
\put(16,28){\line(1,0){8}}%3 to 4
\put(8,36){\line(0,1){8}}%1 to 3
\put(32,36){\line(0,1){8}}%2 to 4
%\put(16,38){$r$}%circle 1
\put(40,58){$2$}%circle 2
%\put(16,14){$r$}%circle 3
%\put(40,14){$r$}%circle 4
\put(6,49){0}%1
\put(30,49){1}%2
\put(6,25){6}%3
\put(30,25){0}%4
\end{picture}\]

\[
\begin{picture}(80,80)
\put(-40,36){$\mapsto$}
\put(0,60){\framebox(20,20)[r]{$1\;\;$}}
\put(0,40){\framebox(20,20)[r]{$2\;\;$}}
\put(0,20){\framebox(20,20)[r]{$3\;\;$}}
\put(0,0){\framebox(20,20)[r]{$4\;\;$}}
\put(20,60){\framebox(20,20)[r]{$5\;\;$}}
\put(20,40){\framebox(20,20)[r]{$6\;\;$}}
\put(20,20){\framebox(20,20)[r]{$7\;\;$}}
\put(20,0){\framebox(20,20)[r]{$8\;\;$}}
\put(40,60){\framebox(20,20)[r]{$9\;\;$}}
\put(40,40){\framebox(20,20)[r]{$10\;$}}
\put(40,20){\framebox(20,20)[r]{$11\;$}}
\put(40,0){\framebox(20,20)[r]{$12\;$}}
\put(60,60){\framebox(20,20)[r]{$13\;$}}
\put(60,40){\framebox(20,20)[r]{$14\;$}}
\put(60,20){\framebox(20,20)[r]{$15\;$}}
\put(60,0){\framebox(20,20)[r]{$16\;$}}

%vertical lines up
\put(10,50){\line(0,1){20}}%2
%\put(10,30){\line(0,1){20}}%3
\put(10,10){\line(0,1){20}}%4
\put(30,50){\line(0,1){20}}%6
\put(30,30){\line(0,1){20}}%7
%\put(30,10){\line(0,1){20}}%8
%
%\put(50,50){\line(0,1){20}}%10
\put(50,30){\line(0,1){20}}%11
\put(50,10){\line(0,1){20}}%12
\put(70,50){\line(0,1){20}}%14
%\put(70,30){\line(0,1){20}}%15
\put(70,10){\line(0,1){20}}%16

%Horizontal lines to the right
%\put(10,70){\line(1,0){20}}%1
%\put(10,50){\line(1,0){20}}%2
\put(10,30){\line(1,0){20}}%3
%\put(10,10){\line(1,0){20}}%4
%
%\put(30,70){\line(1,0){20}}%5
%\put(30,50){\line(1,0){20}}%6
%\put(30,30){\line(1,0){20}}%7
\put(30,10){\line(1,0){20}}%8
%
%\put(50,70){\line(1,0){20}}%9
%\put(50,50){\line(1,0){20}}%10
%\put(50,30){\line(1,0){20}}%11
%\put(50,10){\line(1,0){20}}%12
\end{picture}
\qquad \qquad
\begin{picture}(64,80)
\put(8,64){\circle{16}}%circle 1
\put(32,64){\circle{16}}%circle 2
\put(56,64){\circle{16}}%circle 3
\put(8,40){\circle{16}}%circle 4
\put(32,40){\circle{16}}%circle 5
\put(56,40){\circle{16}}%circle 6
\put(8,16){\circle{16}}%circle 7
\put(32,16){\circle{16}}%circle 8
\put(56,16){\circle{16}}%circle 9

\put(16,64){\line(1,0){8}}%1 to 2
\put(40,64){\line(1,0){8}}%2 to 3
\put(16,40){\line(1,0){8}}%4 to 5
\put(40,40){\line(1,0){8}}%5 to 6
\put(16,16){\line(1,0){8}}%6 to 7
\put(40,16){\line(1,0){8}}%7 to 8

\put(8,48){\line(0,1){8}}%4 to 1
\put(32,48){\line(0,1){8}}%5 to 2
\put(56,48){\line(0,1){8}}%6 to 3
\put(8,24){\line(0,1){8}}%7 to 4
\put(32,24){\line(0,1){8}}%8 to 5
\put(56,24){\line(0,1){8}}%9 to 6

\put(6,61){0}%1
\put(30,61){1}%2
\put(54,61){7}%3
\put(6,37){6}%4
\put(30,37){0}%5
\put(54,37){1}%6
\put(6,13){2}%7
\put(30,13){6}%8
\put(54,13){0}%9
\end{picture}
\qquad \qquad
\begin{picture}(40,80)
\put(8,52){\circle{16}}%circle 1
\put(32,52){\circle{16}}%circle 2
\put(8,28){\circle{16}}%circle 3
\put(32,28){\circle{16}}%circle 4
\put(16,52){\line(1,0){8}}%1 to 2
\put(16,28){\line(1,0){8}}%3 to 4
\put(8,36){\line(0,1){8}}%1 to 3
\put(32,36){\line(0,1){8}}%2 to 4
%\put(16,38){$r$}%circle 1
\put(40,58){$2$}%circle 2
%\put(16,14){$r$}%circle 3
%\put(40,14){$r$}%circle 4
\put(6,49){0}%1
\put(30,49){1}%2
\put(6,25){6}%3
\put(30,25){0}%4
\end{picture}\]

\[
\begin{picture}(80,80)
\put(-40,36){$\mapsto$}
\put(0,60){\framebox(20,20)[r]{$1\;\;$}}
\put(0,40){\framebox(20,20)[r]{$2\;\;$}}
\put(0,20){\framebox(20,20)[r]{$3\;\;$}}
\put(0,0){\framebox(20,20)[r]{$4\;\;$}}
\put(20,60){\framebox(20,20)[r]{$5\;\;$}}
\put(20,40){\framebox(20,20)[r]{$6\;\;$}}
\put(20,20){\framebox(20,20)[r]{$7\;\;$}}
\put(20,0){\framebox(20,20)[r]{$8\;\;$}}
\put(40,60){\framebox(20,20)[r]{$9\;\;$}}
\put(40,40){\framebox(20,20)[r]{$10\;$}}
\put(40,20){\framebox(20,20)[r]{$11\;$}}
\put(40,0){\framebox(20,20)[r]{$12\;$}}
\put(60,60){\framebox(20,20)[r]{$13\;$}}
\put(60,40){\framebox(20,20)[r]{$14\;$}}
\put(60,20){\framebox(20,20)[r]{$15\;$}}
\put(60,0){\framebox(20,20)[r]{$16\;$}}

%vertical lines up
\put(10,50){\line(0,1){20}}%2
%\put(10,30){\line(0,1){20}}%3
\put(10,10){\line(0,1){20}}%4
\put(30,50){\line(0,1){20}}%6
\put(30,30){\line(0,1){20}}%7
%\put(30,10){\line(0,1){20}}%8
%
\put(50,50){\line(0,1){20}}%10
\put(50,30){\line(0,1){20}}%11
\put(50,10){\line(0,1){20}}%12
\put(70,50){\line(0,1){20}}%14
%\put(70,30){\line(0,1){20}}%15
\put(70,10){\line(0,1){20}}%16

%Horizontal lines to the right
%\put(10,70){\line(1,0){20}}%1
%\put(10,50){\line(1,0){20}}%2
\put(10,30){\line(1,0){20}}%3
%\put(10,10){\line(1,0){20}}%4
%
%\put(30,70){\line(1,0){20}}%5
%\put(30,50){\line(1,0){20}}%6
%\put(30,30){\line(1,0){20}}%7
\put(30,10){\line(1,0){20}}%8
%
%\put(50,70){\line(1,0){20}}%9
%\put(50,50){\line(1,0){20}}%10
%\put(50,30){\line(1,0){20}}%11
%\put(50,10){\line(1,0){20}}%12
\end{picture}
\qquad \qquad
\begin{picture}(64,80)
\put(8,64){\circle{16}}%circle 1
\put(32,64){\circle{16}}%circle 2
\put(56,64){\circle{16}}%circle 3
\put(8,40){\circle{16}}%circle 4
\put(32,40){\circle{16}}%circle 5
\put(56,40){\circle{16}}%circle 6
\put(8,16){\circle{16}}%circle 7
\put(32,16){\circle{16}}%circle 8
\put(56,16){\circle{16}}%circle 9

\put(16,64){\line(1,0){8}}%1 to 2
\put(40,64){\line(1,0){8}}%2 to 3
\put(16,40){\line(1,0){8}}%4 to 5
\put(40,40){\line(1,0){8}}%5 to 6
\put(16,16){\line(1,0){8}}%6 to 7
\put(40,16){\line(1,0){8}}%7 to 8

\put(8,48){\line(0,1){8}}%4 to 1
\put(32,48){\line(0,1){8}}%5 to 2
\put(56,48){\line(0,1){8}}%6 to 3
\put(8,24){\line(0,1){8}}%7 to 4
\put(32,24){\line(0,1){8}}%8 to 5
\put(56,24){\line(0,1){8}}%9 to 6

\put(6,61){0}%1
\put(30,61){0}%2
\put(54,61){0}%3
\put(6,37){6}%4
\put(30,37){0}%5
\put(54,37){1}%6
\put(6,13){2}%7
\put(30,13){6}%8
\put(54,13){0}%9
\end{picture}
\qquad \qquad
\begin{picture}(40,80)
\put(8,52){\circle{16}}%circle 1
\put(32,52){\circle{16}}%circle 2
\put(8,28){\circle{16}}%circle 3
\put(32,28){\circle{16}}%circle 4
\put(16,52){\line(1,0){8}}%1 to 2
\put(16,28){\line(1,0){8}}%3 to 4
\put(8,36){\line(0,1){8}}%1 to 3
\put(32,36){\line(0,1){8}}%2 to 4
%\put(16,38){$r$}%circle 1
\put(40,58){$2$}%circle 2
%\put(16,14){$r$}%circle 3
%\put(40,14){$r$}%circle 4
\put(6,49){0}%1
\put(30,49){1}%2
\put(6,25){6}%3
\put(30,25){0}%4
\end{picture}\]

\[
\begin{picture}(80,80)
\put(-40,36){$\mapsto$}
\put(0,60){\framebox(20,20)[r]{$1\;\;$}}
\put(0,40){\framebox(20,20)[r]{$2\;\;$}}
\put(0,20){\framebox(20,20)[r]{$3\;\;$}}
\put(0,0){\framebox(20,20)[r]{$4\;\;$}}
\put(20,60){\framebox(20,20)[r]{$5\;\;$}}
\put(20,40){\framebox(20,20)[r]{$6\;\;$}}
\put(20,20){\framebox(20,20)[r]{$7\;\;$}}
\put(20,0){\framebox(20,20)[r]{$8\;\;$}}
\put(40,60){\framebox(20,20)[r]{$9\;\;$}}
\put(40,40){\framebox(20,20)[r]{$10\;$}}
\put(40,20){\framebox(20,20)[r]{$11\;$}}
\put(40,0){\framebox(20,20)[r]{$12\;$}}
\put(60,60){\framebox(20,20)[r]{$13\;$}}
\put(60,40){\framebox(20,20)[r]{$14\;$}}
\put(60,20){\framebox(20,20)[r]{$15\;$}}
\put(60,0){\framebox(20,20)[r]{$16\;$}}

%vertical lines up
\put(10,50){\line(0,1){20}}%2
%\put(10,30){\line(0,1){20}}%3
\put(10,10){\line(0,1){20}}%4
\put(30,50){\line(0,1){20}}%6
\put(30,30){\line(0,1){20}}%7
%\put(30,10){\line(0,1){20}}%8
%
\put(50,50){\line(0,1){20}}%10
\put(50,30){\line(0,1){20}}%11
\put(50,10){\line(0,1){20}}%12
\put(70,50){\line(0,1){20}}%14
%\put(70,30){\line(0,1){20}}%15
\put(70,10){\line(0,1){20}}%16

%Horizontal lines to the right
%\put(10,70){\line(1,0){20}}%1
%\put(10,50){\line(1,0){20}}%2
%\put(10,30){\line(1,0){20}}%3
%\put(10,10){\line(1,0){20}}%4
%
%\put(30,70){\line(1,0){20}}%5
%\put(30,50){\line(1,0){20}}%6
%\put(30,30){\line(1,0){20}}%7
\put(30,10){\line(1,0){20}}%8
%
%\put(50,70){\line(1,0){20}}%9
%\put(50,50){\line(1,0){20}}%10
%\put(50,30){\line(1,0){20}}%11
%\put(50,10){\line(1,0){20}}%12
\end{picture}
\qquad \qquad
\begin{picture}(64,80)
\put(8,64){\circle{16}}%circle 1
\put(32,64){\circle{16}}%circle 2
\put(56,64){\circle{16}}%circle 3
\put(8,40){\circle{16}}%circle 4
\put(32,40){\circle{16}}%circle 5
\put(56,40){\circle{16}}%circle 6
\put(8,16){\circle{16}}%circle 7
\put(32,16){\circle{16}}%circle 8
\put(56,16){\circle{16}}%circle 9

\put(16,64){\line(1,0){8}}%1 to 2
\put(40,64){\line(1,0){8}}%2 to 3
\put(16,40){\line(1,0){8}}%4 to 5
\put(40,40){\line(1,0){8}}%5 to 6
\put(16,16){\line(1,0){8}}%6 to 7
\put(40,16){\line(1,0){8}}%7 to 8

\put(8,48){\line(0,1){8}}%4 to 1
\put(32,48){\line(0,1){8}}%5 to 2
\put(56,48){\line(0,1){8}}%6 to 3
\put(8,24){\line(0,1){8}}%7 to 4
\put(32,24){\line(0,1){8}}%8 to 5
\put(56,24){\line(0,1){8}}%9 to 6

\put(6,61){0}%1
\put(30,61){0}%2
\put(54,61){0}%3
\put(6,37){7}%4
\put(30,37){0}%5
\put(54,37){1}%6
\put(6,13){1}%7
\put(30,13){6}%8
\put(54,13){0}%9
\end{picture}
\qquad \qquad
\begin{picture}(40,80)
\put(8,52){\circle{16}}%circle 1
\put(32,52){\circle{16}}%circle 2
\put(8,28){\circle{16}}%circle 3
\put(32,28){\circle{16}}%circle 4
\put(16,52){\line(1,0){8}}%1 to 2
\put(16,28){\line(1,0){8}}%3 to 4
\put(8,36){\line(0,1){8}}%1 to 3
\put(32,36){\line(0,1){8}}%2 to 4
%\put(16,38){$r$}%circle 1
\put(40,58){$2$}%circle 2
%\put(16,14){$r$}%circle 3
%\put(40,14){$r$}%circle 4
\put(6,49){0}%1
\put(30,49){1}%2
\put(6,25){6}%3
\put(30,25){0}%4
\end{picture}\]

\[
\begin{picture}(80,80)
\put(-40,36){$\mapsto$}
\put(0,60){\framebox(20,20)[r]{$1\;\;$}}
\put(0,40){\framebox(20,20)[r]{$2\;\;$}}
\put(0,20){\framebox(20,20)[r]{$3\;\;$}}
\put(0,0){\framebox(20,20)[r]{$4\;\;$}}
\put(20,60){\framebox(20,20)[r]{$5\;\;$}}
\put(20,40){\framebox(20,20)[r]{$6\;\;$}}
\put(20,20){\framebox(20,20)[r]{$7\;\;$}}
\put(20,0){\framebox(20,20)[r]{$8\;\;$}}
\put(40,60){\framebox(20,20)[r]{$9\;\;$}}
\put(40,40){\framebox(20,20)[r]{$10\;$}}
\put(40,20){\framebox(20,20)[r]{$11\;$}}
\put(40,0){\framebox(20,20)[r]{$12\;$}}
\put(60,60){\framebox(20,20)[r]{$13\;$}}
\put(60,40){\framebox(20,20)[r]{$14\;$}}
\put(60,20){\framebox(20,20)[r]{$15\;$}}
\put(60,0){\framebox(20,20)[r]{$16\;$}}

%vertical lines up
\put(10,50){\line(0,1){20}}%2
\put(10,30){\line(0,1){20}}%3
\put(10,10){\line(0,1){20}}%4
\put(30,50){\line(0,1){20}}%6
\put(30,30){\line(0,1){20}}%7
%\put(30,10){\line(0,1){20}}%8
%
\put(50,50){\line(0,1){20}}%10
\put(50,30){\line(0,1){20}}%11
\put(50,10){\line(0,1){20}}%12
\put(70,50){\line(0,1){20}}%14
%\put(70,30){\line(0,1){20}}%15
\put(70,10){\line(0,1){20}}%16

%Horizontal lines to the right
%\put(10,70){\line(1,0){20}}%1
%\put(10,50){\line(1,0){20}}%2
%\put(10,30){\line(1,0){20}}%3
%\put(10,10){\line(1,0){20}}%4
%
%\put(30,70){\line(1,0){20}}%5
%\put(30,50){\line(1,0){20}}%6
%\put(30,30){\line(1,0){20}}%7
\put(30,10){\line(1,0){20}}%8
%
%\put(50,70){\line(1,0){20}}%9
%\put(50,50){\line(1,0){20}}%10
%\put(50,30){\line(1,0){20}}%11
%\put(50,10){\line(1,0){20}}%12
\end{picture}
\qquad \qquad
\begin{picture}(64,80)
\put(8,64){\circle{16}}%circle 1
\put(32,64){\circle{16}}%circle 2
\put(56,64){\circle{16}}%circle 3
\put(8,40){\circle{16}}%circle 4
\put(32,40){\circle{16}}%circle 5
\put(56,40){\circle{16}}%circle 6
\put(8,16){\circle{16}}%circle 7
\put(32,16){\circle{16}}%circle 8
\put(56,16){\circle{16}}%circle 9

\put(16,64){\line(1,0){8}}%1 to 2
\put(40,64){\line(1,0){8}}%2 to 3
\put(16,40){\line(1,0){8}}%4 to 5
\put(40,40){\line(1,0){8}}%5 to 6
\put(16,16){\line(1,0){8}}%6 to 7
\put(40,16){\line(1,0){8}}%7 to 8

\put(8,48){\line(0,1){8}}%4 to 1
\put(32,48){\line(0,1){8}}%5 to 2
\put(56,48){\line(0,1){8}}%6 to 3
\put(8,24){\line(0,1){8}}%7 to 4
\put(32,24){\line(0,1){8}}%8 to 5
\put(56,24){\line(0,1){8}}%9 to 6

\put(6,61){0}%1
\put(30,61){0}%2
\put(54,61){0}%3
\put(6,37){0}%4
\put(30,37){0}%5
\put(54,37){1}%6
\put(6,13){1}%7
\put(30,13){6}%8
\put(54,13){0}%9
\end{picture}
\qquad \qquad
\begin{picture}(40,80)
\put(8,52){\circle{16}}%circle 1
\put(32,52){\circle{16}}%circle 2
\put(8,28){\circle{16}}%circle 3
\put(32,28){\circle{16}}%circle 4
\put(16,52){\line(1,0){8}}%1 to 2
\put(16,28){\line(1,0){8}}%3 to 4
\put(8,36){\line(0,1){8}}%1 to 3
\put(32,36){\line(0,1){8}}%2 to 4
%\put(16,38){$r$}%circle 1
\put(40,58){$2$}%circle 2
%\put(16,14){$r$}%circle 3
%\put(40,14){$r$}%circle 4
\put(6,49){0}%1
\put(30,49){1}%2
\put(3,25){13}%3
\put(30,25){0}%4
\end{picture}\]

\[
\begin{picture}(80,80)
\put(-40,36){$\mapsto$}
\put(0,60){\framebox(20,20)[r]{$1\;\;$}}
\put(0,40){\framebox(20,20)[r]{$2\;\;$}}
\put(0,20){\framebox(20,20)[r]{$3\;\;$}}
\put(0,0){\framebox(20,20)[r]{$4\;\;$}}
\put(20,60){\framebox(20,20)[r]{$5\;\;$}}
\put(20,40){\framebox(20,20)[r]{$6\;\;$}}
\put(20,20){\framebox(20,20)[r]{$7\;\;$}}
\put(20,0){\framebox(20,20)[r]{$8\;\;$}}
\put(40,60){\framebox(20,20)[r]{$9\;\;$}}
\put(40,40){\framebox(20,20)[r]{$10\;$}}
\put(40,20){\framebox(20,20)[r]{$11\;$}}
\put(40,0){\framebox(20,20)[r]{$12\;$}}
\put(60,60){\framebox(20,20)[r]{$13\;$}}
\put(60,40){\framebox(20,20)[r]{$14\;$}}
\put(60,20){\framebox(20,20)[r]{$15\;$}}
\put(60,0){\framebox(20,20)[r]{$16\;$}}

%vertical lines up
\put(10,50){\line(0,1){20}}%2
\put(10,30){\line(0,1){20}}%3
\put(10,10){\line(0,1){20}}%4
\put(30,50){\line(0,1){20}}%6
\put(30,30){\line(0,1){20}}%7
%\put(30,10){\line(0,1){20}}%8
%
\put(50,50){\line(0,1){20}}%10
\put(50,30){\line(0,1){20}}%11
\put(50,10){\line(0,1){20}}%12
\put(70,50){\line(0,1){20}}%14
\put(70,30){\line(0,1){20}}%15
\put(70,10){\line(0,1){20}}%16

%Horizontal lines to the right
%\put(10,70){\line(1,0){20}}%1
%\put(10,50){\line(1,0){20}}%2
%\put(10,30){\line(1,0){20}}%3
%\put(10,10){\line(1,0){20}}%4
%
%\put(30,70){\line(1,0){20}}%5
%\put(30,50){\line(1,0){20}}%6
%\put(30,30){\line(1,0){20}}%7
\put(30,10){\line(1,0){20}}%8
%
%\put(50,70){\line(1,0){20}}%9
%\put(50,50){\line(1,0){20}}%10
%\put(50,30){\line(1,0){20}}%11
%\put(50,10){\line(1,0){20}}%12
\end{picture}
\qquad \qquad
\begin{picture}(64,80)
\put(8,64){\circle{16}}%circle 1
\put(32,64){\circle{16}}%circle 2
\put(56,64){\circle{16}}%circle 3
\put(8,40){\circle{16}}%circle 4
\put(32,40){\circle{16}}%circle 5
\put(56,40){\circle{16}}%circle 6
\put(8,16){\circle{16}}%circle 7
\put(32,16){\circle{16}}%circle 8
\put(56,16){\circle{16}}%circle 9

\put(16,64){\line(1,0){8}}%1 to 2
\put(40,64){\line(1,0){8}}%2 to 3
\put(16,40){\line(1,0){8}}%4 to 5
\put(40,40){\line(1,0){8}}%5 to 6
\put(16,16){\line(1,0){8}}%6 to 7
\put(40,16){\line(1,0){8}}%7 to 8

\put(8,48){\line(0,1){8}}%4 to 1
\put(32,48){\line(0,1){8}}%5 to 2
\put(56,48){\line(0,1){8}}%6 to 3
\put(8,24){\line(0,1){8}}%7 to 4
\put(32,24){\line(0,1){8}}%8 to 5
\put(56,24){\line(0,1){8}}%9 to 6

\put(6,61){0}%1
\put(30,61){0}%2
\put(54,61){0}%3
\put(6,37){0}%4
\put(30,37){0}%5
\put(54,37){0}%6
\put(6,13){1}%7
\put(30,13){6}%8
\put(54,13){0}%9
\end{picture}
\qquad \qquad
\begin{picture}(40,80)
\put(8,52){\circle{16}}%circle 1
\put(32,52){\circle{16}}%circle 2
\put(8,28){\circle{16}}%circle 3
\put(32,28){\circle{16}}%circle 4
\put(16,52){\line(1,0){8}}%1 to 2
\put(16,28){\line(1,0){8}}%3 to 4
\put(8,36){\line(0,1){8}}%1 to 3
\put(32,36){\line(0,1){8}}%2 to 4
%\put(16,38){$r$}%circle 1
\put(40,58){$2$}%circle 2
%\put(16,14){$r$}%circle 3
%\put(40,14){$r$}%circle 4
\put(6,49){0}%1
\put(30,49){0}%2
\put(3,25){13}%3
\put(30,25){0}%4
\end{picture}\]

We may rearrange the above ordering using the Yang-Baxter relations (\ref{triple}) and (\ref{commute}), and doing so does not change the value of the product. We do this to obtain the following ordering.

\small
\[\begin{array}{l} f_{1,2}f_{1,3}f_{2,3}f_{1,4}f_{2,4}f_{3,4}\left(f_{1,5}f_{1,6}\right)f_{1,7}f_{2,5}\left(f_{2,6}f_{2,7}\right)f_{3,5}f_{3,6}f_{3,7}f_{4,5}f_{4,6}f_{4,7}f_{5,6}f_{5,7}f_{6,7}f_{8,9}f_{1,9}f_{2,9}f_{3,9}\\
f_{4,9}f_{5,9}f_{6,9}f_{7,9}f_{8,10}f_{1,10}f_{2,10}f_{3,10}f_{4,10}\left(f_{5,10}f_{6,10}\right)f_{7,10}f_{9,10}f_{8,11}\left(f_{1,11}f_{2,11}\right)f_{3,11}f_{4,11}f_{5,11}\left(f_{6,11}f_{7,11}\right)\\
f_{9,11}f_{10,11}f_{1,8}f_{1,12}\left(f_{2,8}f_{2,12}\right)\left(f_{3,8}f_{4,8}\right)f_{3,12}f_{4,12}f_{5,8}f_{5,12}f_{6,8}f_{6,12}\left(f_{7,8}f_{7,12}\right)f_{8,12}f_{9,12}f_{10,12}f_{11,12}f_{1,13}\\
f_{1,14}f_{2,13}f_{2,14}f_{3,13}f_{3,14}f_{4,13}f_{4,14}f_{5,13}f_{5,14}f_{6,13}f_{6,14}f_{7,13}f_{7,14}f_{8,13}f_{8,14}\left(f_{9,13}f_{9,14}\right)f_{10,13}f_{10,14}f_{11,13}\\
f_{11,14}f_{12,13}f_{12,14}f_{13,14}f_{1,15}f_{2,15}f_{3,15}f_{4,15}\left(f_{5,15}f_{6,15}\right)f_{7,15}f_{8,15}f_{9,15}\left(f_{10,15}f_{11,15}\right)f_{12,15}f_{13,15}f_{14,15}\\
\left(f_{1,16}f_{2,16}\right)f_{3,16}f_{4,16}f_{5,16}\left(f_{6,16}f_{7,16}\right)f_{8,16}f_{9,16}f_{10,16}\left(f_{11,16}f_{12,16}\right)f_{13,16}f_{14,16}f_{15,16}
\end{array}\]
\large

This ordering is shared by the previous six diagrams, hence their associated products all have the same value at $z_1 = \cdots = z_n$. We now complete the chain, starting with our last diagram.

\[
\begin{picture}(80,80)
%\put(-40,36){$\mapsto$}
\put(0,60){\framebox(20,20)[r]{$1\;\;$}}
\put(0,40){\framebox(20,20)[r]{$2\;\;$}}
\put(0,20){\framebox(20,20)[r]{$3\;\;$}}
\put(0,0){\framebox(20,20)[r]{$4\;\;$}}
\put(20,60){\framebox(20,20)[r]{$5\;\;$}}
\put(20,40){\framebox(20,20)[r]{$6\;\;$}}
\put(20,20){\framebox(20,20)[r]{$7\;\;$}}
\put(20,0){\framebox(20,20)[r]{$8\;\;$}}
\put(40,60){\framebox(20,20)[r]{$9\;\;$}}
\put(40,40){\framebox(20,20)[r]{$10\;$}}
\put(40,20){\framebox(20,20)[r]{$11\;$}}
\put(40,0){\framebox(20,20)[r]{$12\;$}}
\put(60,60){\framebox(20,20)[r]{$13\;$}}
\put(60,40){\framebox(20,20)[r]{$14\;$}}
\put(60,20){\framebox(20,20)[r]{$15\;$}}
\put(60,0){\framebox(20,20)[r]{$16\;$}}

%vertical lines up
\put(10,50){\line(0,1){20}}%2
\put(10,30){\line(0,1){20}}%3
\put(10,10){\line(0,1){20}}%4
\put(30,50){\line(0,1){20}}%6
\put(30,30){\line(0,1){20}}%7
%\put(30,10){\line(0,1){20}}%8
%
\put(50,50){\line(0,1){20}}%10
\put(50,30){\line(0,1){20}}%11
\put(50,10){\line(0,1){20}}%12
\put(70,50){\line(0,1){20}}%14
\put(70,30){\line(0,1){20}}%15
\put(70,10){\line(0,1){20}}%16

%Horizontal lines to the right
%\put(10,70){\line(1,0){20}}%1
%\put(10,50){\line(1,0){20}}%2
%\put(10,30){\line(1,0){20}}%3
%\put(10,10){\line(1,0){20}}%4
%
%\put(30,70){\line(1,0){20}}%5
%\put(30,50){\line(1,0){20}}%6
%\put(30,30){\line(1,0){20}}%7
\put(30,10){\line(1,0){20}}%8
%
%\put(50,70){\line(1,0){20}}%9
%\put(50,50){\line(1,0){20}}%10
%\put(50,30){\line(1,0){20}}%11
%\put(50,10){\line(1,0){20}}%12
\end{picture}
\qquad \qquad
\begin{picture}(64,80)
\put(8,64){\circle{16}}%circle 1
\put(32,64){\circle{16}}%circle 2
\put(56,64){\circle{16}}%circle 3
\put(8,40){\circle{16}}%circle 4
\put(32,40){\circle{16}}%circle 5
\put(56,40){\circle{16}}%circle 6
\put(8,16){\circle{16}}%circle 7
\put(32,16){\circle{16}}%circle 8
\put(56,16){\circle{16}}%circle 9

\put(16,64){\line(1,0){8}}%1 to 2
\put(40,64){\line(1,0){8}}%2 to 3
\put(16,40){\line(1,0){8}}%4 to 5
\put(40,40){\line(1,0){8}}%5 to 6
\put(16,16){\line(1,0){8}}%6 to 7
\put(40,16){\line(1,0){8}}%7 to 8

\put(8,48){\line(0,1){8}}%4 to 1
\put(32,48){\line(0,1){8}}%5 to 2
\put(56,48){\line(0,1){8}}%6 to 3
\put(8,24){\line(0,1){8}}%7 to 4
\put(32,24){\line(0,1){8}}%8 to 5
\put(56,24){\line(0,1){8}}%9 to 6

\put(6,61){0}%1
\put(30,61){0}%2
\put(54,61){0}%3
\put(6,37){0}%4
\put(30,37){0}%5
\put(54,37){0}%6
\put(6,13){1}%7
\put(30,13){6}%8
\put(54,13){0}%9
\end{picture}
\qquad \qquad
\begin{picture}(40,80)
\put(8,52){\circle{16}}%circle 1
\put(32,52){\circle{16}}%circle 2
\put(8,28){\circle{16}}%circle 3
\put(32,28){\circle{16}}%circle 4
\put(16,52){\line(1,0){8}}%1 to 2
\put(16,28){\line(1,0){8}}%3 to 4
\put(8,36){\line(0,1){8}}%1 to 3
\put(32,36){\line(0,1){8}}%2 to 4
%\put(16,38){$r$}%circle 1
\put(40,58){$2$}%circle 2
%\put(16,14){$r$}%circle 3
%\put(40,14){$r$}%circle 4
\put(6,49){0}%1
\put(30,49){0}%2
\put(3,25){13}%3
\put(30,25){0}%4
\end{picture}\]

\[
\begin{picture}(80,80)
\put(-40,36){$\mapsto$}
\put(0,60){\framebox(20,20)[r]{$1\;\;$}}
\put(0,40){\framebox(20,20)[r]{$2\;\;$}}
\put(0,20){\framebox(20,20)[r]{$3\;\;$}}
\put(0,0){\framebox(20,20)[r]{$4\;\;$}}
\put(20,60){\framebox(20,20)[r]{$5\;\;$}}
\put(20,40){\framebox(20,20)[r]{$6\;\;$}}
\put(20,20){\framebox(20,20)[r]{$7\;\;$}}
\put(20,0){\framebox(20,20)[r]{$8\;\;$}}
\put(40,60){\framebox(20,20)[r]{$9\;\;$}}
\put(40,40){\framebox(20,20)[r]{$10\;$}}
\put(40,20){\framebox(20,20)[r]{$11\;$}}
\put(40,0){\framebox(20,20)[r]{$12\;$}}
\put(60,60){\framebox(20,20)[r]{$13\;$}}
\put(60,40){\framebox(20,20)[r]{$14\;$}}
\put(60,20){\framebox(20,20)[r]{$15\;$}}
\put(60,0){\framebox(20,20)[r]{$16\;$}}

%vertical lines up
\put(10,50){\line(0,1){20}}%2
\put(10,30){\line(0,1){20}}%3
\put(10,10){\line(0,1){20}}%4
\put(30,50){\line(0,1){20}}%6
\put(30,30){\line(0,1){20}}%7
%\put(30,10){\line(0,1){20}}%8
%
\put(50,50){\line(0,1){20}}%10
\put(50,30){\line(0,1){20}}%11
\put(50,10){\line(0,1){20}}%12
\put(70,50){\line(0,1){20}}%14
\put(70,30){\line(0,1){20}}%15
\put(70,10){\line(0,1){20}}%16

%Horizontal lines to the right
%\put(10,70){\line(1,0){20}}%1
%\put(10,50){\line(1,0){20}}%2
%\put(10,30){\line(1,0){20}}%3
%\put(10,10){\line(1,0){20}}%4
%
%\put(30,70){\line(1,0){20}}%5
%\put(30,50){\line(1,0){20}}%6
%\put(30,30){\line(1,0){20}}%7
%\put(30,10){\line(1,0){20}}%8
%
%\put(50,70){\line(1,0){20}}%9
%\put(50,50){\line(1,0){20}}%10
%\put(50,30){\line(1,0){20}}%11
%\put(50,10){\line(1,0){20}}%12
\end{picture}
\qquad \qquad
\begin{picture}(64,80)
\put(8,64){\circle{16}}%circle 1
\put(32,64){\circle{16}}%circle 2
\put(56,64){\circle{16}}%circle 3
\put(8,40){\circle{16}}%circle 4
\put(32,40){\circle{16}}%circle 5
\put(56,40){\circle{16}}%circle 6
\put(8,16){\circle{16}}%circle 7
\put(32,16){\circle{16}}%circle 8
\put(56,16){\circle{16}}%circle 9

\put(16,64){\line(1,0){8}}%1 to 2
\put(40,64){\line(1,0){8}}%2 to 3
\put(16,40){\line(1,0){8}}%4 to 5
\put(40,40){\line(1,0){8}}%5 to 6
\put(16,16){\line(1,0){8}}%6 to 7
\put(40,16){\line(1,0){8}}%7 to 8

\put(8,48){\line(0,1){8}}%4 to 1
\put(32,48){\line(0,1){8}}%5 to 2
\put(56,48){\line(0,1){8}}%6 to 3
\put(8,24){\line(0,1){8}}%7 to 4
\put(32,24){\line(0,1){8}}%8 to 5
\put(56,24){\line(0,1){8}}%9 to 6

\put(6,61){0}%1
\put(30,61){0}%2
\put(54,61){0}%3
\put(6,37){0}%4
\put(30,37){0}%5
\put(54,37){0}%6
\put(6,13){1}%7
\put(30,13){7}%8
\put(54,13){0}%9
\end{picture}
\qquad \qquad
\begin{picture}(40,80)
\put(8,52){\circle{16}}%circle 1
\put(32,52){\circle{16}}%circle 2
\put(8,28){\circle{16}}%circle 3
\put(32,28){\circle{16}}%circle 4
\put(16,52){\line(1,0){8}}%1 to 2
\put(16,28){\line(1,0){8}}%3 to 4
\put(8,36){\line(0,1){8}}%1 to 3
\put(32,36){\line(0,1){8}}%2 to 4
%\put(16,38){$r$}%circle 1
\put(40,58){$2$}%circle 2
%\put(16,14){$r$}%circle 3
%\put(40,14){$r$}%circle 4
\put(6,49){0}%1
\put(30,49){0}%2
\put(6,25){0}%3
\put(30,25){0}%4
\end{picture}\]

\[
\begin{picture}(80,80)
\put(-40,36){$\mapsto$}
\put(0,60){\framebox(20,20)[r]{$1\;\;$}}
\put(0,40){\framebox(20,20)[r]{$2\;\;$}}
\put(0,20){\framebox(20,20)[r]{$3\;\;$}}
\put(0,0){\framebox(20,20)[r]{$4\;\;$}}
\put(20,60){\framebox(20,20)[r]{$5\;\;$}}
\put(20,40){\framebox(20,20)[r]{$6\;\;$}}
\put(20,20){\framebox(20,20)[r]{$7\;\;$}}
\put(20,0){\framebox(20,20)[r]{$8\;\;$}}
\put(40,60){\framebox(20,20)[r]{$9\;\;$}}
\put(40,40){\framebox(20,20)[r]{$10\;$}}
\put(40,20){\framebox(20,20)[r]{$11\;$}}
\put(40,0){\framebox(20,20)[r]{$12\;$}}
\put(60,60){\framebox(20,20)[r]{$13\;$}}
\put(60,40){\framebox(20,20)[r]{$14\;$}}
\put(60,20){\framebox(20,20)[r]{$15\;$}}
\put(60,0){\framebox(20,20)[r]{$16\;$}}

%vertical lines up
\put(10,50){\line(0,1){20}}%2
\put(10,30){\line(0,1){20}}%3
\put(10,10){\line(0,1){20}}%4
\put(30,50){\line(0,1){20}}%6
\put(30,30){\line(0,1){20}}%7
\put(30,10){\line(0,1){20}}%8
\put(50,50){\line(0,1){20}}%10
\put(50,30){\line(0,1){20}}%11
\put(50,10){\line(0,1){20}}%12
\put(70,50){\line(0,1){20}}%14
\put(70,30){\line(0,1){20}}%15
\put(70,10){\line(0,1){20}}%16

%Horizontal lines to the right
%\put(10,70){\line(1,0){20}}%1
%\put(10,50){\line(1,0){20}}%2
%\put(10,30){\line(1,0){20}}%3
%\put(10,10){\line(1,0){20}}%4
%
%\put(30,70){\line(1,0){20}}%5
%\put(30,50){\line(1,0){20}}%6
%\put(30,30){\line(1,0){20}}%7
%\put(30,10){\line(1,0){20}}%8
%
%\put(50,70){\line(1,0){20}}%9
%\put(50,50){\line(1,0){20}}%10
%\put(50,30){\line(1,0){20}}%11
%\put(50,10){\line(1,0){20}}%12
\end{picture}
\qquad \qquad
\begin{picture}(64,80)
\put(8,64){\circle{16}}%circle 1
\put(32,64){\circle{16}}%circle 2
\put(56,64){\circle{16}}%circle 3
\put(8,40){\circle{16}}%circle 4
\put(32,40){\circle{16}}%circle 5
\put(56,40){\circle{16}}%circle 6
\put(8,16){\circle{16}}%circle 7
\put(32,16){\circle{16}}%circle 8
\put(56,16){\circle{16}}%circle 9

\put(16,64){\line(1,0){8}}%1 to 2
\put(40,64){\line(1,0){8}}%2 to 3
\put(16,40){\line(1,0){8}}%4 to 5
\put(40,40){\line(1,0){8}}%5 to 6
\put(16,16){\line(1,0){8}}%6 to 7
\put(40,16){\line(1,0){8}}%7 to 8

\put(8,48){\line(0,1){8}}%4 to 1
\put(32,48){\line(0,1){8}}%5 to 2
\put(56,48){\line(0,1){8}}%6 to 3
\put(8,24){\line(0,1){8}}%7 to 4
\put(32,24){\line(0,1){8}}%8 to 5
\put(56,24){\line(0,1){8}}%9 to 6

\put(6,61){0}%1
\put(30,61){0}%2
\put(54,61){0}%3
\put(6,37){0}%4
\put(30,37){0}%5
\put(54,37){0}%6
\put(6,13){0}%7
\put(30,13){0}%8
\put(54,13){0}%9
\end{picture}
\qquad \qquad
\begin{picture}(40,80)
\put(8,52){\circle{16}}%circle 1
\put(32,52){\circle{16}}%circle 2
\put(8,28){\circle{16}}%circle 3
\put(32,28){\circle{16}}%circle 4
\put(16,52){\line(1,0){8}}%1 to 2
\put(16,28){\line(1,0){8}}%3 to 4
\put(8,36){\line(0,1){8}}%1 to 3
\put(32,36){\line(0,1){8}}%2 to 4
%\put(16,38){$r$}%circle 1
\put(40,58){$2$}%circle 2
%\put(16,14){$r$}%circle 3
%\put(40,14){$r$}%circle 4
\put(6,49){0}%1
\put(30,49){0}%2
\put(6,25){0}%3
\put(30,25){0}%4
\end{picture}\]

Finally, these last three diagrams share the following ordering and so have the same value at $z_1 = \cdots = z_n$;

\small
\[\begin{array}{l} f_{1,2}f_{1,3}f_{2,3}f_{1,4}f_{2,4}f_{3,4}f_{1,5}f_{2,5}f_{3,5}f_{4,5}\left(f_{1,6}f_{2,6}\right)f_{3,6}f_{4,6}f_{5,6}f_{1,7}\left(f_{2,7}f_{3,7}\right)f_{4,7}f_{5,7}f_{6,7}f_{1,8}f_{2,8}\left(f_{3,8}f_{4,8}\right)\\
f_{5,8}f_{6,8}f_{7,8}f_{1,9}\left(f_{1,10}f_{1,11}\right)f_{1,12}f_{2,9}f_{2,10}\left(f_{2,11}f_{2,12}\right)f_{3,9}f_{3,10}f_{3,11}f_{3,12}f_{4,9}f_{4,10}f_{4,11}f_{4,12}\left(f_{5,9}f_{5,10}\right)\\
f_{5,11}f_{5,12}f_{6,9}\left(f_{6,10}f_{6,11}\right)f_{6,12}f_{7,9}f_{7,10}\left(f_{7,11}f_{7,12}\right)f_{8,9}f_{8,10}f_{8,11}f_{8,12}f_{9,10}f_{9,11}f_{9,12}f_{10,11}f_{10,12}f_{11,12}\\
f_{1,13}f_{2,13}f_{3,13}f_{4,13}f_{5,13}f_{6,13}f_{7,13}f_{8,13}f_{10,13}f_{9,13}f_{11,13}f_{12,13}f_{1,14}f_{2,14}f_{3,14}f_{4,14}f_{5,14}f_{6,14}f_{7,14}f_{8,14}\\
\left(f_{9,14}f_{10,14}\right)f_{11,14}f_{12,14}f_{13,14}f_{1,15}f_{2,15}f_{3,15}f_{4,15}\left(f_{5,15}f_{6,15}\right)f_{7,15}f_{8,15}f_{9,15}\left(f_{10,15}f_{11,15}\right)f_{12,15}f_{13,15}\\
f_{14,15}\left(f_{1,10}f_{1,16}\right)f_{2,16}f_{3,16}f_{4,16}f_{5,16}\left(f_{6,16}f_{7,16}\right)f_{8,16}f_{9,16}f_{10,16}\left(f_{11,16}f_{12,16}\right)f_{13,16}f_{14,16}f_{15,16}
\end{array}\]
\large

In particular the last diagram is the column diagram used in the column fusion procedure. We know the value of its associated product is the diagonal matrix element at $z_1 = \cdots = z_n$, and hence this is the value in the limit for all the diagrams in the chain.
\end{example}

We now go some way to proving the opposite implication of Conjecture \ref{ribbonconjecture}, that is to say, if the defining diagram is invalid then the associated product is not the diagonal matrix element in the limit $z_1 = \cdots = z_n$.

Under the conditions imposed in the row, column and hook fusion procedures we found we had a removable singularity along the line $z_1 = \cdots = z_n$. For rational polynomials, a singularity is removable exactly when there is an (improper) cancellation in the numerator and the denominator. We now have the following immediate consequence:

\begin{proposition}If the diagram $\Lambda$ contains any $2\times2$ square of type \emph{$\bigtypecentre{8}$} to \emph{$\bigtype{14}$} then the product $F_\Lambda(z_1, \dots, z_n)$ has a pole at $z_1 = \cdots = z_n$. \end{proposition}

\newpage
\begin{proof} If $(p,q)$ is a singularity of degree 1 and of type $\bigtypecentre{8}$ to $\bigtype{14}$ then $z_p = z_q$ and there is no possible cancellation of $(z_p - z_q)$ on the subspace $\mathcal{P}_\Lambda$, hence the singularity is not removable, but is instead a pole. \end{proof}

Let us now consider the other kinds of invalid diagram, diagrams that contain a square of type $\bigtype{15}$ or contain a union excluded by Definition \ref{validdiagrams}. Taking explicit examples, we identified a coefficient likely to have a pole at $z_1 = \cdots = z_n$, and then used the software packages GAP and Maple to show if this was the case. These examples support the claim that the associated products $F_\Lambda(z_1, \dots, z_n)$ of invalid diagrams have a pole at $z_1 = \cdots = z_n$. We give a couple of such examples below, for details see Appendix \ref{appendix b}.

Consider the column tableau, $\Lambda^c$, of the Young diagram $\lambda = (3,3)$. We proceed to calculate the coefficient of one of the elements in the expansion of $F_{\Lambda^c}$, and show this coefficient has a pole in the limit $z_1 = \cdots = z_n$. By observation we identify the permutation (1 4 3 6) as one such possible permutation.

For comparison let us first consider the subspace defined in the row fusion procedure by the diagram

\[
\begin{picture}(60,40)
\put(0,20){\framebox(20,20)[r]{$1\;\;$}}
\put(0,0){\framebox(20,20)[r]{$2\;\;$}}
\put(20,20){\framebox(20,20)[r]{$3\;\;$}}
\put(20,0){\framebox(20,20)[r]{$4\;\;$}}
\put(40,20){\framebox(20,20)[r]{$5\;\;$}}
\put(40,0){\framebox(20,20)[r]{$6\;$}}

\put(10,30){\line(1,0){20}}%1
\put(10,10){\line(1,0){20}}%2
\put(30,30){\line(1,0){20}}%5
\put(30,10){\line(1,0){20}}%6
\end{picture}\]

Maple calculates the coefficient of (1 4 3 6) to be

\[\frac{2z_1^3 + 10z_1^2 + 12z_1 - 6z_1^2z_2  - 20z_1z_2 + 6z_1z_2^2 -2z_2^3 +10z_2^2- 12z_2 -1}{(z_1-z_2+1)^3(z_1-z_2+2)^2}.\] 

So on the subspace where $z_1=z_3=z_5$, $z_2=z_4=z_6$ and $z_1 - z_2 \notin \mathbb{Z}$ there has been a cancellation of the term $(z_1-z_2)$ in the numerator and denominator, hence the singularity is removable and the coefficient of (1 4 3 6) tends to -1/4 in the limit $z_1 = z_2$.

Now consider an invalid diagram that contains a square of type $\bigtype{15}$

\[
\begin{picture}(60,40)
\put(0,20){\framebox(20,20)[r]{$1\;\;$}}
\put(0,0){\framebox(20,20)[r]{$2\;\;$}}
\put(20,20){\framebox(20,20)[r]{$3\;\;$}}
\put(20,0){\framebox(20,20)[r]{$4\;\;$}}
\put(40,20){\framebox(20,20)[r]{$5\;\;$}}
\put(40,0){\framebox(20,20)[r]{$6\;$}}

\put(10,30){\line(1,0){20}}%1
\put(10,10){\line(1,0){20}}%2
\end{picture}\]

We also calculate the coefficient of (1 4 3 6) on this subspace. The numerator of which is too large to reproduce here. However, importantly, the denominator contains the term $(z_1-z_6)$ and so the singularity is not removable.

%\newpage
Finally we consider one more invalid diagram, a diagram that contains the union \type{3}--\type{5}

\[
\begin{picture}(60,40)
\put(0,20){\framebox(20,20)[r]{$1\;\;$}}
\put(0,0){\framebox(20,20)[r]{$2\;\;$}}
\put(20,20){\framebox(20,20)[r]{$3\;\;$}}
\put(20,0){\framebox(20,20)[r]{$4\;\;$}}
\put(40,20){\framebox(20,20)[r]{$5\;\;$}}
\put(40,0){\framebox(20,20)[r]{$6\;$}}

\put(10,30){\line(1,0){20}}%1
\put(30,10){\line(1,0){20}}%6
\end{picture}\]

We use Maple again to calculate the coefficient of (1 4 3 6) on this subspace. Here the denominator contains the term $(z_1-z_4)$, which again shows the singularity is not removable.

%% file: TheHookFusionProcedureForHeckeAlgebras.tex
\section{The hook fusion procedure for Hecke algebras}

The Hecke algebra of type $A$ is a quantum deformation of the symmetric group algebra. We adapt our hook fusion procedure of Chapter 3 to this new setting to construct an element that generates irreducible representations of this deformation. This is also the content of the author's second paper \cite{Gr2}.

\subsection{Hecke algebras and affine Hecke algebras}

We begin by defining the Hecke and the affine Hecke algebra of type $A$, for a description of these algebras of other types see \cite{C3}. Let $H_n$ be the finite 
dimensional Hecke algebra over the
field $\mathbb{C}(q)$ of rational functions in $q$,
with the generators $T_1, \ldots , T_{n-1}$ and the relations 
\begin{equation}\label{Hecke1}
(T_i-q)(T_i+q^{-1})=0;
\end{equation}
\begin{equation}\label{Hecke2}
T_i T_{i+1} T_i=T_{i+1} T_i T_{i+1};
\end{equation}
\begin{equation}\label{Hecke3}
T_i T_j=T_j T_i,
\quad
j\neq i, i+1
\end{equation}
for all possible indices $i$ and $j$.

These generators are invertible in $H_n$ since
\begin{equation}\label{Tinverse}
T^{-1}_i=T_i-q+q^{-1}
\end{equation}
due to (\ref{Hecke1}).

As usual let $\sigma_i=(i, i+1)$ be 
the adjacent transposition in the symmetric group $S_n$. 
Take any element $\sigma \in S_n$
and choose a reduced decomposition $\sigma=\sigma_{i_1}\ldots \sigma_{i_l}$ and set $T_\sigma=T_{i_1}\ldots T_{i_l}$.
This element of the algebra $H_n$ does not depend on the
choice of reduced decomposition of $\sigma$ due to (\ref{Hecke2}) and
(\ref{Hecke3}).  We denote the element of maximal length in
% the symmetric group 
$S_n$ by $\sigma_0$, but will write $T_0$ instead of 
$T_{\sigma_0}$ for  short. The elements $T_\sigma$ form a basis of $H_n$ 
as a vector space over the field $\mathbb{C}(q)$. We will also use the 
basis in $H_n$  formed by the elements $T_\sigma^{-1}$.

If $q$ is not a root of unity then the $\mathbb{C}(q)$-algebra $H_n$ is semisimple, for a short proof see \cite[Section 4]{GU}. 
Like irreducible representations of the symmetric group $S_n$, the simple ideals of $H_n$ are
labeled by partitions $\lambda$ of $n$, 
In this chapter we 
will construct
a certain non-zero element $F_\Lambda\in H_n$ for any standard tableau $\Lambda$ of shape $\lambda$. The element $F_\Lambda$ is related to the $q$-analogue
of the Young symmetrizer in the group ring $\mathbb{C} S_n$
constructed in \cite{Gy}. Under left multiplication by
the elements of $H_n$, the left ideal $H_nF_\Lambda\subset H_n$ is an 
irreducible $H_n$-module. The $H_n$-modules $V_\lambda$
for different partitions $\lambda$ are pairwise non-equivalent, see
Corollary \ref{C3.6}. Note, at $q=1$, the algebra $H_n(q)$ specialises
to the group ring $\mathbb{C}S_n$, where $T_\sigma$ becomes the permutation $\sigma \in S_n$.  The $H_n(q)$-module
$V_\lambda$ then specialises to the irreducible representation of 
$S_n$, corresponding to the partition $\lambda$.

%%%%%%%%%%%%%%%%%%%%%%%%%%%%%%%%%%%%%%%%%%%%%%%%%%%%%%%%%%%%%%%%%%%%%

Let us also consider the \emph{affine Hecke algebra} $\widehat{H}_n$ is the $\mathbb{C}(q)$-algebra generated by the elements $T_1, \dots, T_{n-1}$ and the pairwise commuting invertible elements $Y_1, \dots, Y_n$ subject to the relations (\ref{Hecke1})--(\ref{Hecke3}) and
\begin{eqnarray}\label{affinehecke} T_i Y_j &=& Y_j T_i, \quad j\neq i, i+1, \\
T_iY_iT_i &=& Y_{i+1}. \end{eqnarray}

By definition, the affine Hecke algebra $\widehat{H}_n$ contains $H_n$ as a subalgebra. And for any nonzero $z \in \mathbb{C}(q)$, one can define a homomorphism $\pi_z : \widehat{H}_n \to H_n$, identical on the subalgebra $H_n \subset \widehat{H}_n$, defined by $\pi_z(Y_1) = z$. It is completely defined due to (\ref{affinehecke}). By pulling any irreducible $H_n$-module $V$ back through the homomorphism $\pi_z$ we obtain an irreducible module over the algebra $\widehat{H}_n$, called an \emph{evaluation module} at $z$ and denoted $V(z)$.

Specifically, setting $z=1$ we have the homomorphism $\pi_1 : \widehat{H}_n \to H_n$, identical on the subalgebra $H_n \subset \widehat{H}_n$, defined by $\pi_1(Y_1) =1$. Denote by $X_i$ the
expression \[ \pi_1(Y_i) = T_{i-1}\ldots T_1 T_1\ldots T_{i-1}.\] 
Using the relations (\ref{Hecke2}), (\ref{Hecke3}), and (\ref{Tinverse})
one can check that 
and that the elements $X_1, \ldots , X_n$ are
pairwise commuting and invertible in $H_n$.  The elements $X_1, \ldots , X_n$ are called
the \emph{Murphy elements} \cite{M2} of the Hecke algebra
$H_n$.

%%%%%%%%%%%%%%%%%%%%%%%%%%%%%%%%%%%%%%%%%%%%%%%%%%%%%%%%%%%%%%%%%

To motivate the study of $H_n$-modules of
hook shape let us consider certain identities of representations of the Hecke algebra. Again we can think of the following
identities as dual to the Jacobi-Trudi identities.

If $\lambda = (\lambda_1, \dots , \lambda_k)$ such that $n =
\lambda_1 + \dots + \lambda_k$ then we have the following
decomposition of the induced representation of the tensor product
of modules corresponding to the rows of $\lambda$;
\begin{equation}\label{qinducedidentity} \textrm{Ind}_{H_{\lambda_1} \times H_{\lambda_2} \times \cdots
\times H_{\lambda_k}}^{H_n } V_{(\lambda_1)} \otimes
V_{(\lambda_2)} \otimes \cdots \otimes V_{(\lambda_k)} \cong
\bigoplus_{\mu} (V_\mu)^{\oplus K_{\mu \lambda}},\end{equation} where the sum
is over all partitions of $n$. Here $V_{(\lambda_i)}$ denotes the
trivial representation of $H_{\lambda_i}$, that sends generators $T_i$ to $q$. Again, the coefficients
$K_{\mu \lambda}$ are Kostka
numbers, with $K_{\lambda \lambda} =
1.$

Similarly we have the equivalent identity for columns, \begin{equation}\label{qinducedidentity2}
\textrm{Ind}_{H_{\lambda'_1} \times H_{\lambda'_2} \times \cdots
\times H_{\lambda'_l}}^{H_n } V_{(1^{\lambda'_1})} \otimes
V_{(1^{\lambda'_2})} \otimes \cdots \otimes V_{(1^{\lambda'_l})}
\cong \bigoplus_{\mu} (V_{\mu})^{\oplus K_{\mu' \lambda'}},\end{equation}
where $l$ is the number of columns of $\lambda$. In this case
$V_{(1^{\lambda'_i})}$ is the alternating representation of
$H_{\lambda'_i}$.

However, we are interested in the following identity related to the the Giambelli identity and the principal hooks of $\lambda$.
\begin{equation}\label{qinducedidentity3} \textrm{Ind}_{H_{h_1} \times H_{h_2} \times \cdots
\times H_{h_d}}^{H_n} V_{(\alpha_1 | \beta_1)} \otimes
V_{(\alpha_2 | \beta_2)} \otimes \cdots \otimes V_{(\alpha_d |
\beta_d)} \cong \bigoplus_{\mu} (V_{\mu})^{\oplus D_{\mu
\lambda}},\end{equation} where $h_i$ is the length of the $i^{\mbox{\scriptsize th}}$
principal hook, and the sum is over all partitions of $n$. This is
a decomposition of the induced representation of the tensor
product of modules of hook shape. Further these hooks are the
principal hooks of $\lambda$. The coefficients, $D_{\mu \lambda}$,
are non-negative integers, and in particular $D_{\lambda \lambda}
=1$. The hook fusion procedure for Hecke algebras will provide a way of singling out the component $V_\lambda$ in this decomposition.

%%%%%%%%%%%%%%%%%%%%%%%%%%%%%%%%%%%%%%%%%%%%%%%%%%%%%%%%%%%%%%%%%%%%%%%%%%%%%%%%%

\subsection{The hook fusion procedure for Hecke algebras}

For each $i=1, \ldots , n-1$ introduce the $H_n$-valued rational function in 
two variables $a \neq 0$, $a \neq b \in\mathbb{C}(q)$
\begin{equation}\label{q-smallf}
F_i(a, b)=T_i+\frac{q-q^{-1}}{a^{-1}b - 1}.
\end{equation}

Now introduce $n$ variables $z_1, \ldots , z_n\in\mathbb{C}(q)$.
Equip the set of all pairs $(i, j)$
where $1\leqslant i<j\leqslant n$, with the following ordering.
The pair $(i, j)$ precedes another pair
$(i',j')$ % from the same set
if $j<j'$, or if $j=j'$ but $i<i'$. Call this the \emph{reverse-lexicographic} ordering.
Take the ordered product
\begin{equation}\label{q-bigf}
F_\Lambda(z_1, \ldots , z_n) = \prod_{(i,j)}^{\rightarrow}\ F_{j-i}
\bigl( q^{2c_i(\Lambda)}z_i, q^{2c_j(\Lambda)}z_j\bigr)
\end{equation}
over this set. Consider the product (\ref{q-bigf})
as a rational function taking values in $H_n$,
of the variables $z_1, \ldots , z_n$.
Note, if $i$ and $j$ sit on
the same diagonal in the tableau $\Lambda$, then $F_{j-i}
\bigl( q^{2c_i(\Lambda)}z_i, q^{2c_j(\Lambda)}z_j\bigr)$ has a pole at $z_i = z_j \neq 0$.

Let $\mathcal{R}_\Lambda$ be the vector subspace in $\mathbb{C}(q)^n$
consisting of all tuples $(z_1, \dots , z_n)$ such that $z_i =
z_j$ whenever the numbers $i$ and $j$ appear in the same row of
the tableau $\Lambda$.

As a direct calculation using (\ref{Hecke1}) and (\ref{Hecke2}) shows,
these functions satisfy 
\begin{equation}\label{q-triple}
F_i(a, b)F_{i+1}(a, c)F_i(b, c)=
F_{i+1}(b, c)F_i(a, c)F_{i+1}(a, b).
\end{equation}
Due to (\ref{Hecke3}) these rational functions also satisfy the relations
\begin{equation}\label{q-commute}
F_i(a, b) F_j(c, d)=F_j(c, d) F_i(a, b);
\qquad
j\neq i, i+1.
\end{equation}

Using (\ref{q-triple}) and (\ref{q-commute}) we may reorder the
product $F_\Lambda(z_1, \dots , z_n)$ such that each singularity
is contained in an expression known to be regular at $z_1 = z_2 =
\dots = z_n \neq 0$, \cite{N2}. It is by this method that it was shown
that the restriction of the rational function $F_\Lambda(z_1,
\dots ,z_n)$ to the subspace $\mathcal{R}_\Lambda$ is regular at
$z_1 = z_2 = \cdots = z_n \neq 0$. Furthermore the following theorem was proved;

\begin{theorem}[Nazarov] \label{q-MNtheorem} Restriction to $\mathcal{R}_\Lambda$ of the
rational function $F_\Lambda(z_1, \dots , z_n)$ is regular at $z_1
= z_2 = \cdots = z_n \neq 0$ and has value $F_\Lambda \in H_n$. The left ideal generated by this element is irreducible, and the $H_n$-modules for different partitions $\lambda$ are pairwise non-equivalent.
\end{theorem}

So, on the subspace $\mathcal{R}_\Lambda$, if $z_i / z_j \notin
q^\mathbb{Z}$ when $i$ and $j$ are in different rows of $\Lambda$,
the induced module (\ref{qinducedidentity}) may be realised as the left ideal in
$H_n$ generated by $F_\Lambda(z_1, \dots, z_n)$. Hence Theorem \ref{q-MNtheorem} provides a way of singling out the 
irreducible component $V_\lambda$ after taking the limit $z_1 = z_2 = \cdots = z_n \neq 0$.

Similarly, we may form another expression for $F_\Lambda$ by
considering the subspace in $\mathbb{C}(q)^n$ consisting of all
tuples $(z_1, \dots , z_n)$ such that $z_i = z_j$ whenever the
numbers $i$ and $j$ appear in the same column of the tableau
$\Lambda$ \cite{N1}. Again, the irreducible representation $V_{\lambda}$ that
appears in the decomposition (\ref{qinducedidentity2}) may be formed from the ideal generated by
$F_{\Lambda}(z_1, \dots, z_n)$ when $z_1 = z_2 = \cdots = z_n \neq 0$.

Let $\mathcal{H}_\Lambda$ be the vector subspace in $\mathbb{C}(q)^n$
consisting of all tuples $(z_1, \dots , z_n)$ such that $z_i =
z_j$ whenever the numbers $i$ and $j$ appear in the same principal
hook of the tableau $\Lambda$. We will prove the following
theorem.

\begin{theorem}\label{q-fulltheorem} Restriction to $\mathcal{H}_\Lambda$ of the
rational function $F_\Lambda(z_1, \dots , z_n)$ is regular at $z_1
= z_2 = \cdots = z_n \neq 0$ and has value $F_\Lambda \in H_n$. The left ideal generated by this element is irreducible, and the $H_n$-modules for different partitions $\lambda$ are pairwise non-equivalent. The element is the same as the element in Theorem \ref{q-MNtheorem}.
\end{theorem}

In particular, this hook fusion procedure can be used to form
irreducible representations of $H_n$ corresponding to Young
diagrams of hook shape using only one auxiliary parameter, $z$. By
taking this parameter to be 1 we find that no parameters are
needed for diagrams of hook shape. Therefore if $\nu$ is a
partition of hook shape, and $N$ a standard tableau of shape $\nu$, we have
\begin{equation}\label{q-bigfhook} F_N = F_N(z) =
\prod_{(p,q)}^\rightarrow F_{j-i}(q^{2c_i(\Lambda)}, q^{c_j(\Lambda)}), \end{equation} with
the pairs $(i, j)$ in the product ordered reverse-lexicographically.

\begin{example} Let $N = \young(13,2)$ \phantom{x} as before. Then on the subspace $\mathcal{H}_N$ we have,
\begin{eqnarray*}F_N &=& \left(T_1 + \frac{q-q^{-1}}{q^{-2}-1} \right)\left(T_2 + \frac{q-q^{-1}}{q^2-1}\right)\left(T_1 + \frac{q-q^{-1}}{q^4-1}\right) \\ &\phantom{=}& \\
&=& T_1T_2T_1 + \frac{1}{q(q^2+1)} T_1T_2 - q T_2T_1 - \frac{1}{q^2+1} T_1 - \frac{1}{q^2 + 1} T_2 + \frac{q}{q^2+1}.
\end{eqnarray*}

When $q=1$ this element becomes \[(1\phantom{x}3) + \frac{(1\phantom{x}2\phantom{x}3)}{2} - (1\phantom{x}3\phantom{x}2) - \frac{(1\phantom{x}2)}{2} - \frac{(2\phantom{x}3)}{2} + \frac{1}{2}1.\]

Notice this element is not the diagonal matrix element described in Example 2.6, but rather it is the diagonal matrix element multiplied on the right by the longest element $\sigma_0 \in S_n$, which, in this case, is the element $\sigma_1\sigma_2\sigma_1 = (1\phantom{x}3) \in S_3$.
\end{example}

%%%%%%%%%%%%%%%%%%%%%%%%%%55

On the subspace $\mathcal{H}_\Lambda$, if $z_i / z_j \notin
q^\mathbb{Z}$ when $i$ and $j$ are in different principal hooks of
$\Lambda$, the induced module (\ref{qinducedidentity}) may be realised as the
left ideal in
$H_n$ generated by $F_\Lambda(z_1, \dots, z_n)$. The irreducible representation $V_\lambda$ appears in the
decomposition of this induced module with coefficient 1, and is
the ideal of $H_n$ generated by $F_\Lambda(z_1, \dots ,
z_n)$ when $z_1 = z_2 = \cdots = z_n \neq 0$. So, as promised, our hook fusion procedure 
will provide a way of singling out the irreducible component
$V_\lambda$ from the this induced module.

%%%%%%%%%%%%%%%%%%%%%%%%%%%%%%%%%%%%%%%%%%%%%%%%%%%%%%%%%%%%%%%%%%%%%

Consider (\ref{q-bigf}) as a rational function of the variables
$z_1, \dots , z_n$ with values in $H_n$. Using the substitution \begin{equation}\label{substitution} w_i = q^{c_i(\Lambda)}z_i, \end{equation} the factor
$F_{i}(w_a, w_b)$ has a pole at $z_a = z_b$ if and
only if the numbers $a$ and $b$ stand on the same diagonal of a
tableau $\Lambda$. So, as in Chapter 3, let us again call the pair $(a, b)$ a singularity, and the corresponding term $F_{i}(w_a,w_b)$ a singularity term, or simply a singularity.

Let $a$ and $b$ be in the same principal hook of $\Lambda$. If $a$ and $b$
are next to one another in the column of the hook then, on
$\mathcal{H}_\lambda$, $F_{i}(w_a,w_b) = T_{i} - q$. Since \begin{equation}\label{P-}(T_i - q)^2 = -(q+q^{-1})(T_i-q)\end{equation}
then $\frac{-1}{q+q^{-1}} F_{i}(w_a,w_b)$ is an
idempotent. Denote this idempotent $P^-_i$.

Similarly, if $a$ and $b$ are next to one another in
the same row of the hook then $F_{i}(w_a,w_b)= T_{i} +
q^{-1}$. And since \begin{equation}\label{P+}(T_i + q^{-1})^2 = (q+q^{-1})(T_i+q^{-1})\end{equation} then
$\frac{1}{q+q^{-1}} F_{i}(w_a,w_b)$ is an idempotent. Denote this idempotent $P^+_i$.\\

We also have
\begin{equation}\label{q-inverse} F_i(a, b) F_i(b, a)=1-\frac{(q-q^{-1})^2 ab}{(a-b)^2}.
\end{equation}
Therefore, if the contents $c_a(\Lambda)$ and $c_b(\Lambda)$ differ by a number
greater than one, then the factor $F_{i}(w_a,w_b)$
is invertible in $H_n$ when $z_a = z_b \neq 0$ for all values of $q$.

The presence of singularity terms in the product $F_\Lambda(z_1,
\dots , z_n)$ mean this product may or may not be regular on the
vector subspace of $\mathcal{H}_\lambda$ consisting of all tuples
$(z_1, \dots , z_n)$ such that $z_1 = z_2 = \cdots = z_n \neq 0$. Using
the following lemma, we will be able to show that $F_\Lambda(z_1,
\dots , z_n)$ is indeed regular on this subspace.

\begin{lemma}\label{q-regular} Restriction of the rational function $F_i(a, b)F_{i+1}(a, c)F_i(b, c)$ to the set of
$(a, b, c)$ such that $a=q^{\pm 2}b$, 
is regular at $a = c \neq 0$. \end{lemma}
\begin{proof} If we write $F_{i+1}(a,c)$ in full, one may expand the left hand side of (\ref{q-triple}) and obtain the sum 
\[
F_i(a, b)T_{i+1}F_i(b, c)+
\frac{q-q^{-1}}{a^{-1}c-1}F_i(a, b)F_i(b, c).
\]
Here the restriction to $a=q^{\pm 2}b$ of the first summand is
evidently regular at $a=c$. After the substitution 
$b=q^{\mp2}a$, the second summand takes the form
\[
\frac{q-q^{-1}}{a^{-1}c-1}
\bigl( T_i\mp q^{\pm1}\bigr)
\biggl(
T_i+\frac{q-q^{-1}}{q^{\pm2}a^{-1} c-1}
\biggr)
=
\frac{q-q^{-1}}{a^{-1} c-q^{\mp2}}
( q^{\pm1}\mp T_i).
\]
The rational function of $a, c$ at the right hand side
of the last displayed equality is also evidently regular at $a=c$. \end{proof}

In particular, if the middle term on the left hand side of (\ref{q-triple}) is a singularity
and the other two terms are an appropriate idempotent and triple term, then this three term product is regular at $z_1 = z_2 =
\dots = z_n \neq 0$. we may now prove the first statement of Theorem
\ref{q-fulltheorem}.

\begin{proposition}\label{q-jimtheorem1} The restriction of the rational function
$F_\Lambda (z_1, \dots , z_n)$ to the subspace
$\mathcal{H}_\lambda$ is regular at $z_1 = z_2 = \cdots = z_n \neq 0$.
\end{proposition}

\begin{proof}
Consider any standard tableau $\Lambda'$ obtained from the tableau $\Lambda$
by an adjacent transposition of its entries, say by $\sigma_k\in S_n$.
Using the relations (\ref{q-triple}) and (\ref{q-commute}), we derive
the equality of rational functions in the variables $z_1, \ldots , z_n$
\[
F_\Lambda(z_1, \ldots , z_n)F_{n-k}
\bigl( 
q^{2c_{k+1}(\Lambda)}z_{k+1}, q^{2c_k(\Lambda)}z_k
\bigr)
=
\]
\begin{equation}\label{q-jimtheorem1equation}
F_k
\bigl( 
q^{2c_k(\Lambda)}z_k, q^{2c_{k+1}(\Lambda)}z_{k+1}
\bigr)
F_{\Lambda^{\prime}}(z'_1, \ldots ,z'_n),
\end{equation}
where the sequence of variables $(z'_1, \ldots ,z'_n)$\ is obtained from
the sequence $(z_1, \ldots , z_n)$ by exchanging the terms $z_k$ and 
$z_{k+1}$. Observe that
\[
(z'_1, \ldots ,z'_n)\in\mathcal{H}_{\Lambda'}
\quad\Leftrightarrow\quad
(z_1, \ldots , z_n)\in\mathcal{H}_\Lambda.
\]
Also observe that here $| c_k(\Lambda)-c_{k+1}(\Lambda)|\geqslant2$
because the tableaux $\Lambda$ and $\Lambda'$ are standard.
Therefore the functions 
\[
F_k
\bigl( 
q^{2c_k(\Lambda)}z_k, q^{2c_{k+l}(\Lambda)}z_{k+1}
\bigr)
\ \quad\textrm{and}\ \quad
F_{n-k}
\bigl( 
q^{2c_{k+1}(\Lambda)}z_{k+1}, q^{2c_k(\Lambda)}z_k
\bigr)
\]
appearing in the equality (\ref{q-jimtheorem1equation}),
are regular at $z_k=z_{k+1} \neq 0$.
Moreover, their values at $z_k=z_{k+1} \neq 0$ are invertible
in the algebra $H_n$, see the relation (\ref{q-inverse}). 
Due to these two observations, the equality (\ref{q-jimtheorem1equation})
shows that Proposition \ref{q-jimtheorem1} is equivalent to its counterpart for
the tableau $\Lambda'$ instead of $\Lambda$.

Let us take the hook tableau $\Lambda^\circ$ of shape $\lambda$, as described in Section \ref{hookfusionprocedure}.  
There is a chain $\Lambda,\Lambda', \ldots ,\Lambda^\circ$ of standard tableaux
of the same shape $\lambda$, such that each subsequent tableau in the 
chain is
obtained from the previous one by an adjacent transposition of the 
entries.
Due to the above argument, it now suffices to prove Proposition \ref{q-jimtheorem1} 
only in the case $\Lambda=\Lambda^\circ$.

We will prove the statement by reordering the factors of the
product \\$F_{\Lambda^\circ} (z_1, \dots , z_n)$, using relations
(\ref{q-triple}) and (\ref{q-commute}), in such a way that each
singularity is part of a triple which is regular at $z_1 = z_2 =
\dots = z_n \neq 0$, and hence the whole of $F_{\Lambda^\circ} (z_1, \dots ,
z_n)$ will be manifestly regular.

Define $g(a,b)$ to be the following; \[ g(a,b) = \left\{
\begin{array}{ccc}
  F_{b-a} (w_a, w_b) & \textrm{if} & a<b \\
  1 & \textrm{if} & a \geqslant b
\end{array} \right. \]
where $w_i$ is the substitution (\ref{substitution}).

Now, let us divide the diagram $\lambda$ into two parts,
consisting of those boxes with positive contents and those with
non-positive contents as in Figure \ref{steps}. Consider the
entries of the $i^{\mbox{\small th}}$ column of the hook tableau
$\Lambda^\circ$ of shape $\lambda$ that lie below the steps. If $u_1,
u_2, \dots , u_k$ are the entries of the $i^{\mbox{\small th}}$ column
below the steps, we define
\begin{equation}\label{q-cproduct} C_i = \prod_{j=1}^n g(u_1 , j)
g(u_2 , j) \cdots g(u_k , j).
\end{equation} Now consider the entries of the $i^{\mbox{\small th}}$ row of $\Lambda^\circ$ that lie
above the steps. If $v_1, v_2, \dots , v_l$ are the entries of the
$i^{\mbox{\small th}}$ row above the steps, we define
\begin{equation}\label{q-rproduct} R_i = \prod_{j=1}^n g(v_1 , j)
g(v_2 , j) \cdots g(v_l , j).
\end{equation}

Our choice of the hook tableau was such that the following is
true; if $d$ is the number of principal hooks of $\lambda$ then by
relations (\ref{q-triple}) and (\ref{q-commute}) we may reorder the
factors of $F_{\Lambda^\circ} (z_1, \dots , z_n)$ such that \[ F_{\Lambda^\circ}
(z_1, \dots , z_n) = \prod_{i=1}^d C_i R_i . \]

Now, each singularity $(a,b)$ has its corresponding term $F_{b-a}
(w_a,w_b)$ contain in some product $C_i$ or $R_i$.
This singularity term will be on the immediate left of the
 term $F_{b-a-1} (w_{a+1}, w_b)$. Also, this ordering has been chosen such
that the product of factors to the left of any such singularity in
$C_i$ or $R_i$ is divisible on the right by $F_{b-a-1}(w_a,w_{b+1})$.\\
Therefore we can replace the pair $F_{b-a}
(w_a,w_b)F_{b-a-1} (w_{a+1}, w_b)$ in $C_i$ by the triple
\[ P_{b-a-1}^- F_{b-a}
(w_a,w_b)F_{b-a-1} (w_{a+1}, w_b), \] where
$P_{b-a-1}^- = \frac{-1}{q+q^{-1}} F_{b-a-1}(w_a,w_{a+1})$ is the
idempotent (\ref{P-}). Divisibility on the right by $F_{b-a-1}(w_a,w_{a+1})$ means the addition of the idempotent has no
effect on the value of the product $C_i$. Similarly, in the product $R_i$ we can replace the pair by \[ P_{b-a-1}^+ F_{b-a}
(w_a,w_b)F_{b-a-1} (w_{a+1}, w_b). \] 
By Lemma \ref{q-regular}, the above triples are regular at $z_1 =
z_2 = \cdots = z_n \neq 0$, and therefore, so too are the products $C_i$
and $R_i$, for all $1 \leqslant i \leqslant d$. Moreover, this
means $F_{\Lambda^\circ} (z_1, \dots , z_n)$ is regular at $z_1 = z_2 =
\dots = z_n \neq 0$.
\end{proof}

{\addtocounter{definition}{1} \bf Example \thedefinition .} As an
example consider the hook tableau of the Young diagram $\lambda =
(3,3,2)$, as shown in Example 3.2.

In the original reverse-lexicographic ordering the product $F_{\Lambda^{\circ}}(z_1,
\dots , z_n)$ is written as;
\begin{normalsize}
\[
\begin{array}{rl}
   F_{\Lambda^\circ} (z_1, \dots , z_n) = & F_1(w_1,w_2)F_2(w_1,w_3)F_1(w_2,w_3)F_3(w_1,w_4)F_2(w_2,w_4)F_1(w_3,w_4)\\&
   F_4(w_1,w_5)F_3(w_2,w_5)F_2(w_3,w_5)F_1(w_4,w_5)F_5(w_1,w_6)F_4(w_2,w_6)\\&F_3(w_3,w_6)F_2(w_4,w_6)F_1(w_5,w_6)F_6(w_1,w_7)F_5(w_2,w_7)F_4(w_3,w_7)\\&F_3(w_4,w_7)F_3(w_4,w_7)F_2(w_5,w_7)F_1(w_6,w_7)F_7(w_1,w_8)F_6(w_2,w_8)\\&F_5(w_3,w_8)F_4(w_4,w_8)F_3(w_5,w_8)F_2(w_6,w_8)F_1(w_7,w_8)\\
\end{array}
\]\end{normalsize}
we may now reorder this product into the form below using
relations (\ref{q-triple}) and (\ref{q-commute}) as described in the
above proposition. The terms bracketed are the singularity terms
with their appropriate triple terms.
\begin{normalsize}\[
\begin{array}{rl}
   F_{\Lambda^\circ} (z_1, \dots , z_n) = &
   F_1(w_1,w_2)F_2(w_1,w_3)F_1(w_2,w_3)F_3(w_1,w_4)F_2(w_2,w_4)F_1(w_3,w_4)\\&F_4(w_1,w_5)F_3(w_2,w_5)F_2(w_3,w_5)\left(F_5(w_1,w_6)F_4(w_2,w_6)\right)F_3(w_3,w_6)\\&F_6(w_1,w_7)\left(F_5(w_2,w_7)F_4(w_3,w_7)\right)F_7(w_1,w_8)F_6(w_2,w_8)F_5(w_3,w_8)\\&
   \cdot F_1(w_4,w_5)F_2(w_4,w_6)F_1(w_5,w_6)F_3(w_4,w_7)F_2(w_5,w_7)\\&\left(F_4(w_4,w_8)F_3(w_5,w_8)\right) \cdot F_1(w_6,w_7)F_2(w_6,w_8)F_1(w_7,w_8)\\
\end{array}
\]\end{normalsize}
We may now add the appropriate idempotents to these singularity-triple term pairs to form triples. Since each of these triples are regular at $z_1 = z_2 = \dots
= z_n$ then so too is the whole of $F_{\Lambda^\circ} (z_1, \dots , z_n)$.
{\nolinebreak \hfill \rule{2mm}{2mm}

\quad

Therefore, due to the above proposition an element $F_\Lambda \in
H_n$ can now be defined as the value of $F_\Lambda (z_1,
\dots , z_n)$ at $z_1 = z_2 = \cdots = z_n \neq 0$. Note that for $n=1$ we have $F_\Lambda=1$. For any $n\geqslant1$, 
take the expansion of the element $F_\Lambda\in H_n$ in the basis of the 
elements $T_\sigma$ where $\sigma$ is ranging over $S_n$.

\begin{proposition}\label{T_0coeff} The coefficient in $F_\Lambda\in H_n$ of the element $T_0$ is $1$.
\end{proposition}
\begin{proof} Expand the product (\ref{q-bigf})
as a sum of the elements $T_\sigma$ with coefficients
from the field of rational functions of $z_1, \ldots , z_n$;
these functions take values in $\mathbb{C}(q)$. 
The decomposition in $S_n$
with ordering of the pairs $(i, j)$ as in (\ref{q-bigf})
\[
\sigma_0=\prod_{(i,j)}^{\longrightarrow}\ \sigma_{j-i}
\]
is reduced, hence the coefficient at $T_0=T_{\sigma_0}$ in the
expansion of (\ref{q-bigf}) is $1$. By the definition of $F_\Lambda$,
then the coefficient of $T_0$ in $F_\Lambda$  must be also $1$
\end{proof}

In particular this shows that $F_\Lambda \neq 0$ for any nonempty
diagram $\lambda$. Let us now denote by $\varphi$ the involutive
antiautomorphism of the algebra $H_n$ over the field $\mathbb{C}(q)$, defined by
$\varphi (T_i) = T_i$ for every $i \in 1, \dots , n-1$.

\newpage
\begin{proposition}\label{q-varphi} The element $F_\Lambda T^{-1}_0$ is $\varphi$-invariant. \end{proposition}
\begin{proof} Any element of the algebra $H_n$ of the form
$F_i(a, b)$ is $\varphi$-invariant. 
Hence applying the antiautomorphism $\varphi$ to an element of 
$H_n$ the form (\ref{q-bigf}) just reverses the ordering of
the factors corresponding to the pairs $(i, j)$.
Using the relations (\ref{q-triple}) and (\ref{q-commute}),
we can rewrite the reversed product as
\[
\prod_{(i,j)}^{\longrightarrow}\ F_{ n-j+i}
\bigl( q^{2c_i(\Lambda)}z_i, q^{2c_j(\Lambda)}z_j\bigr)
\]
where the pairs $(i, j)$ are again ordered as in (\ref{q-bigf}).
But due to (\ref{Hecke2}) and (\ref{Hecke3}),
we also have the identity in the algebra $H_n$
\[
F_{ n-i}(a,b)T_0=T_0 F_i(a, b).
\]
This identity along with the equality 
$\varphi(T_0)=T_0$ implies that
any value of the function $F_\Lambda(z_1, \ldots , z_n)T^{-1}_0$ is 
$\varphi$-invariant. So is the element $F_\Lambda T^{-1}_0 \in H_n$
\end{proof}

\begin{proposition}\label{q-stripping} If $\lambda= (\alpha_1,
\alpha_2, \dots, \alpha_d | \beta_1, \beta_2, \dots , \beta_d)$
and \\$\mu = (\alpha_{k+1}, \alpha_{k+2}, \dots , \alpha_d |
\beta_{k+1}, \beta_{k+2}, \dots , \beta_d)$, then $F_{\Lambda^\circ} = P
\cdot F_{\mathrm{M}^\circ}$, for some element $P \in H_n$.
\end{proposition}
\begin{proof}
Here the shape $\mu$ is obtained by removing the first $k$
principal hooks of $\lambda$. Let $x$ be last entry in the $k^{\mbox{\small th}}$ row
of the hook tableau of shape $\lambda$. By the ordering described in
Proposition \ref{q-jimtheorem1},
\[ F_{\Lambda^\circ}(z_1, \dots ,
z_n) = \prod_{i=1}^k C_iR_i \cdot F_{\mathrm{M}^\circ}(z_{x+1}, \dots,
z_{n}),\]
where $C_i$ and $R_i$ are defined by (\ref{q-cproduct}) and (\ref{q-rproduct}).\\
Since all products $C_i$ and $R_i$ are regular at $z_1 = z_2 =
\dots = z_n \neq 0$, Proposition \ref{q-jimtheorem1} then gives us the
required statement. \end{proof}

%In any given ordering of $F_\Lambda(z_1, \dots , z_n)$, we want a
%singularity term to be placed next to an appropriate triple term
%such that we may then form a regular triple. In that case we will
%say these two terms are `tied'. However, proving the divisibilities described in the next two propositions require some pairs
%to be `untied', in which case we must form a new ordering. This is
%the content of the following proofs. Some explicit examples will
%then given in Example 2.10 below.

To prove the following divisibilities we require to form new singularity and triple term pairs, in which case we must form a new ordering. We will give explicit examples below.

\begin{proposition}\label{q-jimtheorem2} Suppose the numbers $u < v$ stand next to each
other in the same column of the hook tableau $\Lambda^\circ$ of shape
$\lambda$. First, let $s$ be the last entry in the row containing
$u$. If $c_v < 0$ then the element $F_{\Lambda^\circ} \in H_n$
is divisible on the left %and on the right 
by $F_{u}(q^{2c_u(\Lambda^\circ)}, q^{2c_v(\Lambda^\circ)}) =
T_{u} - q$. If $c_v \geqslant 0$ then the element $F_{\Lambda^\circ} \in
H_n$ is divisible on the left by the product
\[ \prod_{i = u, \dots, s}^\leftarrow \left( \prod_{j= s+1, \dots,
v}^\rightarrow F_{i+j-s-1}(q^{2c_i(\Lambda^\circ)}, q^{2c_j(\Lambda^\circ)}) \right) \]
\end{proposition}

\begin{proof}
%By Proposition \ref{q-varphi}, the divisibility of $F_\lambda$ by
%the element $1- (uv)$ on the left is equivalent to the
%divisibility by the same element on the right. Let us prove
%divisibility on the left.
%
%By Proposition \ref{q-stripping}, if $F_\mu$ is divisible
%on the right by $f_{uv}(c_u, c_v)$, or $f_{uv}(c_u, c_v)$ followed
%by some invertible terms, then so too will $F_\lambda$. 
%If $F_\Mu$ is divisible on the left by $F_{k}(q^{2c_u(\Lambda)}, q^{2c_v(\Lambda)})$,
%or $F_{v-u}(q^{2c_u(\Lambda)}, q^{2c_v(\Lambda)})$ preceded by some invertible terms, then $F_\Mu$ 
%is divisible on the right by $F_{v-u}(q^{2c_u(\Lambda)}, q^{2c_v(\Lambda)})$,
%or $F_{v-u}(q^{2c_u(\Lambda)}, q^{2c_v(\Lambda)})$ followed by some invertible terms.
%Hence, by Proposition \ref{q-stripping}, the product $F_\Lambda = P \cdot F_\Mu$ is also divisible on the right by $F_{v-u}(q^{2c_u(\Lambda)}, q^{2c_v(\Lambda)})$,
%or $F_{v-u}(q^{2c_u(\Lambda)}, q^{2c_v(\Lambda)})$ followed by some invertible terms. And so we can say $F_\Lambda$ is divisible on the left by $F_{v-u}(q^{2c_u(\Lambda), 2c_v(\Lambda))$,
%or $F_{v-u}(q^{2c_u(\Lambda)}, q^{2c_v(\Lambda)})$ preceded by some invertible terms. Hence we only need to prove the statement for
%$(u,v)$ such that $u$ is in the first row or first column of
%$\Lambda$.
%
Let $\lambda$ and $\mu$ be as in Proposition \ref{q-stripping} with $\lambda$ a partition of $n$ and $\mu$ a partition of $m$. If $F_\mathrm{M}(z_1, \dots, z_m)$ is divisible on the left by $F_i(w_a,w_b)$ then $F_\mathrm{M}(z_{x+1}, \dots, z_n)$ is divisible on the left by $F_i(w_{a+x},w_{b+x})$. Then, by Proposition \ref{q-varphi}, $F_\mathrm{M}(z_{x+1}, \dots, z_n)$ is divisible on the right by $F_{m-i}(w_{a+x},w_{b+x})$. Therefore, by Proposition \ref{q-stripping}, $F_{\Lambda^\circ}(z_1, \dots, z_n)$ is divisible on the right by $F_{m-i}(w_{a+x},w_{b+x})$ and so $F_{\Lambda^\circ}(z_1, \dots, z_n)$ is divisible on the left by \\$F_{n-m+i}(w_{a+x}, w_{b+x}) = F_{i+x}(w_{a+x},w_{b+x})$. Hence we only need to prove the statement for $(u,v)$ such that $u$ is in the first row or first column of $\Lambda^\circ$.

Let $r$ be the last entry in the first column of $\Lambda^\circ$, $s$
the last entry in the first row of $\Lambda^\circ$, and $t$ the last
entry in the second column of $\Lambda^\circ$, as shown in the proof of Proposition \ref{jimtheorem2}.

%\begin{figure}[h] \label{jimtheorem2fig}
%
%\begin{normalsize}
%\begin{center}
%\begin{picture}(275,200)
%\begin{small}
%\put(0,175){\framebox(25,25)[c]{ 1 }}
%\put(25,175){\framebox(25,25)[c]{ $r + 1$ }}
%\put(50,175){\framebox(25,25)[c]{ $r+2$ }}
%\put(0,150){\framebox(25,25)[c]{ 2 }}
%\put(25,150){\framebox(25,25)[c]{ $s+1$ }}
%\put(50,150){\framebox(25,25)[c]{ $t+1$ }}
%
%
%
%\put(0,25){\line(0,1){125}} \put(25,25){\line(0,1){125}}
%\put(0,0){\framebox(25,25)[c]{ $r$ }}
%
%\put(75,200){\line(1,0){50}} \put(75,175){\line(1,0){50}}
%\put(125,175){\framebox(25,25)[c]{ $u$ }}
%\put(150,200){\line(1,0){100}} \put(150,175){\line(1,0){100}}
%\put(250,175){\framebox(25,25)[c]{ $s$ }}
%
%\put(50,100){\line(0,1){50}} \put(25,75){\framebox(25,25)[c]{ $t$
%}}
%
%\put(75,150){\line(1,0){50}} \put(125,150){\framebox(25,25)[c]{
%$v$ }} \put(150,150){\line(1,0){50}} \put(200,150){\line(0,1){25}}
%
%\linethickness{1.5pt} \put(0,200){\line(1,0){25}}
%\put(25,175){\line(1,0){25}} \put(50,150){\line(1,0){25}}
%\put(25,175){\line(0,1){25}} \put(50,150){\line(0,1){25}}
%
%\put(125,150){\line(0,1){50}}\put(150,150){\line(0,1){50}}
%\put(125,200){\line(1,0){25}}\put(125,150){\line(1,0){25}}
%
%\put(95,185){$\dots$} \put(170,185){$\dots$}\put(220,185){$\dots$}
%\put(95,160){$\dots$} \put(170,160){$\dots$}
%
%\put(10,120){$\vdots$} \put(35,
%120){$\vdots$}\put(10,80){$\vdots$} \put(10,40){$\vdots$}
%
%\end{small}
%\end{picture}
%\end{center}
%\end{normalsize}
% \caption{The first two principal hooks of the hook tableau $\Lambda^\circ$}
% \end{figure}

We now continue this proof by considering three cases and showing
the appropriate divisibility in each.

\emph{(i)} \quad  If $c_v < 0$ (i.e. $u$ and $v$ are in the first
column of $\Lambda^\circ$) then $v = u+1$ and $F_{\Lambda^\circ}(z_1, \dots ,
z_n)$ can be written as $F_{\Lambda^\circ}(z_1, \dots, z_n) = F_{u}(w_u, w_v) \cdot F$. \\
Starting with $F_{\Lambda^\circ}(z_1, \dots , z_n)$ written in the
ordering described in Proposition \ref{q-jimtheorem1} and simply
moving the term $F_{v-u}(w_u, w_v)$ to the left
results in all the singularity terms in the product $F$ remaining
tied to the same triple terms as originally described in that
ordering, and the index of $F_{v-u}(w_u, w_v)$ increases from $u-v$ to $u$. Therefore we may still form regular triples for
each singularity in $F$, and hence $F$ is regular at $z_1 = z_2 = \cdots = z_n \neq 0$.\\
So by considering this expression for $F_{\Lambda^\circ}(z_1, \dots ,
z_n)$ at $z_1 = z_2 = \cdots = z_n \neq 0$ we see that $F_{\Lambda^\circ}$ will be
divisible on the left by $F_{u}(q^{2c_u(\Lambda^\circ)}, q^{2c_v(\Lambda^\circ)}) = T_{v-u} - q$.

\emph{(ii)} \quad  If $c_v = 0$ then $v=s+1$, and $F_{\Lambda^\circ}(z_1,
\dots , z_n)$ can be written as \[ F_{\Lambda^\circ}(z_1, \dots , z_n) =
\prod_{i = u, \dots , s}^\leftarrow F_{i}(w_i, w_{s+1}) \cdot F' . \] Again, starting with the ordering
described in Proposition \ref{q-jimtheorem1}, this results in  all
the singularity terms in the product $F'$ remaining tied to the
same triple terms as originally described in that ordering. Hence
$F'$ is regular at $z_1 = z_2 = \cdots = z_n \neq 0$. And so $F_{\Lambda^\circ}$
is divisible on the left by
\[ \prod_{i = u, \dots , s}^\leftarrow F_{i}(q^{2c_i(\Lambda^\circ)}, q^{2c_{s+1}(\Lambda^\circ)}).
\]

\emph{(iii)} \quad  If $c_v > 0$ (i.e. $v$ is above the steps)
then $F_{v-u}(w_u, w_v)$ is tied to the singularity
$F_{v-u+1}(w_{u-1}, w_v)$ as a triple term. To
show divisibility by $F_{v-u}(w_u, w_v)$ in this case
we need an alternative expression for $F_{\Lambda^\circ}(z_1, \dots ,
z_n)$ that is regular when $z_1 = z_2 = \cdots = z_n \neq 0$. Define a
permutation $\tau$ as follows,

\[ \tau = \prod_{i = u, \dots, s}^\rightarrow \left( \prod_{j= s+1, \dots,
v}^\leftarrow (i j) \right)
\phantom{XXXXXXXXXXXXXXXXXXXXXXXXXXXXXXX}
\]
\[ = \left(
\textrm{ \scriptsize $\begin{array}{ccccccccccccccccc}
                   1 & 2 & \dots & u-1 & u & u+1 & \dots &  & \dots & \dots  &  & \dots & v-1 & v & v+1 & \dots & n \\
                   1 & 2 & \dots & u-1 & s+1 & s+2 & \dots & v-1 & v & u & u+1 & \dots & s-1 & s & v+1 & \dots & n \\
                 \end{array}$ } \right) \]

From the definition of $C_1$ in (\ref{q-cproduct}) we now define
$C'_1 = \psi_\tau C_1$, where $\psi_\tau$ is a homomorphism such that
\[ \psi_\tau F_{j-i}(w_i,w_j) = F_{j-i}(w_i,w_{\tau j}). \]

%\newpage
%For the rest of this proof we will simply write $f_{ij}$
%instead of $f_{ij}(z_i + c_i, z_j + c_j)$. 

Define $R'_1$ as,
\newpage
\begin{eqnarray*}
R'_1 & = & \prod_{i= r+2, \dots, u-1}^\leftarrow \left(
\prod_{j=s+1, \dots ,v}^\rightarrow F_{i+j-s-r-1}(w_i,w_j) \right)  \\
 && \times  \prod_{i=
s+1, \dots, t-1}^\rightarrow \left( \prod_{j=i+1, \dots
,t}^\rightarrow F_{j-i+1}(w_i,w_j) \right) \cdot \left( \prod_{j=s+1, \dots
,t}^\leftarrow F_{t+1-j}(w_{r+1},w_j) \right)  \\
   && \times \left( \prod_{j=t+1, \dots
,v}^\rightarrow F_{j-s}(w_{r+1},w_j) \right) \cdot \prod_{i= r+1,
\dots, s-1}^\rightarrow \left( \prod_{j=i+1, \dots ,s}^\rightarrow
F_{j-i+v-s}(w_i,w_j) \right)  \\
&& \times \prod_{j= v+1, \dots, n}^\rightarrow
\left( \prod_{i=r+1, \dots ,s}^\rightarrow F_{j-i}(w_i,w_j) \right).
\end{eqnarray*}

Finally, define $C'_2$ as, \[ C'_2 = \prod_{j= t+1, \dots,
n}^\rightarrow \left( \prod_{i=s+1, \dots ,t}^\rightarrow F_{j-i}(w_i,w_j)
\right). \]

Then, \[ F_\Lambda(z_1, \dots , z_n) = \prod_{i = u, \dots ,
s}^\leftarrow \left( \prod_{j=s+1, \dots, v}^\rightarrow F_{i+j-s-1}(w_i,w_j)
\right) \cdot C'_1 R'_1 C'_2 R_2 \cdot \prod_{i=3}^d C_iR_i, \]
where $d$ is the number of principal hooks of $\lambda$.

The product $C'_1 R'_1 C'_2 R_2 \cdot \prod C_iR_i$ is regular at
$z_1 = z_2 = \dots z_n \neq 0$ since, as before, for any singularity
$(a,b)$ the terms $F_{i}(w_a,w_b)F_{i-1}(w_{a+1},w_b)$ can be replaced by the triple
$P_{i-1}^\pm F_{i}(w_a,w_b)F_{i-1}(w_{a+1},w_b)$ for some index $i$ -- except in the
expression $R'_1$ where the terms $F_{i-1}(w_a,w_l)F_{i}(w_a,w_b)$ are replaced by
$F_{i-1}(w_a,w_l)F_{i}(w_a,w_b)P_{i-1}^+$, where $l$ is the entry to the
immediate left of $b$. Note that $l = b-1$ when $c_b(\Lambda^\circ) > 1$ and $l =
s+1$ when $c_b(\Lambda^\circ) = 1$. \\
And so by letting $z_1 = z_2 = \dots =z_n \neq 0$ we see that $F_{\Lambda^\circ}$
is divisible on the left by \[ \prod_{i = u, \dots , s}^\leftarrow
\left( \prod_{j=s+1, \dots, v}^\rightarrow F_{i+j-s-1}(q^{2c_i(\Lambda^\circ)}, q^{2c_j(\Lambda)})
\right).
\]
\end{proof}

\newpage
\begin{proposition}\label{q-jimtheorem3} Suppose the numbers $u < v$ stand next to each
other in the same row of the hook tableau $\Lambda^\circ$ of shape
$\lambda$. Let $r$ be the last entry in the column containing $u$.
If $c_u > 0$ then the element $F_{\Lambda^\circ} \in H_n$ is
divisible on the left by $F_{u}(q^{2c_u(\Lambda^\circ)}, q^{c_v(\Lambda^\circ)}) = T_{u} +
q^{-1}$. If $c_u \leqslant 0$ then the element $F_{\Lambda^\circ} \in
H_n$ is divisible on the left by the product
\[ \prod_{i = u, \dots, r}^\leftarrow \left( \prod_{j= r+1, \dots,
v}^\rightarrow F_{i+j-r-1}(q^{c_i(\Lambda^\circ)}, q^{c_j(\Lambda^\circ)}) \right) \]
\end{proposition}
\begin{proof}
Suppose $u$ is to the immediate left of $v$ in some row of
$\Lambda^\circ$. As in the proof of Proposition \ref{q-jimtheorem2}, we
need only consider $(u,v)$ such that $u$ is in the first row or
first column of $\Lambda^\circ$. Let $r$ be the last entry in the first
column of $\Lambda^\circ$ and $s$ the last entry in the first row of
$\lambda$.

As in the proof of Proposition \ref{q-jimtheorem2}, we consider
three cases.

\emph{(i)} \quad  If $c_u > 0$ (i.e. $u$ and $v$ are in the first
row of $\Lambda^\circ$) then $v = u+1$ and $F_{\Lambda^\circ}(z_1, \dots , z_n)$
can be written as $F_{\Lambda^\circ}(z_1, \dots, z_n) = F_{u}(w_u,
w_v) \cdot F$. Singularity terms in the product $F$ remaining
tied to the same triple terms as originally described in
Proposition \ref{q-jimtheorem1}, hence $F$ is regular at $z_1 = z_2
= \cdots = z_n \neq 0$. By considering this expression for $F_{\Lambda^\circ}(z_1,
\dots , z_n)$ at $z_1 = z_2 = \cdots = z_n \neq 0$ we see that $F_{\Lambda^\circ}$
will be divisible on the left by $F_{u}(q^{2c_u(\Lambda^\circ)}, q^{2c_v(\Lambda^\circ)}) = T_u +
q^{-1}$.

\emph{(ii)} \quad  If $c_u = 0$ then $u=1$, $v=r+1$, and
$F_{\Lambda^\circ}(z_1, \dots , z_n)$ can be written as \[ F_{\Lambda^\circ}(z_1,
\dots , z_n) = \prod_{i = 1, \dots , r}^\leftarrow F_{i}(w_i, w_{r+1}) \cdot F' .
\] Again, singularities in $F'$ remain tied to the same triple terms as
originally described in Proposition \ref{q-jimtheorem1}. And so
$F_{\Lambda^\circ}$ is divisible on the left by
\[ \prod_{i = 1, \dots , r}^\leftarrow F_{i}(q^{2c_i(\Lambda^\circ)}, q^{2c_{r+1}(\Lambda^\circ)}).
\]

\emph{(iii)} \quad  If $c_u < 0$ (i.e. $v$ is below the steps)
then $F_{v-u}(w_u, w_v)$ is tied to the singularity
$F_{v-u+1}(w_{u-1}, w_v)$ as a triple term.

%For the rest of this proof we will simply write $f_{ij}$ instead
%of $f_{ij}(z_i + c_i, z_j + c_j)$. 
Define $C''_1$ as,
\begin{eqnarray*}
C''_1 & = & \prod_{i= 2, \dots, u-1}^\leftarrow \left(
\prod_{j=r+1, \dots ,v}^\rightarrow F_{i+j-r-1}(w_i,w_j) \right)  \\
   && \times \prod_{i=
r+1, \dots, s-1}^\rightarrow \left( \prod_{j=i+1, \dots
,s}^\rightarrow F_{j-i+1}(w_i,w_j) \right) \cdot \left( \prod_{j=r+1, \dots
,s}^\leftarrow F_{s+1-j}(w_i,w_j) \right) \\
&& \times F_{s-r+1}(w_1,w_{s+1})  \cdot \left(\prod_{i=r+1, \dots ,s}^\rightarrow F_{s+1-i}(w_i,w_{s+1}) \right) \\
   && \times
\left( \prod_{j= s+2, \dots, v}^\rightarrow F_{j-r}(w_1,w_j)\right) \cdot
    \prod_{i= 1,
\dots, r-1}^\rightarrow \left( \prod_{j=i+1, \dots ,r}^\rightarrow
F_{j-i+v-r}(w_i,w_j) \right) \\
&& \times \prod_{j= v+1, \dots, n}^\rightarrow
\left( \prod_{i=1, \dots ,r}^\rightarrow F_{j-i}(w_i,w_j) \right), \\
\end{eqnarray*}

and define $R''_1$ as, \[ R''_1 = \prod_{j= s+2, \dots,
n}^\rightarrow \left( \prod_{i=r+1, \dots ,s}^\rightarrow F_{j-i}(w_i,w_j)
\right). \]

Then, \[ F_{\Lambda^\circ}(z_1, \dots , z_n) = \prod_{i = u, \dots ,
r}^\leftarrow \left( \prod_{j=r+1, \dots, v}^\rightarrow F_{i+j-r-1}(w_i,w_j)
\right) \cdot C''_1 R''_1 \cdot \prod_{i=2}^d C_iR_i, \] where $d$
is the number of principal hooks of $\lambda$.

The product $C''_1 R''_1 \cdot \prod C_iR_i$ is regular at $z_1 =
z_2 = \dots z_n \neq 0$ since, as before, for any singularity $(a,b)$ the
terms $F_{i}(w_a,w_b)F_{i-1}(w_{a+1},w_b)$ can be replaced by the triple
$P_{i-1}^\pm F_{i}(w_a,w_b)F_{i-1}(w_{a+1},w_b)$ for some index $i$ -- except in the
expression $C''_1$ where $F_{s-r+1}(w_1, w_{s+1})F_{s-r}(w_{r+1}, w_{s+1})$ is replaced by
$P_{s-r}^+ F_{s-r+1}(w_1, w_{s+1})F_{s-r}(w_{r+1}, w_{s+1})$, and $F_{i-1}(w_a,
w_{b-1})F_i(w_a,w_b)$ is replaced by
$F_{i-1}(w_a,w_{b-1})F_i(w_a,w_b)P_{i-1}^-$ for all other singularities $(a, b)$ in $C''_1$.\\
And so by letting $z_1 = z_2 = \dots =z_n \neq 0$ we see that $F_{\Lambda^\circ}$
is divisible on the left by \[ \prod_{i = u, \dots , r}^\leftarrow
\left( \prod_{j=r+1, \dots, v}^\rightarrow F_{i+j-r-1}(q^{2c_i(\Lambda^\circ)}, q^{2c_j(\Lambda^\circ)})
\right).
\]
\end{proof}

Let us now
consider an example that allows us to see how the product
$F_{\Lambda^\circ}(z_1, \dots, z_n)$ is broken down in the proofs of
Proposition \ref{q-jimtheorem2} and Proposition \ref{q-jimtheorem3}.

{\addtocounter{definition}{1} \bf Example \thedefinition .} We
again consider the hook tableau of the Young diagram $\lambda =
(3,3,2)$.

We begin with the product $F_{\Lambda^\circ}(z_1, \dots, z_n)$ in the
ordering described in Proposition \ref{q-jimtheorem1}. %For
%simplicity we again write $f_{pq}$ in place of the term
%$f_{pq}(z_p + c_p, z_q + c_q)$. 
We have also marked out the
singularities in this expansion along with their triple terms, but
no idempotents have yet been added which would form regular
triples.

\begin{normalsize}
\begin{equation}\label{q-example1}
\begin{array}{rl}
   F_{\Lambda^\circ} (z_1, \dots , z_n) = &
   F_1(w_1,w_2)F_2(w_1,w_3)F_1(w_2,w_3)F_3(w_1,w_4)F_2(w_2,w_4)F_1(w_3,w_4)\\&
   F_4(w_1,w_5)F_3(w_2,w_5)F_2(w_3,w_5)(F_5(w_1,w_6)F_4(w_2,w_6))F_3(w_3,w_6)\\&
   F_6(w_1,w_7)(F_5(w_2,w_7)F_4(w_3,w_7))F_7(w_1,w_8)F_6(w_2,w_8)F_5(w_3,w_8)\\&
   \cdot F_1(w_4,w_5)F_2(w_4,w_6)F_1(w_5,w_6)F_3(w_4,w_7)F_2(w_5,w_7)\\&
   (F_4(w_4,w_8)F_3(w_5,w_8)) \cdot F_1(w_6,w_7)F_2(w_6,w_8)F_1(w_7,w_8)\\
\end{array}
\end{equation}
\end{normalsize}

Let $u=4$ and $v=6$ in Proposition \ref{q-jimtheorem2}. Then by that
proposition we may arrange the above product as follows. Notice
since $c_v = 0$ all singularity-triple term pairs remain the same.

\begin{normalsize}
\[
\begin{array}{rl}
   F_{\Lambda^\circ} (z_1, \dots , z_n) = &
   F_5(w_5,w_6)F_4(w_4,w_6) \cdot F_1(w_1,w_2)F_2(w_1,w_3)F_1(w_2,w_3)(F_3(w_1,w_6)\\&
   F_2(w_2,w_6))F_1(w_3,w_6)F_4(w_1,w_4)F_3(w_2,w_4)F_2(w_3,w_4)F_5(w_1,w_5)\\&
   F_4(w_2,w_5)F_3(w_3,w_5)F_6(w_1,w_7)(F_5(w_2,w_7)F_4(w_3,w_7))F_7(w_1,w_8)\\&
   F_6(w_2,w_8)F_5(w_3,w_8) \cdot F_2(w_4,w_5)F_3(w_4,w_7)F_2(w_5,w_7)\\&
   (F_4(w_4,w_8)F_3(w_5,w_8)) \cdot F_1(w_6,w_7)F_2(w_6,w_8)F_1(w_7,w_8)\\
\end{array}
\]
\end{normalsize}
We may now add the appropriate idempotents so that all
singularities remain in regular triples. And so by considering the
product at $z_1 = z_2 = \cdots = z_n \neq 0$ we have that $F_{\Lambda^\circ}$ is
divisible on the left by $T_4 - q$, preceded only by invertible
terms, as desired.

\newpage
Now let $u=3$ and $v=7$ in Proposition \ref{q-jimtheorem3}. Then by
that proposition we may arrange (\ref{q-example1}) as follows.
Singularities in $C_1$ have been marked out with their alternative
triple terms, while all other singularity-triple term pairs remain
the same.
\begin{normalsize}\[
\begin{array}{rl}
   F_{\Lambda^\circ} (z_1, \dots , z_n) = &
   F_3(w_3,w_4)F_4(w_3,w_5)F_5(w_3,w_6)F_6(w_3,w_7) \cdot F_2(w_2,w_4)F_3(w_2,w_5)\\&
   (F_4(w_2,w_6)F_5(w_2,w_7)) \cdot F_2(w_4,w_5) \cdot F_1(w_1,w_5)(F_2(w_1,w_4) \cdot F_3(w_1,w_6))\\& 
   \cdot F_2(w_4,w_6)F_1(w_5,w_6) \cdot F_4(w_1,w_7)\cdot F_5(w_1,w_2)F_6(w_1,w_3)F_5(w_2,w_3)\\& 
   \cdot F_7(w_1,w_8)F_6(w_2,w_8)F_5(w_3,w_8) \cdot F_3(w_4,w_7)F_2(w_5,w_7)\\&
   (F_4(w_4,w_8)F_3(w_5,w_8)) \cdot F_1(w_6,w_7)F_2(w_6,w_8)F_1(w_7,w_8)\\
\end{array}
\]\end{normalsize}
In moving $F_4(w_3,w_7)$ to the left it is untied from the singularity
$F_5(w_2,w_7)$. So we must form new triples which are regular at $z_1 =
z_2 = \cdots = z_n \neq 0$ by the method described in Proposition
\ref{q-jimtheorem3}. Therefore, by considering the product at $z_1 =
z_2 = \cdots = z_n \neq 0$, we have that $F_{\Lambda^\circ}$ is divisible on the
left by $T_6 + q^{-1}$, again preceded only by invertible terms, as
desired. {\nolinebreak \hfill \rule{2mm}{2mm}

\quad

\begin{lemma}\label{q-divisibilitybyadjacenttransposition} Let $\Lambda$ and $\tilde{\Lambda}$ be tableaux of the same shape such that $k= \Lambda(a,b) = \Lambda(a+1,b)-1$ and $\tilde{k}=\tilde{\Lambda}(a,b) = \tilde{\Lambda}(a+1,b)-1$. Then $F_\Lambda \in H_n$ is divisible on the left by $T_k - q$ if and only if $F_{\tilde{\Lambda}} \in H_n$ is divisible on the left by $T_{\tilde{k}} - q$. \end{lemma}

\begin{proof}
Let $\sigma$ be the permutation such that $\tilde{\Lambda}=\sigma\cdot\Lambda$.
% Note that $\sigma$ is not necessarily an adjacent transposition. 
There is a decomposition
$\sigma=\sigma_{i_N}\ldots\sigma_{i_1}$ such that for each $M=1, \ldots , N-1$
the tableau $\Lambda_{ M}=\sigma_{i_M}\ldots\sigma_{i_1}\cdot\Lambda$ is
standard. Note that this decomposition is not necessarily reduced.

Denote $F_\Lambda$ by $F_k(q^{2c_k(\Lambda)},q^{2c_{k+1}(\Lambda)}) \cdot F$. Then, by using the relations
(\ref{q-triple}) and (\ref{q-commute}) and  Proposition \ref{q-jimtheorem1}, we have the following chain of equalities:

\[
F_{\tilde{k}}
\bigl( 
q^{2c_{\tilde{k}}(\tilde{\Lambda})}, q^{2c_{\tilde{k}+1}(\tilde{\Lambda})}
\bigr)
\ \cdot\ 
\prod_{M=1, \ldots , N}^{\longleftarrow}
F_{ n- i_M}
\bigl( 
q^{2 c_{i_M+1}(\Lambda_M)}
, 
q^{2 c_{i_M}(\Lambda_M)}
\bigr) \cdot F =  \]
\[\prod_{M=1, \ldots , N}^{\longleftarrow}
F_{ i_M}
\bigl( 
q^{2 c_{i_M}(\Lambda_M)}
, 
q^{2 c_{i_M+1}(\Lambda_M)}
\bigr)
\ \cdot\ 
F_{ k}
\bigl( q^{2c_k(\Lambda)}, q^{2c_{k+1}(\Lambda)}\bigr) \cdot F =\]
\[
\prod_{M=1, \ldots , N}^{\longleftarrow}
F_{ i_M}
\bigl( 
q^{2 c_{i_M}(\Lambda_M)}
, 
q^{2 c_{i_M+1}(\Lambda_M)}
\bigr)
\ \cdot\ F_\Lambda\ =
\]
\[
F_{\tilde{\Lambda}}\ \cdot
\prod_{M=1, \ldots , N}^{\longleftarrow}
F_{ n- i_M}
\bigl( 
q^{2 c_{i_M+1}(\Lambda_M)}
, 
q^{2 c_{i_M}(\Lambda_M)}
\bigr)
\]

Hence divisibility by $T_k - q$ for $F_\Lambda$ implies its counterpart for the tableau
$\tilde{\Lambda}$ and the index $\tilde{k}$, and vice versa.
Here we also use the equalities
\[
F_{ k}
\bigl( q^{2c_k(\Lambda)}, q^{2c_{k+1}(\Lambda)}\bigr)
=T_k-q,
\]\[
F_{\tilde{k}}
\bigl( 
q^{2c_{\tilde{k}}(\tilde{\Lambda})}, q^{2c_{\tilde{k}+1}(\tilde{\Lambda})}
\bigr)
=T_{\tilde{k}}-q.
\]
\end{proof}

\begin{corollary}\label{q-divisibilitycorollary}
If $k=\Lambda(a, b)$ and $k+1=\Lambda(a+1, b)$
then the element $F_\Lambda\in H_n$ is divisible on the left by $T_k-q$. If $k=\Lambda(a, b)$ and $k+1=\Lambda(a, b+1)$
then the element $F_\Lambda\in H_n$ is divisible on the left by 
$T_k+q^{-1}$.
\end{corollary}

\begin{proof}
Due to Lemma \ref{q-divisibilitybyadjacenttransposition}
it suffices to prove the first part of Corollary \ref{q-divisibilitycorollary} for only one 
standard tableau $\Lambda$ of shape $\lambda$. Therefore, using Proposition \ref{q-jimtheorem2} and taking $\tilde{\Lambda}$ to be the hook tableau $\Lambda^\circ$ of shape $\lambda$ we have shown the first part of Corollary \ref{q-divisibilitycorollary} in the case $c_v(\Lambda) < 0$.

Next let $\Lambda^{\circ}(a,b)= u$, $\Lambda^{\circ}(a+1,b) = v$ and $s$ be the last entry in the row containing $u$. Then for $c_v({\Lambda^\circ}) \geqslant 0$ Proposition \ref{q-jimtheorem2} showed that $F_{\Lambda^\circ} \in H_n$ is divisible on the left by 
\begin{equation}\label{q-divisibilitycorollaryequation2} \prod_{i = u, \dots, s}^\leftarrow \left( \prod_{j= s+1, \dots,
v}^\rightarrow F_{i+j-s-1}(q^{2c_i(\Lambda^\circ)}, q^{2c_j(\Lambda^\circ)}) \right) \end{equation}
Put $k=u+v-s-1$,
this is the value of the index $i+j-s-1$ in (\ref{q-divisibilitycorollaryequation2})
when $i=u$ and $j=v$. Let $\Lambda$ be the tableau
such that $\Lambda^\circ$ is obtained from 
the tableau $\sigma_k\cdot\Lambda$ by the permutation
\[
\prod_{i =  u, \ldots, s}^{\longleftarrow}\,
\biggl(\ 
\prod_{j = s+1, \ldots,  v}^{\longrightarrow}
\sigma_{ i+j-s-1} \biggr).
\]
The tableau $\Lambda$ is standard. Moreover, then
$\Lambda(a, b)=k$ and $\Lambda(a+1, b)=k+1$.
Note that the rightmost factor in the product (\ref{q-divisibilitycorollaryequation2}),
corresponding to $i=u$ and $j=v$, is
\[
F_{u+v-s-1}
\bigl( q^{2c_u(\Lambda^\circ)}, q^{2c_{v}(\Lambda^\circ)}\bigr)
=
T_k - q.
\]
Denote by $F$ the product of all factors in (\ref{q-divisibilitycorollaryequation2})
but the rightmost one. Further, denote by $G$ the product obtained
by replacing each factor in $F$
\[
F_{i+j-s-1}
\bigl(q^{2c_i(\Lambda^\circ)}, q^{2c_j(\Lambda^\circ)}\bigr)
\]
respectively by
\[
F_{n-i-j+s+1}
\bigl(q^{2c_j(\Lambda^\circ)}, q^{2c_i(\Lambda^\circ)}\bigr).
\]
The element $F \in H_n$ is invertible, and we have
\[ F \cdot F_\Lambda = F_{\Lambda^\circ} \cdot G = F \cdot (T_k - q) \cdot (C'_1R'_1C'_2R_2 \prod C_iR_i) \cdot G,\] where the final equality is as described in Proposition \ref{q-jimtheorem2}. Therefore the divisibility of 
the element $F_{\Lambda^\circ}$ on the left by the product (\ref{q-divisibilitycorollaryequation2})
will imply the divisibility of the element $F_\Lambda$ on the left by 
$T_k - q$.

This shows the required divisibility for the tableau $\Lambda = \sigma \sigma_k \cdot \Lambda^{\circ}$. Using Lemma \ref{q-divisibilitybyadjacenttransposition} again concludes the proof of the first part of Corollary \ref{q-divisibilitycorollary}. 

The second part of Corollary \ref{q-divisibilitycorollary} may be shown similarly.
\end{proof}

%========================================================================

\subsection{Generating irreducible representations of $H_n$}

In his thesis of 1974 Hoefsmit wrote down the irreducible seminormal representations for the Hecke algebras of type $A$, $B$ and $D$, \cite{H}. In the sequel we will construct the seminormal basis of type $A$. For every standard tableau $\Lambda$ of shape $\lambda$ 
we have defined an element $F_\Lambda$ of the algebra $H_n$. 
Let us now assign to $\Lambda$ another element of $H_n$,
which will be denoted by $G_\Lambda$. 

Let $\rho\in S_n$ be the permutation such 
that $\Lambda=\rho\cdot\Lambda^\circ$. %that is $\Lambda(a, b)=\rho(\Lambda^\circ(a, b))$ for all possible $a$ and $b$.
For any $j=1, \ldots , n$ take the subsequence of the sequence
$\rho(1), \ldots ,\rho(n)$ consisting of all $i<j$ such that
$\rho^{-1}(i)>\rho^{-1}(j)$. Denote by $\mathcal{A}_j$ the result of reversing
this subsequence. Let $|\mathcal{A}_j|$ be the length of sequence $\mathcal{A}_j$. 
We have a reduced decomposition in the symmetric group $S_n$,
\begin{equation}\label{rho}
\rho\ =
\prod_{j=1,\ldots, n}^{\longrightarrow}
\biggl(\ 
\prod_{k=1,\ldots,|\mathcal{A}_j|}^{\longrightarrow}
\ \sigma_{j-k}
\biggr).
\end{equation}
Let $\sigma_{i_L}\ldots\sigma_{i_1}$ be the product of adjacent 
transpositions at the right hand side of (\ref{rho}). 
For each tail $\sigma_{i_K}\ldots\sigma_{i_1}$ of this product, the image 
$\sigma_{i_K}\ldots\sigma_{i_1}\cdot\Lambda^\circ$ is a standard tableau. 
This can easily be proved by induction on the length $K=L, \ldots , 1$ of the 
tail, see also the proof of Proposition \ref{Heckeidentity} below. Note that for
any 
$i\in\mathcal{A}_j$ and $k\in\{1, \ldots , n-1\}$ the elements of the algebra $H_n$,
\[
F_k
\bigl( q^{2c_i(\Lambda)}, q^{2c_j(\Lambda)}\bigr)
\ \quad\textrm{and}\ \quad
F_k
\bigl( q^{2c_j(\Lambda)}, q^{2c_i(\Lambda)}\bigr)
\]
are well defined and invertible. Indeed, if
$i=\Lambda(a, b)$ and $j=\Lambda(c, d)$ for some $a, b$ and $c, d$ 
then $a<c$ and $b>d$. So $c_i(\Lambda)-c_j(\Lambda)=b-a-d+c\geqslant2$ here.

\begin{proposition}\label{Heckeidentity}
We have the equality in the algebra $H_n$
\[
F_\Lambda\ \cdot
\prod_{j=1,\ldots, n}^{\longrightarrow}
\biggl(\ 
\prod_{k=1,\ldots,|\mathcal{A}_j|}^{\longrightarrow}
\ F_{n-j+k}
\bigl( q^{2c_j(\Lambda)}, q^{2c_i(\Lambda)}\bigr)
\biggr)\ =
\]\[
\prod_{j=1,\ldots, n}^{\longrightarrow}
\biggl(\ 
\prod_{k=1,\ldots,|\mathcal{A}_j|}^{\longrightarrow}
\ F_{j-k}
\bigl( q^{2c_i(\Lambda)}, q^{2c_j(\Lambda)}\bigr)
\biggr)
\ \cdot\ 
F_{\Lambda^{\circ}}
\ \quad\textrm{where}\ \quad 
i=\mathcal{A}_j(k).
\]
\end{proposition}

\begin{proof}
We will proceed by induction on the length
$N=|\mathcal{A}_1|+\ldots+|\mathcal{A}_{n}|$
of the element $\rho\in S_n$. Let $m$ be the minimal of the
indices $j$ such that the sequence $\mathcal{A}_j$ is not empty.
Then we have $\mathcal{A}_{m}(1)=m-1$. Indeed, if $\mathcal{A}_{m}(1)<m-1$ then
$\rho^{-1}(\mathcal{A}_{m}(1))>\rho^{-1}(m-1)$. Then
$\mathcal{A}_{m}(1)\in\mathcal{A}_{m-1}$,
which would contradict to the minimality of $m$.  The tableau 
$\sigma_{m-1}\cdot\Lambda$ is standard, denote it by $\Lambda'$. 
In our proof of Proposition \ref{q-jimtheorem1} we used the equality
(\ref{q-jimtheorem1equation}). Setting $k=m-1$ in that equality and then using 
Proposition \ref{q-jimtheorem1} itself, we obtain the equality in $H_n$
\begin{equation}\label{Heckeidentityequation}
F_\Lambda F_{ n-m+1}
\bigl( q^{2c_m(\Lambda)}, q^{2c_{m-1}(\Lambda)}\bigr)
=
F_{m-1}
\bigl( q^{2c_{m-1}(\Lambda)}, q^{2c_m(\Lambda)}\bigr)
F_{\Lambda^{\prime}}.
\end{equation}

For each index $j=1, \ldots , n$ denote by $\mathcal{A}_j^{\prime}$ the counterpart
of the 
sequence $\mathcal{A}_j$ for the standard tableau $\Lambda'$ instead of $\Lambda$. 
Each of the sequences
$\mathcal{A}_{1}^{\prime}, \ldots ,\mathcal{A}_{ m-2}^{\prime}$ and $\mathcal{A}_m^{\prime}$  
is empty. The sequence $\mathcal{A}_{ m-1}^{\prime}$ is obtained from the 
sequence $\mathcal{A}_m$ by removing its first term $\mathcal{A}_{ m}(1)=m-1$. 
By replacing the terms $m-1$ and $m$, whenever any of them occurs,
respectively by  $m$ and $m-1$ in all the sequences 
$\mathcal{A}_{ m+1}, \ldots ,\mathcal{A}_{n}$ we obtain the sequences  
$\mathcal{A}_{ m+1}^{\prime}, \ldots ,\mathcal{A}_{n}^{\prime}$.

Assume that the Proposition \ref{Heckeidentity} is true for $\Lambda'$ instead of
$\Lambda$. Write the product at the left hand side of the equality to be 
proved in Proposition~\ref{Heckeidentity} as
\[
F_\Lambda 
F_{ n-m+1}
\bigl( q^{2c_m(\Lambda)}, q^{2c_{m-1}(\Lambda)}\bigr)
\hskip9pt\cdot\hskip-7pt
\prod_{k=2,\ldots,|\mathcal{A}_m|}^{\longrightarrow}
\ F_{n-m+k}
\bigl( q^{2c_m(\Lambda)}, q^{2c_i(\Lambda)}\bigr)\hskip6pt\times
\]\[
\prod_{j= m+1,\ldots, n}^{\longrightarrow}
\biggl(\ 
\prod_{k=1,\ldots,|\mathcal{A}_j|}^{\longrightarrow}
\ F_{ n-j+k}
\bigl( q^{2c_j(\Lambda)}, q^{2c_i(\Lambda)}\bigr)
\biggr)
\]
where in the first line $i=\mathcal{A}_{ m}(k)$, 
while in the second line $i=\mathcal{A}_j(k)$.
Using the equality (\ref{Heckeidentityequation}) and the description
of the sequences $\mathcal{A}_{1}^{\prime}, \ldots ,\mathcal{A}_{ n}^{\prime}$
as given above, the latter product can be rewritten as
\[
F_{m-1}
\bigl( q^{2c_{m-1}(\Lambda)}, q^{2c_m(\Lambda)}\bigr)
F_{\Lambda^{\prime}}
\ \ \times
\]\[
\prod_{j=1,\ldots, n}^{\longrightarrow}
\biggl(\ 
\prod_{k=1,\ldots,|\mathcal{A}_j^{\prime}|}^{\longrightarrow}
\ F_{ n-j+k}
\bigl( q^{2c_j(\Lambda')}, q^{2c_i(\Lambda')}\bigr)
\biggr)
\ \quad\textrm{where}\ \quad 
i=\mathcal{A}_{ j}^{\prime}(k). 
\]
By the inductive assumption, this product equals
\[
F_{m-1}
\bigl( q^{2c_{m-1}(\Lambda)}, q^{2c_m(\Lambda)}\bigr)
\hskip9pt\cdot
\prod_{j=1,\ldots, n}^{\longrightarrow}
\biggl(\ 
\prod_{k=1,\ldots,|\mathcal{A}_j^{\prime}|}^{\longrightarrow}
\ F_{j-k}
\bigl( q^{2c_i(\Lambda')}, q^{2c_j(\Lambda')}\bigr)
\biggr)
\]
times $F_{\Lambda^{\!\circ}}$, where we keep to the notation
$i=\mathcal{A}_{ j}^{\prime}(k)$. 
Using the description of the sequences 
$\mathcal{A}_{1}^{\prime}, \ldots ,\mathcal{A}_{ n}^{\prime}$
once again, the last product can be rewritten as at the right hand side 
of the equality to be proved in Proposition \ref{Heckeidentity}
\end{proof}

For any $j=1, \ldots , n$ denote by $\mathcal{B}_j$ the subsequence of the sequence
$\rho(1), \ldots ,\rho(n)$ consisting of all $i<j$ such that
$\rho^{-1}(i)<\rho^{-1}(j)$. Note that we have a reduced decomposition
in the symmetric group $S_n$,
\[
\rho\sigma_0\ =
\prod_{j=1,\ldots, n}^{\longrightarrow}
\biggl(\ 
\prod_{k=1,\ldots,|\mathcal{B}_j|}^{\longrightarrow}
\ \sigma_{j-k}
\biggr)
\]
where $|\mathcal{B}_j|$ is the length of sequence
$\mathcal{B}_j$.
Consider the rational
function taking values in $H_n$, of the variables $z_1, \ldots , z_n$
\[
\prod_{j=1,\ldots, n}^{\longrightarrow}
\biggl(\ 
\prod_{k=1,\ldots,|\mathcal{B}_j|}^{\longrightarrow}
\ F_{j-k}
\bigl( q^{2c_i(\Lambda)}z_i, q^{2c_j(\Lambda)}z_j\bigr)
\biggr)
\ \quad\textrm{where}\ \quad 
i=\mathcal{B}_j(k).
\]
Denote this rational function by $G_\Lambda(z_1, \ldots , z_n)$.
Using induction on the length of the element $\rho\in S_n$ 
as in the proof of Proposition \ref{Heckeidentity}, one can prove that
\[
F_\Lambda(z_1, \ldots , z_n)\ =\  
G_\Lambda(z_1, \ldots , z_n)\ \ \times
\]\[
\prod_{j=1,\ldots, n}^{\longleftarrow}
\biggl(\ 
\prod_{k=1,\ldots,|\mathcal{A}_j|}^{\longleftarrow}
\ F_{ n-j+k}
\bigl( q^{2c_i(\Lambda)}z_i, q^{2c_j(\Lambda)}z_j\bigr)
\biggr)
\ \quad\textrm{where}\ \quad 
i=\mathcal{A}_j(k).
\]
Hence restriction of $G_\Lambda(z_1, \ldots , z_n)$ 
to the subspace $\mathcal{H}_\Lambda\subset\mathbb{C}(q)^{n}$
is regular on the line $z_1 = \cdots = z_n \neq 0$ due to Proposition \ref{q-jimtheorem1}.
The value of that restriction is our element $G_\Lambda\in H_n$ 
by definition. Moreover, then $F_\Lambda$ equals
\[
%F_\Lambda\ =\  
G_\Lambda\ \cdot
\prod_{j=1,\ldots, n}^{\longleftarrow}
\biggl(\ 
\prod_{k=1,\ldots,|\mathcal{A}_j|}^{\longleftarrow}
\ F_{ n-j+k}
\bigl( q^{2c_i(\Lambda)}, q^{2c_j(\Lambda)}\bigr)
\biggr)
\ \quad\textrm{where}\ \quad 
i=\mathcal{A}_j(k).
\]
Using the relation (\ref{q-inverse}),
this factorization of $F_\Lambda$ implies that the left hand side
of the equality in Proposition \ref{Heckeidentity} also equals $G_\Lambda$ times
\[
\prod_{j=1,\ldots, n}
\biggl(\ 
\prod_{k=1,\ldots,|\mathcal{A}_j|}
\biggl(1-
%\frac{(q-q^{-1})^{2}}
%{( q^{ c_i(\La)-c_j(\La)}-q^{ c_j(\La)-c_i(\La)})^{2}}
\frac{(q-q^{-1})^{2}q^{2c_i(\Lambda)+2c_j(\Lambda)}}
{(q^{2c_i(\Lambda)}-q^{2c_j(\Lambda)})^{2}}
\biggr)
\biggr)
\quad\textrm{where}\ \quad 
i=\mathcal{A}_j(k).
\]
By rewriting the factors of this product, 
Proposition \ref{Heckeidentity} yields

\begin{corollary}\label{C3.1}
{\bf\hskip-6pt.\hskip1pt} 
We have the equality in the algebra $H_n$
\[
\prod_{j=1,\ldots, n}
\biggl(\ 
\prod_{k=1,\ldots,|\mathcal{A}_j|}
\biggl(1-
\frac{(q-q^{-1})^{2}}
{( q^{ c_i(\Lambda)-c_j(\Lambda)}-q^{ c_j(\Lambda)-c_i(\Lambda)})^{2}}
\biggr)
\biggr)
\ \cdot\ 
G_\Lambda\ =
\]\[
\prod_{j=1,\ldots, n}^{\longrightarrow}
\biggl(\ 
\prod_{k=1,\ldots,|\mathcal{A}_j|}^{\longrightarrow}
\ F_{j-k}
\bigl( q^{2c_i(\Lambda)}, q^{2c_j(\Lambda)}\bigr)
\biggr)
\ \cdot\ 
F_{\Lambda^{\!\circ}}
\ \quad\textrm{where}\ \quad 
i=\mathcal{A}_j(k).
\]
\end{corollary}

Yet arguing like in the proof of Proposition \ref{T_0coeff},
the definition of $G_\Lambda$ implies

\begin{proposition}\label{P3.2}
{\bf\hskip-6pt.\hskip1pt} 
The element $G_\Lambda$ equals $T_{\rho\sigma_0}$ plus a sum of 
the elements 
$T_\sigma$ with certain non-zero coefficients from $\mathbb{C}(q)$, 
where the length of 
each $\sigma\in S_n$ is less than that of $\rho\sigma_0$. 
\end{proposition}

Note that $G_{\Lambda^\circ}=F_{\Lambda^\circ}$ by definition.
Denote by $V_\lambda$ the left ideal in the algebra $H_n$
generated by the element $F_{\Lambda^\circ}$. Due to Corollary \ref{C3.1}
we have $G_\Lambda\in V_\lambda$ for any standard tableau $\Lambda$ of shape $\lambda$.
Proposition \ref{P3.2} shows that the elements $G_\Lambda\in H_n$ for all
pairwise distinct standard tableaux $\Lambda$ of shape $\lambda$
are linearly independent. The next proposition implies, in particular, that
these elements also span the vector space $V_\lambda$.

For any $k=1, \ldots , n-1$ denote $d_k(\Lambda)=c_k(\Lambda)-c_{k+1}(\Lambda)$.
If the tableau $\sigma_k\Lambda$ is not standard, then the numbers $k$ and
$k+1$ stand next to each other in the same row or in the same
column of $\Lambda$, that is $k+1=\Lambda(a, b+1)$ or $k+1=\Lambda(a+1, b)$ 
for $k=\Lambda(a, b)$. Then we have $d_k(\Lambda)=-1$ or $d_k(\Lambda)=1$
respectively. But if the tableau $\sigma_k\Lambda$ is standard, then we have
$| d_k(\Lambda)|\geqslant2$. 

\begin{proposition}\label{T3.3}
{\bf\hskip-6pt.\hskip1pt} 
For any standard tableau $\Lambda$ and any $k=1, \ldots , n-1$ we have:

\vspace{10pt}

a)

\vspace{-42pt}

\[
T_kG_\Lambda=
\left\{
\begin{array}{rl}
qG_\Lambda&\textrm{\quad if\quad\ }d_k(\Lambda)=-1,
\\[1mm]
- q^{-1}G_\Lambda&\textrm{\quad if\quad\ }d_k(\Lambda)=1;
\end{array}
\right.
\]

\vspace{-5pt}

b)

\vspace{-33pt}

\[
T_kG_\Lambda=
\frac{q-q^{-1}}{1-q^{2d_k(\Lambda)}}G_\Lambda\ +\ 
G_{\sigma_k\Lambda}\ \times
\]\[
\phantom{T_kG_\Lambda=}
\left\{
\begin{array}{ll}
1-
{
\displaystyle 
\frac{(q-q^{-1})^{2}}{(q^{ d_k(\Lambda)}-q^{-d_k(\Lambda)})^{2}}
}
&\textrm{\quad if\quad\ }d_k(\Lambda)\leqslant-2,
\\
1&\textrm{\quad if\quad\ }d_k(\Lambda)\geqslant2.
\end{array}
\right.
\]
\end{proposition}

\begin{proof}
The element $G_\Lambda\in H_n$ is obtained by multiplying $F_\Lambda$
on the right by a certain 
% invertible 
element of $H_n$. Hence Part (a)
of Proposition \ref{T3.3} immediately follows from Corollary \ref{divisibilitycorollary}. 
Now suppose that the tableau $\sigma_k\Lambda$ is standard.
Moreover, suppose that $d_k(\Lambda)\geqslant2$. Using Corollary \ref{C3.1} along with
the relations (\ref{q-triple})~and~(\ref{q-commute}), one can get the equality
\begin{equation}\label{3.1}
\biggl(1-
\frac{(q-q^{-1})^{2}}
{(q^{ d_k(\Lambda)}-q^{-d_k(\Lambda)})^{2}}
\biggr)
G_\Lambda\ =\ F_k
\bigl( q^{2c_k(\Lambda)}, q^{2c_{k+1}(\Lambda)}\bigr)
G_{\sigma_k\Lambda}.
\end{equation}
Using the relation (\ref{q-inverse}), we obtain from (\ref{3.1}) the equality
\[
F_k
\bigl( q^{2c_{k+1}(\Lambda)}, q^{2c_k(\Lambda)}\bigr)
G_\Lambda\ =\ G_{\sigma_k\Lambda}.
\]
The last equality implies Part (b) of Theorem 3.3 in the case
when $d_k(\Lambda)\geqslant2$, see the definition (\ref{q-smallf}).
Exchanging the tableaux $\Lambda$ and $\sigma_k\Lambda$ in (\ref{3.1}), so that
the resulting equality applies in the case when $d_k(\Lambda)\leqslant-2$, 
we prove Part (b) of Proposition \ref{T3.3} in this remaining case. \end{proof}

Thus the elements $G_\Lambda\in H_n$ for all
pairwise distinct standard tableaux $\Lambda$ of shape $\lambda$
form a basis in the vector space $V_\lambda$. This basis
is distinguished due to

\begin{proposition}\label{P3.4}
{\bf\hskip-6pt.\hskip1pt} 
We have $X_iG_\Lambda= q^{2c_i(\Lambda)}G_\Lambda$	
for each $i=1, \ldots , n$.
\end{proposition}

\begin{proof}
We will proceed by induction on $i=1, \ldots , n$. By definition,
$X_1=1$. On the other hand, $c_1(\Lambda)=0$ for any standard tableau
$\Lambda$. Thus Proposition~\ref{P3.4} is true for $i=1$.
Now suppose that Proposition \ref{P3.4} is true for $i=k$ where 
$k<n$. To show that it is also true for $i=k+1$, we will use 
Proposition \ref{T3.3}. Note that $X_{k+1}= T_kX_kT_k$.
If $d_k(\Lambda)=\pm1$, then $T_kX_kT_kG_\Lambda$ equals
\[
\mp q^{\mp1}T_kX_kG_\Lambda
=\mp q^{2c_k(\Lambda)\mp1}T_kG_\Lambda
=q^{2c_k(\Lambda)\mp2}G_\Lambda
=q^{2c_{k+1}(\Lambda)}G_\Lambda
\]
respectively. 
If $d_k(\Lambda)\geqslant2$, then the product $T_kX_kT_kG_\Lambda$ equals
\begin{eqnarray*}
&\phantom{x}& T_kX_k
\biggl(
\frac{q-q^{-1}}{1-q^{2d_k(\Lambda)}}G_\Lambda+G_{\sigma_k\Lambda}
\biggr)
=
\\%\]\[
&\phantom{x}& q^{2c_{k+1}(\Lambda)}
T_k
\biggl(
\frac{q-q^{-1}}{q^{-2d_k(\Lambda)}-1}G_\Lambda+
G_{\sigma_k\Lambda}
\biggr)=
\\%\]\[
&\phantom{x}& q^{2c_{k+1}(\Lambda)}
\biggl(
\frac{q-q^{-1}}{q^{-2d_k(\Lambda)}-1}
\biggl(
\frac{q-q^{-1}}{1-q^{2d_k(\Lambda)}}G_\Lambda+G_{\sigma_k\Lambda}
\biggr)+
\\%\]\[
&\phantom{x}& \frac{q-q^{-1}}{1-q^{-2d_k(\Lambda)}}G_{\sigma_k\Lambda}+ 
\biggl(
1-
\frac{(q-q^{-1})^{2}}{(q^{ d_k(\Lambda)}-q^{-d_k(\Lambda)})^{2}}
\biggr)
G_\Lambda
\biggr)
=
q^{2c_{k+1}(\Lambda)}G_\Lambda.
\end{eqnarray*}

In the case when $d_k(\Lambda)\leqslant-2$,
the proof of the equality $X_{k+1}G_\Lambda=q^{2c_{k+1}(\Lambda)}G_\Lambda$
is similar. \end{proof}
 
Let us now consider the left ideal $V_\lambda\subset H_n$ as 
$H_n$-module. Here the algebra $H_n$ acts via left multiplication.

\begin{corollary}\label{C3.5}
{\bf\hskip-6pt.\hskip1pt} 
The $H_n$-module $V_\lambda$ is irreducible.
\end{corollary}

\begin{proof}
Let $A$ be the commutative subalgebra of $H_n$ generated by the Murphy elements $X_1, \dots, X_n$. For any $H_n$ module $V$ and any $n$-tuple $\textbf{y} = (y_1, \dots, y_n)$ we define a subspace \[ V(\textbf{y}) = \{ v \in V : X_i v = y_i v \textrm{ for all } 1 \leqslant i \leqslant n \}. \] Let $V = V_\lambda$, then the vectors $G_\Lambda\in V_\lambda$, where $\Lambda$ is ranging over the set of all
standard tableaux of the given shape $\lambda$, form an eigenbasis for the
action on $V_\lambda$ of the Murphy elements in $H_n$. So each $G_\Lambda$ belongs to some $V(\textbf{y})$. Moreover, we have the direct sum \[ V_\lambda = \bigoplus_\textbf{y} V(\textbf{y}). \] If $\Lambda \neq M$ then $G_\Lambda \in V(\textbf{y})$ and $G_M \in V(\textbf{z})$ with $\textbf{y} \neq \textbf{z}$ - otherwise $c_i(\Lambda) = c_i(M)$ for all $i$. So $\dim V(\textbf{y}) \leqslant 1$. This means, as an $A$-module, $V_\lambda$ is a direct sum of 1-dimensional non-isomorphic $A$-submodules, namely the non-zero $V(\textbf{y})$.

Let $U$ be a $H_n$-submodule of $V_\lambda$. Then $U$ is also a submodule for $A$ and so contains a non-zero $V(\textbf{y}) = \mathbb{C}G_\Lambda$, for some standard tableau $\Lambda$. However, by Corollary \ref{C3.1} any basis vector $G_\Lambda\in V_\lambda$ 
can be obtained by acting on the element $G_{\Lambda^\circ}\in V_\lambda$ by 
a certain invertible element of $H_n$. So $G_{\Lambda^\circ} \in U$, but since $G_{\Lambda^\circ}$ generates all of $V_\lambda$ we have $U = V_\lambda$, hence $V_\lambda$ is simple. \end{proof}

\begin{corollary}\label{C3.6}
{\bf\hskip-6pt.\hskip1pt} 
The $H_n$-modules $V_\lambda$ for different partitions $\lambda$ of $n$
are pairwise non-equivalent.
\end{corollary}

\begin{proof}
Take any symmetric polynomial $f$ in $n$ variables
over the field $\mathbb{C}(q)$. For
all standard tableaux $\Lambda$ of the same shape $\lambda$, the values
of this polynomial
\begin{equation}\label{fval}
f\bigl(q^{2c_1(\Lambda)}, \ldots , q^{2c_n(\Lambda)}\bigr)\in\mathbb{C}(q)
\end{equation}
are the same.
Hence by Proposition \ref{P3.4}, the element $f(X_1, \ldots , X_n)\in H_n$ acts 
on $V_\lambda$ via multiplication by the scalar (\ref{fval}).
On the other hand, the partition $\lambda$ can be uniquely
restored from the values (\ref{fval}) where the polynomial $f$ varies.
Thus the $H_n$-modules $V_\lambda$ with different
partitions $\lambda$ cannot be equivalent. \end{proof}

%\emph{Remark.}\hskip6pt
%The centre of the algebra $\Hh_n$
%consists of all the Laurent polynomials in the generators $Y_1, \ldots , Y_n$
%which are invariant under permutations of these generators; 
%see for instance \cite[Proposition 3.11]{L}. In particular, the element
%$f(X_1, \ldots , X_n)\in H_n$ is central, as the image
%of a central element of $\Hh_n$ under the homomorphism  $\pi$.
%Moreover, the centre of the algebra $H_n$ coincides with the collection
%of all elements $f(X_1, \ldots , X_n)$ where the symmetric polynomial $f$ 
%varies; cf.\ \cite{J}. However, we do not use any of these
%facts in this section
%\qed
%
For any $k=1, \ldots , n-1$ consider the restriction of the $H_n$-module
$V_\lambda$ to the subalgebra $H_k\subset H_n$. %We use the
%standard embedding $H_k\to H_n$, where $T_i\mapsto T_i$ for 
%each index $i=1, \ldots , k-1$.

\begin{corollary}\label{Vkappa}
The vector $G_\Lambda\in V_\lambda$ belongs to the $H_k$-invariant subspace
in $V_\lambda$, equivalent to the $H_k$-module $V_{\kappa}$ where
the partition $\kappa$ is the shape of the tableau obtained by
removing from $\Lambda$ the entries $k+1, \ldots , n$.
\end{corollary}

\begin{proof}
It suffices to consider the case $k=n-1$ only.
For each index $a$ such that $\lambda_a>\lambda_{a+1}$,
denote by $V_a$ the vector subspace in $V_\lambda$ spanned by those vectors $G_\Lambda$ where $\Lambda(a,\lambda_a)=n$.
By Proposition \ref{T3.3}, the subspace
$V_a$ is preserved by the action 
of the subalgebra $H_{n-1}\subset H_n$ on $V_\lambda$.
Moreover, Proposition \ref{T3.3} shows that the $H_{n-1}$-module $V_a$
is equivalent to $V_{\kappa}$ where the partition $\kappa$ of $n-1$
is obtained by decreasing the $a^{\mbox{\scriptsize th}}$ part of $\lambda$ by $1$. \end{proof}

The properties of the vector $G_\Lambda$ given 
by Corollary \ref{Vkappa} for $k=1, \ldots , n-1$,
determine this vector 
in $V_\lambda$ uniquely up to a non-zero factor from $\mathbb{C}(q)$. 
%These properties can be restated for any irreducible
%$H_n$-module $V$ equivalent to $V_\lambda$.
%Explicit formulas for the action of the generators $T_1, \ldots , T_{n-1}$
%of $H_n$ on the vectors in $V$ determined by these properties,
%are known; cf.\ \cite[Theorem 6.4]{M1}.
% However, our proof of Theorem \ref{T3.3} here
% based on the results of Section 2, 
% is new 
Setting $q=1$, the algebra $H_n$ specializes to the symmetric group ring
$\mathbb{C}S_n$. The element $T_\sigma\in H_n$ then specializes to the
permutation $\sigma\in S_n$ itself. 
The proof of Proposition \ref{q-jimtheorem1} demonstrates that the coefficients
in the expansion of the element $F_\Lambda\in H_n$ relative to the basis
of the elements $T_\sigma$, are regular at $q=1$ as rational functions
of the parameter $q$. Thus the specialization of the element
$F_\Lambda\in H_n$ at $q=1$ is well defined. The same is true for the
element $G_\Lambda\in H_n$, see Corollary \ref{C3.1}.
The specializations at $q=1$ of the basis vectors $G_\Lambda\in V_\lambda$
form the \emph{Young seminormal basis} in the corresponding
irreducible representation of the group $S_n$.
The action of the generators $\sigma_1, \ldots ,\sigma_{n-1}$ of $S_n$ on the 
vectors of the latter basis was first given by \cite[Theorem IV]{Y2}.
For the interpretation of the elements $F_\Lambda$ and $G_\Lambda$ using
representation theory of the affine Hecke algebra $\widehat{H}_n$, 
see \cite[Section 3]{C2}.% and references therein.

%% file: Appendix.tex
\section{Appendix A}\label{appendix a}

This appendix gives the GAP code used to calculate the results of \ref{nonmixedhookresults}. My thanks to Alexander Konovalov and Alexander Hulpke and all those of the GAP support mailing list for their help, and patience!

\small
\begin{verbatim}
###IOTest###

#--------Mixed Hooks Problem-------

############################################################

#--------CancellationFunction-------
#        Added by AK 04/01/05

CancellationOfRationalFunction:=function(r)
local fam, num, den, numpol, denpol, div, gcdnum, gcdden, 
      coeffs, i, anum, aden;
if IsZero(r) then
  return r;
fi;
#
# Step 1. Divide numerator and denominator by their GCD
#
fam:=RationalFunctionsFamily(FamilyObj(1));;
num:=ExtRepNumeratorRatFun(r);
den:=ExtRepDenominatorRatFun(r);
numpol:=PolynomialByExtRep(fam,num);
denpol:=PolynomialByExtRep(fam,den);
div:=Gcd(numpol,denpol);
r:=(numpol/div)/(denpol/div);
#
# Step 2. Divide numerator and denominator by GCD of their
# coefficients and then multiply the result by their fraction
#
num:=ExtRepNumeratorRatFun(r);
den:=ExtRepDenominatorRatFun(r);
numpol:=PolynomialByExtRep(fam,num);
denpol:=PolynomialByExtRep(fam,den);
# making coefficients integer
coeffs:=CoefficientsOfUnivariatePolynomial(numpol);
anum:=Lcm(List( coeffs, i -> DenominatorRat(i)));
if anum<>1 then
  numpol:=numpol*anum;
fi;  
coeffs:=CoefficientsOfUnivariatePolynomial(denpol);
aden:=Lcm(List( coeffs, i -> DenominatorRat(i)));
if aden<>1 then
  denpol:=denpol*aden;
fi;  
gcdnum:=Gcd(CoefficientsOfUnivariatePolynomial(numpol));
gcdden:=Gcd(CoefficientsOfUnivariatePolynomial(denpol));
r:=(numpol/gcdnum)/(denpol/gcdden);
return (aden/anum)*(gcdnum/gcdden)*r;
end;

############################################################

#--------Trace Function-------

TraceOfProduct:=function(a,b)
local subtr, la, lb, i, j, inv, sum, s;
if IsZero(a) then 
  return ZeroCoefficient(a);
elif IsZero(b) then
  return ZeroCoefficient(b);
else
  #
  # Step 1. For each element x from the supp(a) we find its inverse
  # in the supp(b). If such exists, we multiply their coefficients
  # and store them in the list 'subtr'.
  #
  subtr :=[ ZeroCoefficient(a) ];
  la:=CoefficientsAndMagmaElements(a);
  lb:=CoefficientsAndMagmaElements(b);
  Print("Computing ", Length(la)/2, " coefficients \n");
  for i in [1,3..Length(la)-1] do
    Print( (i+1)/2, "\r");
    inv:=la[i]^-1;
    for j in [1,3..Length(lb)-1] do
      if lb[j]=inv then
        Add( subtr, la[i+1]*lb[j+1] );
        break;
      fi;
    od;
    # now we do cancellation on each step
    subtr[Length(subtr)]:=CancellationOfRationalFunction(subtr[Length(subtr)]);
  od;
  #
  # Step 2. Computing the sum of elements from the list 'subtr',
  # performing cancellation on each step.
  #
  sum:=ZeroCoefficient(a);
  Print("Computing the sum of ", Length(subtr), " coefficients \n");
  for s in [ 1 .. Length(subtr) ] do
    Print(s, "\r");
    sum:=sum+subtr[s];
    sum:=CancellationOfRationalFunction(sum);
  od;  
  return sum;
fi;
end;

############################################################

#--------Central Idempotent Function------

CentralIdempotentSymmetricGroup:= function( QG, pi )
local G, OrderG, Emb, n, t, charprm, classprm, pos, chi, idcoeffs,
      ccl, dom, cycletypes, sortedccl, elms, fusions, eC, j, val, k;

# Initialization
if not ( IsFreeMagmaRing( QG ) and
         IsGroup( UnderlyingMagma( QG ) ) and
         IsRationals( LeftActingDomain( QG ) ) ) then
    Error( "The input must be a rational group algebra" );
fi;
G := UnderlyingMagma( QG );
OrderG := Size( G );
Emb := Embedding( G, QG );
if not IsNaturalSymmetricGroup( G ) then
    Error( "<G> must be a symmetric group in its natural representation" );
fi;
n := NrMovedPoints( G );

# Compute the coefficients of the idempotent w.r.t. the character 
# table.
t := CharacterTable( "Symmetric", n );
charprm := CharacterParameters( t );
classprm := ClassParameters( t );
pos := Position( charprm, [ 1, pi ] );
if pos = fail then
    Error( "<pi> must be a partition of ", n );
fi;
chi := ValuesOfClassFunction( Irr( t )[ pos ] );
idcoeffs := ( chi[1] / OrderG ) * chi;

# Identify the classes of the character table and the group.
# (This is not needed if `G' was constructed with `SymmetricGroup'
# but in general the classes may be sorted in a different way.)
ccl := ConjugacyClasses( G );
dom := [ 1 .. n ];
cycletypes := List( ccl, c -> CycleLengths( Representative( c ), dom ) 
);
sortedccl := List( classprm, p -> ccl[ Position( cycletypes, p[2] ) ] 
);

# Form the group ring element.
# Compute the PCIs of QG.
elms:= Elements( G );
fusions:= List( sortedccl,
                c -> List( Elements( c ),
                           x -> PositionSorted( elms, x ) ) );
eC:= 0 * [ 1 .. OrderG ];
for j in [ 1 .. Length( chi ) ] do
    val:= idcoeffs[j];
    for k in fusions[j] do
        eC[k]:= val;
    od;
od;

return ElementOfMagmaRing( ElementsFamily( FamilyObj( QG ) ),
                           0, eC, elms );
end;

############################################################

#--------Function To Compute Eigenvalue----------

IOTest:=function(clambmu,Ilambmu,nu)

local r,z,u,fam,s8,g,QS,y,
#clambmu,Ilambmu,nu,
Flambdawithz,
elementsFlambdawithz,
coeffsFlambdawithz,
coeffevalFlambdawithz,
Flambda,
Fmuwithz,
elementsFmuwithz, 
coeffsFmuwithz,
coeffevalFmuwithz,
Fmu,
ZnuInQS,
Znu,
vnu,
Rlambdamu,
TraceRlambdamuvnu,
rlambdamunu,
num,
den,
numfactors,
denfactors,
elementsRlambdamu,
coeffsRlambdamu,
i,c;

############################################################

#--------The Initial Setup-----------

r:=PolynomialRing(Rationals,["z","u"]);;
z:=IndeterminatesOfPolynomialRing(r)[1];;
u:=IndeterminatesOfPolynomialRing(r)[2];;
fam:=RationalFunctionsFamily(FamilyObj(1));;
s8:=SymmetricGroup(8);;
g:=GroupRing(r,s8);;
QS:= GroupRing(Rationals,s8);;
y:=One(g);;

#clambmu:=[0,1,2,-1,0,1,2,-1];;
#Ilambmu:=[1,1,1,1,1,1,1,1];;
#nu:=[6,2];;

Print("Initial setup completed \n");

############################################################

#--------Calculating Flambda (two hooks only)---------

Print("Computing Flambda... \c");
Flambdawithz:=
(y - y*(3,4)/(z*((Ilambmu[3] - Ilambmu[4])^2) + clambmu[3] - clambmu[4]))*
(y - y*(2,4)/(z*((Ilambmu[2] - Ilambmu[4])^2) + clambmu[2] - clambmu[4]))*
(y - y*(2,3)/(z*((Ilambmu[2] - Ilambmu[3])^2) + clambmu[2] - clambmu[3]))*
(y - y*(1,4)/(z*((Ilambmu[1] - Ilambmu[4])^2) + clambmu[1] - clambmu[4]))*
(y - y*(1,3)/(z*((Ilambmu[1] - Ilambmu[3])^2) + clambmu[1] - clambmu[3]))*
(y - y*(1,2)/(z*((Ilambmu[1] - Ilambmu[2])^2) + clambmu[1] - clambmu[2]));;

elementsFlambdawithz:=Support(Flambdawithz);;
coeffsFlambdawithz:=CoefficientsBySupport(Flambdawithz);;
coeffevalFlambdawithz:=List(coeffsFlambdawithz,i->Value(i,[z],[0]));;
Flambda:=ElementOfMagmaRing(FamilyObj(Zero(g)), Zero(r),
One(r)*coeffevalFlambdawithz, elementsFlambdawithz);;

Unbind(Flambdawithz);;
Unbind(elementsFlambdawithz);;
Unbind(coeffsFlambdawithz);;
Unbind(coeffevalFlambdawithz);;
Print("OK \n");

############################################################

#--------Calculating Fmu (two hooks only)----------

Print("Computing Fmu... \c");
Fmuwithz:= 
(y - y*(7,8)/(z*((Ilambmu[7] - Ilambmu[8])^2) + clambmu[7] - clambmu[8]))*
(y - y*(6,8)/(z*((Ilambmu[6] - Ilambmu[8])^2) + clambmu[6] - clambmu[8]))*
(y - y*(6,7)/(z*((Ilambmu[6] - Ilambmu[7])^2) + clambmu[6] - clambmu[7]))*
(y - y*(5,8)/(z*((Ilambmu[5] - Ilambmu[8])^2) + clambmu[5] - clambmu[8]))*
(y - y*(5,7)/(z*((Ilambmu[5] - Ilambmu[7])^2) + clambmu[5] - clambmu[7]))*
(y - y*(5,6)/(z*((Ilambmu[5] - Ilambmu[6])^2) + clambmu[5] - clambmu[6]));;

elementsFmuwithz:=Support(Fmuwithz);;
coeffsFmuwithz:=CoefficientsBySupport(Fmuwithz);;
coeffevalFmuwithz:=List(coeffsFmuwithz,i->Value(i,[z],[0]));;
Fmu:=ElementOfMagmaRing(FamilyObj(Zero(g)), Zero(r), One(r)*coeffevalFmuwithz,
elementsFmuwithz);;

Unbind(Fmuwithz);;
Unbind(elementsFmuwithz);;
Unbind(coeffsFmuwithz);;
Unbind(coeffevalFmuwithz);;
Print("OK \n");

############################################################

#--------Calculating Znu----------

Print("Computing Znu... \c");
ZnuInQS:= CentralIdempotentSymmetricGroup(QS,nu);;
Znu:=ElementOfMagmaRing( 
           FamilyObj(Zero(g)), 
           Zero(r), 
           One(r)*CoefficientsBySupport(ZnuInQS), 
           Support(ZnuInQS));;
Print("OK \n");

############################################################

#--------Calculating Test Eigenvector vnu---------

Print("Computing vnu = Znu*Flambda*Fmu... \c");
vnu:= Znu*Flambda;
Unbind(Flambda);;
Unbind(Znu);;
vnu:= vnu*Fmu;;
Unbind(Fmu);;
Print("OK \n");

############################################################

#--------Calculating Operator Rlambdamu--------

Print("Computing Rlambdamu... \n");
Rlambdamu:= 
(y - y*(4,5)/(u + clambmu[4] - clambmu[5]))*
(y - y*(4,6)/(u + clambmu[4] - clambmu[6]))*
(y - y*(4,7)/(u + clambmu[4] - clambmu[7]))*
(y - y*(4,8)/(u + clambmu[4] - clambmu[8]))*
(y - y*(3,5)/(u + clambmu[3] - clambmu[5]))*
(y - y*(3,6)/(u + clambmu[3] - clambmu[6]))*
(y - y*(3,7)/(u + clambmu[3] - clambmu[7]))*
(y - y*(3,8)/(u + clambmu[3] - clambmu[8]))*
(y - y*(2,5)/(u + clambmu[2] - clambmu[5]))*
(y - y*(2,6)/(u + clambmu[2] - clambmu[6]))*
(y - y*(2,7)/(u + clambmu[2] - clambmu[7]))*
(y - y*(2,8)/(u + clambmu[2] - clambmu[8]))*
(y - y*(1,5)/(u + clambmu[1] - clambmu[5]))*
(y - y*(1,6)/(u + clambmu[1] - clambmu[6]))*
(y - y*(1,7)/(u + clambmu[1] - clambmu[7]))*
(y - y*(1,8)/(u + clambmu[1] - clambmu[8]));;

# AK 04/01/05 ###############################################################
                                                                            #
elementsRlambdamu:=Support(Rlambdamu);;                                     #
coeffsRlambdamu:=CoefficientsBySupport(Rlambdamu);;                         #
Unbind(Rlambdamu);;   
Print("Performing cancellations on ", Length(coeffsRlambdamu), " coefficients \n");
for i in [ 1 .. Length(coeffsRlambdamu) ] do
Print(i, "\r");
c:= coeffsRlambdamu[i];;                                                    #
c:=CancellationOfRationalFunction(c);;
od;                                                                         #
Rlambdamu:=ElementOfMagmaRing( FamilyObj(Zero(g)),                          #
                               Zero(r),                                     #
                               One(r)*coeffsRlambdamu,                      #
                               elementsRlambdamu);;                         #
                                                                            #
#############################################################################

Print("...OK \n");

############################################################

#--------Calculating Eigenvalue rlambdamunu--------

Print("Computing trace of Rlambdamu*vnu... \n");
TraceRlambdamuvnu:= TraceOfProduct(Rlambdamu, vnu);;
Unbind(Rlambdamu);;
Print("...OK \n");
Print("Computing rlambdamunu ... \c");
rlambdamunu:= TraceRlambdamuvnu/CoefficientsAndMagmaElements(vnu)[2];;
rlambdamunu:= CancellationOfRationalFunction( rlambdamunu );;
Print("OK \n");

#Print("rlambdamunu = ", rlambdamunu, "\n");

############################################################

#--------Factorising numerator and denominator--------

Print("Computing factors ... \n");
num:=ExtRepNumeratorRatFun(rlambdamunu);;
den:=ExtRepDenominatorRatFun(rlambdamunu);;
numfactors:= Factors(PolynomialByExtRep(fam,num));;
denfactors:= Factors(PolynomialByExtRep(fam,den));;
# return [numfactors, denfactors];
Print("numfactors = ", numfactors, "\n");
Print("denfactors = ", denfactors, "\n");
AppendTo("output.txt", "numfactors = ", numfactors, "\n");
AppendTo("output.txt", "denfactors = ", denfactors, "\n");
end;

############################################################

#--------End-----------
\end{verbatim}
\large

We then input the necessary parameters, these are the contents of $\lambda$ and $\mu$, a vector indicating whether a box is in the first or second hook of $\lambda$ and $\mu$, and the partition $\nu$. We consider six cases below.

\small
\begin{verbatim}
###input###

Read("IOTest.txt");

AppendTo("output.txt", "############ Test 1: \n");

clambmu:=[0,1,2,-1,0,1,2,-1];;
Ilambmu:=[1,1,1,1,1,1,1,1];;
nu:=[5,3];;
IOTest(clambmu,Ilambmu,nu);

AppendTo("output.txt", "############ Test 2: \n");

clambmu:=[0,1,2,-1,0,1,-1,0];;
Ilambmu:=[1,1,1,1,1,1,1,2];;
nu:=[4,2,2];;
IOTest(clambmu,Ilambmu,nu);

AppendTo("output.txt", "############ Test 3: \n");

clambmu:=[0,1,-1,0,0,1,-1,0];;
Ilambmu:=[1,1,1,2,1,1,1,2];;
nu:=[4,3,1];;
IOTest(clambmu,Ilambmu,nu);

AppendTo("output.txt", "############ Test 4: \n");

clambmu:=[0,1,-1,0,0,1,-1,0];;
Ilambmu:=[1,1,1,2,1,1,1,2];;
nu:=[3,2,2,1];;
IOTest(clambmu,Ilambmu,nu);

AppendTo("output.txt", "############ Test 5: \n");

clambmu:=[0,1,-1,0,0,1,-1,-2];;
Ilambmu:=[1,1,1,2,1,1,1,1];;
nu:=[3,3,1,1];;
IOTest(clambmu,Ilambmu,nu);

AppendTo("output.txt", "############ Test 6: \n");

clambmu:=[0,1,-1,-2,0,1,-1,-2];;
Ilambmu:=[1,1,1,1,1,1,1,1];;
nu:=[2,2,2,1,1];;
IOTest(clambmu,Ilambmu,nu);
\end{verbatim}
\large

Which outputs;

\small
\begin{verbatim}
###output###

############ 

Test 1: 
numfactors = [ u-4, u-1 ]
denfactors = [ u, u+1 ]

############ 

Test 2: 
numfactors = [ -u+3, u+4 ]
denfactors = [ -u+1, u+2 ]

############ 

Test 3: 
numfactors = [ -u+3, u-2, u+2 ]
denfactors = [ -u, u, u+1 ]

############ 

Test 4: 
numfactors = [ -u+2, u+2, u+3 ]
denfactors = [ -u+1, u, u ]

############ 

Test 5: 
numfactors = [ -u+2, u+3 ]
denfactors = [ -u, u+1 ]

############ 

Test 6: 
numfactors = [ u+1, u+4 ]
denfactors = [ u-2, u-1 ]
\end{verbatim}
\large

\section{Appendix B}\label{appendix b}

This appendix describes how we used the software packages GAP and Maple to calculate certain examples. These examples support the claim that the associated products $F_\Lambda(z_1, \dots, z_n)$ of invalid diagrams, as described in Section \ref{ribbonfp}, have a pole at $z_1 = \cdots = z_n$.

Consider the column tableau, $\Lambda^c$, of the Young diagram $\lambda = (3,3)$ below:

\[
\begin{picture}(60,40)
\put(0,20){\framebox(20,20)[r]{$1\;\;$}}
\put(0,0){\framebox(20,20)[r]{$2\;\;$}}
\put(20,20){\framebox(20,20)[r]{$3\;\;$}}
\put(20,0){\framebox(20,20)[r]{$4\;\;$}}
\put(40,20){\framebox(20,20)[r]{$5\;\;$}}
\put(40,0){\framebox(20,20)[r]{$6\;$}}
\end{picture}\]

We will proceed to calculate the coefficient of one of the elements in the expansion of $F_{\Lambda^c}$, and show this coefficient has a pole in the limit $z_1 = \cdots = z_n$. By observation we identify the permuation (1 4 3 6) as one such possible permutation.

If we order the product $F_\Lambda(z_1, \dots, z_n)$ lexicographically, the following code, written for GAP by James Woodward, tells us how (1 4 3 6) is generated in the expansion of the product.

\small
\begin{verbatim}
myList := [(1,2) , (1,3), (1,4), (1,5), (1,6), (2,3), (2,4), (2,5), (2,6),
           (3,4), (3,5), (3,6), (4,5), (4,6), (5,6)];

chosenElement := (1,4,3,6);
Count := 0;

for i in [1..2^Length(myList)-1] do

  g := ();
  myString := "";

  for j in [1..Length(myList)] do
    if Int(i/(2^(j-1))) mod 2 = 1 then 
       g := myList[j] * g; 
       myString := Concatenation(myString , String(myList[j]));
    fi;
  od;

  if g = chosenElement then
    Count := Count + 1;
    Print(Count, " ", g, "=" , myString, "\n");
  fi;

od;

Print("There are ", Count, " ways to do it! \n");
\end{verbatim}
\large

This gives us the following list:

\small
\begin{verbatim}
1 (1,4,3,6)=(1,3)(1,4)(3,6)
2 (1,4,3,6)=(1,2)(1,3)(1,4)(2,3)(3,6)
3 (1,4,3,6)=(1,2)(1,4)(2,3)(2,4)(3,6)
4 (1,4,3,6)=(1,4)(3,4)(3,6)
5 (1,4,3,6)=(1,2)(1,4)(2,4)(3,4)(3,6)
6 (1,4,3,6)=(1,4)(3,6)(4,6)
7 (1,4,3,6)=(1,2)(1,4)(2,4)(3,6)(4,6)
8 (1,4,3,6)=(1,3)(1,4)(3,4)(3,6)(4,6)
9 (1,4,3,6)=(1,2)(1,3)(1,4)(2,3)(3,4)(3,6)(4,6)
10 (1,4,3,6)=(1,2)(1,4)(2,3)(2,4)(3,4)(3,6)(4,6)
11 (1,4,3,6)=(1,3)(1,4)(3,5)(3,6)(5,6)
12 (1,4,3,6)=(1,2)(1,3)(1,4)(2,3)(3,5)(3,6)(5,6)
13 (1,4,3,6)=(1,2)(1,4)(2,3)(2,4)(3,5)(3,6)(5,6)
14 (1,4,3,6)=(1,4)(3,4)(3,5)(3,6)(5,6)
15 (1,4,3,6)=(1,2)(1,4)(2,4)(3,4)(3,5)(3,6)(5,6)
16 (1,4,3,6)=(1,4)(3,5)(3,6)(4,5)(5,6)
17 (1,4,3,6)=(1,2)(1,4)(2,4)(3,5)(3,6)(4,5)(5,6)
18 (1,4,3,6)=(1,3)(1,4)(3,4)(3,5)(3,6)(4,5)(5,6)
19 (1,4,3,6)=(1,2)(1,3)(1,4)(2,3)(3,4)(3,5)(3,6)(4,5)(5,6)
20 (1,4,3,6)=(1,2)(1,4)(2,3)(2,4)(3,4)(3,5)(3,6)(4,5)(5,6)
21 (1,4,3,6)=(1,4)(3,6)(4,5)(4,6)(5,6)
22 (1,4,3,6)=(1,2)(1,4)(2,4)(3,6)(4,5)(4,6)(5,6)
23 (1,4,3,6)=(1,3)(1,4)(3,4)(3,6)(4,5)(4,6)(5,6)
24 (1,4,3,6)=(1,2)(1,3)(1,4)(2,3)(3,4)(3,6)(4,5)(4,6)(5,6)
25 (1,4,3,6)=(1,2)(1,4)(2,3)(2,4)(3,4)(3,6)(4,5)(4,6)(5,6)

There are 25 ways to do it!
\end{verbatim}
\large

We now use Maple to calculate the coefficient of (1 4 3 6) in the expansion of $F_\Lambda(z_1, \dots, z_n)$. We first do this without imposing any conditions on the auxiliary parameters $z_i$.

\small
\begin{verbatim}
># The coefficients of each term in the product

> c12:=(-1)/(z1-z2+1):
> c13:=(-1)/(z1-z3-1):
> c14:=(-1)/(z1-z4):
> c15:=(-1)/(z1-z5-2):
> c16:=(-1)/(z1-z6-1):
> c23:=(-1)/(z2-z3-2):
> c24:=(-1)/(z2-z4-1):
> c25:=(-1)/(z2-z5-3):
> c26:=(-1)/(z2-z6-2):
> c34:=(-1)/(z3-z4+1):
> c35:=(-1)/(z3-z5-1):
> c36:=(-1)/(z3-z6):
> c45:=(-1)/(z4-z5-2):
> c46:=(-1)/(z4-z6-1):
> c56:=(-1)/(z5-z6+1):

># The coefficients of each word that makes (1,4,3,6)

> w1:=c13*c14*c36:
> w2:=c12*c13*c14*c23*c36:
> w3:=c12*c14*c23*c24*c36:
> w4:=c14*c34*c36:
> w5:=c12*c14*c24*c34*c36:
> w6:=c14*c36*c46:
> w7:=c12*c14*c24*c36*c46:
> w8:=c13*c14*c34*c36*c46:
> w9:=c12*c13*c14*c23*c34*c36*c46:
> w10:=c12*c14*c23*c24*c34*c36*c46:
> w11:=c13*c14*c35*c36*c56:
> w12:=c12*c13*c14*c23*c35*c36*c56:
> w13:=c12*c14*c23*c24*c35*c36*c56:
> w14:=c14*c34*c35*c36*c56:
> w15:=c12*c14*c24*c34*c35*c36*c56:
> w16:=c14*c35*c36*c45*c56:
> w17:=c12*c14*c24*c35*c36*c45*c56:
> w18:=c13*c14*c34*c35*c36*c45*c56:
> w19:=c12*c13*c14*c23*c34*c35*c36*c45*c56:
> w20:=c12*c14*c23*c24*c34*c35*c36*c45*c56:
> w21:=c14*c36*c45*c46*c56:
> w22:=c12*c14*c24*c36*c45*c46*c56:
> w23:=c13*c14*c34*c36*c45*c46*c56:
> w24:=c12*c13*c14*c23*c34*c36*c45*c46*c56:
> w25:=c12*c14*c23*c24*c34*c36*c45*c46*c56:

># The coefficient of (1436) in the product

> simplify(w1+w2+w3+w4+w5+w6+w7+w8+w9+w10+w11+w12
              +w13+w14+w15+w16+w17+w18+w19+w20+w21+w22+w23+w24+w25);
\end{verbatim}
\large

The result is a rational polynomial in $z_1, z_2, \dots, z_6$. The numerator of which is too large to reproduce here. The denominator however is calculated to be \begin{eqnarray*} & & (z_1 - z_2 +1)(z_1 - z_4)(z_2-z_4-1)(z_3-z_5-1)(z_3-z_6)(z_4-z_5-2)\\ && \times (z_5-z_6+1)(z_1-z_3-1)(z_3-z_4+1)(z_2-z_3-2)(z_4-z_6-1).\end{eqnarray*}

The terms $(z_1-z_4)$ and $(z_3-z_6)$ in the denominator show we have a non-removable singularity at $z_1 = \cdots = z_6$ here.

For comparison let us first consider the subspace defined in the row fusion procedure by the diagram

\[
\begin{picture}(60,40)
\put(0,20){\framebox(20,20)[r]{$1\;\;$}}
\put(0,0){\framebox(20,20)[r]{$2\;\;$}}
\put(20,20){\framebox(20,20)[r]{$3\;\;$}}
\put(20,0){\framebox(20,20)[r]{$4\;\;$}}
\put(40,20){\framebox(20,20)[r]{$5\;\;$}}
\put(40,0){\framebox(20,20)[r]{$6\;$}}

\put(10,30){\line(1,0){20}}%1
\put(10,10){\line(1,0){20}}%2
\put(30,30){\line(1,0){20}}%5
\put(30,10){\line(1,0){20}}%6
\end{picture}\]

We use Maple to calculate the coefficient of (1 4 3 6) on this subspace.

\small
\begin{verbatim}
> z3:=z1:
> z5:=z1:
> z4:=z2:
> z6:=z2:

> simplify(w1+w2+w3+w4+w5+w6+w7+w8+w9+w10+w11+w12
           +w13+w14+w15+w16+w17+w18+w19+w20+w21+w22+w23+w24+w25);

              1                /    3        2       2                      
------------------------------ \2 z1  + 10 z1  - 6 z1  z2 + 12 z1 - 20 z1 z2
             3               2                                              
(z1 - z2 + 1)  (-z2 + z1 + 2)                                               

            2                   3        2\
   + 6 z1 z2  - 1 - 12 z2 - 2 z2  + 10 z2 /
\end{verbatim}
\large

i.e. the coefficient is

\[\frac{2z_1^3 + 10z_1^2 - 6z_1^2z_2 + 12z_1 - 20z_1z_2 + 6z_1z_2^2 - 12z_2-1-2z_2^3+10z_2^2}{(z_1-z_2+1)^3(-z_2+z_1+2)^2}.\] 

So on the subspace where $z_1=z_3=z_5$ and $z_2=z_4=z_6$ the coefficient tends to -1/4 in the limit $z_1 = z_2$.

Now consider an invalid diagram, for example a diagram that contains a square of type $\bigtype{15}$

\[
\begin{picture}(60,40)
\put(0,20){\framebox(20,20)[r]{$1\;\;$}}
\put(0,0){\framebox(20,20)[r]{$2\;\;$}}
\put(20,20){\framebox(20,20)[r]{$3\;\;$}}
\put(20,0){\framebox(20,20)[r]{$4\;\;$}}
\put(40,20){\framebox(20,20)[r]{$5\;\;$}}
\put(40,0){\framebox(20,20)[r]{$6\;$}}

\put(10,30){\line(1,0){20}}%1
\put(10,10){\line(1,0){20}}%2
\end{picture}\]

Then the coefficient of (1 4 3 6) is given by

\small
\begin{verbatim}
> z3:=z1:
> z6:=z4:

> c:=simplify(w1+w2+w3+w4+w5+w6+w7+w8+w9+w10+w11+w12
                  +w13+w14+w15+w16+w17+w18+w19+w20+w21+w22+w23+w24+w25);

     /    2   2       2   2       2   3     3   3          3     
c := \3 z1  z5  + 4 z1  z6  + 5 z2  z5  - z2  z5  + 5 z1 z5  - z1

            3             2                2           2              
   + 2 z1 z2  z5 z6 + 5 z1  z2 z5 z6 - 8 z1  - z1 z2 z5  z6 + 10 z1 z2

       2         2     2   2          2      2     2      2     3   2   
   - z1  z2 z5 z6  - z1  z5  z6 - 6 z1  z2 z6  + z1  z2 z5  - z1  z6  z2

          2      2       3   2     3   2           2       3          3      
   - z1 z2  z5 z6  - 9 z2  z5  + z2  z6  + 20 z1 z5  + 7 z2  z6 + 6 z2  z5 z6

          2          2           2   2          3          2      2       3
   + z1 z2  z5 - 7 z5  z2 z6 - z5  z6  + 3 z1 z5  z6 + 3 z5  z2 z6  - 3 z1 

          3   2           2           3                   3          2      
   - z1 z2  z6  - 10 z1 z2  z5 z6 - z1  z5 z6 - 2 z1 z2 z5  z6 - 6 z1  z5 z6

         2                3             3          2   2           2   2   
   - 8 z1  z2 z6 + 4 z1 z2  z5 - 4 z1 z2  z6 + 8 z1  z2  z6 - 11 z1  z2  z5

       2   2           3   2                2   2        2   2              
   - z1  z2  z5 z6 + z1  z6  z5 + 2 z1 z2 z5  z6  - z1 z2  z5  z6 + 20 z1 z5

          2                                3                2   2     2   3   
   + 32 z2  z5 - 14 z1 z6 + 14 z2 z6 - 3 z5  z2 z6 - 3 z1 z5  z6  + z2  z5  z6

       2   2   2        2                        2        3                 
   - z2  z5  z6  - 17 z2  z6 + 4 z5 z6 - 15 z1 z2  - 14 z2  z5 + 16 z1 z2 z6

          2   2     2   3        2   2   2     2      3        2      
   - 11 z1  z2  + z1  z5  z6 - z1  z5  z6  - z1  z2 z5  + 18 z1  z2 z5

         2        2                2       2   2          2      2
   - 7 z1  z5 - z1  z6 - 8 z1 z5 z6  + 2 z2  z5  z6 - 6 z2  z5 z6 

            2   2     3      2          2   3           2   2       2   
   + 3 z1 z2  z6  + z2  z5 z6  + 2 z1 z2  z5  + 10 z1 z2  z5  + 4 z5  z6

         3                3          2             2       2              2
   - 6 z5  z2 - 7 z1 z2 z5  + 6 z1 z5  z6 - 3 z5 z6  - 2 z6  + 10 z5 z2 z6 

         3             2   3          3              2        2   
   + 5 z1  z2 z5 + 2 z1  z5  + 7 z1 z2  - 27 z1 z2 z5  - 24 z5  z2

         3           3          2      2                        2      
   - 3 z1  z2 z6 + z5  z6 + 2 z1  z2 z5  z6 - 6 z1 z5 z6 - 11 z2  z5 z6

                                        3   2                           2     
   + 4 z5 z2 z6 + 16 z1 z2 z5 z6 - z1 z2  z5  - 24 z5 z2 + 6 z1 z2 z5 z6  + z6

          2   2          2          2       2   2                     2   3   
   + 26 z2  z5  - 4 z1 z6  + 7 z2 z6  - 5 z2  z6  - 20 z1 z2 z5 - 2 z1  z2  z6

       3   2          3          3          3          3   2     4   2
   + z1  z2  z6 + 4 z1  z2 - 6 z1  z5 + 2 z1  z6 - 2 z1  z5  + z2  z5 

       4          4          2   3        4     3   2       3   2        4   
   - z2  z6 + 2 z2  z5 + 2 z1  z2  - z1 z2  - z1  z2  + 2 z1  z6  - z1 z2  z5

          2          2   3           4        3   2          2   2   2
   + 17 z1  z2 + 2 z1  z2  z5 + z1 z2  z6 - z1  z5  z6 + 2 z1  z2  z6 

       2   2   2     4           3   2        3      2\//                  
   - z1  z2  z5  - z2  z5 z6 - z1  z2  z5 + z1  z2 z5 / \(z2 - z5 - 2) (-z5

                                                       2                  
   + z6 - 1) (z1 - z5 - 1) (-z2 + z1 + 2) (z1 - z2 + 1)  (z2 - z6 - 1) (z1

        \
   - z6)/
\end{verbatim}
\large

Importantly, the denominator contains the term $(z_1-z_6)$ and so the singularity is not removable.

Finally we consider one more invalid diagram, a diagram that contains the union \type{3}--\type{5}

\[
\begin{picture}(60,40)
\put(0,20){\framebox(20,20)[r]{$1\;\;$}}
\put(0,0){\framebox(20,20)[r]{$2\;\;$}}
\put(20,20){\framebox(20,20)[r]{$3\;\;$}}
\put(20,0){\framebox(20,20)[r]{$4\;\;$}}
\put(40,20){\framebox(20,20)[r]{$5\;\;$}}
\put(40,0){\framebox(20,20)[r]{$6\;$}}

\put(10,30){\line(1,0){20}}%1
\put(30,10){\line(1,0){20}}%6
\end{picture}\]

Then the coefficient of (1 4 3 6) is given by

\small
\begin{verbatim}
> z3:=z1:
> z6:=z4:

> c:=simplify(w1+w2+w3+w4+w5+w6+w7+w8+w9+w10+w11+w12
                    +w13+w14+w15+w16+w17+w18+w19+w20+w21+w22+w23+w24+w25);

     /     3             3             3      2       3      2       3   2   
c := \-6 z1  z2 z4 + 6 z1  z2 z5 + 2 z1  z2 z5  - 2 z1  z4 z5  + 4 z1  z4  z5

         3   3       3          2                         2   2   
   - 2 z1  z4  + 3 z4  - 8 z4 z5  z2 + 16 z4 z5 z2 + 12 z4  z5  z2

          2                3          3          2          3             3   
   + 22 z4  z5 z2 - 6 z4 z5  z2 + 3 z4  z5 + 7 z4  z2 - 6 z4  z5 z2 - 6 z4  z2

            3       3          2                   2           2   
   + 3 z4 z5  - 8 z5  z2 - 8 z1  + 12 z2 z1 + 15 z1  z2 - 12 z2  z1

         2   2     2           2              2                    
   - 6 z1  z2  + z1  z4 - 19 z1  z2 z4 + 11 z2  z1 z4 + 10 z2 z1 z4

         2      2       2      2              2       3           2   
   - 4 z1  z2 z4  + 2 z2  z1 z4  - 10 z2 z1 z4  + 3 z2  z1 - 20 z2  z4

         3           2   2       2   2       3          2           2   
   + 5 z2  z4 + 10 z1  z4  + 3 z2  z4  - 8 z2  z5 - 7 z1  z5 + 32 z2  z5

          2   2          3             2      3       2         
   + 12 z1  z2  z4 - 6 z2  z1 z4 + 4 z1  z2 z4  + 2 z1  z2 z4 z5

                          2      3             3       2   2   2
   + 36 z2 z1 z4 z5 - 2 z2  z1 z4  + 2 z2 z1 z4  - 4 z1  z2  z4 

         3      2                      2             2              2   2   
   + 2 z2  z1 z4  - 16 z2 z1 z5 + 22 z1  z2 z5 - 8 z2  z1 z5 - 12 z1  z2  z5

         3   2       2   3       2   3                   2         
   - 2 z2  z4  - 4 z1  z4  + 2 z2  z4  - 16 z1 z4 - 10 z2  z1 z4 z5

          2      2          2      2                2          2   2      
   - 10 z1  z2 z4  z5 + 8 z2  z1 z4  z5 - 2 z2 z1 z4  z5 + 8 z1  z2  z4 z5

         3                3              2              3             2   2   
   - 4 z2  z1 z4 z5 + 6 z2  z1 z5 - 32 z2  z4 z5 + 10 z2  z4 z5 + 4 z1  z4  z5

         2   2                            2           2                     3
   - 2 z2  z4  z5 - 10 z1 z4 z5 - 19 z1 z4  z5 - z1 z4  + 20 z1 z5 + 3 z1 z4 

             2        2   2              2       2      2           2   
   + 20 z1 z5  + 32 z2  z5  - 32 z2 z1 z5  + 8 z2  z1 z5  + 11 z1 z5  z4

         3   2              2                2   2             3   
   - 8 z2  z5  + 10 z2 z1 z5  z4 + 8 z2 z1 z5  z4  - 4 z2 z1 z5  z4

         2   2                 3       2      3           2   2          3   
   - 8 z2  z5  z4 - 10 z2 z1 z5  + 4 z2  z1 z5  - 12 z1 z5  z4  + 6 z1 z5  z4

         2   2   2       2   3             3       2   3       3   3
   - 4 z2  z5  z4  + 2 z2  z5  z4 + 5 z1 z5  + 8 z2  z5  - 2 z2  z5 

         3   2           2      2             3          2      2       2   2
   + 4 z2  z5  z4 - 10 z2  z1 z5  z4 + 6 z1 z4  z5 + 2 z1  z2 z5  + 3 z1  z5 

         2      2          3      2       2   2   2        2      
   + 8 z1  z2 z5  z4 + 2 z2  z1 z5  - 4 z1  z2  z5  - 10 z1  z5 z4

         2   2          2   2   2             3          3   2   
   - 2 z1  z5  z4 - 4 z1  z5  z4  - 4 z2 z1 z4  z5 - 2 z2  z4  z5

         2   3          2   3          2   3       2      3          
   + 2 z1  z4  z5 + 2 z2  z4  z5 + 2 z1  z5  - 2 z1  z2 z5  - z1 + z4

         3             3       3          3          3   2       3   2
   - 2 z1  z4 z5 - 3 z1  + 3 z1  z4 - 6 z1  z5 + 4 z1  z4  - 2 z1  z5 

                        2       2   2        2                     2
   + 12 z4 z5 + 12 z4 z5  - 6 z4  z5  - 15 z4  z5 + 20 z4 z2 - 8 z4 

          2                     3      2       3                2   3   
   - 32 z5  z2 - 32 z5 z2 + 2 z1  z2 z4  - 4 z1  z2 z4 z5 + 2 z1  z5  z4

         3   \/((-z5 + z4 - 1) (z4 - z5 - 2) (z1 - z5 - 1) (z1 - z4 + 1) (z2
   + 3 z1  z2/                                                              

   - z4 - 1) (-z2 + z1 + 2) (z1 - z4) (z1 - z2 + 1))
\end{verbatim}
\large

Here the denominator contains the term $(z_1-z_4)$, this shows the singularity is not removable.

%% file: References.tex
%\section{References}

%\addcontentsline{toc}{}{\textbf{References \qquad\qquad\qquad\qquad\qquad\qquad\qquad\qquad\qquad\qquad\qquad\qquad}}
\addtocontents{toc}{\vspace{12pt} \textbf{References \qquad\qquad\qquad\qquad\qquad\qquad\qquad\qquad\qquad\qquad\qquad\quad 132}}